\documentclass[a4paper,10pt]{article}

\usepackage[T1]{fontenc}
\usepackage{lmodern}
\usepackage{microtype}
\usepackage[margin=20mm]{geometry}
\usepackage{amsmath,amssymb,amsthm,mathtools}
\usepackage{mathrsfs}
\usepackage{array}
\usepackage{booktabs}
\usepackage{enumitem}
\usepackage{needspace}
\usepackage{xcolor}
\usepackage{tikz}
\usetikzlibrary{arrows.meta,positioning}
\usepackage{caption}
\usepackage[authoryear,square]{natbib}
\setcitestyle{aysep={},citesep={;}}
\usepackage[colorlinks=true,linkcolor=blue!55!black,citecolor=blue!55!black,urlcolor=blue!55!black]{hyperref}
\let\cite\citep
\let\citet\citep

\newtheorem{theorem}{Theorem}[section]
\newtheorem{proposition}[theorem]{Proposition}
\newtheorem{lemma}[theorem]{Lemma}
\newtheorem{corollary}[theorem]{Corollary}
\newtheorem{conjecture}[theorem]{Conjecture}
\theoremstyle{definition}
\newtheorem{definition}[theorem]{Definition}
\theoremstyle{remark}
\newtheorem{remark}[theorem]{Remark}
\newlist{romanenu}{enumerate}{1}
\setlist[romanenu]{label=\textup{\roman*)},ref=\roman*,
  leftmargin=0pt,itemindent=0pt,labelsep=.5em,align=right}

\newcommand{\Q}{\mathbb Q}
\newcommand{\C}{\mathbb C}
\newcommand{\Acon}{A_{\mathrm{con}}}
\newcommand{\Spec}{\operatorname{Spec}}
\newcommand{\Spf}{\operatorname{Spf}}
\newcommand{\Der}{\operatorname{Der}}
\newcommand{\End}{\operatorname{End}}
\newcommand{\Stab}{\operatorname{Stab}}
\newcommand{\GL}{\operatorname{GL}}
\newcommand{\Aut}{\operatorname{Aut}}
\newcommand{\Soc}{\operatorname{Soc}}
\newcommand{\Com}{\mathrm{Com}}
\newcommand{\Sm}{\mathrm{Sm}}
\newcommand{\cent}{\mathrm{cent}}
\newcommand{\su}{\mathsf u}
\newcommand{\sv}{\mathsf v}
\newcommand{\sw}{\mathsf w}
\newcommand{\sq}{\mathsf q}
\newcommand{\mail}[1]{\href{mailto:#1}{\texttt{#1}}}

\title{Semisimple Deformations of Contraction Algebras\\
for Smooth Threefold Flops}
\author{Shutao Bi\thanks{Email: \mail{shutaobi@mail.sdu.edu.cn}. Homepage: \url{https://shutaobi.github.io/}.}}
\date{October 2026}

\numberwithin{equation}{section}
\hypersetup{
 pdftitle={Semisimple Deformations of Contraction Algebras for Smooth Threefold Flops},
 pdfauthor={Shutao Bi},
 pdfsubject={Deformations of contraction algebras, schemes of associative algebra laws, and Gopakumar--Vafa invariants},
 bookmarksnumbered=true, bookmarksopen=false}

\begin{document}
\maketitle

\begin{abstract}
Building on the work of Bryan--Katz--Leung, Hua--Toda and Kawamata, we study deformations of contraction algebras and the local geometry of the scheme of associative algebra laws.  Let \(A\) be the basic contraction algebra of a completed projective flopping contraction of a smooth complex threefold. We construct a finite free polynomial family that becomes a constant semisimple algebra away from the origin, with the original curve projectors fixed.  For each restricted positive root \(\beta=(\beta_i)\), with entries indexed by the retained vertices, the split generic fibre has \(n_\beta\) matrix blocks of size \(\sum_i\beta_i\), and the labelled projector \(e_i\) has rank \(\beta_i\) on each block. Here \(n_\beta\) is the corresponding genus-zero Gopakumar--Vafa invariant.  If the geometric generic fibre of a finite free associative formal deformation of \(A\) is semisimple, it has the same matrix factors, and any lifts of the original projectors have the same ranks on those factors.  The scheme of associative algebra laws is smooth at the point defined by \(A\), and its unique irreducible component through that point is the closure of the corresponding semisimple orbit.

The proof uses the curve divisors described by Katz--Morrison and Nabijou--Wemyss, presentations of the universal contraction quotient and its cocenter, and the orbit method of Jelisiejew--Smoktunowicz.  The presentations also prove the geometric Hua--Toda flatness assertions: in the geometric families considered here, the relative contraction algebra, its centre and its cocenter are finite free over arbitrary complete Noetherian local complex bases.

On a smooth open neighbourhood of \(A\) in the scheme of algebra laws, the centres and cocenters form dual vector bundles compatible with arbitrary base change.  For the type \(D\) algebras \(D_{n,2}\), \(n\ge2\), the centre morphism is not smooth.  These are counterexamples to the smoothness assertion in Conjecture~1.4(2) of Jelisiejew--Smoktunowicz.

\end{abstract}

\clearpage
\begingroup \small
\tableofcontents
\endgroup
\clearpage

\section{Introduction}
\label{sec:introduction}

The geometric families of \citet[Theorem~4.2]{HuaToda} and \citet[Theorem~3.1(2)--(3)]{KawamataNC} have semisimple generic fibres. We ask whether the contraction algebra determines the matrix factors of every geometrically semisimple generic fibre of a finite free associative formal deformation. When the exceptional locus has several curves, we also retain the projectors indexed by the original curves and ask whether their ranks in these matrix factors are determined.

For a product of matrix algebras over a field, the dimension of its cocenter is the number of factors. We compute \(HH_0(A)=A/[A,A]\) for the contraction algebra \(A\). Together with the orbit method of \citet[Section~1 and Theorem~2.14]{JS}, this calculation determines the possible semisimple geometric generic fibres and proves that \(A\) is a smooth point of the scheme of associative algebra laws. On a smooth open neighbourhood of this point, the centres of the parametrized algebras form a vector bundle compatible with base change.

Let \(f:X\to\Spec R\) be the completion at the exceptional image point of a simple flopping contraction of a smooth complex threefold, with exceptional curve \(C\cong\mathbb P^1\). The contraction algebra of \citet{DonovanWemyss} represents the noncommutative deformation functor of \(\mathcal O_C(-1)\). We use the opposite-algebra convention corresponding to right actions on the deformed sheaf. For the basic Van den Bergh generator \(R\oplus N\), put
\[
 \Lambda=\End_R(R\oplus N)^{\mathrm{op}},\qquad
 A=\Acon=\Lambda/\Lambda e_0\Lambda,
\]
where \(e_0\) projects onto the free summand. The ideal generated by this projector consists of the endomorphisms that factor through finite sums of the free summand.

Bryan, Katz and Leung deform the neighbourhood of \(C\) to a family with isolated \((-1,-1)\)-curves and identify the number of such curves in each homology class with the corresponding Hilbert-scheme multiplicity \cite[Section~2, Proposition~2.3 and Lemma~2.8]{BKL}. For a contraction of length \(\ell\), these multiplicities are the genus-zero GV invariants \(n_1,\ldots,n_\ell\); \(n_j\) counts curves of class \(j[C]\) \cite{BKL,Toda}. Katz's local DT--GV construction recovers the same multiplicities for contractible curves \cite[Proposition~3.3 and Corollary~3.4]{KatzGV}. The invariants \(n_j\) define the semisimple algebra
\begin{equation}
\label{eq:general-gv-target}
 S_{\mathrm{GV}}=\prod_{j=1}^{\ell}M_j(\C)^{\times n_j}.
\end{equation}
Toda's width formula gives \(\dim_\C A=\sum_jj^2n_j=\dim_\C S_{\mathrm{GV}}\) \cite[Theorem~1.1]{Toda}.

For finite-dimensional algebras over an algebraically closed field \(k\) of characteristic zero, we say that \(A\) \emph{deforms to} \(S\) if a finite free unital \(k[[s]]\)-algebra has special fibre \(A\) and geometric generic fibre \(S\otimes_k\overline{k((s))}\). Jelisiejew--Smoktunowicz prove that this deformation notion agrees with the polynomially split and orbit-closure formulations \cite[Theorem~2.20]{JS}. In the polynomially split formulation, the family is defined over \(k[u]\) and becomes the constant algebra \(S\) after inverting \(u\). Thus, if the geometric generic fibre is semisimple, \(S\) is a product of matrix algebras over \(k\), and its block sizes and multiplicities agree with those of that fibre.

\begin{conjecture}[Wemyss; cf.\ {\cite[Questions~1.1--1.2]{JS}}]
\label{conj:wemyss-gv}
The algebra \(A\) deforms to \(S_{\mathrm{GV}}\), and every complex semisimple algebra to which \(A\) deforms is isomorphic to \(S_{\mathrm{GV}}\).
\end{conjecture}

\subsection*{Semisimple generic fibres and smoothness of algebra laws}

For the general statements, let \(f_0:X_0\to Y_0\) be the completion at an exceptional image point of a projective flopping contraction of a smooth complex threefold. Let \(A\) be its contraction algebra. Fix the normalized basic Van den Bergh generator and its labelled nonfree projectors \((e_i)_{i\in J}\). A general surface section gives an ADE diagram with retained subset \(J\).  We use the simultaneous ADE resolutions of \citet[Theorem~1]{KatzMorrison}, together with the restricted-root formulation and perturbation construction of \citet[Theorem~2.6 and Section~2.6]{NabijouWemyss}.  A nonzero restricted positive root \(\beta=(\beta_i)_{i\in J}\) is the vector of retained coefficients of a positive root; distinct proportional vectors remain distinct.

Let \(n_\beta\) be the corresponding genus-zero GV invariant, using the curve-counting description in \cite[Section~3.1]{NabijouWemyss}, and set
\[
 s_\beta=\sum_i\beta_i,\qquad
 N_{\mathrm{GV}}=\sum_\beta n_\beta,\qquad
 E=\bigoplus_{i\in J}\C e_i.
\]

A marking is a fixed embedding \(E\hookrightarrow A\), with each projector assigned to its labelled exceptional curve.  The associated Peirce spaces are the vector spaces \(e_iAe_j\); together they decompose the underlying vector space of \(A\).  On a matrix factor, the marking is described by the ranks of the orthogonal projectors \(e_i\), whose sum is the identity of that factor.

Put \(d=\dim_\C A\). Write \(X_d\) for the scheme of associative laws on a fixed \(d\)-dimensional vector space with fixed identity, and let \(G=\Stab_{\GL_d}(1)\) be the group of allowed changes of basis. The fixed-idempotent scheme \(X_E\) also fixes the Peirce spaces \(e_iAe_j\) and the labelled vectors \(e_i\); denote its change-of-basis group by \(G_E\).  Orbits in these schemes represent isomorphism classes with the indicated data fixed.

\begin{theorem}[Semisimple generic fibres and smoothness]
\label{thm:intro-marked-semisimple}
For the basic contraction algebra \(A\) of \(f_0\), put
\[
 S_E=\prod_\beta
 \End_\C\left(\bigoplus_{i\in J}\C^{\beta_i}\right)^{\times n_\beta}.
\]
\begin{enumerate}[label=\textup{(\roman*)}]
\item There is a finite free marked \(\C[u]\)-algebra \(\mathscr P\) with \(\mathscr P/(u)\cong A\) and \(\mathscr P[u^{-1}]\cong S_E\otimes_\C\C[u,u^{-1}]\). On each \(\beta\)-factor, \(e_i\) has rank \(\beta_i\).
\item Every finite free formal associative deformation of \(A\) over \(\C[[t]]\) with geometrically semisimple generic fibre has, on that geometric generic fibre and with respect to any lifts of the original projectors, the matrix factors and ranks prescribed by \(S_E\). Forgetting the projectors, the corresponding complex semisimple algebra is \(\prod_m M_m(\C)^{\times N_m}\), where \(N_m=\sum_{s_\beta=m}n_\beta\).
\item \(\dim HH_0(A)=\dim Z(A)=N_{\mathrm{GV}}\).  Both \(X_d\) and \(X_E\) are smooth at \(A\).  In each scheme, the unique irreducible component through \(A\) is the closure of the corresponding semisimple orbit.
\end{enumerate}
\end{theorem}

For one curve, \(\beta\) consists of the single coefficient \(j\), so the theorem proves Conjecture~\ref{conj:wemyss-gv}. Theorem~\ref{thm:marked-semisimple-existence} constructs the family with split generic fibre \(S_E\), and Theorem~\ref{thm:marked-general-uniqueness} proves that all geometrically semisimple generic fibres have these matrix factors and projector ranks. The theorem determines the geometric generic fibre with these projectors; the integral families can still differ.

\subsection*{Relation to earlier work}

For a single curve, Hua--Toda recover the genus-zero GV invariants from the contraction algebra by a noncommutative Donaldson--Thomas formula \cite[Theorem~4.4 and Corollary~4.6]{HuaToda}.  We ask in addition which semisimple algebras can occur as geometric generic fibres of associative deformations.  Hua--Toda's flatness conjecture \cite[Conjecture~4.3]{HuaToda} concerns the relative algebra, centre and cocenter in geometric families.  We treat these three objects over complete Noetherian local bases in Section~\ref{sec:hua-toda-flatness}.

Kawamata computes the multiplicity with which each original generator summand contributes to a generic matrix factor \cite[Theorem~3.1(3) and its proof]{KawamataNC}. These multiplicities are the projector ranks in the generic matrix factors.  Our construction through the universal quotient makes the comparisons over formal classifying arcs explicit and produces a polynomially split family with fixed projectors.  These steps are proved in Section~\ref{sec:marked-semisimple}.

Explicit algebraic approaches include the deformation of a nine-dimensional contraction algebra to \(M_2(\C)\times\C^5\) by \citet[Section~4]{BowmanPuljicSmoktunowicz}, and Smoktunowicz's study of strongly flat deformations \cite{SmoktunowiczFlat}. In their February 2026 preprint, Joachim Jelisiejew and Agata Smoktunowicz use the local normal forms of \citet{BrownWemyss} to prove semisimple existence and uniqueness for types A and D \cite[Theorem~1.3]{JS}. In their October 2026 preprint, Agata Smoktunowicz and Alicja Smoktunowicz construct flat deformations to explicit semisimple algebras for a family of type E contraction algebras. Their result concerns completed Jacobi algebras of potentials \(x^3+xy^3+p(y)\), where \(p\in\C[y]\) has order at least four, under the nilpotence and dimension hypotheses of \cite[Theorem~7.1]{SmoktunowiczSmoktunowicz}.

For a flop with a single exceptional curve, Karmazyn and Pan determine the dimension of the centre and the unique semisimple deformation type. They also prove that, in geometric families over complete Noetherian local bases, the relative contraction algebra, its centre and its cocenter are free \cite[Theorems~B and~C]{KP26}. These results overlap with ours for one curve. Their Proposition~D also gives a counterexample to smoothness of the centre morphism for the Laufer flop. We prove the flatness and uniqueness results for multi-curve flops, where uniqueness determines both the matrix factors and the rank of each labelled curve idempotent on each factor. Our proofs combine the orbit criterion of Jelisiejew and Smoktunowicz with the Lie-theoretic description of the universal contraction quotients. The latter gives the cocenter dimension needed to compare the tangent-space and orbit dimensions. Relative Hochschild cochains give the corresponding tangent calculation with every projector fixed.

This manuscript was developed independently of the preprint by Karmazyn and Pan. Preliminary versions have been publicly available on the \href{https://shutaobi.github.io/}{author's academic homepage} since October~2026, with a \href{https://github.com/Shutaobi/shutaobi.github.io/commits/main/thesis/Wemyss_Semi-Simple_deformation_on_GV_algebras.pdf}{version history on GitHub}. The author shared the manuscript with Joachim Jelisiejew and Agata Smoktunowicz on 20~October~2026 and with Shengyong Pan on 21~October~2026.

\subsection*{Flatness for several exceptional curves}

For several curves, the universal contraction quotient with the chosen curve projectors, and its cocenter, admit presentations by square matrices. A geometric family pulls these matrices back to \(B[[w]]\), where \(B\) is the deformation base and \(w\) is the surface-section parameter. We prove that the resulting cokernels are finite free over \(B\), including when this complete Noetherian local ring has nilpotents. Their ranks are the orders in \(w\) of the determinants after reduction modulo the maximal ideal of \(B\). The generic fibre need not be semisimple.

\begin{theorem}[Geometric Hua--Toda flatness]
\label{thm:intro-hua-toda-flatness}
Let \(f_B\) be a flat deformation of \(f_0\) over a complete Noetherian local complex algebra \((B,\mathfrak m,\C)\).  Form the relative contraction quotient \(\mathscr A_B\) using a lift of the fixed marked basic generator. Retain all points of its support when passing to fibres. Then the algebra, cocenter, and centre are finite free, with ranks
\[
 \operatorname{rank}_B\mathscr A_B=\sum_\beta n_\beta s_\beta^2,
 \qquad
 \operatorname{rank}_BHH_0(\mathscr A_B/B)
 =\operatorname{rank}_BZ(\mathscr A_B)=N_{\mathrm{GV}}.
\]
The centre and cocenter commute with every commutative scalar extension \(B\to B'\), and every symmetrizing trace of the central algebra lifts to a perfect relative trace.  The corresponding local statement holds over Noetherian algebraic or analytic bases when the relative quotient is finite over the base.
\end{theorem}

We obtain the algebra presentation from retained corners of the Etingof--Rains central extension \cite[Theorem~3.2]{EtingofRains}.  For the cocenter, we combine the full-algebra trace results of \citet[Theorem~3.1 and Corollary~3.5]{EtingofLatourRains} with a vertex-deletion argument.  The resulting Theorem~\ref{thm:all-corner-trace-excision} identifies the corner cocenter with the kernel of the map from the cocenter of the full central extension to that of the deleted subdiagram. The dimension of the corner cocenter is therefore the difference between these two cocenter dimensions. Pulling back the algebra and cocenter presentations gives the ranks for the relative quotient defined in Remark~\ref{rem:hua-toda-support}. Flatness and local base change are proved in Theorem~\ref{thm:multicurve-hua-toda} and Corollary~\ref{cor:hua-toda-local-bases}.

\subsection*{Consequences for centres}

The main theorem gives a smooth open neighbourhood of \(A\) in the scheme of associative algebra laws.  On this neighbourhood, centres and cocenters vary in vector bundles of constant rank.  The vector-bundle statement applies to arbitrary finite free associative families through \(A\).

We then ask whether a prescribed deformation of the commutative algebra \(Z(A)\) can be obtained as the centre of a deformation of \(A\). The tangent map of the centre morphism detects first-order obstructions to this lifting problem. For a commutative algebra \(Z\), let \(T^1_{\mathrm{com}}(Z)\) be the space of first-order commutative deformations, computed by Harrison cohomology \cite{Harrison}.

\Needspace{12\baselineskip}
\begin{theorem}[The centre deformation problem]
\label{thm:intro-centre-morphism}
For \(A\) as above, with \(d=\dim A\), the following hold.
\begin{enumerate}[label=\textup{(\roman*)}]
\item \(\widehat{\mathcal O}_{X_d,A}\cong \C[[t_1,\ldots,t_{d(d-1)-d+N_{\mathrm{GV}}}]]\); in particular \(A\) is unobstructed.
\item Near \(A\) in \(X_d\), the centre and cocenter of the tautological algebra are dual vector bundles of rank \(N_{\mathrm{GV}}\), compatible with arbitrary base change.  Every finite free associative deformation of \(A\) over a Noetherian local complex base inherits these properties and a lift of any chosen symmetrizing trace of the closed fibre.
\item The centre morphism is smooth at an algebra law \(B\) in this neighbourhood if and only if \(\chi_B:HH^2(B)\to T^1_{\mathrm{com}}(Z(B))\) is surjective. This holds on a dense open subset containing the semisimple laws.
\item For \(D_{n,2}\), \(n\ge2\), a polynomial family of centres inside the smoothable component has first-order class outside the image of \(\chi\).  Thus the centre morphism is not smooth at \(D_{n,2}\). For \(D_{2,2}\), nonlifting also holds without a chosen identification of the special centre.
\end{enumerate}
\end{theorem}

These consequences establish unobstructedness and local constancy of the centre dimension in \cite[Conjecture~1.4]{JS} for the contractions considered here.  The type \(D\) algebras \(D_{n,2}\), \(n\ge2\), give counterexamples to the smoothness assertion in part~(2).  The centre morphism remains non-smooth even after restricting its target to the smoothable component of commutative algebra laws.

\subsection*{Outline of the proof}

The orbit argument requires dimensions of the algebra and its cocenter. We compute these dimensions by restricting universal presentations to classifying arcs.  Write \(C_J\) for the universal contraction quotient with the projectors indexed by \(J\) retained, and \(R_J\) for its parameter ring.  Let \(P_\beta=0\) be the divisor along which a curve of restricted class \(\beta\) appears.  The zeroth Fitting ideals of the following modules are generated by products of these divisor equations.  The exponent of each divisor equation is as follows:
\[
\begin{array}{c|c|l}
\text{module}&\text{exponent of }P_\beta&\text{information recorded}\\ \hline
 e_iC_Je_j&\beta_i\beta_j&\text{individual Peirce dimensions}\\
 C_J&s_\beta^2&\text{dimension of the matrix algebra}\\
 HH_0(C_J/R_J)&1&\text{number of simple blocks}
\end{array}
\]
Along the classifying arc, the factor \(P_\beta\) has order \(n_\beta\). The rank of each pulled-back module is therefore the sum of these orders multiplied by the exponents in its determinant. For a single curve these exponents are \(j^2\) for the algebra and \(1\) for its cocenter.

For one curve, the comparison with the original GV invariant uses the relative local Hilbert ring. Its closed fibre has finite length, and the embedded-deformation functor has a complete obstruction space of the same dimension as its tangent space. These properties make the relative Hilbert ring finite flat. The points of its reduced geometric generic fibre correspond to the isolated \((-1,-1)\)-curves (Proposition~\ref{prop:conserved-multiplicities-gv}). For a marked configuration, the same local calculation compares formal classifying arcs with the analytic arc used in the Nabijou--Wemyss counting convention (Corollary~\ref{cor:multicurve-cartan-gram}).

The rank formulas give the dimension equality used in the uniqueness proof. The geometric construction identifies the multiplication and the projectors in the marked semisimple product. We perturb the classifying arc in the simultaneous-resolution families of Katz--Morrison and Bryan--Katz--Leung \cite{KatzMorrison,BKL}.  The determinant degrees of the transported generator summands on each exceptional curve determine their multiplicities in the corresponding matrix block, following Hua--Toda and Kawamata \cite[Theorem~4.2 and its proof]{HuaToda} \cite[Theorem~3.1(3) and its proof]{KawamataNC}.

\begin{figure}[htbp]
\centering
\begin{tikzpicture}[>={Stealth[length=2mm]},x=15mm,y=15mm,
   dot/.style={circle,fill,inner sep=1.5pt}]
\draw[thick,blue!55!black] (-2.45,0) -- (2.45,0)
      node[anchor=west,font=\footnotesize] {\(D_1\)};
\draw[thick,red!60!black] (-1.25,1.1875) -- (1.55,-1.4725);
\node[red!60!black,font=\footnotesize,anchor=south east] at (-1.2,1.16) {\(D_2\)};
\draw[very thick,domain=-2.05:2.05,samples=80,smooth]
      plot (\x,{0.5*\x*\x});
\node[font=\footnotesize,anchor=south west] at (2.0,2.02) {\(f_0\)};
\draw[very thick,densely dashed,black!65,domain=-2.05:2.05,samples=80,smooth]
      plot (\x,{0.5*\x*\x-0.72});
\node[font=\footnotesize,anchor=north west,black!65] at (2.0,1.30) {\(f_0+sv\)};
\node[dot] at (0,0) {};
\node[font=\footnotesize,anchor=south east,inner sep=2pt] at (0,0.02) {\(0\)};
\node[dot,blue!55!black] at (-1.2,0) {};
\node[dot,blue!55!black] at ( 1.2,0) {};
\node[blue!55!black,font=\footnotesize,anchor=north east,inner sep=3pt]
     at (-1.18,-0.04) {\(M_1\)};
\node[blue!55!black,font=\footnotesize,anchor=north west,inner sep=3pt]
     at ( 1.22,-0.04) {\(M_1\)};
\node[dot,red!60!black] at (0.5805,-0.5515) {};
\node[red!60!black,font=\footnotesize,anchor=west,inner sep=3pt]
     at (0.63,-0.60) {\(M_2\)};
\end{tikzpicture}
\caption{A schematic slice of the marked base with two curve divisors. The central classifying arc is tangent to \(D_1\); the dashed arc is obtained by translating the solid arc by the constant vector \(sv\), directed downwards in this slice. The translation separates the tangency into transverse intersections.  An intersection with \(D_j\) gives the local factor \(M_j\), while the conserved intersection number determines how many such factors occur in the generic product of matrix algebras.  This figure depicts the one-curve divisor-count construction of \citet[Section~2.1, Proposition~2.3, and Lemma~2.8]{BKL}.}
\label{fig:morsification}
\end{figure}
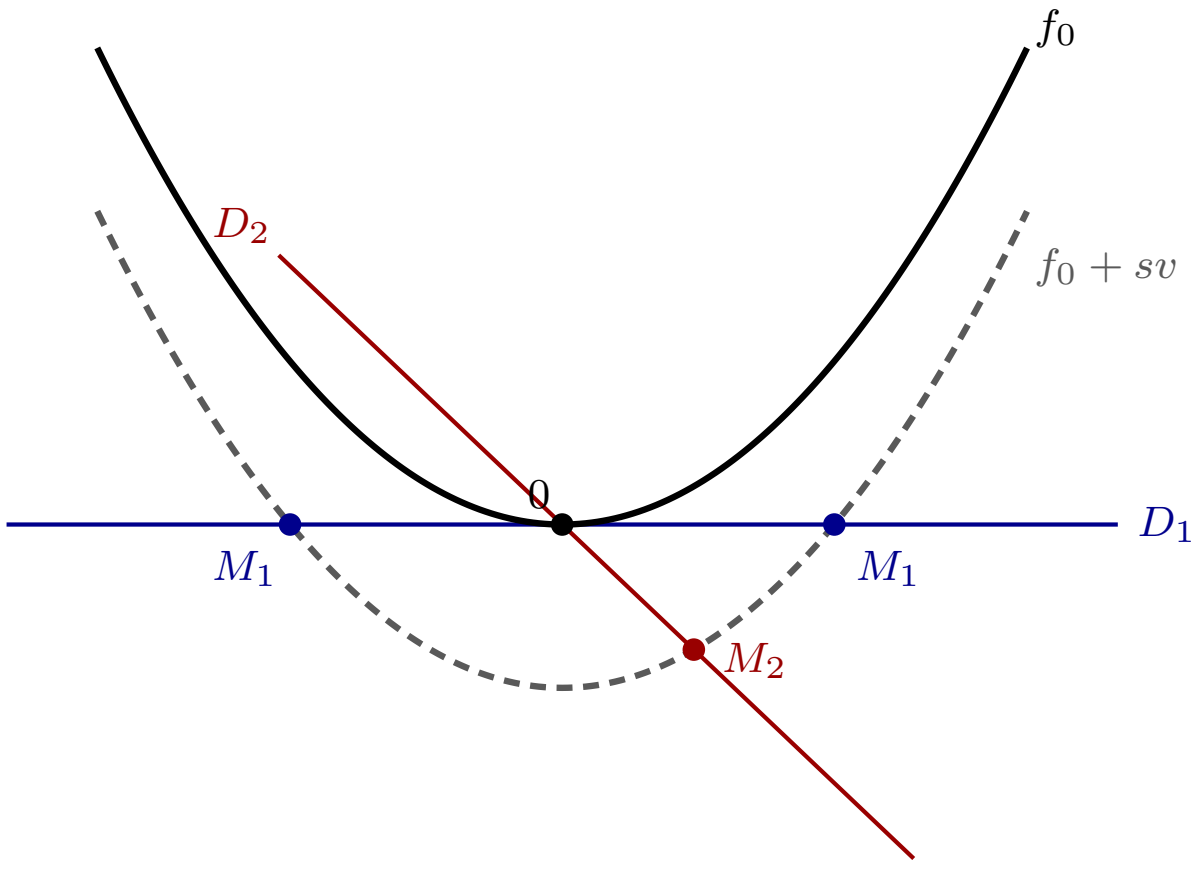

August proves that \(A\) is symmetric \cite[Theorem~2.6]{August}. The bimodule sequence of Jasso--Keller--Muro gives four-periodicity \cite[Theorem~3.2.5]{JKM}. Linckelmann's Tate duality \cite[equation~(2.11)]{Linckelmann} then gives \(\dim HH^1(A)=\dim HH^2(A)\). The cocenter dimension makes the tangent dimension equal to the dimension of the constructed semisimple orbit.  This equality forces smoothness and a unique component through \(A\), and hence determines the possible semisimple geometric generic fibres.  Relative Hochschild cochains give the corresponding tangent formula with every projector fixed, which proves the full marked statement.

For the one-curve cocenter formula, Subsection~\ref{subsec:cocenter-statements-reduction} states the universal presentation and performs the determinant--length reduction. The following three subsections prove the presentation over \(\widehat{\mathcal O}_{\mathfrak h,0}\), identify the geometric algebras, and descend the presentation to \(R_H\). The \(E_8\) relations for length five, rank certificate and tail estimate are collected in Appendix~\ref{app:e8-executable-certificate}, using the results of Etingof--Rains and Etingof--Latour--Rains \cite{EtingofRains,EtingofLatourRains}.  The exact sequence in Theorem~\ref{thm:all-corner-trace-excision} subsequently computes the cocenter dimension for arbitrary retained subsets.  The appendices record the exact root enumeration, the \(E_8\) calculation, and the \(D_{2,2}\) tangent data; the source distribution includes their verification programs.

\paragraph{Conventions and organization.}
All algebras are associative and unital.  All parameter rings and deformation bases in the geometric constructions, including local Hilbert algebras, are commutative.  The geometry is over \(\C\) unless a field extension or a split local model over another field is specified; the general orbit criterion works over an algebraically closed field of characteristic zero.  We write \(A^e=A\otimes A^{\mathrm{op}}\), \(HH^i(A)=\operatorname{Ext}_{A^e}^i(A,A)\), and \(D=\operatorname{Hom}_\C(-,\C)\).  The symbol \([A,A]\) denotes the linear span of commutators, and \(HH_0(A)=A/[A,A]\). The parameters \(w,s,z\) denote respectively the surface slice, the external deformation, and the central-extension element.  In the multi-curve construction \(\mathscr E\) denotes the central extension; \(E=\bigoplus_i\C e_i\) denotes the separable algebra of labelled idempotents.

\Needspace{12\baselineskip}
Section~\ref{sec:general-existence} constructs a family with semisimple generic fibre for one curve. Section~\ref{sec:general-uniqueness} computes the cocenter and applies the orbit criterion to prove uniqueness of the matrix factors. Section~\ref{sec:hua-toda-flatness} begins with a criterion for freeness after pullback and the lifting of symmetrizing traces.  It then constructs the universal presentations and identifies their geometric pullbacks, first for one curve and then for several curves. Section~\ref{sec:marked-semisimple} determines the semisimple geometric generic fibres with the curve projectors fixed and proves smoothness of the corresponding schemes of algebra laws. Section~\ref{sec:centre-morphism} derives the consequences for centres in associative families and studies the centre-lifting obstruction.

\section{Constructing a semisimple deformation for one curve}
\label{sec:general-existence}

We construct a semisimple deformation for a simple one-curve flop from its marked surface deformation and a fixed tilting generator.

Katz--Morrison's length classification and simultaneous partial resolution are recalled in \cite[Theorems~1.1--1.2]{Karmazyn}. Bryan--Katz--Leung describe the curve-divisor count in a one-parameter analytic deformation \cite[Section~2.1, Proposition~2.3, and Lemma~2.8]{BKL}, and Toda relates the resulting GV multiplicities to the contraction-algebra dimension \cite[Section~2.4 and Theorem~1.1]{Toda}. At each transverse intersection of the perturbed classifying arc with a curve divisor, the completed contraction quotient of the generator is a matrix algebra. The intersection multiplicity of the original arc with that divisor gives the number of such factors. Their size is the degree of the determinant of the nonfree summand on the corresponding exceptional curve.

Throughout this section the ground field is \(\C\). Let
\[
 f\colon X\longrightarrow Y=\Spec R
\]
be the completion at the contraction point of a simple flopping contraction of a smooth complex threefold.  Thus \(\operatorname{Exc}(f)=C\cong\mathbb P^1\), the scheme-theoretic length \(\ell\) belongs to \(\{1,\ldots,6\}\), and the contraction algebra is denoted by \(\Acon\).  We write \(n_j\) for the genus-zero Gopakumar--Vafa invariant in class \(j[C]\).

\Needspace{10\baselineskip}
\begin{theorem}
\label{thm:general-semisimple-existence}
Let \(f\) be any completed simple one-curve flop of a smooth complex threefold, of length \(\ell\), with contraction algebra \(\Acon\) and genus-zero GV invariants \(n_j\).  There exists a finite free \(\C[[s]]\)-algebra \(\mathscr A\) such that
\[
 \mathscr A/(s)\cong \Acon.
\]
After a finite splitting extension \(K'/\C((s))\),
\begin{equation}
  \mathscr A\otimes_{\C[[s]]}K'
 \cong\prod_{j=1}^{\ell}M_j(K')^{\times n_j}.
 \label{eq:general-semisimple-product}
\end{equation}
\end{theorem}

\begin{definition}
\label{def:marked-partial-resolution-data}
Let \((\Phi,*)\) be the primitive marked ADE pair associated to a general hyperplane section, with retained vertex \(*\).  The root system \(\Phi\) has complex Cartan space \(\mathfrak h\), Weyl group \(W\), and simple roots \(\Delta=\{\alpha_i\}\), which determine the positive roots \(\Phi^+\).  The marking fixes the retained simple root \(\alpha_*\) and the corresponding summand projector in the endomorphism algebra.  For a root \(\alpha=\sum_i c_i(\alpha)\alpha_i\), the marked coefficient is \(c_*(\alpha)\); the Levi roots are those with \(c_*(\alpha)=0\).  Put
\[
 W_0=\langle s_i:i\ne *\rangle,\qquad H=\mathfrak h/W_0,
\]
and let \(\pi_W:H\to\mathfrak h/W\) be the induced finite map to the semiuniversal ADE base.  Choose affine coordinates on \(H\) centred at the origin.

If \(\theta\) is the highest root, the length classification gives
\[
 \ell=c_*(\theta)
\]
\cite{KatzMorrison}.  The curve divisors and symbolic thickenings are indexed by \(1\le j\le\ell\) \cite[Proposition~2.3 and Lemma~2.4]{BKL}. Modulo the Levi root lattice, \([\alpha]=c_*(\alpha)[\alpha_*]\).  Lemma~\ref{lem:marked-degree-specialization} shows that the corresponding exceptional curve has class \(c_*(\alpha)[C]\), by computing line-bundle degrees on that curve.
\end{definition}

\begin{figure}[htbp]
\centering
\begin{tikzpicture}[
  every node/.style={inner sep=0pt},
  nd/.style={circle,draw,black,fill=white,minimum size=4.4pt,inner sep=0pt},
  mk/.style={circle,draw,black,fill=black,minimum size=4.4pt,inner sep=0pt},
  lb/.style={font=\footnotesize},
  x=6.2mm,y=6.2mm]

\def\rowlabel#1#2{\node[lb,anchor=east] at (-1.15,#1) {#2};}
\def\rowlevi#1#2{\node[lb,anchor=west] at (8.4,#1) {#2};}

\rowlabel{0}{\(\ell=1\)}
\node[mk] (a) at (0,0) {};
\rowlevi{0}{\(A_1\): Levi \(\varnothing\)}

\rowlabel{-1.6}{\(\ell=2\)}
\node[nd] (d1) at (0,-1.6) {};  \node[mk] (d2) at (1,-1.6) {};
\node[nd] (d3) at (2,-1.6) {};  \node[nd] (d4) at (1,-2.5) {};
\draw (d1)--(d2)--(d3); \draw (d2)--(d4);
\rowlevi{-1.6}{\(D_4\): Levi \(A_1^{3}\)}

\rowlabel{-3.9}{\(\ell=3\)}
\node[nd] (e1) at (0,-3.9) {}; \node[nd] (e2) at (1,-3.9) {};
\node[mk] (e3) at (2,-3.9) {}; \node[nd] (e4) at (3,-3.9) {};
\node[nd] (e5) at (4,-3.9) {}; \node[nd] (e6) at (2,-4.8) {};
\draw (e1)--(e2)--(e3)--(e4)--(e5); \draw (e3)--(e6);
\rowlevi{-3.9}{\(E_6\): Levi \(A_1\sqcup A_2^{2}\)}

\rowlabel{-6.2}{\(\ell=4\)}
\node[nd] (f1) at (0,-6.2) {}; \node[nd] (f2) at (1,-6.2) {};
\node[mk] (f3) at (2,-6.2) {}; \node[nd] (f4) at (3,-6.2) {};
\node[nd] (f5) at (4,-6.2) {}; \node[nd] (f6) at (5,-6.2) {};
\node[nd] (f7) at (2,-7.1) {};
\draw (f1)--(f2)--(f3)--(f4)--(f5)--(f6); \draw (f3)--(f7);
\rowlevi{-6.2}{\(E_7\): Levi \(A_2\sqcup A_1\sqcup A_3\)}

\rowlabel{-8.5}{\(\ell=5\)}
\node[nd] (g1) at (0,-8.5) {}; \node[nd] (g2) at (1,-8.5) {};
\node[nd] (g3) at (2,-8.5) {}; \node[mk] (g4) at (3,-8.5) {};
\node[nd] (g5) at (4,-8.5) {}; \node[nd] (g6) at (5,-8.5) {};
\node[nd] (g7) at (6,-8.5) {}; \node[nd] (g8) at (2,-9.4) {};
\draw (g1)--(g2)--(g3)--(g4)--(g5)--(g6)--(g7); \draw (g3)--(g8);
\rowlevi{-8.5}{\(E_8\): Levi \(A_3\sqcup A_4\)}

\rowlabel{-10.8}{\(\ell=6\)}
\node[nd] (h1) at (0,-10.8) {}; \node[nd] (h2) at (1,-10.8) {};
\node[mk] (h3) at (2,-10.8) {}; \node[nd] (h4) at (3,-10.8) {};
\node[nd] (h5) at (4,-10.8) {}; \node[nd] (h6) at (5,-10.8) {};
\node[nd] (h7) at (6,-10.8) {}; \node[nd] (h8) at (2,-11.7) {};
\draw (h1)--(h2)--(h3)--(h4)--(h5)--(h6)--(h7); \draw (h3)--(h8);
\rowlevi{-10.8}{\(E_8\): Levi \(A_4\sqcup A_1\sqcup A_2\)}
\end{tikzpicture}
\caption{The six primitive marked ADE pairs.  The filled vertex is the retained vertex \(\alpha_*\); the highest root has coefficient \(\ell=c_*(\theta)\) at this vertex, and the Levi type is the diagram left after deleting it.  These are the primitive one-curve markings in Bryan--Katz--Leung's description of the curve divisors \cite[Section~2.1 and Proposition~2.3]{BKL}.  For lengths \(2,3,4,6\), the retained vertex is the branching vertex used in the corner computations of Etingof--Latour--Rains.  The corner for the \(E_8\) marking of length five is computed in Proposition~\ref{prop:e8-exact-modular}.}
\label{fig:marked-pairs}
\end{figure}
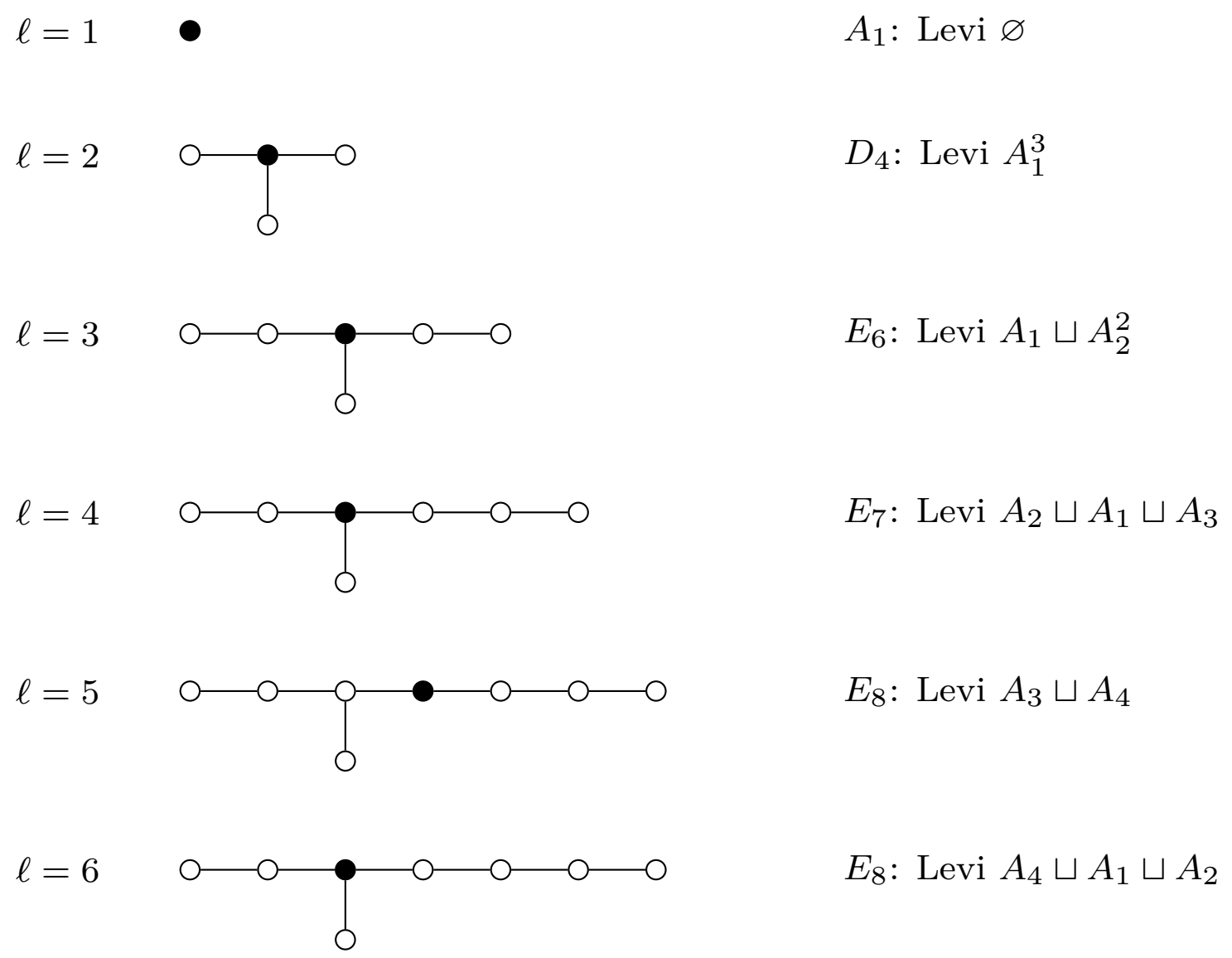

A root with retained coefficient \(j\) corresponds to the curve class \(j[C]\). The quotient \(H=\mathfrak h/W_0\) of the Cartan space by the Levi Weyl group is the deformation base of the marked partial resolution. The six primitive cases are shown in Figure~\ref{fig:marked-pairs}.

\begin{lemma}
\label{lem:marked-analytic-classifying-arc}
After replacing the source and target of the contraction by sufficiently small analytic neighbourhoods, there is a holomorphic map
\[
 f_0\colon\mathbb D_w\longrightarrow H,\qquad f_0(0)=0,
\]
for which the given analytic contraction is isomorphic, as a marked contraction, to the pullback of the simultaneous partial resolution over \(H\).  The completion of \(f_0\) at \(w=0\) is a formal classifying arc
\[
 \nu_0\colon\Spf\C[[w]]\longrightarrow H,
\]
whose pullback is isomorphic to the original completed marked contraction.
\end{lemma}

\begin{proof}
Choose a general analytic hyperplane parameter \(w\) through the contraction point.  The central fibre of the resulting hyperplane-section family is the marked partial resolution of the ADE surface section.  The semiuniversal deformation spaces of the full resolution, the marked partial resolution, and the ADE surface are, respectively,
\[
 \mathfrak h,
 \qquad H=\mathfrak h/W_0,
 \qquad \mathfrak h/W.
\]
The compatible analytic classifying maps and their curve-divisor geometry are established by Bryan--Katz--Leung \cite[Section~2, preprint pp.~6--8]{BKL}, and the local universal-flop formulation used here is stated in \cite[Theorem~1.2]{Karmazyn}.  These results identify the resolved threefold neighbourhood with the Cartesian pullback of the simultaneous marked partial resolution along the analytic classifying map \(f_0\).

The classifying map \(f_0\) completes at \(w=0\) to the formal arc \(\nu_0\colon\Spf\C[[w]]\to H\). Recovering the completed contraction also requires completion of the pulled-back target at its central point, in the surface coordinates as well as the arc coordinate. Write \(\mathbb H_\ell=\mathcal O(H)=\Gamma(H,\mathcal O_H)\) for the coordinate ring of the affine parameter space \(H=\mathfrak h/W_0\), namely its ring of regular algebraic functions. Its completion at the maximal ideal of \(0\in H\) is \(\widehat{\mathcal O}_{H,0}\).

Use the affine ADE deformation family \(\mathcal X^{\mathrm{ADE}}\to\mathfrak h/W\) constructed in \cite[Section~3, Theorem~1]{KatzMorrison}, following the classifying-map diagram of \cite[Section~2, preprint pp.~6--7]{BKL}. Pulling this family back to \(H\) yields the affine target of the chosen simultaneous partial resolution:
\[
 \mathcal X_\ell
 =\mathcal X^{\mathrm{ADE}}\times_{\mathfrak h/W}H
 =\Spec\mathcal R_\ell.
\]
This description, including the affineness of \(\mathcal X_\ell\), is recorded explicitly in \cite[Section~2.3, preprint p.~10, properties~(2) and~(5)]{Karmazyn}. The coordinate ring \(\mathcal R_\ell\) is a finitely generated \(\mathbb H_\ell\)-algebra. Choose algebra generators \(z_1,\ldots,z_N\) that vanish at the central contraction point, together with polynomials \(G_1,\ldots,G_q\) generating their defining ideal of relations. This presents the algebra as
\[
 \mathcal R_\ell\cong
 \mathbb H_\ell[z_1,\ldots,z_N]/(G_1,\ldots,G_q).
\]

The formal arc \(\nu_0\) induces a ring homomorphism \(\nu_0^*:\mathbb H_\ell\to\C[[w]]\). Forming the ordinary affine pullback of \(\mathcal X_\ell\) along \(\nu_0\) gives the following coordinate ring
\[
 \mathcal R_\ell\otimes_{\mathbb H_\ell,\nu_0^*}\C[[w]]
 \cong
 \C[[w]][z_1,\ldots,z_N]/
       (G_1(\nu_0(w),z),\ldots,G_q(\nu_0(w),z)).
\]
Here \(G_i(\nu_0(w),z)\) indicates that \(\nu_0^*\) is applied to each coefficient of \(G_i\), leaving the variables \(z_1,\ldots,z_N\) unchanged.

\Needspace{5\baselineskip}
In this ordinary pullback over \(\C[[w]]\), the central contraction point has maximal ideal \((w,z_1,\ldots,z_N)\). Completion of the target algebra at this maximal ideal gives the complete local ring
\[
 \C[[w,z_1,\ldots,z_N]]/
       (G_1(\nu_0(w),z),\ldots,G_q(\nu_0(w),z)).
\]
Completing the analytic identification above identifies this ring with the complete local ring of the original target at the contraction point. Pulling back the projective family \(\mathcal Y_\ell\to\mathcal X_\ell\) to the spectrum of this completed ring then recovers the original completed marked contraction. The comparison of the tilting generators below identifies their contraction quotients.
\end{proof}

Let \(S_0\) be the central partial-resolution surface over \(0\in H\), and write \(\mathcal I_C\subset\mathcal O_{S_0}\) for the ideal sheaf of \(C\) in that surface.

\begin{proposition}
\label{prop:general-source-package}
For the classifying arc of Lemma~\ref{lem:marked-analytic-classifying-arc}, the simultaneous partial resolution has the following properties.
\begin{enumerate}[label=\textup{(\roman*)}]
\item The base \(H\) is smooth affine space, the quotient \(\pi_W\colon H\to\mathfrak h/W\) is finite, and the universal simultaneous partial resolution
\[
 \pi_\ell\colon\mathcal Y_\ell\longrightarrow
 \mathcal X_\ell=\Spec\mathcal R_\ell
\]
is projective with source and target flat over \(H\).  The fibres of \(\pi_\ell\) have dimension at most one and
\[
 R\pi_{\ell *}\mathcal O_{\mathcal Y_\ell}
 \cong\mathcal O_{\mathcal X_\ell}.
\]
The universal flopping algebra is finite over its centre \(\mathcal R_\ell\).
\item The curve discriminant is a union of precisely \(\ell\) reduced irreducible divisors
\[
 D^{\mathrm{curv}}=D_1\cup\cdots\cup D_\ell.
\]
At a general point of \(D_j\), the partial-resolution surface contains one smooth \((-2)\)-curve in class \(j[C]\).  To specialize a chosen curve, retain a positive root \(\alpha=\sum_i c_i\alpha_i\), with \(c_*=j\), and its hyperplane \(H_\alpha\subset\mathfrak h\), before taking the quotient by \(W_0\). If \(E_i\) are the simple exceptional curves on the full central resolution, the effective root divisor \(E_\alpha=\sum_i c_iE_i\) has a projective flat deformation over this chosen \(H_\alpha\), with central fibre \(E_\alpha\). At a point where only \(\pm\alpha\) vanish, its fibre maps isomorphically to the corresponding curve in the partial resolution.  This parameterization keeps the root choice.  Independently, write \(C_j\subset S_0\) for the symbolic thickening defined by \(\mathcal I_C^{(j)}\).
\item A sufficiently small generic perturbation of \(f_0\), with no intersection with the curve divisors on the boundary of the chosen disk, meets \(D_j\) in \(k_j\) points, counted with multiplicity, where \(k_j\) is the length of the local Hilbert scheme at \(C_j\).
\item Toda's pure thickening is
\[
 \mathcal O_{C^{(j)}}=(\mathcal O_{S_0}/\mathcal I_C^j)/Q,
\]
where \(Q\) is the maximal zero-dimensional subsheaf of \(\mathcal O_{S_0}/\mathcal I_C^j\), and \(n_j\) is the length of the local Hilbert ring of the original threefold at \([C^{(j)}]\).  A generic geometric deformation has exactly \(n_j\) disjoint \((-1,-1)\)-curves in class \(j[C]\).
\item The contraction algebra satisfies
\[
 \dim_{\C}\Acon=
 \sum_{j=1}^{\ell}j^2n_j.
\]
\end{enumerate}
\end{proposition}

\begin{proof}
\emph{(i).} Use the Katz--Morrison family \(\pi_\ell:\mathcal Y_\ell\to\mathcal X_\ell\) specified in Lemma~\ref{lem:marked-analytic-classifying-arc} and recalled in \cite[Section~2.3]{Karmazyn}. Karmazyn's reflection-group description identifies \(H=\mathfrak h/W_0\) as an affine space and \(\pi_W:H\to\mathfrak h/W\) as a finite morphism. Karmazyn also records projectivity of \(\pi_\ell\) and flatness of \(\mathcal Y_\ell\) and \(\mathcal X_\ell\) over \(H\). Moreover, \cite[Section~3.3, first paragraph, and Assumption~3.1]{Karmazyn} states that \(\pi_\ell\) has fibres of dimension at most one and satisfies \(R\pi_{\ell *}\mathcal O_{\mathcal Y_\ell} \cong\mathcal O_{\mathcal X_\ell}\). Let \(\mathcal T=\mathcal O_{\mathcal Y_\ell}\oplus\mathcal T'\) be the basic tilting bundle defining the universal flopping algebra \(\mathcal A_\ell=\End_{\mathcal Y_\ell}(\mathcal T)^{\mathrm{op}}\) in \cite[Section~3.3]{Karmazyn}. Properness of \(\pi_\ell\) makes \(\pi_{\ell *}\mathcal E nd(\mathcal T)\) coherent on the affine scheme \(\mathcal X_\ell\), which proves that \(\mathcal A_\ell\) is a finite \(\mathcal R_\ell\)-module.

To identify the centre \(Z(\mathcal A_\ell)\) with \(\mathcal R_\ell\), we pass to the dual generator \(\mathcal P=\mathcal T^\vee\). By \cite[Propositions~3.2.5--3.2.7]{VanDenBergh}, \(\mathcal P\) is globally generated and contains \(\mathcal O_{\mathcal Y_\ell}\) as a direct summand. Duality identifies \(\mathcal A_\ell\) with \(\End_{\mathcal Y_\ell}(\mathcal P)\).

Any central endomorphism of \(\mathcal P\) commutes with the projector onto the summand \(\mathcal O_{\mathcal Y_\ell}\), so it preserves this summand and acts on it as multiplication by some \(r\in\Gamma(\mathcal Y_\ell,\mathcal O_{\mathcal Y_\ell})=\mathcal R_\ell\). Because the endomorphism also commutes with every homomorphism \(\mathcal O_{\mathcal Y_\ell}\to\mathcal P\), it acts as multiplication by \(r\) on all global sections of \(\mathcal P\). Global generation of \(\mathcal P\) therefore identifies this endomorphism with \(r\,\mathrm{id}_{\mathcal P}\). Conversely, multiplication by any \(r\in\mathcal R_\ell\) is central in \(\mathcal A_\ell\). Thus \(Z(\mathcal A_\ell)=\mathcal R_\ell\), which completes~(i).

\emph{(ii).} The Katz--Morrison description of the curve discriminant, in the restricted-root form of \cite[Theorem~2.6]{NabijouWemyss}, identifies \(D^{\mathrm{curv}}\) with the images in \(H\) of the non-Levi root hyperplanes. In the primitive one-curve case, the nonzero marked coefficients are \(1,\ldots,\ell\), and the roots of each coefficient form one \(W_0\)-orbit \cite[Proposition~2.3 and Lemma~2.4]{BKL}. Let \(D_j\subset H\) be the image of the root hyperplanes with marked coefficient \(j\). Lemma~\ref{lem:coefficient-orbits} verifies that these roots form a single \(W_0\)-orbit, so \(D_j\) is irreducible. Lemma~\ref{lem:root-divisors-and-bad-locus} shows that the descended root equation has generic multiplicity one, so \(D_j\) is reduced.

To identify the class of the corresponding curve, fix a positive root \(\alpha=\sum_i c_i\alpha_i\) with \(c_*=j\), and consider its root hyperplane \(H_\alpha\subset\mathfrak h\). Lemma~\ref{lem:chosen-root-divisor-family} constructs a projective flat family \(\mathcal E_\alpha\to H_\alpha\) on the full simultaneous resolution, whose central fibre is \(E_\alpha=\sum_i c_iE_i\) as a scheme. Projectivity, flatness, and this central-fibre identification are preserved by base change. At any point of \(H_\alpha\) where only \(\pm\alpha\) vanish, no Levi root vanishes. Thus the contraction from the full resolution to the marked partial resolution maps the exceptional curve isomorphically onto its image (Remark~\ref{rem:root-family-full-resolution}). Along this family, Lemma~\ref{lem:marked-degree-specialization} computes line-bundle degrees to identify the class of the image curve as \(j[C]\).

\emph{(iii) and (iv).} Lemma~\ref{lem:symbolic-pure-thickening} identifies the symbolic thickening \(C_j\) with Toda's pure thickening \(C^{(j)}\), so the two subschemes have the same local Hilbert ring. Katz compares the local Hilbert component at \(C_j\) with the stable-sheaf component containing \(\mathcal O_{C_j}\), and identifies the Hilbert length with the genus-zero GV invariant of class \(j[C]\) \cite[Lemma~3.2 and Proposition~3.3]{KatzGV}. Thus \(k_j=n_j\).

We next compare this Hilbert length with the intersection count. Let \(P_j\) be a reduced equation of \(D_j\), and put \(m_j=\operatorname{ord}_{w=0}P_j(f_0(w))\) for the intersection multiplicity of the classifying arc with \(D_j\). Lemma~\ref{lem:arc-not-in-discriminant} shows that \(0<m_j<\infty\). In the affine coordinates chosen on \(H\), fix a translation vector \(v\in\C^{\dim H}\) as in Lemma~\ref{lem:good-constant-translate-expanded}. The translated arc \(w\mapsto f_0(w)+v\) meets the curve divisors transversely at smooth points and avoids them on the boundary of the disk. For the family \(F(s,w)=f_0(w)+sv\), Proposition~\ref{prop:conserved-multiplicities-gv} identifies the number of generic roots of \(P_j(F(s,w))\) with the original Hilbert length \(k_j\). Thus \(m_j=k_j=n_j\).

Finally, for any sufficiently small transverse perturbation of \(f_0\) with no boundary intersection, the total intersection multiplicity with \(D_j\) remains \(m_j=n_j\). A generic perturbation avoids the higher-codimension root strata of Lemma~\ref{lem:root-divisors-and-bad-locus}. At each transverse intersection with \(D_j\), Lemma~\ref{lem:root-to-odp} yields a split ordinary double point whose chosen small resolution contains one \((-1,-1)\)-curve of class \(j[C]\). The curve-discriminant description in~(ii) accounts for every complete curve in the deformed family. This proves~(iii) and~(iv).

\emph{(v).} The dimension formula is provided by Toda's width theorem, which relates the contraction algebra completed at the original contraction point to the Katz invariants \(n_j\) of class \(j[C]\) \cite[Theorem~1.1 and Sections~2.2--2.4]{Toda}. Since the noncommutative width is \(\dim_\C\Acon\), the formula gives \(\dim_\C\Acon=\sum_{j=1}^{\ell}j^2n_j\). Passing to the opposite algebra, if required by the endomorphism convention, leaves this dimension unchanged.
\end{proof}

\subsection{Root divisors and a deformation into ordinary double points}
\label{subsec:full-marked-morsification}

For each positive marked coefficient, the image under the quotient map of the corresponding root hyperplanes forms a divisor on the deformation base \(H\). Bryan--Katz--Leung constructed these divisors in an analytic family \cite[Section~2.1 and Proposition~2.3]{BKL}, building on Katz--Morrison's root-hyperplane description of simultaneous ADE resolutions \cite[Theorem~1]{KatzMorrison}. We choose a constant translation of the classifying arc so that the translated arc meets the divisors transversely. We then compare its analytic intersections with the roots over the formal generic field.

Let \(q\colon\mathfrak h\to H\) be the quotient map by the Levi Weyl group \(W_0\), and for each root \(\alpha\) let \(H_\alpha\subset\mathfrak h\) be its root hyperplane. Because the reflections generating \(W_0\) preserve the marked coefficient \(c_*(\alpha)\), the group \(W_0\) permutes the non-Levi roots of any fixed marked coefficient. For any such \(W_0\)-orbit \(\mathcal O\), the square-free product of linear forms
\[
 p_{\mathcal O}=\prod_{\alpha\in\mathcal O}\alpha
\]
is \(W_0\)-invariant and therefore descends to a polynomial equation on \(H\). By Lemma~\ref{lem:coefficient-orbits}, the positive roots carrying marked coefficient \(j\) constitute a single \(W_0\)-orbit \(\mathcal O_j\). Consequently, the product \(p_{\mathcal O_j}\) descends to a polynomial \(P_j\in\mathcal O(H)\) satisfying \(q^*P_j=p_{\mathcal O_j}\). The polynomial \(P_j\) cuts out the reduced irreducible divisor \(D_j\) introduced in Proposition~\ref{prop:general-source-package}(ii). We set
\[
 D^{\mathrm{curv}}=\bigcup_{j=1}^{\ell}D_j;
\]
The curve discriminant \(D^{\mathrm{curv}}\) is the locus of parameters in \(H\) over which the partial-resolution fibres contain a complete exceptional curve. We will choose the constant translation below so that the deformed arc avoids the images in \(H\) of all intersections of distinct root hyperplanes. Consequently, at every point where the translated arc meets \(D^{\mathrm{curv}}\), every lift to \(\mathfrak h\) lies on a unique root hyperplane.

\Needspace{21\baselineskip}
\subsubsection{Root divisors and transverse local models}

\begin{lemma}
\label{lem:coefficient-orbits}
For each of the six primitive marked ADE pairs and each \(1\le j\le\ell\), the set
\[
 \mathcal O_j:=\{\alpha\in\Phi^+:c_*(\alpha)=j\}
\]
is a single \(W_0\)-orbit.  Its cardinalities are:

\begin{center}
\begin{tabular}{r|c|l}
\toprule
length \(\ell\) & marked pair &
  \(|\mathcal O_1|,\ldots,|\mathcal O_\ell|\) \\
\midrule
1 & \(A_1\)                    & \(1\) \\
2 & \(D_4\)                    & \(8,\,1\) \\
3 & \(E_6\)                    & \(18,\,9,\,2\) \\
4 & \(E_7\)                    & \(24,\,18,\,8,\,3\) \\
5 & \(E_8\), \(c_*(\theta)=5\)      & \(40,\,30,\,20,\,10,\,4\) \\
6 & \(E_8\), \(c_*(\theta)=6\)      & \(30,\,30,\,20,\,15,\,6,\,5\) \\
\bottomrule
\end{tabular}
\end{center}
In particular \(\sum_j|\mathcal O_j|\) is the number of positive roots with \(c_*(\alpha)>0\), written \(N_p\) in Subsection~\ref{subsec:universal-marked-cocenter}, and each \(D_j\) is irreducible.
\end{lemma}

\begin{proof}
Bryan--Katz--Leung associate a curve-divisor component to each positive marked coefficient \cite[Proposition~2.3]{BKL}. Because simple reflections across Levi roots leave the marked coefficient \(c_*\) invariant, the subgroup \(W_0\) acts on each set \(\mathcal O_j\); the explicit reflection calculations in Appendix~\ref{app:general-certificates} confirm that this action is transitive with the cardinalities listed in the table. Transitivity also proves irreducibility. If the descended product \(P_j\) had a nonconstant proper factor, the pullback of that factor would be a proper nonempty \(W_0\)-invariant subproduct of \(p_{\mathcal O_j}\) up to a scalar. Its factors would yield a proper nonempty \(W_0\)-stable subset of \(\mathcal O_j\), contradicting transitivity. The divisor \(D_j\) is therefore irreducible.

For reducedness, Lemma~\ref{lem:root-divisors-and-bad-locus} below proves that the quotient map \(q\colon\mathfrak h\to H\) is \'etale at a general point of each non-Levi root hyperplane. The descended equation \(P_j\) therefore defines \(D_j\) with multiplicity one generically, establishing reducedness.
\end{proof}

\begin{lemma}
\label{lem:root-divisors-and-bad-locus}
Let \(\alpha\) be a non-Levi root.  The quotient \(\mathfrak h\to H=\mathfrak h/W_0\) is \'etale at a general point of \(H_\alpha\).  Consequently the divisor descended from every \(W_0\)-orbit of non-Levi root hyperplanes is generically reduced.

There is a closed subset \(\Sigma^{\mathrm{bad}}\subset D^{\mathrm{curv}}\) of codimension at least two in \(H\), containing:
\begin{enumerate}[label=\textup{(\roman*)}]
\item the singular loci of all curve-discriminant components;
\item intersections of distinct curve-discriminant components;
\item the images of intersections \(H_\alpha\cap H_\beta\), where \(\alpha\) is non-Levi and \(\beta\ne\pm\alpha\).
\end{enumerate}
At every point of \(\mathfrak h\) above \(D_j\setminus\Sigma^{\mathrm{bad}}\), exactly one root pair \(\{\pm\alpha\}\) vanishes.
\end{lemma}

\begin{proof}
To prove the first assertion, choose a point  \[
 x\in H_\alpha\setminus
 \bigcup_{\beta\ne\pm\alpha}H_\beta.
\]
The stabilizer of \(x\) in the full Weyl group \(W\) is the order-two subgroup \(\{1,s_\alpha\}\). Because the root \(\alpha\) is non-Levi, the reflection \(s_\alpha\) does not belong to \(W_0\), so the stabilizer of \(x\) in the Levi subgroup \(W_0\) is trivial. In characteristic zero, the quotient map \(q\colon\mathfrak h\to H\) is therefore \'etale at \(x\). Consequently, the square-free product of linear forms over the \(W_0\)-orbit of \(\alpha\) descends to a polynomial on \(H\) that is reduced at the image point \(q(x)\).

To construct \(\Sigma^{\mathrm{bad}}\), recall that distinct root hyperplanes in \(\mathfrak h\) intersect in codimension two. Because the finite quotient map \(q\) preserves the complex dimension of irreducible closed subsets, the images of these intersections specified in~\textup{(iii)} have codimension at least two in \(H\). Since the curve divisors \(D_j\) are reduced and distinct, both their singular loci in~\textup{(i)} and their pairwise intersections in~\textup{(ii)} likewise have codimension at least two in the smooth variety \(H\). Defining \(\Sigma^{\mathrm{bad}}\) as the union of the closed loci~\textup{(i)--(iii)} thus yields the stated codimension bound. By property~\textup{(iii)}, every preimage in \(\mathfrak h\) of a point in \(D_j\setminus\Sigma^{\mathrm{bad}}\) lies on exactly one root hyperplane. In type \(A_1\), there are no intersections of distinct root hyperplanes, so \(\Sigma^{\mathrm{bad}}\) is empty.
\end{proof}

\begin{lemma}
\label{lem:quadratic-root-coordinate}
Let \(\bar x\in D_j\setminus\Sigma^{\mathrm{bad}}\), and choose \(x\in H_\alpha\subset\mathfrak h\) with \(q(x)=\bar x\).  There are completed local coordinates \((\zeta,z_2,\ldots,z_r)\) at \(\bar x\) and \((t,z_2,\ldots,z_r)\) at its image in \(\mathfrak h/W\) such that
\[
 \widehat{\mathcal O}_{H,\bar x}
  =\C[[\zeta,z_2,\ldots,z_r]],\qquad
 \widehat{\mathcal O}_{\mathfrak h/W,\pi_W(\bar x)}
  =\C[[t,z_2,\ldots,z_r]],
\]
and \(\pi_W^*t=\zeta^2\), while \(D_j=(\zeta=0)\).
\end{lemma}

\begin{proof}
At the chosen lift, the \(W_0\)-stabilizer is trivial, whereas the full Weyl stabilizer is the order-two subgroup generated by \(s_\alpha\). This stabilizer description holds at any point where exactly one positive root vanishes. To verify this, write a complex point of \(\mathfrak h\) as \(a+ib\) with \(a,b\) real and choose a generic real scalar \(c\). Avoiding the finitely many linear conditions in this choice ensures that the vector \(a+cb\) has the same vanishing roots and the same stabilizer.

For a real vector \(v\) in the closed dominant chamber, a reduced-word expansion of \(v-wv\) expresses the difference as a sum of positive roots whose coefficients are given by its nonnegative simple-root pairings. If \(wv=v\), all these coefficients must vanish, showing that the stabilizer is generated by the reflections across the chamber walls containing it. Here the only such root pair is \(\{\pm\alpha\}\). Moreover \(s_\alpha\notin W_0\), since the retained fundamental weight is fixed by \(W_0\) but moved by this non-Levi reflection.

Because the \(W_0\)-stabilizer of \(x\) is trivial, the quotient morphism \(q\colon\mathfrak h\to H\) is \'etale at \(x\), inducing an isomorphism of complete local rings \(\widehat{\mathcal O}_{H,\bar x}\cong\widehat{\mathcal O}_{\mathfrak h,x}\). The full Weyl quotient map \(\mathfrak h\to\mathfrak h/W\) factors through \(q\), so that \(\widehat{\mathcal O}_{\mathfrak h/W,\pi_W(\bar x)}\) is identified with the \(\langle s_\alpha\rangle\)-invariant subring of \(\widehat{\mathcal O}_{\mathfrak h,x}\). We choose local coordinates \((\zeta,z_2,\ldots,z_r)\) on \(\mathfrak h\) centered at \(x\) with \(H_\alpha=(\zeta=0)\), such that the reflection \(s_\alpha\) acts by \(\zeta\mapsto-\zeta\) on the normal coordinate while fixing \(z_2,\ldots,z_r\) along \(H_\alpha\). The invariant completed local ring is therefore \(\C[[\zeta^2,z_2,\ldots,z_r]]\). Setting \(t=\zeta^2\) identifies the completed local ring of the target as \(\C[[t,z_2,\ldots,z_r]]\) with \(\pi_W^*t=\zeta^2\), while the descended curve divisor \(D_j\) is cut out in the completed local ring at \(\bar x\) by \(\zeta=0\).
\end{proof}

\begin{lemma}
\label{lem:chosen-root-divisor-family}
Let \(\mathcal Z\to\mathfrak h\) be the standard full simultaneous resolution of an ADE surface, with the central exceptional curves \(C_i\) identified with the simple roots \(\alpha_i\). For a positive root \(\alpha=\sum_i c_i\alpha_i\), put \(H_\alpha=\{\alpha=0\}\subset\mathfrak h\). There is a closed subscheme
\[
 \mathcal E_\alpha\subset
 \mathcal Z\times_{\mathfrak h}H_\alpha
\]
which is projective and flat over \(H_\alpha\), whose central fibre is exactly the effective Cartier divisor \(\sum_i c_iC_i\), with its scheme structure. At a point where the only vanishing roots are \(\{\pm\alpha\}\), its fibre is the reduced exceptional \((-2)\)-curve. Formation of the displayed family by base change preserves projectivity, flatness, and this central fibre. In particular, for an arc \(\operatorname{Spec}R\to H_\alpha\) from a complete DVR, taking the closed point to \(0\in\mathfrak h\) and the generic point to a point where only \(\pm\alpha\) vanish, the pullback is a projective flat curve over \(R\) with these special and generic fibres.
\end{lemma}

\begin{proof}
We adopt the fixed marking and standard simultaneous resolution of \cite[Theorem~1 and Sections~4--5]{KatzMorrison}. All page numbers in this proof refer to the cited arXiv version \href{https://arxiv.org/abs/alg-geom/9202002v1}{\texttt{alg-geom/\allowbreak9202002v1}}. The compatible standard partial models are those of \cite[Section~2.4, diagram~(2.E)]{NabijouWemyss}. 

Under a fixed marking, Katz--Morrison's Theorem~1(3) ensures fibrewise effectivity. To construct a family of divisors over the chosen root hyperplane, we form the scheme-theoretic closure of the generic exceptional curve in types \(A\) and \(D\), and take the zero divisor of a section of a line bundle on a projective compactification in type \(E\). In types \(A\) and \(D\), this construction works directly on the ambient simultaneous resolution rather than using the root-family equations printed on page~24.

\smallskip
\noindent\emph{A closure criterion.}
Let \(S\) be a smooth integral complex variety, and let \(q:T\to S\) be a smooth morphism with geometrically integral, nonproper surface fibres, factored as a projective morphism \(\pi:T\to X\) followed by \(X\to S\). Let \(\sigma:S\to X\) be a section that is a closed immersion. For a geometrically integral curve \(C_\eta\subset\pi^{-1}(\sigma(\eta))\) over the generic point \(\eta\in S\), let \(E\subset T\) be the reduced scheme-theoretic closure of \(C_\eta\).

The closure \(E\) is an integral closed subscheme of codimension one in the regular scheme \(T\), and is therefore an effective Cartier divisor on \(T\). Because the closed subscheme \(\pi^{-1}(\sigma(S))\subset T\) contains \(C_\eta\), its defining ideal vanishes on the reduced scheme-theoretic closure \(E\), so that \(E\) lies scheme-theoretically inside \(\pi^{-1}(\sigma(S))\). As \(\pi\) is projective, the resulting map \(E\to S\) is projective.

No geometric fibre \(E_b\) can coincide with the full surface \(T_b\), because \(E_b\) is proper while \(T_b\) is nonproper. A local equation for \(E\) therefore restricts to a nonzero regular function on each integral surface \(T_b\), acting as a nonzerodivisor there. By the fibrewise Cartier criterion \cite[Tag~\href{https://stacks.math.columbia.edu/tag/062Y}{062Y}, part~(3)]{Stacks}, \(E\) is a relative effective Cartier divisor, flat over \(S\). Because its image in \(S\) is closed and contains the generic point \(\eta\), every fibre is nonempty. Moreover, being an effective Cartier divisor on a smooth surface, every fibre has no embedded points.

To identify the multiplicities of the central fibre \(E_0\), we compare line-bundle degrees between the generic and central fibres. For any projective flat curve family \(E\to S\) and any line bundle \(L\) on \(E\), the Euler characteristics \(\chi(E_b,L_b)\) and \(\chi(E_b,\mathcal O_{E_b})\) are locally constant on \(S\). This follows, for example, by representing their derived pushforwards locally by perfect complexes and taking alternating ranks \cite[Tags~\href{https://stacks.math.columbia.edu/tag/0DI4}{0DI4} and~\href{https://stacks.math.columbia.edu/tag/0DJT}{0DJT}]{Stacks}. The difference of these Euler characteristics measures a weighted sum of degrees. Indeed, if a fibre is the Cartier divisor \(D=\sum m_iC_i\) on a smooth surface, then
\[
 \chi(D,L|_D)-\chi(D,\mathcal O_D)
       =\sum_i m_i\deg(L|_{C_i}).                 \tag{*}
\]
Adding one component at a time gives \(0\to\mathcal O_{C_i}(-D')\to \mathcal O_{D'+C_i}\to\mathcal O_{D'}\to0\). Tensoring by \(L\), subtracting Euler characteristics, and using Riemann--Roch on \(C_i\) proves \((*)\) by induction. This argument allows repeated components.

\smallskip
\noindent\emph{Type \(A_{n-1}\).}
For the \(A_{n-1}\) family \(XY=\prod_{a=1}^n(Z+t_a)\), we use the graph-closure resolution constructed on Katz--Morrison page~20, with projective coordinates \([u_l:v_l]\) for \(1\leq l\leq n-1\). The projection \(\pi\colon T\to X\) to the affine surface family is projective. The total space \(T\) is smooth, and its fibres are the smooth minimal resolutions of the ADE surface fibres; by miracle flatness and the smooth-fibre criterion, \(q\colon T\to\mathfrak h\) is a smooth family of integral surfaces. Restricting this family to the chosen root hyperplane \(H=\{t_i=t_j\}\) with \(i<j\), we set \(X_H=X\times_{\mathfrak h}H\) and \(T_H=T\times_{\mathfrak h}H\). Because each surface fibre maps surjectively onto the two-dimensional affine surface, it is nonproper.

The polynomial section
\[
 \sigma:H\longrightarrow X_H,\qquad
 (X,Y,Z)=(0,0,-t_i)
\]
takes values in the singular locus of each surface fibre. Since \(X_H\to H\) is affine, it is separated, and its section \(\sigma\) is a closed immersion \cite[Tag~\href{https://stacks.math.columbia.edu/tag/01KT}{01KT}]{Stacks}. Over the generic point \(\eta\in H\), the value \(\sigma(\eta)\) is the unique ordinary double point of the fibre \(X_\eta\). Let \(E\) be the reduced scheme-theoretic closure of the unique reduced exceptional \(\mathbb P^1\) lying over \(\sigma(\eta)\). Applying the closure criterion shows that \(E\) is projective and flat over \(H\), and that every fibre is an effective Cartier divisor.

To compute the multiplicity of \(C_l\) in the central fibre, let \(L_l\) be the pullback of \(\mathcal O_{\mathbb P^1}(1)\) from the \(l\)-th factor of the graph. On the central exceptional configuration, the coordinate description on page~20 gives \(\deg(L_l|_{C_a})=\delta_{la}\), because only the \(a\)-th coordinate varies on \(C_a\), parametrizing that curve isomorphically.

On the generic exceptional curve above \(t_i=t_j\), the factors \(l<i\) and \(l\geq j\) are constant. The factors with \(i\leq l<j\) are related by multiplication by nonzero products of \(t_a-t_i\), so each parametrizes that \(\mathbb P^1\) isomorphically. Thus the generic \(L_l\)-degree equals one precisely for \(i\leq l<j\), and vanishes otherwise. Since the central fibre \(E_0\) is supported on the central exceptional configuration, write \(E_0=\sum m_aC_a\). Applying the degree-invariance formula \((*)\) gives \(m_a=1\) for \(i\leq a<j\) and zero otherwise. This matches exactly the positive-root divisor. For type \(A_1\), the hyperplane \(H\) is a point, and the same calculation applies with \(E=C_1\).

\smallskip
\noindent\emph{Type \(D_n\).}
We follow the two-stage simultaneous resolution on Katz--Morrison pages~21--22. The first stage blows up the ideal \((X-P,Y-S)\) with homogeneous coordinates \([u:v]\), while the second stage resolves the affine chart \(v\ne0\) via the \(A\)-type graph construction. Here the polynomials \(P,Q,S,G\) are defined exactly as in Katz--Morrison equations~(4)--(7) and~(12):
\[
 f(U)=\prod_a(U+t_a)=U P(-U^2)+Q(-U^2),\quad
 g(Z)=\prod_a(Z+t_a^2),\quad
 \gamma=\prod_a t_a,
\]
and the affine surface family is given by
\[
 X^2+Y^2Z-\frac{g(Z)-\gamma^2}{Z}+2\gamma Y=0,
\]
where the quotient is polynomial. Each blowup stage is projective. Smoothness of the total family follows from its smooth total space and smooth two-dimensional fibres via the miracle-flatness argument applied in type \(A\). The following homogeneous coordinate ratios are obtained from the ambient construction on page~22 and follow directly from the \(A\)-type graph construction:
\[
 [u_l:v_l]
  =[v^{n-l-2}(v^2Z+u^2):\prod_{a=l+1}^n(u+t_av)]
       \quad(1\leq l\leq n-2),                 \tag{D1}
\]
\[
 [u_{n-1}:v_{n-1}]
   =[v^2Z+u^2:v(u+t_nv)].                      \tag{D2}
\]
These expressions hold wherever the entries do not vanish simultaneously. On the locus \(v\ne0\), the graph construction resolves all indeterminacies. Near \(v=0\) (where \(u\ne0\)), the second entry of~(D1) and the first entry of~(D2) are invertible, so these ratios extend the graph morphisms regularly across \(v=0\) and agree on the overlap. Consequently, these ratios define global morphisms to \(\mathbb P^1\). Let \(L_0\) be the pullback of \(\mathcal O_{\mathbb P^1}(1)\) along \([u:v]\), and for each \(1\leq l\leq n-1\) let \(L_l\) be the pullback of \(\mathcal O_{\mathbb P^1}(1)\) along \([u_l:v_l]\).

For a root hyperplane \(H=\{t_j=\varepsilon t_k\}\), where \(j<k\) and \(\varepsilon\in\{1,-1\}\), put
\[
 \sigma(t)=\left(0,\ \varepsilon\prod_{a\ne j,k}t_a,\ -t_k^2\right).
                                                               \tag{D3}
\]
Substituting \(\sigma(t)\) into the affine surface equation and its partial derivatives confirms that \(\sigma\colon H\to X_H\) is a polynomial section taking values in the singular locus of the surface fibres; this identity may be verified on the dense locus \(t_k\ne0\), where \(-t_k^2\) is a double root of \(g(Z)\), and then extends to all of \(H\). Here too \(X_H\to H\) is affine, so the section is a closed immersion. Over the generic point \(\eta\in H\), the value \(\sigma(\eta)\) is the unique ordinary double point of the surface fibre \(X_\eta\). We let \(E\subset T_H\) be the reduced scheme-theoretic closure of the unique reduced exceptional curve \(E_\eta\) lying over \(\sigma(\eta)\); by the closure criterion, \(E\to H\) is projective and flat, with Cartier fibres.

The family is now projective and flat. To identify its central fibre \(E_0=\sum m_aC_a\), we compute the degrees of \(L_0,\ldots,L_{n-1}\) on the central components \(C_1,\ldots,C_n\). The explicit descriptions on Katz--Morrison page~23 show that for \(1\leq a<n\),
\[
 \deg(L_0|_{C_a})=0,\qquad
 \deg(L_l|_{C_a})=\delta_{la}\quad(1\leq l<n).
\]
The remaining central curve \(C_n\) is parametrized by the first blowup coordinate \([u:v]\). Setting \(t_a=0\) in~(D1)--(D2) and cancelling common factors on the dense locus \(uv\ne0\) shows that the remaining coordinate maps on \(C_n\) restrict to
\[
 [u_l:v_l]=[v^{n-l-2}:u^{n-l-2}]\ (l\leq n-2),\qquad
 [u_{n-1}:v_{n-1}]=[u:v].
\]
The resulting degrees on \(C_n\) are therefore
\[
 \deg(L_0|_{C_n})=1,\quad
 \deg(L_l|_{C_n})=n-l-2\ (l\leq n-2),\quad
 \deg(L_{n-1}|_{C_n})=1.                       \tag{D4}
\]
In particular, this degree matrix determines every coefficient of an exceptional divisor from its \(L_0,L_1,\ldots,L_{n-1}\) degrees.

When \(\varepsilon=1\) (so that \(\alpha=e_j-e_k\)), the first blowup has constant coordinate \([u:v]=[-t_j:1]\) along the generic exceptional curve \(E_\eta\), and the second stage reduces to the type-\(A\) graph resolution for the repeated factor \(U+t_j=U+t_k\). Thus \(\deg(L_0|_{E_\eta})=0\), while \(\deg(L_l|_{E_\eta})=1\) for \(j\leq l<k\) and \(0\) otherwise; applying \((*)\) together with (D4) yields \(E_0=C_j+\cdots+C_{k-1}\).

When \(\varepsilon=-1\), choose generic parameters on \(H\), so that \(t_k\ne0\) and no other squares \(t_a^2\) coincide. Since \(f\) has no repeated factors at the generic point \(\eta\in H\), the second-stage \(A\)-resolution is an isomorphism on the generic fibre. The first blowup resolves the ordinary double point \(\sigma(\eta)\), and its exceptional curve \(E_\eta\) maps isomorphically onto the first \(\mathbb P^1\) factor \([u:v]\). For example, \(f(t_k)=f(-t_k)=0\) implies \(P(-t_k^2)=Q(-t_k^2)=0\), and hence \(Y=S\) on (D3), so both defining linear relations of that blowup vanish on this whole \(\mathbb P^1\). Its reduced exceptional curve is therefore that \(\mathbb P^1\), with \(L_0\)-degree one.

In (D1), substitute \(Z=-t_k^2\) and cancel precisely those factors of \((u+t_jv),(u+t_kv)\) that occur in the denominator. The resulting homogeneous numerator and denominator share no common zeros. The resulting degrees are
\[
 d_l=\deg(L_l|_{E_\eta})=
 \begin{cases}
 n-l-2,&l<j,\\
 n-l-1,&j\leq l<k,\\
 n-l,&k\leq l\leq n-2,
 \end{cases}
 \qquad(1\leq l\leq n-2),                    \tag{D5}
\]
where empty index ranges are omitted. The final coordinate ratio is obtained from (D2), which gives
\[
 d_{n-1}=\begin{cases}1,&k=n,\\2,&k<n.\end{cases}        \tag{D6}
\]
There is no common zero at \(v=0\), because the denominator of (D1) is nonzero there and the numerator of (D2) is nonzero there.

Write \(E_0=\sum m_aC_a\). Equations \((*)\), (D4)--(D6) give
\[
 m_n=1,\qquad m_l=d_l-(n-l-2)\ (l\leq n-2),\qquad
 m_{n-1}=d_{n-1}-1.                           \tag{D7}
\]
In the standard \(D_n\) convention where \(\alpha_l=e_l-e_{l+1}\) for \(l<n\) and \(\alpha_n=e_{n-1}+e_n\), these are exactly the coefficients of \(e_j+e_k\). When \(k\leq n-2\), they are one from \(j\) through \(k-1\), two from \(k\) through \(n-2\), and one at both fork vertices. When \(k=n-1\), they are one from \(j\) through \(n\). When \(k=n\), they are one from \(j\) through \(n-2\) and at \(n\), while vanishing at \(n-1\). For \((j,k)=(n-1,n)\) this is just \(C_n\).

An effective Cartier divisor on a smooth surface is determined by its component multiplicities. The calculation therefore gives the ideal \(\mathcal O(-\sum c_aC_a)\) of the central fibre, including its scheme structure. The Cartier condition excludes embedded points.

At any point where only \(\pm\alpha\) vanish, the fibre is \(mC\) for the sole exceptional \((-2)\)-curve \(C\). The generic fibre is \(\mathbb P^1\), so flatness gives \(\chi(\mathcal O_{E_b})=1\). The filtration of \(\mathcal O_{mC}\) gives \(\chi(\mathcal O_{mC})=\sum_{r=0}^{m-1}(2r+1)=m^2\). Consequently \(m=1\) at every such point, not only at a general point of \(H_\alpha\).

\smallskip
\noindent\emph{Types \(E_6,E_7,E_8\).}
For \(E_n\) (\(n=6,7,8\)), we construct \(\mathcal E_\alpha\) as the zero divisor of a section of a line bundle on a projective compactification. Katz--Morrison construct a smooth projective family of rational surfaces
\[
 p:\overline{\mathcal Z}\longrightarrow\mathfrak h,
 \qquad \mathcal Z=\overline{\mathcal Z}\setminus\mathcal B,
\]
by successively blowing up the ordered sections on the smooth locus of a fixed cuspidal cubic in \(\mathbb P^2\times\mathfrak h\) \cite[Section~5, arXiv version pp.~25--26]{KatzMorrison}. Here \(p\) is smooth and projective, and \(\mathcal B_b\) is an integral anticanonical curve on each fibre. The anticanonical class is nef and has square \(9-n>0\).

The ordered sections continue to define sections after each blowup: blowing up a section on the relative cubic leaves its strict transform isomorphic to the cubic, since that section is a Cartier divisor on the cubic. Thus this construction provides global line bundles for the pullback of a line in \(\mathbb P^2\) and for the total transforms of the ordered exceptional divisors.

We use these line bundles to extend the marked class \(\sum_i c_i[C_i]\) to a line bundle \(L_\alpha\) on \(\overline{\mathcal Z}\), and restrict it to the source of \(p_\alpha:\overline{\mathcal Z}_{H_\alpha}\to H_\alpha\). The coefficients of the marked class in the blowup basis are fixed integers, so the same line bundle is defined over the entire family. By Katz--Morrison's Theorem~1(3), with the marking fixed immediately before that theorem, this class is effective on the fibre over every closed complex point of \(H_\alpha\). The ordered blowup basis identifies the same class on all fibres. By upper semicontinuity for \(L_\alpha\) along the projective flat morphism \(p_\alpha\), the locus \(\{b\in H_\alpha : h^0(L_{\alpha,b})\geq1\}\) is closed in \(H_\alpha\). Because this locus contains all closed complex points of the finite-type complex scheme \(H_\alpha\), it equals \(H_\alpha\). Since cohomology on a proper scheme commutes with ground-field extension, \(L_{\alpha,b}\) is effective on every geometric fibre. On each fibre, the class of \(L_\alpha\) has self-intersection \(-2\) and canonical pairing zero.

The cited arXiv version contains a misprint in the last row of Table~4 on page~26. That row prints the upper limit \(7\) twice, whereas the stated \(E_8\) cubic configuration requires \(8\). With distinct indices, the printed class \(q_7=3e_0-2e_{i_1}-\sum_{a=2}^{7}e_{i_a}\) has \(q_7^2=-1\) and \(K\cdot q_7=-1\), so is not a root. The corrected class and equation are
\[
 q_8=3e_0-2e_{i_1}-\sum_{a=2}^{8}e_{i_a},\qquad
 2t_{i_1}+\sum_{a=2}^{8}t_{i_a}=0.
\]
They satisfy \(q_8^2=-2\), \(K\cdot q_8=0\), and the equation is exactly \(q_8^\perp\), by the distinguished coordinates on page~11. This is an explicit correction justified by the lattice and the adjacent eight-point description, not an author-issued erratum. The proof uses Theorem~1(3), not the erroneous printed row.

There is exactly one effective divisor in this class on every geometric fibre. Indeed, if its complete linear system had a nonzero moving part \(M\), nefness of \(-K\) would give \((-K)\cdot M=0\). The intersection form in the blowup basis is \(\operatorname{diag}(1,-1,\ldots,-1)\), and \((-K)^2=9-n>0\), so its restriction to \((-K)^\perp\) is negative definite. The numerical class of a nonzero effective divisor is nonzero, since it has positive intersection with an ample divisor. Thus \(M^2<0\). Two general members of the moving part have no common irreducible component and therefore have nonnegative intersection, a contradiction. This rules out a moving part. It follows that \(h^0(L_{\alpha,b})=1\). Serre duality gives \(h^2(L_{\alpha,b})=0\), because an effective divisor of class \(K-L_{\alpha,b}\) would have negative intersection \(-(9-n)\) with the nef class \(-K\). Riemann--Roch on this rational surface gives
\[
 \chi(L_{\alpha,b})
 =1+\tfrac12\bigl((-2)-0\bigr)=0,
 \qquad h^1(L_{\alpha,b})=1.
\]

The fibrewise divisors can now be assembled into a family. Cohomology and base change show that \(M_\alpha=p_{\alpha *}L_\alpha\) is a line bundle whose formation commutes with arbitrary base change on \(H_\alpha\). The flat projective morphism and the line bundle give a perfect complex \(Rp_{\alpha *}L_\alpha\), compatible with arbitrary derived base change \cite[Tags~\href{https://stacks.math.columbia.edu/tag/0DI4}{0DI4} and~\href{https://stacks.math.columbia.edu/tag/0DJT}{0DJT}]{Stacks}. Locally represent this complex by a bounded complex of finite free modules. The constant fibre cohomology dimensions just computed force constant ranks of its differentials. Since \(H_\alpha\) is reduced, the larger minors vanish, and the constant-rank matrices split locally. The cohomology modules of this complex are therefore locally free, and their formation commutes with arbitrary base change. In degree zero the rank is one.

This pushforward line bundle supplies the section defining the family: the evaluation map \(p_\alpha^*M_\alpha\to L_\alpha\) defines a section of \(L_\alpha\otimes p_\alpha^*M_\alpha^{-1}\) that is nonzero on every geometric fibre. Let \(\mathcal E_\alpha\) be its zero scheme. Each fibre is an effective Cartier divisor on a smooth integral surface. The fibrewise Cartier criterion ensures that \(\mathcal E_\alpha\) is a relative effective Cartier divisor, hence flat over \(H_\alpha\) \cite[Tag~\href{https://stacks.math.columbia.edu/tag/062Y}{062Y}]{Stacks}. As a closed subscheme of \(\overline{\mathcal Z}_{H_\alpha}\), it is projective.

It remains to verify that this divisor on the compactification is supported entirely within the open surface family \(\mathcal Z_{H_\alpha}\). Indeed, on each geometric fibre, every irreducible component of \((\mathcal E_\alpha)_b\) has intersection zero with the nef class \(-K\). Because the boundary curve \(\mathcal B_b\) has anticanonical class \(-K\) and self-intersection \((-K)^2=9-n>0\), it cannot occur as a component of \((\mathcal E_\alpha)_b\). Moreover, \(\mathcal B_b\cdot(\mathcal E_\alpha)_b=0\). By positivity of local intersection multiplicities for distinct curves, \(\mathcal B_b\) and \((\mathcal E_\alpha)_b\) have disjoint supports. Consequently, the divisor lies in \(\mathcal Z_{H_\alpha}\) while remaining the same closed divisor on the projective compactification, so it remains projective over \(H_\alpha\).

On the central fibre over \(0\in H_\alpha\), the effective Cartier divisor \(\sum_i c_iC_i\) already belongs to the linear system \(|L_{\alpha,0}|\), so \(h^0(L_{\alpha,0})=1\) forces the scheme-theoretic equality \((\mathcal E_\alpha)_0=\sum_i c_iC_i\). Its defining ideal sheaf is \(\mathcal O(-\sum_i c_iC_i)\), having no embedded points. At any point where only \(\pm\alpha\) vanish, the support of the fibre is the sole exceptional curve. The square \(-2\) of the root class forces its multiplicity to be one.

Finally, projectivity and flatness are preserved under arbitrary base change, as is the property of being a relative effective Cartier divisor \cite[Tag~\href{https://stacks.math.columbia.edu/tag/056Q}{056Q}]{Stacks}. The assertion regarding complete DVRs follows by algebraic base change: given a formal arc into the completion of \(H_\alpha\) at the origin, its coordinates form elements of \(R=\varprojlim R/\mathfrak m^a\) satisfying the linear equation of \(H_\alpha\). Because \(H_\alpha\) is affine, these elements define an actual scheme morphism \(\operatorname{Spec}R\to H_\alpha\). Pulling back \(\mathcal E_\alpha\) along this morphism produces a projective flat divisor over \(R\), whose formal completion coincides with the compatible system of its reductions modulo \(\mathfrak m^a\).
\end{proof}

\begin{remark}
\label{rem:root-family-full-resolution}
A line bundle on the partial family is pulled back to the full family and then restricted to \(\mathcal E_\alpha\). Where exactly one positive non-Levi root vanishes, the morphism from the full resolution to the partial resolution is an isomorphism on the unique exceptional curve, while on the central fibre the components indexed by the deleted vertices are contracted. This is the compatible standard resolution diagram of \cite[Section~2.4 and proof of Theorem~2.6]{NabijouWemyss}. Consequently no preservation of flatness under contraction, and no identification with a symbolic thickening on the partial surface, is needed here.
\end{remark}

\begin{lemma}
\label{lem:marked-degree-specialization}
Let \(\mathscr L\) be a line bundle on the marked partial-resolution family, and let \(E_b\) be the exceptional curve in the partial-resolution fibre over a point \(b\in D_j\setminus\Sigma^{\mathrm{bad}}\).  Then
\[
 \deg_{E_b}\mathscr L_b=j\deg_C\mathscr L_0.
\]
\end{lemma}

\begin{proof}
Choose a point \(x\in H_\alpha\subset\mathfrak h\) mapping to \(b=q(x)\in D_j\setminus\Sigma^{\mathrm{bad}}\), where \(\alpha=\sum_i c_i\alpha_i\) is a root with marked coefficient \(c_*=j\). To compute the degree on the full resolution, pull the line bundle \(\mathscr L\) back from the marked partial resolution to the full simultaneous resolution over \(\mathfrak h\). By Lemma~\ref{lem:chosen-root-divisor-family}, there exists a projective flat family of effective Cartier divisors \(\mathcal E_\alpha\to H_\alpha\) on the full resolution whose central fibre is scheme-theoretically \(C_\alpha=\sum_i c_iC_i\). Restricting this family along the radial line \(\tau\mapsto\tau x\) for \(\tau\in\mathbb C\) yields a flat family of projective curves whose central fibre at \(\tau=0\) is \(C_\alpha\), while its fibre at \(\tau=1\) is the single smooth \((-2)\)-curve associated with the root \(\alpha\).

Because the degree of a line bundle on a flat family of projective curves equals the difference of Euler characteristics \(\chi(L)-\chi(\mathcal O)\), it is constant along the radial line. The central degree is given by
\[
 \deg_{C_\alpha}\mathscr L
 =\sum_i c_i\deg_{C_i}\mathscr L.
\]
We now compare this transported degree with the partial resolution. For each deleted vertex \(i\ne *\), the simple exceptional curve \(C_i\) is contracted to a point by the partial resolution, so the pullback of \(\mathscr L\) has degree zero on \(C_i\). The retained component \(C_*\) maps isomorphically to \(C\), so that with \(c_*=j\), the central degree equals \(j\deg_C\mathscr L_0\). At the parameter point \(x\) (where \(\tau=1\)), no Levi root vanishes, which implies that the contraction morphism from the full resolution to the partial resolution maps the unique exceptional curve at \(\tau=1\) isomorphically onto \(E_b\). The degree at \(\tau=1\) therefore coincides with \(\deg_{E_b}\mathscr L_b\), establishing the equality \(\deg_{E_b}\mathscr L_b=j\deg_C\mathscr L_0\). For a point of \(D_j\setminus\Sigma^{\mathrm{bad}}\) over a field on which the root curve is split, base change the algebraic family of divisors and its radial line to that field. The same Euler-characteristic calculation gives the identical degree formula.
\end{proof}

\begin{lemma}
\label{lem:root-to-odp}
At a point of \(D_j\setminus\Sigma^{\mathrm{bad}}\), the transverse surface deformation is
\[
 xy-z^2=t.
\]
After pullback to the marked partial-resolution base it is
\[
 xy-z^2=\zeta^2,
\]
and hence, after \(\xi_+=z+i\zeta\), \(\xi_-=z-i\zeta\), it is the ordinary double point
\[
 xy-\xi_+\xi_-=0.
\]
If a one-dimensional classifying curve meets \(D_j\) transversely, the pullback of the simultaneous partial resolution to that curve contains one \((-1,-1)\)-curve in numerical class \(j[C]\).
\end{lemma}

\begin{proof}
At a point of \(D_j\setminus\Sigma^{\mathrm{bad}}\), the vanishing root subsystem is \(\{\pm\alpha\}\cong A_1\). Katz--Morrison's explicit simultaneous resolution identifies this fibre configuration with a single \(A_1\) singularity \cite[Theorem~1 and Sections~4--5]{KatzMorrison}. In their hypersurface deformation, there is an independent constant-term coefficient \cite[Tables~2--3]{KatzMorrison}, so differentiating with respect to that coefficient proves that the total space is smooth. This local description holds at every point of \(D_j\setminus\Sigma^{\mathrm{bad}}\).

Because the fibre Hessian is invertible at the node, the critical point extends to a formal critical section over the base. Translating by this section and applying formal Morse elimination in the fibre variables reduces the local equation to \(xy-z^2=g(t,z_2,\ldots,z_r)\). Smoothness of the total space ensures that \(dg\ne0\) at the point. Since the discriminant locus is supported on \(t=0\) by Lemma~\ref{lem:quadratic-root-coordinate}, we obtain \(g=tu\), where \(dg\ne0\) forces \(u\) to be a unit. We scale \(t\) by \(tu\) and \(\zeta\) by \(\zeta\sqrt{u(\zeta^2,z_2,\ldots,z_r)}\). This coordinate change preserves the square relation \(t=\zeta^2\) and normalizes the family to
\[
 xy-z^2=t,\qquad t=\zeta^2.
\]

At the chosen lift \(x\in H_\alpha\), no Levi root vanishes, so the contraction from the full resolution to the marked partial resolution changes no curve in this neighbourhood. Formal versality of the resolved \(A_1\) surface, whose resolved parameter \(\sigma\) squares to the smoothing parameter \(\sigma^2\) \cite[Section~2.3]{Karmazyn}, identifies this local simultaneous resolution with one of the two standard \(A_1\) resolutions, normalized by \(\sigma=\pm\zeta\). Using \((z+i\zeta)(z-i\zeta)=z^2+\zeta^2\), the pulled-back equation becomes \(xy-\xi_+\xi_-=0\).

On a classifying curve meeting the curve divisor transversely, \(\zeta\) restricts to a local uniformizing parameter. Using \(\zeta\) as the curve coordinate absorbs higher-order terms in the parametrization, recovering the completed ordinary double point (ODP) equation in the threefold total space. The simultaneous partial resolution provides one of the two small resolutions of this ODP, so its exceptional curve has normal bundle \(\mathcal O_{\mathbb P^1}(-1)^{\oplus2}\). Finally, the degree comparison of Lemma~\ref{lem:marked-degree-specialization} identifies the curve's numerical class as \(j[C]\).
\end{proof}

\subsubsection{Perturbing the classifying arc}

We choose one deformation of the classifying arc that meets all curve divisors transversely.  Bryan--Katz--Leung use an analytic generic perturbation for which the curve intersections remain proper over the parameter disk \cite[Section~2.1 and Lemma~2.8]{BKL}.  For our constant translation of the classifying arc, the divisor equations remain nonzero on the boundary of a fixed disk. The argument principle then preserves the total intersection multiplicity with each divisor.

\begin{lemma}
\label{lem:arc-not-in-discriminant}
The image of \(f_0\) is contained in no \(D_j\).  After shrinking to a closed disk \(\overline{\mathbb D}_\rho\),
\[
 f_0^{-1}(D^{\mathrm{curv}})\cap\overline{\mathbb D}_\rho=\{0\},
 \qquad f_0(\partial\mathbb D_\rho)\cap D^{\mathrm{curv}}=\varnothing.
\]
If \(P_j\) is a reduced equation of \(D_j\), then
\[
 0<m_j:=\operatorname{ord}_{w=0}P_j(f_0(w))<\infty.
\]
\end{lemma}

\begin{proof}
If the image of the classifying arc \(f_0\) were contained in \(D_j\), then for every nearby nonzero parameter \(w\), the corresponding partial-resolution surface would contain a complete exceptional curve, by the full curve-discriminant description \cite[Theorem~2.6]{NabijouWemyss}. Pullback along \(f_0\) identifies these curves with exceptional curves of the original threefold contraction at nonzero values of \(w\). The original contraction has only the exceptional curve \(C\), over \(w=0\), which is a contradiction.

This noncontainment ensures that every reduced equation \(P_j\) pulls back to a nonzero holomorphic function whose zero at the origin is isolated. Because there are only finitely many indices \(j\), a single shrinking of the disk yields the boundary nonvanishing statement and establishes the finite, positive vanishing orders.
\end{proof}

\begin{lemma}
\label{lem:symbolic-pure-thickening}
\emergencystretch=1em\relax
The Bryan--Katz--Leung symbolic thickening \(C_j\) equals Toda's pure one-dimensional thickening \(C^{(j)}\).  Consequently their Hilbert multiplicities satisfy \(k_j=n_j\).
\end{lemma}

\begin{proof}
We work in a stalk \(\mathcal S\) of the normal surface \(S_0\) at a point of \(C\), letting \(\mathfrak p\subset\mathcal S\) denote the height-one prime ideal defining \(C\). By definition, the \(j\)-th symbolic power is
\[
 \mathfrak p^{(j)}=\mathfrak p^j\mathcal S_{\mathfrak p}\cap\mathcal S.
\]
The quotient module \(\mathfrak p^{(j)}/\mathfrak p^j\) is supported only at closed points and therefore has finite length over \(\mathcal S\). Because the symbolic power \(\mathfrak p^{(j)}\) is \(\mathfrak p\)-primary, the quotient ring \(\mathcal S/\mathfrak p^{(j)}\) has \(\mathfrak p\) as its unique associated prime, so it has no nonzero finite-length submodule. In the short exact sequence
\[
 0\longrightarrow\mathfrak p^{(j)}/\mathfrak p^j
  \longrightarrow\mathcal S/\mathfrak p^j
  \longrightarrow\mathcal S/\mathfrak p^{(j)}\longrightarrow0,
\]
the submodule \(\mathfrak p^{(j)}/\mathfrak p^j\) is therefore the unique maximal finite-length submodule of \(\mathcal S/\mathfrak p^j\), which identifies it with Toda's maximal zero-dimensional subsheaf \(Q\). Setting \(I=\mathcal I_C\), these local identifications give the sheaf equality \((\mathcal O/I^j)/Q=\mathcal O/I^{(j)}\), so the pure thickening and symbolic thickening coincide, along with their local Hilbert rings. Finally, Katz's direct identification of this isolated Hilbert component with the stable-sheaf GV component gives \(k_j=n_j\) \cite[Lemma~3.2 and Proposition~3.3]{KatzGV}.
\end{proof}

\begin{lemma}
\label{lem:parametric-boundary-transversality}
Let \(M\) be a compact smooth manifold with boundary, let \(\mathcal U\) and \(N\) be smooth manifolds, and let \(S_1,\ldots,S_m\subset N\) be smooth locally closed submanifolds.  If a smooth map \(\Xi\colon M\times\mathcal U\to N\) and its boundary restriction are transverse to every \(S_i\), then the parameters for which both slice maps are simultaneously transverse to all \(S_i\) form a residual dense subset of \(\mathcal U\).  If \(\dim_{\mathbb R}M=2\), such a slice misses every stratum of real codimension at least four, and its boundary misses every stratum of real codimension at least two.
\end{lemma}

\begin{proof}
For a fixed stratum \(S\subset N\), consider its preimage \(Z=\Xi^{-1}(S)\subset M\times\mathcal U\). Because \(\Xi\) is transverse to \(S\), the incidence set \(Z\) is a smooth submanifold with boundary \(Z_\partial=\Xi_\partial^{-1}(S)\), where \(\Xi_\partial\) denotes the restriction of \(\Xi\) to \(\partial M\times\mathcal U\). A parameter point \(u\in\mathcal U\) is a regular value of the projection \(Z\to\mathcal U\) if and only if the corresponding slice map \(\Xi_u=\Xi(-,u)\colon M\to N\) is transverse to \(S\).

To see this, express the normal differential along the incidence manifold as \(A\,dx+B\,du\). Total transversality means that the combined images of \(A\) and \(B\) span the full normal space to the stratum, while surjectivity of the projection differential for \(Z\to\mathcal U\) is equivalent to the inclusion \(\operatorname{im}B\subseteq\operatorname{im}A\). Combining these conditions is equivalent to surjectivity of \(A\) alone, which is precisely the condition for slice transversality. The same equivalence holds for the boundary restriction.

When an incidence manifold has dimension strictly smaller than the parameter space, its image has measure zero by the dimension estimate for smooth maps. In all remaining dimensions, Sard's theorem \cite[Theorems~4.1 and 7.2]{Sard} guarantees that the set of critical values has measure zero in local parameter charts. Exhaust the closed critical loci in the incidence manifolds by countably many compact sets. Their images are closed sets of measure zero, hence nowhere dense. Removing these critical values leaves a residual dense set of regular values in the parameter manifold.

For boundary charts, apply Sard's theorem to a local smooth extension across the boundary, and apply it separately to the boundary incidence projection. Intersecting the resulting residual dense sets for the finitely many strata preserves density and residuality. Finally, for a transverse slice of real dimension two, the preimages of strata of real codimension at least four are empty by dimension counting, while its one-dimensional boundary likewise misses all strata of real codimension at least two.
\end{proof}

\begin{lemma}
\label{lem:good-constant-translate-expanded}
There is an arbitrarily small vector \(v\in H\cong\C^{\dim H}\) such that \(f_v(w)=f_0(w)+v\) is transverse to every smooth codimension-one curve-discriminant stratum, avoids \(\Sigma^{\mathrm{bad}}\), and has no curve-discriminant point on \(\partial\mathbb D_\rho\).  The vector can be chosen so that, for some \(\epsilon>0\),
\[
 f_0(w)+sv\notin D^{\mathrm{curv}}
 \quad\text{for }w\in\partial\mathbb D_\rho,
 \quad |s|\le1+\epsilon.
\]
\end{lemma}

\begin{proof}
Decompose the curve discriminant into finitely many smooth locally closed algebraic strata, each either contained in or disjoint from \(\Sigma^{\mathrm{bad}}\). We construct such a stratification by taking the smooth loci of the irreducible components away from their mutual intersections and repeating the procedure on the remaining closed subsets. Because the omitted closed sets have strictly smaller dimension at each stage, this stratification process terminates after finitely many steps. The strata lying outside \(\Sigma^{\mathrm{bad}}\) have complex codimension one, whereas all strata within \(\Sigma^{\mathrm{bad}}\) have complex codimension at least two. The translation evaluation map
\[
 \overline{\mathbb D}_\rho\times H\longrightarrow H,
 \qquad (w,\xi)\longmapsto f_0(w)+\xi
\]
and its boundary restriction are smooth submersions. Lemma~\ref{lem:parametric-boundary-transversality} therefore yields a residual dense set of translation parameters for which the disk avoids the strata in \(\Sigma^{\mathrm{bad}}\) (which have real codimension at least four), while its boundary avoids every curve-discriminant stratum.

The image \(f_0(\partial\mathbb D_\rho)\) is disjoint from the closed curve discriminant, so the absolute values of the pulled-back reduced local equations achieve a strictly positive minimum on \(\partial\mathbb D_\rho\). Uniform continuity on the compact boundary supplies a ball of translation vectors that preserve this boundary nonvanishing on the boundary of the disk. Choose a transverse parameter in a strictly smaller ball. Scalar multiples of this vector for parameters in a slightly enlarged disk remain in the original ball, giving the required boundary condition.
\end{proof}

\subsubsection{Root counts and the formal generic fibre}

Because the translation preserves boundary nonvanishing, the total intersection multiplicity with each curve divisor \(D_j\) remains constant. Weierstrass polynomials describe these intersections. Their roots over the formal generic field locate the points where the local contraction quotients are matrix algebras.

\begin{lemma}
\label{lem:weierstrass-root-polynomials}
Fix \(v\) as in Lemma~\ref{lem:good-constant-translate-expanded} and put \(F(s,w)=f_0(w)+sv\) and \(G_j(s,w)=P_j(F(s,w))\).  On a connected \(s\)-disk \(U\) containing \(0\) and \(1\), the zero scheme of \(G_j(s,\cdot)\) is finite flat of degree \(m_j\) and is represented by a monic polynomial
\[
 W_j(s,T)\in\mathcal O(U)[T],\qquad
 \deg_TW_j=m_j,\qquad W_j(0,T)=T^{m_j}.
\]
\end{lemma}

\begin{proof}
Because \(G_j(s,w)\ne0\) on the product boundary \(U\times\partial\mathbb D_\rho\), the argument-principle integral
\[
 \frac{1}{2\pi i}\int_{\partial\mathbb D_\rho}
 \frac{\partial_wG_j(s,w)}{G_j(s,w)}\,dw
\]
is a continuous, integer-valued function of \(s\in U\), and is therefore constant and equal to \(m_j\), its value at \(s=0\). For each integer \(r\ge1\), inserting the factor \(w^r\) into the integrand yields a weighted integral holomorphic in \(s\) that evaluates to the \(r\)-th power sum of the zeros of \(G_j(s,\cdot)\) in \(\mathbb D_\rho\), counted with multiplicity. Via Newton's identities, these power sums uniquely determine holomorphic elementary symmetric functions on \(U\). We define \(W_j(s,T)\in\mathcal O(U)[T]\) to be the monic polynomial in \(T\) of degree \(m_j\) with the coefficients so obtained.

For each parameter \(s\), the roots of \(W_j\), counted with multiplicity, coincide with the zeros of \(G_j(s,\cdot)\) in the \(w\)-disk. To identify the zero scheme throughout the family, use the following division formula. For a holomorphic function \(h(s,w)\), its remainder modulo \(W_j\) is given by
\[
 R_h(s,w)=\frac{1}{2\pi i}\int_{|z|=\rho}
 \frac{h(s,z)}{W_j(s,z)}
 \frac{W_j(s,z)-W_j(s,w)}{z-w}\,dz.
\]
The remainder \(R_h\) is a polynomial in \(w\) of degree strictly less than \(m_j\), and Cauchy's formula yields the decomposition \(h=W_jQ_h+R_h\) with \(Q_h\) holomorphic. Specializing to \(h=G_j\), the remainder vanishes at every root with its multiplicity, forcing \(R_{G_j}=0\). The quotient \(Q_{G_j}\) is nowhere vanishing, hence a unit on \(U\times\mathbb D_\rho\), because any additional zero would alter the total root count.

Because \(G_j\) differs from \(W_j\) by a unit factor, the two functions cut out the same analytic zero scheme. The pushforward of its structure sheaf to the parameter disk is finite free with module basis \(1,w,\ldots,w^{m_j-1}\). Boundary nonvanishing guarantees that this support is proper over the parameter disk. At the central parameter \(s=0\), all \(m_j\) zeros coalesce at the origin, yielding \(W_j(0,T)=T^{m_j}\).
\end{proof}

\begin{lemma}
\label{lem:good-to-formal-generic}
Put \(B=\C[[s,w]]\) and \(K=\C((s))\). The polynomials \(W_j\in K[T]\) are separable and pairwise coprime. The classifying map on \(\Spec B[1/s]\) avoids \(\Sigma^{\mathrm{bad}}\). At each geometric root of \(W_j\), this map meets the smooth part of \(D_j\) transversely and meets no other curve divisor.
\end{lemma}

\begin{proof}
At \(s=1\), the translated classifying arc intersects each curve divisor \(D_j\) transversely and avoids all intersections between distinct divisors. The polynomials \(W_j(1,T)\) therefore possess only simple zeros and share no common roots, which implies that their discriminants and pairwise resultants are nonzero at \(s=1\). As nonzero holomorphic functions of \(s\in U\), their Taylor series at \(s=0\) are not identically zero; while these series may be divisible by positive powers of \(s\) and thus not units in the power series ring, they become units upon passing to the formal generic field \(K=\C((s))\). This proves that the polynomials \(W_j\in K[T]\) are separable and pairwise coprime.

To exclude \(\Sigma^{\mathrm{bad}}\), put \(W=\prod_jW_j\), a polynomial of degree \(M\). Its zero scheme has coordinate algebra \(\mathcal O(U)[w]/(W)\), which is a free module of rank \(M\) over \(\mathcal O(U)\). Choose finitely many polynomial generators for the defining ideal of \(\Sigma^{\mathrm{bad}}\). Their compositions with \(F\) admit polynomial remainders modulo \(W\) via the contour division established above.

At \(s=1\), these remainders generate the unit ideal in \(\C[w]/(W(1,w))\), because the roots are simple and none maps into the bad locus. In the module basis \(1,w,\ldots,w^{M-1}\), form the matrix whose columns are these remainders together with their products by \(w,\ldots,w^{M-1}\). The image of this matrix is precisely the ideal they generate. Because that ideal is the unit ideal at \(s=1\), the matrix admits a maximal minor \(a(s)\) that is nonzero there. The adjugate of the corresponding square submatrix shows that multiplication by \(a(s)\) annihilates the quotient algebra modulo the ideal generated by these remainders.

Because \(a(s)\) is not identically zero, its Taylor series at the origin is nonzero, and hence \(a(s)\) becomes a unit in \(K\). The quotient module therefore vanishes after inverting \(s\). Since \(\Sigma^{\mathrm{bad}}\subset D^{\mathrm{curv}}\), every point mapping into \(\Sigma^{\mathrm{bad}}\) must lie in the zero scheme of \(W\). Consequently, no root of the formal generic fibre maps into the bad locus. Finally, at a smooth point of \(D_j\), separability is equivalent to the condition \(\partial_wG_j\ne0\), which is precisely transversality of the slice map.
\end{proof}

The following lemma will give the rank of the complete local Hilbert algebra. Under its hypotheses, a bound on the dimension of a complete obstruction space ensures that the algebra is finite free over \(\C[[s]]\), with rank equal to its closed-fibre length.

\begin{lemma}
\label{lem:complete-obstruction-flatness}
Let \(R=\C[[s]]\), and let \(\mathscr H\) be a complete Noetherian commutative local \(R\)-algebra with residue field \(\C\).  Suppose
\[
 k=\dim_\C\mathscr H/(s)<\infty,\qquad k>0,\qquad
 a=\dim_\C\frac{\mathfrak m_{\mathscr H}}
 {\mathfrak m_{\mathscr H}^{\,2}+s\mathscr H}.
\]
Assume that the functor of local continuous \(R\)-homomorphisms \(\mathscr H\to A\), on commutative local Artinian \(R\)-algebras with residue field \(\C\), has a finite-dimensional complete obstruction space \(V\) with \(\dim_\C V\le a\).  Explicitly, for each small extension \(A'\twoheadrightarrow A\), whose kernel \(N\) is annihilated by the maximal ideal of \(A'\), the obstruction lies in \(V\otimes_\C N\), commutes with linear pushouts of \(N\), and vanishes exactly when the given object lifts. Then \(\mathscr H\) is finite free of rank \(k\) over \(R\). It is a one-dimensional complete intersection.
\end{lemma}

\begin{proof}
Lift a basis of the relative cotangent space to elements \(x_1,\ldots,x_a\in\mathfrak m_{\mathscr H}\). Together with \(s\), these elements generate the maximal ideal by Nakayama's lemma. Successive approximation in the complete local ring yields a minimal relative presentation
\[
 B=R[[x_1,\ldots,x_a]]\twoheadrightarrow\mathscr H,
 \qquad J=\ker(B\to\mathscr H),\quad
 \mathfrak m=(s,x_1,\ldots,x_a).
\]
Relative minimality means that \(J\bmod s\subset(x_1,\ldots,x_a)^2\), so every partial derivative \(\partial f/\partial x_i\) for \(f\in J\) lies in the maximal ideal \(\mathfrak m\). The Artin--Rees lemma supplies an integer \(q\) with \(J\cap\mathfrak m^q\subset\mathfrak mJ\) \cite[Lemma~10.51.2, Tag~00IN]{Stacks}. Now consider the small extension
\[
 B/(\mathfrak mJ+\mathfrak m^q)
 \longrightarrow B/(J+\mathfrak m^q),
 \qquad\text{kernel }J/\mathfrak mJ.
\]
Let \(\mathrm{ob}\in V\otimes_\C(J/\mathfrak mJ)\) denote the obstruction to lifting the canonical quotient map \(\mathscr H\to B/(J+\mathfrak m^q)\) across this extension. We claim that the coefficients of \(\mathrm{ob}\) in a basis of \(V\) span the relation space \(J/\mathfrak mJ\). If these coefficients did not span, there would exist a nonzero linear functional \(\lambda\colon J/\mathfrak mJ\twoheadrightarrow\C\) annihilating each coefficient. Pushing out the small extension along \(\lambda\) would then produce a small extension with one-dimensional kernel \(\C\) whose obstruction is the image of \(\mathrm{ob}\), which vanishes by construction. The homomorphism \(\mathscr H\to B/(J+\mathfrak m^q)\) would therefore lift to a local \(R\)-algebra homomorphism into this pushed-out ring.

Any such lift must fix \(s\) and shift each generator \(x_i\) by an element \(\delta_i\) in the square-zero kernel. By Taylor expansion, each defining relation \(f\in J\) evaluates to \(f(s,x)+\sum_i(\partial f/\partial x_i)\delta_i\). Because each \(\partial f/\partial x_i\) lies in \(\mathfrak m\) and \(\mathfrak m\) annihilates the square-zero kernel, the correction terms vanish identically, leaving the image of \(f\) unchanged. Since \(\lambda\ne0\), some relation in \(J\) has a nonzero image under \(\lambda\), which directly contradicts the fact that any homomorphism from \(\mathscr H=B/J\) must annihilate every relation. Hence the coefficients of \(\mathrm{ob}\) must span \(J/\mathfrak mJ\), establishing the dimension bound \(\dim_\C J/\mathfrak mJ\le\dim_\C V\le a\).

By Nakayama's lemma, the minimal number of generators of \(J\) is \(r=\dim_\C J/\mathfrak mJ\le a\). Because \(\mathscr H/(s)\) is a nonzero Artinian ring, the ideal \((J,s)\) is \(\mathfrak m\)-primary, with height equal to \(\dim B=a+1\). By Krull's height theorem, the ideal \((J,s)\) generated by the \(r\) generators of \(J\) together with \(s\) can have height at most \(r+1\) \cite[Lemma~10.60.12, Tag~0BBZ]{Stacks}. This yields the reverse inequality \(a+1\le r+1\), which forces \(r=a\).

Choose a minimal set of generators \(f_1,\ldots,f_a\) of \(J\). The sequence \(f_1,\ldots,f_a,s\) forms a full system of parameters in the regular local ring \(B\), and therefore constitutes a regular sequence \cite[Lemmas~10.106.3 and~10.104.2, Tags~00NQ and~02JN]{Stacks}. Consequently, \(\mathscr H=B/J\) is Cohen--Macaulay of dimension one, and the parameter \(s\) is a nonzerodivisor.

Because the closed quotient is Artinian, the \(s\)-adic and maximal-ideal topologies are cofinal, so \(\mathscr H\) is complete in the \(s\)-adic topology. Lift a basis modulo \(s\). Successive reduction expresses each element as a linear combination of these lifts with remainders divisible by increasing powers of \(s\). Completeness allows the coefficients to be summed. These lifts therefore generate \(\mathscr H\) as an \(R\)-module. Because \(s\) is a nonzerodivisor, this finitely generated \(R\)-module is torsion-free. Over the DVR \(R\), it is finite free, with rank given by the closed-fibre dimension \(k\).
\end{proof}

\begin{proposition}
\label{prop:conserved-multiplicities-gv}
For \(1\le j\le\ell\), let \(m_j=\operatorname{ord}_{w=0}P_j(f_0(w))\) be the intersection multiplicity of the original classifying arc with the reduced curve divisor \(D_j\), as in Lemma~\ref{lem:arc-not-in-discriminant}. Let \(n_j\) be the genus-zero GV invariant of the original flop in class \(j[C]\), defined as the local Hilbert-scheme length at the pure thickening \(C^{(j)}\) in Proposition~\ref{prop:general-source-package}(iv) \cite[Definition~2.4]{Toda}. Then
\begin{equation}
\label{eq:m_j-n_j}
 m_j=n_j.
\end{equation}
Thus the polynomial \(W_j\) of Lemma~\ref{lem:weierstrass-root-polynomials} for the translated formal family has exactly \(n_j\) roots over a splitting field.
\end{proposition}

\begin{proof}
Write \(X_{\mathrm{orig}}\) for the original smooth threefold representative before target completion. We compare \(k_j\), the length of the local ring \(\mathcal O_{\operatorname{Hilb}(X_{\mathrm{orig}}),[C_j]}\) introduced in Proposition~\ref{prop:general-source-package}(iii), with the number of exceptional curves of class \(j[C]\) in the geometric generic fibre of the prescribed constant translation. Lemma~\ref{lem:symbolic-pure-thickening} identifies \(C_j\) with the pure thickening \(C^{(j)}\), and hence gives \(k_j=n_j\). To prove \(m_j=n_j\), we construct a finite free \(\C[[s]]\)-algebra \(\mathscr H_j\) prorepresenting embedded deformations of \(C_j\) in the translated family, show that its closed fibre has length \(k_j=n_j\), and show that its geometric generic fibre is reduced of degree \(m_j\).

We set \(R_0=\C[[s]]\), \(D=R_0[[w]]=\C[[s,w]]\), and \(K=\C((s))\). Pulling back the algebraic simultaneous partial resolution along the translation morphism \(F(s,w)=f_0(w)+sv\) produces the following projective morphism prior to target completion:
\[
 \pi_D:X_D\longrightarrow Y_D=\Spec Q_D.
\]
Here the coordinate algebra \(Q_D=\mathcal R_\ell\otimes_{\mathbb H_\ell}D\) is of finite type over the complete Noetherian ring \(D\), and \(X_D\) is flat over \(D\) (hence flat over \(R_0\)) by Proposition~\ref{prop:general-source-package}(i). Completing \(Y_D\) at the maximal ideal of the central contraction point and pulling back \(X_D\) recovers the completed translated family. We define the central threefold fibre \(X_c:=X_D\times_{\Spec R_0}\Spec\C\), the structure sheaf \(F_j:=\mathcal O_{C_j}\), and the ideal sheaf \(I_j:=\ker(\mathcal O_{X_c}\to F_j)\).

\emph{The relative Hilbert functor.} Let \(M=\operatorname{Hilb}_{X_D/Y_D}\) denote the relative Hilbert space, viewed as an algebraic space over \(R_0\). It is locally of finite presentation over \(Y_D\), since \(\pi_D\) is separated and of finite presentation \cite[Proposition~99.9.4, Tag~0D01]{Stacks}. An object over an algebra \(A\) consists of a target morphism \(y_A:\Spec A\to Y_D\) together with a proper flat closed subscheme of \(X_D\times_{Y_D,y_A}\Spec A\) \cite[Situation~99.9.1, Tag~0CZY, and the following paragraph in Section~99.9, Tag~0CZX]{Stacks}. Choose an \'etale scheme chart \(U\to M\) containing a \(\C\)-point \(u_j\) lying over \([C_j]\), and define
\[
 \mathscr H_j=\widehat{\mathcal O}_{U,u_j}.
\]
This algebra is a complete Noetherian local \(R_0\)-algebra. By formal \'etaleness of the chart, \(\mathscr H_j\) prorepresents the functor of embedded deformations of \([C_j]\) over local Artinian \(R_0\)-algebras with this central point fixed.

This relative Hilbert functor also parametrizes all embedded proper deformations in \(X_D\times_{R_0}\Spec A\), without requiring their target point to be fixed a priori. Indeed, Katz proves that \(H^0(F_j)=\C\) and \(H^q(F_j)=0\) for all \(q>0\) \cite[Lemma~3.2]{KatzGV}. For any flat sheaf deformation \(F_A\) with proper support, induction over small Artin extensions gives \(H^0(F_A)\cong A\) and \(H^q(F_A)=0\) for all \(q>0\). Specifically, any small extension kernel \(N\) yields the long exact cohomology sequence associated with \(0\to F_j\otimes N\to F_A\to F_{A/N}\to0\). Lifting the central generator \(1\) and comparing module lengths shows that the natural map \(A\to H^0(F_A)\) is an isomorphism. The evaluation homomorphism \(\mathcal O_{X_A}\to F_A\) is surjective by Nakayama's lemma. Furthermore, altering the chosen generator of \(H^0(F_A)\) by a unit factor preserves its kernel ideal.

Moreover, \(\operatorname{End}(F_A)=H^0(F_A)=A\), so the module action of \(Q_D\) uniquely determines the moving section \(y_A\). These constructions commute with Artin base change and invert the forgetful functor from relative Hilbert closed subschemes to sheaf deformations.

Restricting to \(s=0\), these deformations correspond to embedded deformations in the original \emph{threefold}. Each deformed structure sheaf is annihilated by a power of the central target ideal, because that action vanishes on the closed sheaf and successive powers land in higher powers of the maximal ideal of the Artinian base. Consequently, the finite target neighbourhoods agree before and after full completion. Comparing with the original contraction thus yields the isomorphism
\[
 \mathscr H_j/(s)\cong
 \widehat{\mathcal O}_{\operatorname{Hilb}(X_{\mathrm{orig}}),[C_j]},
 \qquad
 \dim_\C\mathscr H_j/(s)=k_j=n_j.
\]
Here Lemma~\ref{lem:symbolic-pure-thickening} identifies the symbolic thickening \(C_j\) with the pure thickening defining \(n_j\), guaranteeing that their local Hilbert lengths coincide. The local ring \(\mathscr H_j/(s)\) is Artinian by Katz's isolated-component theorem.

\emph{Complete obstructions for the prescribed ambient family.} Let \(A'\twoheadrightarrow A\) be a small extension of local Artinian \(R_0\)-algebras with kernel \(N\) annihilated by \(\mathfrak m_{A'}\), and let \(0\to K_A\to\mathcal O_{X_A}\to F_A\to0\) represent a flat quotient. Because \(\mathcal O_{X_A}\) and \(F_A\) are flat over \(A\), the ideal sheaf \(K_A\) is likewise \(A\)-flat. Let \(\widetilde K_A\subset\mathcal O_{X_{A'}}\) be the full preimage of \(K_A\). Flatness of \(X_D\) over \(R_0\) implies that the ambient sheaf \(\mathcal O_{X_{A'}}\) is flat over \(A'\), which induces the canonical extension sequence
\[
 0\longrightarrow F_A\otimes_A N
 \longrightarrow
 \widetilde K_A/(K_A\otimes_A N)
 \longrightarrow K_A\longrightarrow0.
\]
The extension class of this sequence in \(\operatorname{Ext}^1_{X_A}(K_A,F_A\otimes_A N)\) is the standard obstruction to lifting flat quotients \cite[Remark~99.7.9, Tag~0CZU]{Stacks}. Its vanishing is equivalent to the existence of a flat quotient lift: in one direction, a splitting yields the lifted ideal sheaf by preimage in \(\widetilde K_A\), whose intersection with \(\mathcal O_{X_A}\otimes_A N\) is \(K_A\otimes_A N\). Consequently, the quotient satisfies the nilpotent flatness criterion of \cite[Lemma~91.11.2, Tag~08MQ]{Stacks}.

Conversely, any flat quotient lift defines such a splitting. This obstruction theory commutes with linear pushouts of \(N\). Because any lifted quotient is coherent and supported over a finite central target thickening, its support is proper; hence it belongs to the relative Hilbert functor defined above. The general quotient formalism of \cite[Situation~99.7.1, Tag~082M]{Stacks} does not require \(X_D/R_0\) to be of finite type.

Derived restriction sends the flat module \(K_A\) to \(I_j\), so the natural obstruction space is \(V_j=\operatorname{Ext}^1_{X_c}(I_j,F_j)\), while differences between lifts are parametrized by the tangent space \(T_j=\operatorname{Hom}_{X_c}(I_j,F_j)\) \cite[Remark~99.7.10, Tag~0CZV]{Stacks}. In the long exact sequence obtained from \(0\to I_j\to\mathcal O_{X_c}\to F_j\to0\), evaluation at the section \(1\) identifies \(\operatorname{Hom}(F_j,F_j)\) with \(H^0(F_j)\). Since \(H^1(F_j)=H^2(F_j)=0\), the long exact sequence yields the natural isomorphisms
\[
 T_j\cong\operatorname{Ext}^1(F_j,F_j),\qquad
 V_j\cong\operatorname{Ext}^2(F_j,F_j).
\]

To compare the tangent and obstruction dimensions, we apply duality to a complex with proper support. Let \(X\) be the original smooth quasi-projective threefold representative, and choose a projective compactification \(\overline X\). The proper support of \(F_j\) is closed in \(\overline X\) and disjoint from the boundary \(\overline X\setminus X\). The extension by zero \(\overline F_j\) is therefore coherent and perfect on \(\overline X\), vanishing on an open neighbourhood of the boundary. Proper duality for \(\overline F_j\) applies even when \(\overline X\) is singular \cite[Lemmas~48.27.1(5) and~48.27.4, Tags~0FVV and~0FVY]{Stacks}. By open restriction of upper shriek \cite[Lemma~48.17.2, Tag~0G4J]{Stacks}, the dualizing complex of \(\overline X\) restricts to that of the smooth open subscheme \(X\). Because \(X\) immerses into projective space, the smooth fundamental-class formula identifies this restriction with \(\omega_X[3]\) \cite[Lemma~48.29.2, Tag~0E9Z]{Stacks}.

The canonical bundle of the Gorenstein target is trivial near the closed point, so crepancy trivializes \(\omega_X\) on a neighbourhood of \(\operatorname{Supp}(F_j)\). Setting \(E=R\mathcal{H}om(\overline F_j,\overline F_j)\), perfect duality and tensor symmetry identify \(E^\vee\) with \(E\). Because \(E\) is supported where \(\omega_X\) is trivial, proper duality gives the duality isomorphism \(\mathbb H^i(\overline X,E)^\vee\cong\mathbb H^{3-i}(\overline X,E)\). In particular, we obtain \(\operatorname{Ext}^2(F_j,F_j)\cong \operatorname{Ext}^1(F_j,F_j)^\vee\).

These self-\(\operatorname{Ext}\) groups are invariant under completion at the contraction point of the target. Indeed, the derived pushforward of the self-\(\mathcal Hom\) complex has finite-length cohomology supported at the contraction point; faithfully flat completion of the base commutes with this pushforward by derived base change \cite[Lemma~36.22.5, Tag~08IB]{Stacks}. On the central fibre \(X_c\), perfection near the support descends from the faithfully flat completed local rings, which are isomorphic to the completed regular local rings of the original smooth threefold. The identical self-\(\mathcal Hom\) comparison therefore transfers the dimension equality to \(X_c\), establishing \(\dim V_j=\dim T_j=:a\).

We now apply Lemma~\ref{lem:complete-obstruction-flatness}: the relative tangent space is \(T_j\), the complete obstruction space is \(V_j\), their dimensions agree, and the closed fibre has length \(k_j\). Hence \(\mathscr H_j\) is finite free over \(R_0\), of rank \(k_j\).

\emph{The generic fibre and its specialization.} It remains to identify and count the geometric generic points of this finite free algebra. To obtain a family over its spectrum, we first return to the local Hilbert chart. The completion homomorphism \(\mathcal O_{U,u_j}\to\mathscr H_j\) is an isomorphism: \(\mathcal O_{U,u_j}\) injects into its completion by the Krull intersection theorem, forming an \(R_0\)-submodule of the finite module \(\mathscr H_j\). This submodule is finite and closed in the \(s\)-adic topology, and is also dense by cofinality with the maximal-ideal topology, so the two rings coincide \cite[Lemma~10.51.4, Tag~00IP, and Lemma~10.97.1, Tag~00MA]{Stacks}. The universal Hilbert family on the chart therefore pulls back to \(\Spec\mathscr H_j\). At the generic-fibre points of \(\Spec\mathscr H_j\), the relative tangent spaces are those of the relative Hilbert functor.

The target morphism \(\Spec\mathscr H_j\to Y_D\) extends uniquely to the completion of \(Y_D\) at its central point: letting \(\mathfrak n\) denote the central target ideal, \(\mathfrak n^N\mathscr H_j\subset s\mathscr H_j\) for sufficiently large \(N\). Base change along this completion preserves the family.

Now work over a common finite splitting field \(\C((u))/K\), with valuation ring \(\C[[u]]\). The Weierstrass and transversality lemmas provide exactly \(m_j\) distinct roots of \(W_j\). Lemma~\ref{lem:root-to-odp} gives a split \((-1,-1)\)-curve at each root, with reduced contraction fibre \(\mathbb P^1\). These constitute all complete curves, because the curve discriminant accounts for the entire locus where complete curves exist.

To extend each generic curve over the valuation ring, first extend its image in the target. Form the generic base from \(D[1/s]\) before completing at any moving point. By Lemma~\ref{lem:multicurve-critical-support}, the image of each such curve in \(Y_D\) lies in the relative critical locus, which is finite over \(R_0\). This finiteness over \(R_0\) ensures that all coordinates of the image point are integral over \(R_0\), extending uniquely to a target section over \(\C[[u]]\) that specializes to the original contraction point.

For each of these \(m_j\) generic curves \(E\), let \(\overline E\) be its scheme-theoretic closure in \(X_D\times_{Y_D}\Spec\C[[u]]\). Because its defining ideal sheaf is \(u\)-saturated, its structure sheaf is \(u\)-torsion-free and therefore flat over the discrete valuation ring \(\C[[u]]\). Saturation also ensures that this flat closure is unique. Fix a relatively very ample line bundle \(L\) from the universal model, and set \(h=\deg_C L>0\). The marked degree formula gives \(\deg_E L=jh\), so constancy of Hilbert polynomials forces \(P_{\overline E_0}(t)=jh\,t+1\).

The special fibre \(Z=\overline E_0\) is scheme-theoretically contained in the central fibre \(C_\ell\subset S_0\) \cite[Lemma~3.2(ii)]{KatzGV}. The pure one-dimensional quotient of \(\mathcal O_Z\) has fundamental cycle \(j[C]\), and hence equals \(\mathcal O_{C_j}\): on the normal surface, height-one primary ideals are determined by their multiplicities and are the symbolic powers. Writing \(T\) for the zero-dimensional kernel, Katz's identity \(\chi(\mathcal O_{C_j})=1\) yields
\[
 1=\chi(\mathcal O_Z)=1+\operatorname{length}T.
\]
The kernel \(T\) therefore has length zero, establishing the scheme-theoretic equality \(Z=C_j\) as closed subschemes. The resulting Hilbert morphism lifts through the chosen chart over the strictly henselian ring \(\C[[u]]\), defining a generic point of \(\Spec\mathscr H_j\).

Conversely, every geometric generic point of \(\Spec\mathscr H_j\) parametrizes a closed subscheme of a contraction fibre having Hilbert polynomial \(jh\,t+1\). That fibre contains a complete curve and is the reduced integral \(\mathbb P^1\) over one of the roots. Because its only one-dimensional closed subscheme is the curve itself, the object is precisely the corresponding exceptional curve.

Two points of this local Hilbert scheme representing the same curve define the same saturated flat closure and the same \'etale lift with prescribed closed point, and must therefore coincide. Thus each embedded exceptional curve corresponds to a unique geometric generic point of the local Hilbert scheme.

It remains to show that the geometric generic fibre of \(\Spec\mathscr H_j\) is reduced. For each split \((-1,-1)\)-curve \(E\), Lemma~\ref{lem:multicurve-odp-sheaf-contribution} proves that \(\operatorname{Ext}^1(\mathcal O_E,\mathcal O_E)=0\). Combining this vanishing with \(H^0(\mathcal O_E)=\C((u))\) and \(H^1(\mathcal O_E)=0\) ensures that the embedded Hilbert tangent space at \([E]\) is zero. Every geometric generic local ring is therefore a field by Nakayama's lemma, which proves that the generic fibre is reduced. Its degree is therefore computed by counting its geometric points, yielding
\[
 n_j=k_j
 =\operatorname{rank}_{R_0}\mathscr H_j
 =\#\{\text{generic exceptional curves of class }j[C]\}
 =m_j.
\]

This equality establishes the curve count in the prescribed translation \(F\). For any sufficiently small transverse analytic perturbation avoiding curve-divisor intersections on the boundary of the chosen disk, the argument principle preserves the reduced-divisor intersection number \(m_j\). This confirms the original local Hilbert count for that perturbation as well.
\end{proof}

\Needspace{15\baselineskip}
\begin{proposition}
\label{prop:full-marked-formal-morsification}
There is a finite extension \(K'/K\) with valuation ring \(\C[[u]]\), where \(K'=\C((u))\) and \(s=u^e\) for some \(e\ge1\).  Completing the contraction at each generic exceptional image gives a small resolution of a split ordinary double point. Moreover,
\[
 W_j(u^e,w)=\prod_{a=1}^{n_j}(w-w_{j,a}(u)),
 \qquad w_{j,a}(u)\in u\C[[u]],
\]
with all roots pairwise distinct, including roots belonging to different coefficients.  The pullback of the universal partial resolution has exactly \(n_j\) pairwise disjoint \((-1,-1)\)-curves in class \(j[C]\), and no other complete curves.
\end{proposition}

\begin{proof}
By Lemma~\ref{lem:good-to-formal-generic}, the polynomials \(W_j\) are separable and pairwise coprime. At each geometric root of \(W_j\), the classifying map \(F\) meets \(D_j\) transversely, and its image lies outside the bad locus \(\Sigma^{\mathrm{bad}}\). By Lemma~\ref{lem:root-divisors-and-bad-locus}, every lift of this image point under the quotient map \(q:\mathfrak h\to H\) lies on a unique root hyperplane. Passing to a finite field extension, all the polynomials \(W_j\) split into linear factors and all these lifts become defined over the same field. Scalar extension preserves both transversality of the intersection and the vanishing of a unique root pair. The proof of Lemma~\ref{lem:root-to-odp} relies on the explicit ADE hypersurface equations, the invertibility of the fibre Hessian, and a change of smoothing parameter that absorbs an invertible unit. All three ingredients remain valid after scalar extension and identify the target singularity at each exceptional image as a geometrically ordinary double point. The formal coordinate recursions require only the inverse Hessian and characteristic-zero integer denominators, while splitting the quadratic tangent cones requires adjoining only finitely many square roots.

Pass to a further finite field extension that supplies these square roots together with the local rulings at all finitely many points. To verify that the formal coordinate transformations converge and are defined over this field, write the completed hypersurface equation as \(q(x)+r_d(x)+\) terms of degree at least \(d+1\), where \(q\) is a nondegenerate quadratic form, \(d\ge3\), and \(r_d\) is homogeneous of degree \(d\). Letting \(A\) denote the Hessian matrix of \(q\), the substitution
\[
 x\longmapsto x-\frac{1}{d}A^{-1}\nabla r_d
\]
eliminates the homogeneous term \(r_d\) of degree \(d\), because Euler's homogeneous relation yields \(x^{\mathsf t}\nabla r_d=d r_d\), while all newly created terms have degree at least \(d+1\). Successive substitutions converge in the maximal-ideal topology to a formal automorphism having identity linear part, reducing the completed equation to its quadratic part. Diagonalizing this quadratic form and adjoining square roots of its nonzero coefficients along with a square root of \(-1\) yields the split ordinary double point equation \(a_1a_2-b_1b_2=0\). The two small resolutions correspond to its two rulings.

Denote the resulting splitting field by \(K'\). The integral closure of \(\C[[s]]\) in \(K'\) is a finite complete local DVR with residue field \(\C\) \cite[Remark~15.113.6 and Lemmas~10.161.8, 10.153.4, Tags~09E8, 032L and~04GH]{Stacks}. Choose a uniformizing parameter \(t\). Successively extracting residues yields a unique \(t\)-adic expansion, identifying this valuation ring with \(\C[[t]]\). Write the ramified parameter as \(s=t^e a(t)\) with \(a(t)\) a unit. Selecting an \(e\)-th root of the nonzero constant term of \(a(t)\) and applying the binomial series to the remaining unit factor defines \(a(t)^{1/e}\). Then \(u=t a(t)^{1/e}\) is a new uniformizer satisfying \(s=u^e\).

Each root \(w_{j,a}(u)\) of \(W_j\) is integral because \(W_j\) is a monic polynomial, and the central specialization \(W_j(0,T)=T^{n_j}\) forces its constant coefficient to vanish. Hence \(w_{j,a}(u)\in u\C[[u]]\). Because the polynomials \(W_j\) are separable and pairwise coprime, all roots \(w_{j,a}(u)\) are pairwise distinct. At each root, Lemma~\ref{lem:root-to-odp} produces a single \((-1,-1)\)-curve of class \(j[C]\). Because these curves lie in distinct \(w\)-fibres, they are pairwise disjoint.

We finally check that these are all the complete curves. Because the completed target scheme is affine, every complete curve in the source must contract to a point in the target, and therefore lies entirely in the exceptional locus. In particular, the parameter coordinate \(w\) is constant along any such curve. By definition, the curve discriminant \(D^{\mathrm{curv}}\) is the locus of parameter values whose simultaneous partial resolutions contain a complete exceptional curve.

For each parameter point in the formal generic disk lying on the curve discriminant, its image in \(H\) belongs to a unique reduced divisor \(D_j\) and avoids \(\Sigma^{\mathrm{bad}}\). The unique exceptional curve lying over that point has already been accounted for in the root count above, which proves that the list of complete curves is exhaustive.
\end{proof}

\begin{remark}
\label{rem:generic-fibre-ring}
The relative family retains all moving roots over \(\C[[s,w]]\). Its generic fibre is formed over \(\C[[s,w]][1/s]\) and finite extensions of that ring, as in Lemma~\ref{lem:good-to-formal-generic}. Over \(\C((u))[[w]]\), the element \(w-w_{j,a}(u)\) is a unit whenever \(w_{j,a}\ne0\).  Thus one must localize the finite relative family before making any further completion; completing first at \(w=0\) would discard the points defined by the nonzero roots of the \(W_j\).
\end{remark}

\subsection{Tilting generators and the local matrix blocks}
\label{subsec:full-universal-tilting}

We know the exceptional curves of the generic fibre and their classes \(j[C]\). We now compute the completed contraction quotient at each curve. The pulled-back universal bundle remains a perverse projective generator. At a split ordinary double point, its positive dual decomposes into free summands and copies of the unique nonfree indecomposable projective in the positive perverse heart. The determinant degree counts these nonfree copies. Their stable endomorphism algebra is a matrix algebra of that size.

\subsubsection{The universal generator and base change}

Karmazyn constructs the simultaneous partial resolution and its universal algebra \cite[Sections~2.3 and 3.3]{Karmazyn}.  We use the normalized basic generator on this family given by Van den Bergh's construction \cite[Propositions~3.2.5--3.2.7]{VanDenBergh}. The universal length-\(\ell\) partial resolution carries a basic zero-perverse generator and its positive dual
\[
 \mathcal T=\mathcal O\oplus\mathcal T',
 \qquad
 \mathcal P=\mathcal T^\vee=\mathcal O\oplus\mathcal P'.
\]
Set
\[
 \mathcal A_\ell=\End_{\mathcal Y_\ell}(\mathcal T)^{\mathrm{op}}
 \cong\End_{\mathcal Y_\ell}(\mathcal P),
 \qquad e_0=[\mathcal O].
\]
The bundle \(\mathcal P\) belongs to the \((-1)\)-perverse coherent heart, while \(\mathcal T\) belongs to the zero-perverse coherent heart. Duality preserves the free projector and replaces the endomorphism algebra by its opposite.

By the universal-flop construction \cite[Sections~2.3 and 3.3]{Karmazyn}, the universal source \(\mathcal Y_\ell\), target \(\mathcal X_\ell=\Spec\mathcal R_\ell\), centre, and endomorphism algebra \(\mathcal A_\ell\) are flat over the marked base \(H=\Spec\mathbb H_\ell\), and the target \(\mathcal X_\ell\) is affine of finite type over \(H\). We will use this flatness over \(H\) to establish Tor-independence and ordinary endomorphism base change for every Noetherian classifying algebra.

On the universal family, the positive nonfree summand \(\mathcal P'\) is the normalized maximal extension
\begin{equation}
 0\longrightarrow\mathcal O^{\oplus(\ell-1)}
 \longrightarrow\mathcal P'
 \longrightarrow\mathcal L\longrightarrow0,
 \qquad \deg_C\mathcal L=1.
 \label{eq:universal-maximal-extension}
\end{equation}
Consequently, \(\operatorname{rk}\mathcal P'=\ell\) and \(\deg_C\det\mathcal P'=1\). These numerical normalizations follow directly from the specification of the basic generator. Indeed, Karmazyn's graded choice completes at the central target point to the basic generator \cite[Remark~3.5]{Karmazyn}. Van den Bergh's complete-local classification identifies its nonfree positive summand as having determinant degree one on the exceptional curve \(C\) \cite[Proposition~3.5.4 and Theorem~3.5.5]{VanDenBergh}, and this summand has rank \(\ell\) \cite[Remark~3.11]{Karmazyn}.

The maximal extension sequence \eqref{eq:universal-maximal-extension} can be constructed directly on the universal family itself. Put \(E=\mathcal P'\), and let \(\mathcal J_C\) be the ideal sheaf of the reduced central curve. By structure-sheaf cohomology vanishing on the contraction, \(H^1(\mathcal J_C)=0\). Because \(E\) is globally generated and the contraction has cohomological dimension one, we likewise have \(H^1(\mathcal J_C E)=0\). The global sections of \(E\) therefore surject onto its sections on the central curve \(C\).

On the curve \(C\), the globally generated restriction \(E|_C\) has rank \(\ell\) and degree one, which forces the splitting \(E|_C\cong\mathcal O_C(1)\oplus\mathcal O_C^{\ell-1}\). Choosing \(\ell-1\) general sections defines a line-bundle quotient of degree one, and the vanishing \(H^1(\mathcal O_C(-1))=0\) splits the resulting extension. The line subbundle \(\mathcal O_C(1)\) is invariant under the grading action. Because the multiplicative group \(\mathbb G_m\) is linearly reductive, the induced short exact sequence of global sections on \(C\) admits an equivariant splitting. A weight basis of the complementary trivial bundle on the curve then lifts to homogeneous global sections of \(E\).

The degeneracy locus of these lifted sections is closed and \(\mathbb G_m\)-invariant. If this locus were nonempty, its image under the projective contraction would contain the vertex of the positively graded affine target, contradicting the linear independence of the lifted sections on \(C\). Consequently, these sections define a bundle injection \(\mathcal O^{\oplus(\ell-1)}\to E\) whose cokernel is the line bundle \(\mathcal L=\det E\) everywhere on the universal family.

Finally, the vanishing \(H^1(E^\vee)=0\) shows that the extension classes generate \(H^1(\mathcal L^{-1})\). Upon completion, these extension classes are minimal because the nonfree positive summand is indecomposable \cite[Lemma~3.5.2]{VanDenBergh}.
\begin{lemma}
\label{lem:arbitrary-pullback-tor-independent}
Let \(B\) be a Noetherian \(\mathbb H_\ell\)-algebra, and form
\[
\begin{array}{ccc}
 \mathcal Y_B&\xrightarrow{g}&\mathcal Y_\ell\\
 \pi_B\downarrow&&\downarrow\pi\\
 \mathcal X_B&\xrightarrow{h}&\mathcal X_\ell.
\end{array}
\]
This square is Tor-independent.  Consequently
\begin{equation}
 R\pi_{B*}Lg^*E\simeq Lh^*R\pi_*E
 \label{eq:derived-base-change-expanded}
\end{equation}
for every perfect \(E\) on \(\mathcal Y_\ell\), and \(R\pi_{B*}\mathcal O_{\mathcal Y_B}=\mathcal O_{\mathcal X_B}\).
\end{lemma}

\begin{proof}
Write \(C\to A\to Q\) for the homomorphisms of local rings of the marked base \(H\), the target \(\mathcal X_\ell\), and the source \(\mathcal Y_\ell\), respectively. By the flatness of the universal family over the base, both \(A\) and \(Q\) are flat over \(C\). For any \(C\)-algebra \(B\), we have
\[
 Q\otimes_A^L(A\otimes_C B)
 \simeq Q\otimes_C^LB
 \simeq Q\otimes_CB.
\]
The displayed Cartesian square is therefore Tor-independent. The morphism \(\pi\) is projective, hence quasi-compact and quasi-separated, while the perfect complex \(E\) has quasi-coherent cohomology. The derived base-change theorem of \cite[Lemma~36.22.5, Tag~\href{https://stacks.math.columbia.edu/tag/08IB}{08IB}]{Stacks} applies to this square and gives the canonical comparison map
\[
 Lh^*R\pi_*E\longrightarrow R\pi_{B*}Lg^*E
\]
as an isomorphism, giving \eqref{eq:derived-base-change-expanded}. Specializing to \(E=\mathcal O\) and applying the identity \(R\pi_*\mathcal O=\mathcal O\) on the universal family yields the final assertion.
\end{proof}

\begin{lemma}
\label{lem:perverse-generator-survives}
The pullbacks \(\mathcal P_B=g^*\mathcal P\) and \(\mathcal T_B=g^*\mathcal T\) are projective generators of the \((-1)\)- and zero-perverse hearts of \(\pi_B\) respectively, and they induce the corresponding equivalences on bounded coherent derived categories.
\end{lemma}

\begin{proof}
We verify Van den Bergh's projective-generator criterion on the pulled-back family. Write \(\mathcal X_\ell=\Spec\mathcal R_\ell\). Because \(\mathcal R_\ell\) is an algebra of finite type over \(\mathbb H_\ell\), the pullback target scheme \(\mathcal X_B=\Spec(\mathcal R_\ell\otimes_{\mathbb H_\ell}B)\) is Noetherian. The pulled-back contraction \(\pi_B\) is projective, its fibres have dimension at most one, and derived base change from Lemma~\ref{lem:arbitrary-pullback-tor-independent} yields \(R\pi_{B*}\mathcal O=\mathcal O\). These checks confirm the geometric hypotheses of the criterion over the Noetherian base, with no requirement of smoothness or reducedness on \(B\).

The category \(\mathfrak V\) in this criterion consists of globally generated bundles \(M\) with \(R^1\pi_*(M^\vee)=0\). A finite evaluation surjection for the globally generated bundle \(\mathcal P\) remains surjective under pullback. Because \(\mathcal P\in\mathfrak V\), the pushforward on the universal family is concentrated in degree zero: \(R\pi_*\mathcal T=(\pi_*\mathcal T)[0]\). Applying Lemma~\ref{lem:arbitrary-pullback-tor-independent} gives
\[
 R\pi_{B*}\mathcal T_B\simeq Lh^*(\pi_*\mathcal T).
\]
The complex on the left side lies in \(D^{\ge0}\), while the derived pullback of a sheaf on the right side lies in \(D^{\le0}\). Both sides are therefore concentrated in degree zero, giving \(R^1\pi_{B*}\mathcal T_B=0\).

Pullback likewise preserves the relative ampleness of \(\det\mathcal P\), global generation, and the distinguished direct summand \(\mathcal O\). Van den Bergh's projective-generator criterion therefore applies to the pulled-back family \cite[Propositions~3.2.5--3.2.7]{VanDenBergh}, and the derived equivalences follow from \cite[Corollary~3.2.8]{VanDenBergh}.
\end{proof}

\begin{proposition}
\label{prop:ordinary-endomorphism-base-change}
There is a canonical isomorphism of marked algebras
\begin{equation}
 \End_{\mathcal Y_B}(\mathcal T_B)^{\mathrm{op}}
 \cong\mathcal A_\ell\otimes_{\mathbb H_\ell}B.
 \label{eq:ordinary-end-base-change}
\end{equation}
Moreover, \(\operatorname{Ext}^i_{\mathcal Y_B}(\mathcal T_B,\mathcal T_B)=0\) for every \(i>0\).  The isomorphism preserves \(e_0\), and the bundle pullbacks preserve ranks and determinants. The algebra identification remains valid after central localization and Noetherian completion.
\end{proposition}

\begin{proof}
Lemma~\ref{lem:arbitrary-pullback-tor-independent} establishes derived base change for perfect complexes on the classifying pullback. We choose the perfect complex to be the endomorphism bundle \(\mathcal E nd(\mathcal T)\). By the tilting property on the universal family, the derived pushforward is concentrated in degree zero, giving
\[
 R\pi_*\mathcal E nd(\mathcal T)
  =\pi_*\mathcal E nd(\mathcal T)[0].
\]
Set \(M=\Gamma(\mathcal Y_\ell,\mathcal E nd(\mathcal T)) =\mathcal A_\ell^{\mathrm{op}}\), and write \(\mathcal X_\ell=\Spec\mathcal R_\ell\). Both \(\mathcal R_\ell\) and \(M\) are flat over the marked base ring \(\mathbb H_\ell\) on the universal family \cite[Section~3.3, the defining paragraph for \(\mathcal A_\ell\)]{Karmazyn}, so
\[
 M\otimes_{\mathcal R_\ell}^{L}
       (\mathcal R_\ell\otimes_{\mathbb H_\ell}B)
 \simeq M\otimes_{\mathbb H_\ell}^{L}B
 \simeq M\otimes_{\mathbb H_\ell}B.
\]
Consequently, the derived pushforward of the pulled-back endomorphism bundle remains concentrated in degree zero. Taking global sections over the affine target scheme establishes the algebra isomorphism \eqref{eq:ordinary-end-base-change} along with the vanishing of all higher self-Ext groups.

In degree zero, the comparison homomorphism is induced by pullback of endomorphisms, so it preserves composition and the identity. Pullback likewise preserves the inclusion and retraction of the distinguished free summand, and consequently preserves its projector \(e_0\). Furthermore, rank and determinant commute with bundle pullback.

Finally, because the contraction morphism \(\pi_B\colon\mathcal Y_B\to\mathcal X_B\) is proper, the endomorphism algebra \(\End_{\mathcal Y_B}(\mathcal T_B)\) is a finitely generated module over the Noetherian coordinate algebra \(\mathcal O(\mathcal X_B)\). Both the localization of this Noetherian target ring and the subsequent ideal-adic completion of the localized ring are flat ring extensions. Applying flat base change to the coherent pushforward sheaf \(\pi_{B*}\mathcal E nd(\mathcal T_B)\) identifies the endomorphism algebra of the pulled-back generator with the corresponding ordinary scalar extension. Because the algebra is finite over the target ring, completion identifies this scalar extension with the completed endomorphism algebra. Flat base change likewise preserves the vanishing of all higher self-Ext groups and the distinguished idempotent projector \(e_0\).
\end{proof}

The affine target ring is also the centre of this endomorphism algebra. Indeed, any central endomorphism of \(\mathcal P_B=\mathcal O\oplus\mathcal P'_B\) commutes with the free projector and is therefore diagonal, with a scalar \(r\) on \(\mathcal O\).  Commuting with every map \(\mathcal O\to\mathcal P'_B\) forces the other diagonal block to act by \(r\) on all global sections.  Global generation then makes that block \(r\) times the identity. 

Hence the centre is \(\Gamma(\mathcal Y_B,\mathcal O)=\Gamma(\mathcal X_B,\mathcal O)\). This argument also applies after the stated localizations and completions, and duality gives the same centre for the zero-perverse convention.

Write \(\Lambda_B=\mathcal A_\ell\otimes_{\mathbb H_\ell}B\). Since the endomorphism comparison preserves \(e_0\), the relative contraction quotient is \(\Lambda_B/\Lambda_Be_0\Lambda_B\).

\begin{proposition}
\label{prop:literal-special-basic-nccr}
Let \(B=\C[[s,w]]\) carry a classifying map whose restriction modulo \(s\) is the central arc.  Localize \(\Lambda_B\) at the central contraction point of its centre and complete there. The reduction \(\widehat{\Lambda_B}/(s)\) is isomorphic to the basic noncommutative crepant resolution (NCCR) of the original flop, preserving the free idempotent \(e_0\).
\end{proposition}

\begin{proof}
Because the endomorphism algebra \(\Lambda_B\) is a finite module over its Noetherian centre, completion at the maximal ideal of the central contraction point commutes with reduction modulo \(s\). Lemma~\ref{lem:marked-analytic-classifying-arc} then identifies the geometric closed fibre at \(s=0\) with the original completed flop contraction.

Lemma~\ref{lem:perverse-generator-survives} ensures that the reduced positive bundle remains a projective generator in the perverse coherent heart. On the completed local one-curve flop, Van den Bergh's classification establishes that there are exactly two indecomposable positive perverse projectives, namely the trivial bundle \(\mathcal O\) and the maximal extension bundle \(\mathcal M\), with \(\operatorname{rk}\mathcal M=\ell\) and \(\deg_C\det\mathcal M=1\) \cite[Proposition~3.5.4 and Theorem~3.5.5]{VanDenBergh}.

We write the reduced nonfree positive summand as a direct sum \(\mathcal O^{\oplus a}\oplus\mathcal M^{\oplus b}\). Comparing the determinant degrees on \(C\) forces \(b=1\), and subsequently comparing vector bundle ranks yields \(\ell=a+\ell\), so that \(a=0\). The reduced positive generator is therefore isomorphic to \(\mathcal O\oplus\mathcal M\). Taking duals recovers the chosen zero-perverse basic generator of the original flop and preserves its distinguished free projector \(e_0\).
\end{proof}

\subsubsection{The generator at an intersection with \texorpdfstring{\(D_j\)}{D j}}

We work over the finite splitting extension field \(K'=\C((u))\) of Proposition~\ref{prop:full-marked-formal-morsification}, over which all local models split. Choose a point \(\lambda\) lying over a transverse intersection of the geometric generic classifying arc with the curve divisor \(D_j\). We write \(\widehat Y_\lambda\to\Spec\widehat R_\lambda\) for the completed small resolution at \(\lambda\) and \(C_\lambda\cong\mathbb P^1_{K'}\) for its isolated exceptional curve.

\begin{lemma}
\label{lem:split-odp-projectives}
Let \(F\) be any characteristic-zero field and put
\[
 R_F=F[[a_1,a_2,b_1,b_2]]/(a_1a_2-b_1b_2),\qquad I=(a_1,b_1).
\]
On the small resolution \(\pi_Z\colon Z=\operatorname{Proj}_{R_F}\bigoplus_{n\ge0}I^n\to\Spec R_F\), let \(\mathcal Q=\mathcal O_Z(1)\).  Then \(\mathcal O_Z\oplus\mathcal Q\) is the basic positive perverse projective generator.  Every positive perverse projective vector bundle is a direct sum of \(\mathcal O_Z\) and \(\mathcal Q\), and
\[
 \operatorname{rk}\mathcal Q=1,\qquad
 \deg_C\mathcal Q=1,\qquad \pi_{Z*}\mathcal Q=I.
\]
These assertions remain true after completion along the exceptional curve.
\end{lemma}

\begin{proof}
Write \(\theta=b_1/a_1\) on the \(a_1\)-chart and \(\theta^{-1}=a_1/b_1\) on the \(b_1\)-chart. After removing the Rees algebra torsion and completing along the exceptional curve, the coordinate charts and their transition gluing are given by
\[
 \widehat U_{a_1}=\Spf F[\theta][[a_1,b_2]],
 \quad b_1=a_1\theta,\quad a_2=\theta b_2,
\]
\[
 \widehat U_{b_1}=\Spf F[\theta^{-1}][[b_1,a_2]],
 \quad a_1=\theta^{-1}b_1,\quad b_2=\theta^{-1}a_2.
\]
Thus \(Z\) is the split Atiyah small resolution. Its exceptional curve \(C\) is cut out by \(a_1=b_2=0\) on \(\widehat U_{a_1}\) and by \(b_1=a_2=0\) on \(\widehat U_{b_1}\), identifying the curve as \(C\cong\mathbb P^1_F\).

The transition function of the tautological quotient line bundle \(\mathcal Q\) is that of \(\mathcal O_{\mathbb P^1}(1)\), which establishes the restriction \(\mathcal Q|_C=\mathcal O_C(1)\). Because the quotient ring satisfies \(R_F/I\cong F[[a_2,b_2]]\), the depth lemma shows that \(I\) is a maximal Cohen--Macaulay (MCM) module and is therefore reflexive.

By Hartogs extension across the codimension-two curve \(C\) on the smooth threefold resolution \(Z\), the direct image \(\pi_{Z*}\mathcal Q\) is likewise reflexive. Because \(\pi_{Z*}\mathcal Q\) coincides with \(I\) away from the singular vertex, we obtain \(\pi_{Z*}\mathcal Q=I\).

To verify the perverse-generator conditions, recall that the normal bundle is \(N_{C/Z}\cong\mathcal O_C(-1)^{\oplus2}\). On the \(n\)-th infinitesimal graded piece of the formal neighbourhood of \(C\), the restrictions of \(\mathcal O_Z\) and \(\mathcal Q^{-1}\) are given respectively by
\[
 \operatorname{Sym}^nN_{C/Z}^{\vee}\cong
 \mathcal O_C(n)^{\oplus(n+1)}
 \quad\text{and}\quad
 \mathcal O_C(n-1)^{\oplus(n+1)}.
\]
A direct computation on the two affine charts shows that \(H^1(\mathbb P^1_F,\mathcal O(k))=0\) for all \(k\ge-1\); consequently, the first cohomology groups of both graded pieces vanish for all degrees \(n\ge0\). The chart transition relations identify the defining ideal sheaf of \(C\) with \(\mathfrak m_{R_F}\mathcal O_Z\), so that the infinitesimal thickenings of \(C\) match the formal thickenings in the theorem on formal functions. Induction on the thickening order together with \cite[Theorem~30.20.5, Tag~\href{https://stacks.math.columbia.edu/tag/02OC}{02OC}]{Stacks} yields the higher direct image vanishings
\[
 R^1\pi_{Z*}\mathcal O_Z=R^1\pi_{Z*}\mathcal Q^{-1}=0.
\]
The algebraic small resolution of the affine quadratic cone is the total space of the vector bundle \(\mathcal O_{\mathbb P^1}(-1)^{\oplus2}\), and the formal charts describe its completion along the zero section. The line bundle \(\mathcal Q\) is relatively ample and globally generated. Combined with the direct image vanishings established above and the presence of the free direct summand \(\mathcal O_Z\), these properties satisfy Van den Bergh's projective-generator criterion \cite[Propositions~3.2.5--3.2.7]{VanDenBergh}, both on the algebraic scheme over \(\Spec R_F\) and on the formal completion along the exceptional curve \(C\).

For the classification of projectives, the endomorphism algebra of \(\mathcal O_Z\oplus\mathcal Q\) is finite over the complete local ring \(R_F\) and is semiperfect. The two summands have local endomorphism rings \(\End(\mathcal O_Z)=R_F=\End(\mathcal Q)\), and are distinguished by their degrees on \(C\). The two summand projectors form a complete set of pairwise inequivalent primitive idempotents. Projectives over a semiperfect ring decompose as direct sums of the corresponding indecomposable projectives, so the tilting equivalence yields the decomposition of every projective in the perverse heart.

Both summands and their idempotents are defined over the ground field \(F\). Because \(Z\) is projective over the complete Noetherian ring \(R_F\), Grothendieck's existence theorem identifies coherent sheaves on \(Z\) with coherent sheaves on its completion along \(C\) \cite[Lemma~30.24.3, Tag~0885]{Stacks}. Formal functions preserves the cohomological vanishings defining the perverse hearts, and full faithfulness preserves their Hom spaces. Hence this equivalence preserves perverse projectives and the displayed decomposition on the formal completion.
\end{proof}

\begin{proposition}
\label{prop:coefficient-j-projective-decomposition}
Let \(\mathcal Q\) be the unique nonfree indecomposable positive perverse projective on the completed ODP small resolution, normalized by \(\operatorname{rk}\mathcal Q=1\) and \(\deg\det\mathcal Q|_{C_\lambda}=1\).  Then
\begin{equation}
 \widehat{\mathcal P}'_\lambda
 \cong\mathcal O^{\oplus(\ell-j)}\oplus\mathcal Q^{\oplus j}.
 \label{eq:coefficient-j-projective-decomposition}
\end{equation}
Let \(N_\lambda\) denote the pushforward of \(\widehat{\mathcal P}'_\lambda\) to the completed target. Then
\[
 N_\lambda\cong
 \widehat R_\lambda^{\oplus(\ell-j)}\oplus I_\lambda^{\oplus j},
\]
where \(I_\lambda\) is the rank-one Atiyah MCM module.
\end{proposition}
\begin{proof}
Lemma~\ref{lem:perverse-generator-survives} shows that the pulled-back positive bundle remains projective in the perverse heart. Proposition~\ref{prop:ordinary-endomorphism-base-change} identifies its endomorphism algebra after field specialization and completion, preserving direct summands, ranks, and determinants. Because the roots are split over \(K'\), the local model is the split ordinary double point of Lemma~\ref{lem:split-odp-projectives} with \(F=K'\). Applying that classification lemma and the Krull--Schmidt theorem yields the decomposition \(\widehat{\mathcal P}'_\lambda\cong \mathcal O^{\oplus a}\oplus\mathcal Q^{\oplus b}\). Comparing vector bundle ranks on both sides gives \(a+b=\ell\). Taking determinants commutes with pullback, and Lemma~\ref{lem:marked-degree-specialization} computes the specialized degree as
\[
 \deg_{C_\lambda}\det\widehat{\mathcal P}'_\lambda
 =j\deg_C\det\mathcal P'=j.
\]
On the other hand, the direct sum decomposition \(\widehat{\mathcal P}'_\lambda\cong \mathcal O^{\oplus a}\oplus\mathcal Q^{\oplus b}\) has determinant degree \(b\) on \(C_\lambda\). Equating these degrees gives \(b=j\), which in turn implies \(a=\ell-j\). Finally, Lemma~\ref{lem:split-odp-projectives} identifies the pushforward of \(\mathcal Q\) with the Atiyah module \(I_\lambda\).
\end{proof}

\subsubsection{The local contraction quotient}

The local calculation now reduces to endomorphisms of the copies of the Atiyah module in the preceding decomposition, modulo maps through free modules. For any module \(M\) over a commutative local ring \(R\), we denote by \(\underline{\End}_R(M)\) the stable endomorphism ring, defined as the quotient of \(\End_R(M)\) by the two-sided ideal of endomorphisms that factor through finite free \(R\)-modules.

\begin{lemma}
\label{lem:atiyah-stable-endomorphisms}
Let
\[
 R_{\mathrm{odp}}=F[[a_1,a_2,b_1,b_2]]/(a_1a_2-b_1b_2),\qquad I=(a_1,b_1),
\]
where \(F\) has characteristic zero.  Then
\[
 \underline{\End}_{R_{\mathrm{odp}}}(I)\cong F,
 \qquad
 \underline{\End}_{R_{\mathrm{odp}}}(I^{\oplus j})\cong M_j(F).
\]
\end{lemma}

\begin{proof}
Write \(R=R_{\mathrm{odp}}\). Because \(I\) is a rank-one reflexive module over the normal domain \(R\), its endomorphisms are multiplications by fractions. At every height-one prime ideal of \(R\), the module \(I\) is locally free of rank one; consequently, the colon ideal satisfies \((I:I)\subseteq\bigcap_{\operatorname{ht}p=1}R_p=R\). Thus \(\End_R(I)=R\).

To identify the subspace of endomorphisms factoring through a free module, consider the fractional ideal \(I^{-1}=\{q\in\operatorname{Frac}R:qI\subseteq R\}\). We compute
\[
 I^{-1}=R+R(b_2/a_1).
\]
The right side lies in \(I^{-1}\) because \((b_2/a_1)a_1=b_2\) and \((b_2/a_1)b_1=a_2\). Conversely, let \(qI\subseteq R\) and define \(r=qa_1\in R\). Then \(b_1r=a_1(qb_1)\in a_1R\). In \(R/(a_1)=F[[a_2,b_1,b_2]]/(b_1b_2)\), the annihilator of \(b_1\) is \((b_2)\). Hence \(r=a_1r_0+b_2r_1\) for elements \(r_0,r_1\in R\), so \(q=r_0+(b_2/a_1)r_1\) and proves the formula for \(I^{-1}\).

The endomorphisms of \(I\) that factor through a finite free \(R\)-module form the product ideal
\[
 II^{-1}=(a_1,b_1)+(b_2,a_2)=(a_1,a_2,b_1,b_2)=\mathfrak m.
\]
The stable endomorphism algebra is therefore \(\underline{\End}_R(I)=R/\mathfrak m=F\). For the direct sum \(I^{\oplus j}\), both endomorphisms and maps factoring through free modules are matrices of the corresponding maps for \(I\). Their quotient algebra is therefore isomorphic to \(M_j(F)\).
\end{proof}

\begin{proposition}
\label{prop:idempotent-resolved-mcm-comparison}
For the completed ODP small resolution, pushforward gives an isomorphism
\begin{equation}
 \End_{\widehat Y_\lambda}(\widehat{\mathcal P}_\lambda)
 \cong
 \End_{\widehat R_\lambda}
 \bigl(\widehat R_\lambda^{\oplus(\ell-j+1)}
       \oplus I_\lambda^{\oplus j}\bigr)
 \label{eq:resolved-mcm-comparison}
\end{equation}
which preserves the projector onto the distinguished free summand.
\end{proposition}

\begin{proof}
The completed ordinary double point \(\Spec\widehat R_\lambda\) is a normal complete local Gorenstein threefold scheme, and the source of its small resolution is regular. The contraction morphism \(\widehat\pi\colon\widehat Y_\lambda\to\Spec\widehat R_\lambda\) is projective, birational, and crepant, having fibres of dimension at most one and satisfying \(R\widehat\pi_*\mathcal O_{\widehat Y_\lambda}=\widehat R_\lambda\). By Lemma~\ref{lem:perverse-generator-survives}, the positive vector bundle \(\widehat{\mathcal P}_\lambda\) is a perverse projective generator. The inclusion \(\widehat{\mathcal P}_\lambda\in\mathfrak V\) guarantees the vanishing \(R^1\widehat\pi_*(\widehat{\mathcal P}_\lambda^\vee)=0\), while global generation along with the fibre dimension bound yields \(R^1\widehat\pi_*\widehat{\mathcal P}_\lambda=0\). Van den Bergh's endomorphism comparison theorem therefore identifies the endomorphism algebra on the resolution with the endomorphism algebra of the pushed-forward maximal Cohen--Macaulay module \cite[Lemma~3.2.9 and Proposition~3.2.10]{VanDenBergh}. By naturality of the direct image functor, this isomorphism maps the projector onto the trivial line bundle summand \(\mathcal O_{\widehat Y_\lambda}\) to the projector onto the distinguished free summand \(\widehat R_\lambda\).
\end{proof}

\begin{theorem}
\label{thm:entire-completed-j-factor}
At every point over a transverse intersection with \(D_j\) produced by the formal marked deformation of the classifying arc, the completion of the relative contraction quotient is the matrix algebra \(M_j(K')\).  The maximal ideal of the completed target ring \(\widehat R_\lambda\) acts trivially on this factor.
\end{theorem}

\begin{proof}
After the arbitrary classifying pullback, both central localization and Noetherian completion are flat base changes. Consequently, for any module-finite algebra \(\Lambda\) with an idempotent \(e_0\),
\[
 \widehat{\Lambda/\Lambda e_0\Lambda}
 \cong\widehat\Lambda/\widehat\Lambda e_0\widehat\Lambda.
\]
By Proposition~\ref{prop:idempotent-resolved-mcm-comparison}, the completed quotient is the endomorphism algebra of \(\widehat R^{\oplus(\ell-j+1)}\oplus I^{\oplus j}\) modulo the ideal generated by one distinguished free projector.

Every other free projector factors through that copy, so the quotient is the stable endomorphism algebra of \(I^{\oplus j}\). Lemma~\ref{lem:atiyah-stable-endomorphisms} gives \(M_j(K')\).  The completed target ring acts through its residue field, proving the last assertion.
\end{proof}

\subsection{Assembling a finite free deformation}
\label{subsec:full-general-assembly}

We know the exceptional curves of the generic fibre and the completed contraction quotient at each curve. The root multiplicities come from the Bryan--Katz--Leung divisor count and the GV interpretation used above \cite[Lemma~2.8]{BKL}. Toda's formula gives the dimension of the closed fibre \cite[Theorem~1.1]{Toda}. We now construct a single finite relative quotient whose generic fibre decomposes into these local factors.

Put \(B=\C[[s,w]]\), and let \(\nu:\Spec B\to H\) be the completion of the constant-translation family \(F(s,w)=f_0(w)+sv\) from Lemmas~\ref{lem:good-constant-translate-expanded} and \ref{lem:weierstrass-root-polynomials}.  Set
\[
 R_B=\mathcal R_\ell\otimes_{\mathbb H_\ell}B.
\]
Let \(\mathfrak n\subset R_B\) be the maximal ideal of the central contraction point (it contains \(s\)), and define the central completion
\[
 S_\nu=\widehat{(R_B)_{\mathfrak n}}.
\]
The contraction over this complete local target is
\[
 \pi_\nu:\mathcal Y_\nu
 :=\mathcal Y_\ell\times_{\Spec\mathcal R_\ell}\Spec S_\nu
 \longrightarrow\Spec S_\nu.
\]
Here the map \(\Spec S_\nu\to\Spec\mathcal R_\ell\) comes from the composite \(\mathcal R_\ell\to R_B\to S_\nu\). The universal source and target are flat over the marked base. After pullback to \(B\), localization and completion of the affine target are flat, and their pullbacks to the source remain flat. Thus \(\mathcal Y_\nu\) and \(\Spec S_\nu\) are flat over \(B\), and hence over \(\C[[s]]\). The fibre of \(\pi_\nu\) at \(s=0\) is the original completed flop.

Let \(g_\nu:\mathcal Y_\nu\to\mathcal Y_\ell\) denote the canonical projection, and let \(\mathcal T=\mathcal O_{\mathcal Y_\ell}\oplus\mathcal T'\) be the universal tilting bundle fixed in Subsection~\ref{subsec:full-universal-tilting}. We form the pulled-back vector bundle \(\mathcal T_\nu=g_\nu^*\mathcal T=\mathcal O_{\mathcal Y_\nu}\oplus g_\nu^*\mathcal T'\) on \(\mathcal Y_\nu\). By Proposition~\ref{prop:ordinary-endomorphism-base-change}, the opposite of its endomorphism algebra is canonically isomorphic to the scalar extension of \(\mathcal A_\ell\) over \(S_\nu\), and this isomorphism preserves the distinguished idempotent projector \(e_0=[\mathcal O_{\mathcal Y_\nu}]\) onto the trivial line bundle. Passing to the quotient by the two-sided ideal generated by \(e_0\) defines the relative algebras
\begin{equation}
 \Lambda_\nu=\mathcal A_\ell\otimes_{\mathcal R_\ell}S_\nu
              \cong\End_{\mathcal Y_\nu}(\mathcal T_\nu)^{\mathrm{op}},
 \qquad
 \mathscr A=\Lambda_\nu/\Lambda_\nu e_0\Lambda_\nu.
 \label{eq:expanded-relative-contraction-algebra}
\end{equation}

\begin{proposition}
\label{prop:general-special-fibre-identification}
There is a marked algebra isomorphism
\[
 \mathscr A/(s)\cong \Acon.
\]
\end{proposition}

\begin{proof}
By Proposition~\ref{prop:literal-special-basic-nccr}, specializing modulo \(s\) yields an isomorphism \(\Lambda_\nu/(s)\cong\Lambda\) as marked algebras preserving the distinguished projector \(e_0\), where \(\Lambda\) is the selected basic NCCR of the original flop. Right exactness of the tensor product with \(\C[[s]]/(s)\) over \(\C[[s]]\) then gives
\[
 \mathscr A/(s)
 \cong
 \frac{\Lambda}{\Lambda e_0\Lambda}
 =\Acon.
\]
\end{proof}

\begin{lemma}
\label{lem:topological-finiteness-expanded}
Let \(R_0=\C[[s]]\), and define \(S\) as follows:
\[
 S=R_0[[x_1,\ldots,x_m]]/J.
\]
Then \(S\) is complete Noetherian.  Let \(M\) be finite over \(S\).  If \(M/sM\) is finite-dimensional over \(\C\), then \(M\) is finite over \(R_0\).
\end{lemma}

\begin{proof}
Choose elements \(m_1,\ldots,m_d\in M\) lifting a \(\C\)-basis of the closed fibre module \(M/sM\). For any given element \(m\in M\), successive approximation modulo powers of \(sM\) produces the truncation expansions
\[
 m=\sum_{i=1}^d\left(\sum_{a=0}^{q-1}c_{ia}s^a\right)m_i+s^qm^{(q)},
 \qquad c_{ia}\in\C.
\]
The compatible coefficient truncations define the formal power series \(c_i(s)=\sum_{a\ge0}c_{ia}s^a\) are well-defined elements of the complete discrete valuation ring \(R_0=\C[[s]]\). For each truncation order \(q\ge1\), the tail of each series \(c_i(s)\) beyond degree \(q-1\) is divisible by \(s^q\), which implies the containment
\[
 m-\sum_i c_i(s)m_i\in s^qM.
\]
Here \(S\) is a Noetherian local ring with \(s\in\mathfrak m_S\), and \(M\) is finite over \(S\), Krull's intersection theorem guarantees that \(\bigcap_{q\ge1}s^qM\subseteq\bigcap_{q\ge1}\mathfrak m_S^qM=0\). Consequently, the remainder difference \(m-\sum_{i=1}^d c_i(s)m_i\) vanishes identically. The lifts \(m_1,\ldots,m_d\) therefore generate \(M\) as a finite \(R_0\)-module.
\end{proof}

\begin{proposition}
\label{prop:relative-algebra-finite-over-dvr}
The algebra \(\mathscr A\) is a finite \(\C[[s]]\)-module.
\end{proposition}

\begin{proof}
The universal target is of finite type over \(H\), so its pullback \(\Spec R_B\) is of finite type over \(B=\C[[s,w]]\). Completion at the central contraction point therefore gives a presentation
\[
 S_\nu\cong\C[[s,x_1,\ldots,x_m]]/J
\]
for some \(m\) and some ideal \(J\). By Proposition~\ref{prop:general-source-package}(i), \(\mathcal A_\ell\) is finite over \(\mathcal R_\ell\). Base change and completion make \(\Lambda_\nu\) finite over \(S_\nu\), so its quotient \(\mathscr A\) is also finite over \(S_\nu\).

By Proposition~\ref{prop:general-special-fibre-identification}, the closed fibre \(\mathscr A/(s)\cong \Acon\) is finite-dimensional over \(\C\). Applying the topological finiteness criterion of Lemma~\ref{lem:topological-finiteness-expanded} to the finite \(S_\nu\)-module \(\mathscr A\), we conclude that \(\mathscr A\) is finitely generated as a module over \(\C[[s]]\).
\end{proof}

\begin{lemma}
\label{lem:fullness-away-from-roots}
On the geometric generic fibre of the classifying family, the relative contraction quotient vanishes away from the finitely many points over transverse intersections with the curve divisors. It is a finite algebra supported on exactly those closed points.
\end{lemma}

\begin{proof}
The root-divisor description identifies \(D^{\mathrm{curv}}\) with the locus of parameters over which the partial-resolution contraction contains a complete curve \cite[Theorem~2.6]{NabijouWemyss}. Away from \(D^{\mathrm{curv}}\), any positive-dimensional fibre of the projective contraction would contain a complete curve, which implies that the contraction is quasi-finite, and hence finite, outside \(D^{\mathrm{curv}}\). Because \(R\pi_{\ell *}\mathcal O_{\mathcal Y_\ell} =\mathcal O_{\mathcal X_\ell}\), the base-changed contraction is an isomorphism wherever it is finite.

In Proposition~\ref{prop:full-marked-formal-morsification}, the classifying arc meets \(D^{\mathrm{curv}}\) solely at the specified transverse roots, and the source scheme contains a single exceptional curve over each intersection point. Away from these points, every direct summand of the pushed-forward tilting bundle is locally free, and after shrinking becomes finite free.

Every endomorphism of the pushed-forward tilting bundle on this locus factors through a direct sum of copies of the distinguished rank-one free summand. Consequently, the idempotent \(e_0\) is full and the relative contraction algebra quotient \(\mathscr A\) vanishes identically there.

To determine the moving support after completing the target at its central contraction point, we compute the zeroth Fitting ideal over the base. Put \(R_H=\widehat{\mathcal O}_{H,0}\), and let \(\mathcal C_H\) denote the universal contraction quotient completed at the central target point. Proposition~\ref{prop:arbitrary-base-classifying-arcs} identifies \(\mathscr A\) with the ordinary pullback of \(\mathcal C_H\) along the formal arc \(\nu\). Furthermore, Proposition~\ref{prop:weighted-contraction-presentation} supplies a square presentation matrix of \(\mathcal C_H\) over \(R_H\) whose determinant equals an invertible unit times \(\prod_jP_j^{j^2}\). Commuting the zeroth Fitting ideal with base change yields
\[
 \mathscr A\cong\mathcal C_H\otimes_{R_H}B,
 \qquad
 \operatorname{Fitt}^{B}_0(\mathscr A)
 =\left(\prod_jG_j(s,w)^{j^2}\right).
\]

The proofs of these results rely on universal endomorphism base change and the explicit algebra presentation over \(\widetilde R=\widehat{\mathcal O}_{\mathfrak h,0}\), and do not use the semisimple-existence theorem being proved here. For the present application, the external deformation base is \(\C[[s]]\), over which source and target are flat as verified above.

After passing to the finite splitting extension, each root defines a \(K'\)-rational closed point of \(\Spec(B[1/s]\otimes_{\C((s))}K')\). At a root of \(W_j\), the local uniformizing parameter is \(w-w_{j,a}\), and the zeroth Fitting ideal over the base has vanishing order \(j^2>0\). The localization of \(\mathscr A\) is therefore nonzero there.

Because the free projector is full at every point where the contraction is an isomorphism, this nonzero support must lie over an exceptional curve. At each lift of the corresponding parameter to \(\mathfrak h\), exactly one positive root vanishes, so the local geometric model contains precisely one exceptional curve. Thus the completed quotient is nonzero at each listed root and supported nowhere else.

At each root of \(W_j\), Theorem~\ref{thm:entire-completed-j-factor} identifies the completed local factor with the matrix algebra \(M_j(K')\). By Proposition~\ref{prop:relative-algebra-finite-over-dvr}, the generic quotient is finite-dimensional over \(K'\) and hence finite over its Artinian central image. Its support is therefore precisely the stated finite set of \(K'\)-rational closed points.
\end{proof}

\begin{lemma}
\label{lem:no-infinitesimal-thickening-expanded}
Let \(S\) be a Noetherian \(K'\)-algebra and \(A\) a finite \(S\)-algebra supported at \(K'\)-rational maximal ideals \(p_1,\ldots,p_r\), so \(S/p_i=K'\).  Suppose
\[
 \widehat{A_{p_i}}\cong M_{d_i}(K')
\]
and the maximal ideal of \(\widehat{S_{p_i}}\) acts trivially.  Then
\[
 A\cong\prod_{i=1}^rM_{d_i}(K').
\]
In particular its scheme-theoretic support is reduced.
\end{lemma}

\begin{proof}
The completion map \(S_{p_i}\to\widehat{S_{p_i}}\) is faithfully flat. Since \(p_i\widehat{A_{p_i}}=0\), faithful flatness implies \(p_iA_{p_i}=0\). Thus the algebra before completion, \(A_{p_i}\), already coincides with its completion and is isomorphic to \(M_{d_i}(K')\). Put \(J=\bigcap_{i=1}^r p_i\). The finite \(S\)-module \(JA\) localizes to zero outside \(\{p_1,\ldots,p_r\}\), and it also localizes to zero at each \(p_i\) because \(p_jS_{p_i}=S_{p_i}\) for \(j\ne i\). Thus \(JA=0\).

By the Chinese remainder theorem, \(S/J\cong(K')^r\). The central idempotents of \(S/J\) split \(A\) into the product of its localizations \(A_{p_i}\cong M_{d_i}(K')\).
\end{proof}

\begin{proposition}
\label{prop:exact-general-generic-product}
After the finite splitting extension \(K'/\C((s))\) of Proposition~\ref{prop:full-marked-formal-morsification},
\begin{equation}
 \mathscr A\otimes_{\C[[s]]}K'
 \cong\prod_{j=1}^{\ell}M_j(K')^{\times n_j}.
 \label{eq:expanded-general-generic-product}
\end{equation}
\end{proposition}

\begin{proof}
Proposition~\ref{prop:conserved-multiplicities-gv} establishes that each polynomial \(W_j\) has degree \(n_j\). By Proposition~\ref{prop:full-marked-formal-morsification}, the roots of all polynomials \(W_j\) are pairwise distinct and split over the finite field extension \(K'=\C((u))\), where \(s=u^e\) and each root satisfies \(w_{j,a}(u)\in u\C[[u]]\). These roots define \(N_{\mathrm{GV}}=\sum_j n_j\) distinct \(K'\)-rational points on the generic fibre of \(\Spec\C[[s,w]]\), all of whose closures specialize to the central origin \((s,w)=(0,0)\). Lemma~\ref{lem:fullness-away-from-roots} proves that the scalar extension \(\mathscr A\otimes_{\C[[s]]}K'\) is supported exclusively at these \(K'\)-rational points, while Proposition~\ref{prop:relative-algebra-finite-over-dvr} ensures that it is a finite-dimensional algebra over \(K'\). At each root of \(W_j\), Theorem~\ref{thm:entire-completed-j-factor} identifies the completed local factor with the matrix algebra \(M_j(K')\) and shows that the maximal ideal of the completed target ring annihilates this factor. Lemma~\ref{lem:no-infinitesimal-thickening-expanded} therefore removes any residual infinitesimal thickening at these support points and gives the product decomposition \eqref{eq:expanded-general-generic-product}.
\end{proof}

\begin{lemma}
\label{lem:dvr-torsion-exclusion-expanded}
Let \(M\) be a finite module over the DVR \(\C[[s]]\).  If
\[
 \dim_{\C}M/sM
 =\dim_{\C((s))}M[1/s]=r,
\]
then \(M\) is free of rank \(r\).
\end{lemma}

\begin{proof}
By the structure theorem for finite modules over a DVR,
\[
 M\cong\C[[s]]^{\oplus r}
 \oplus\bigoplus_{h=1}^{m}\C[[s]]/(s^{e_h}).
\]
After tensoring with the fraction field only the free summands survive, giving dimension \(r\). Reduction modulo \(s\) has dimension \(r+m\).  Equality of these dimensions forces \(m=0\).
\end{proof}

\begin{proof}[Proof of Theorem~\ref{thm:general-semisimple-existence}]
We consider the relative contraction quotient algebra \(\mathscr A\) defined in \eqref{eq:expanded-relative-contraction-algebra}. Proposition~\ref{prop:general-special-fibre-identification} identifies its closed fibre with \(\Acon\), while Proposition~\ref{prop:exact-general-generic-product} decomposes its geometric generic fibre into \(n_j\) matrix factors of size \(j\) for each coefficient \(1\le j\le\ell\). Because vector space dimension is invariant under scalar extension from the formal Laurent field \(\C((s))\) to its finite field extension \(K'\), we obtain the generic dimension
\[
 \dim_{\C((s))}\mathscr A[1/s]
 =\sum_{j=1}^{\ell}j^2n_j.
\]
Proposition~\ref{prop:relative-algebra-finite-over-dvr} establishes that \(\mathscr A\) is a finitely generated \(\C[[s]]\)-module. On the other hand, Toda's dimension formula from Proposition~\ref{prop:general-source-package}(v) and the closed-fibre isomorphism yield
\[
 \dim_{\C}\mathscr A/(s)
 =\dim_{\C}\Acon
 =\sum_{j=1}^{\ell}j^2n_j.
\]
Because the closed-fibre and generic-fibre dimensions coincide, Lemma~\ref{lem:dvr-torsion-exclusion-expanded} shows that \(\mathscr A\) is torsion-free and therefore finite free of rank \(\sum_{j=1}^\ell j^2n_j\) over the discrete valuation ring \(\C[[s]]\).
\end{proof}

\section{The cocenter and uniqueness of the semisimple generic fibre}
\label{sec:general-uniqueness}

Section~\ref{sec:general-existence} constructs a deformation with geometric generic fibre \(S_{\mathrm{GV}}\). Thus the contraction algebra lies in the closure of the change-of-basis orbit of \(S_{\mathrm{GV}}\). Is there any other semisimple orbit whose closure contains the contraction algebra? We prove that this orbit closure is the unique irreducible component of the scheme of associative algebra laws through that point. By symmetry and four-periodicity, the comparison of tangent and orbit dimensions reduces to a computation of the cocenter on the marked ADE base.

\subsection{Statement and proof strategy}
\label{subsec:uniqueness-statement}

For a completed simple one-curve flop \(f\colon X\to\Spec R\) of a smooth complex threefold, let \(A=\Acon\) be its contraction algebra of dimension \(d=\dim_{\C}A\), and let \(n_j\) denote the genus-zero Gopakumar--Vafa invariants of the flop. We write the split semisimple algebra and its number of matrix factors as
\[
 S_{\mathrm{GV}}:=\prod_{j=1}^{\ell}M_j(\C)^{\times n_j},\qquad N_{\mathrm{GV}}:=\sum_{j=1}^{\ell}n_j.
\]
\begin{theorem}
\label{thm:general-uniqueness}
For every completed simple one-curve flop of a smooth complex threefold,
\begin{equation}
\label{eq:uniqueness-main}
 \dim_{\C}HH_0(A)=\dim_{\C}Z(A)=N_{\mathrm{GV}}.
\end{equation}
Consequently \(A\) is a smooth point of the fixed-unit associative-law scheme \(X_d\), the orbit closure \(\overline{G\cdot S_{\mathrm{GV}}}\) is the unique irreducible component of \(X_d\) through \(A\), and every complex semisimple algebra to which \(A\) deforms is isomorphic to \(S_{\mathrm{GV}}\).  More generally, if \(\mathscr D\) is a finite free unital \(\C[[u]]\)-algebra with \(\mathscr D/u\mathscr D\cong A\) and geometrically semisimple generic fibre, then \(\mathscr D\otimes_{\C[[u]]}\overline{\C((u))}\cong S_{\mathrm{GV}}\otimes_{\C}\overline{\C((u))}\).
\end{theorem}

The theorem identifies the geometric generic fibre of every finite free associative deformation with geometrically semisimple generic fibre. Section~\ref{sec:general-existence} separately constructs a finite free \(\C[[u]]\)-deformation that realizes it.

The orbit criterion requires equality of the tangent and orbit dimensions. Symmetry and four-periodicity of the contraction algebra, together with centre--cocenter duality, reduce that equality to \(\dim HH_0(A)=N_{\mathrm{GV}}\). Subsection~\ref{subsec:cocenter-statements-reduction} states this formula and derives it from a universal presentation and the comparison along the original classifying arc. The following three subsections establish those inputs. The detailed \(E_8\) computation for length five is in Appendix~\ref{app:e8-executable-certificate}.

\subsection{The orbit criterion}
\label{sec:orbit-criterion}

The orbit analysis of Jelisiejew--Smoktunowicz establishes local closedness and openness for semisimple algebra multiplication laws \cite[Section~1 and Theorem~2.14]{JS}. We require the corresponding results for multiplications with a fixed identity element. Let \(k\) be an algebraically closed field of characteristic zero, let \(V=k\cdot1\oplus W\) be a \(d\)-dimensional \(k\)-vector space with distinguished non-zero vector \(1\), and let \(X_d\) denote the affine scheme of associative multiplications on \(V\) for which \(1\) is the two-sided identity. The scheme \(X_d\) is defined inside \(\operatorname{Hom}_k(V\otimes_k V,V)\) by the quadratic associativity equations and the linear equations fixing the unit. The algebraic subgroup preserving the distinguished identity element is
\[
G=\Stab_{\GL(V)}(1).
\]

For any unital associative \(k\)-algebra \(\Lambda\), the normalized Hochschild cochain spaces are defined by
\[
 C^n_{\mathrm{norm}}(\Lambda,\Lambda)
 =\operatorname{Hom}_k((\Lambda/k\cdot1)^{\otimes n},\Lambda).
\]
We denote by \(Z^n_{\mathrm{norm}}\) and \(B^n_{\mathrm{norm}}\) the spaces of cocycles and coboundaries of this normalized Hochschild complex.

\begin{lemma}
\label{lem:fixed-unit-tangents}
For \(\Lambda\in X_d(k)\),
\[
T_\Lambda X_d=Z^2_{\mathrm{norm}}(\Lambda,\Lambda),\qquad
T_\Lambda(G\cdot\Lambda)=B^2_{\mathrm{norm}}(\Lambda,\Lambda),
\]
and
\[
\operatorname{Lie}\Stab_G(\Lambda)=\Der(\Lambda).
\]
Moreover, \(G\) is smooth, connected, and irreducible of dimension \(d(d-1)\), and its orbits are precisely the unital isomorphism classes.
\end{lemma}

\begin{proof}
A dual-number multiplication \(m+\varepsilon\phi\) on \(V[\varepsilon]/(\varepsilon^2)\) satisfies the unit condition for the distinguished element \(1\) precisely when \(\phi(1,-)=\phi(-,1)=0\), which places \(\phi\) in the normalized cochain space \(C^2_{\mathrm{norm}}(\Lambda,\Lambda)\). Modulo \(\varepsilon^2\), the associativity condition linearizes to the normalized Hochschild cocycle equation \(\delta\phi=0\), identifying the Zariski tangent space of the law scheme with these normalized cocycles. Differentiating an infinitesimal change of basis \(g=\operatorname{id}_V+\varepsilon f\) with \(f(1)=0\) yields the normalized Hochschild coboundaries as its tangent image. The kernel of this differential consists of all endomorphisms satisfying the Leibniz rule and annihilating the unit, which is the derivation space \(\Der(\Lambda)\).

Relative to the direct-sum decomposition \(k\cdot1\oplus W\), every element of \(G\) takes the unique block form
\[
\begin{pmatrix}1&\lambda\\0&h\end{pmatrix},
\qquad \lambda\in W^*,\quad h\in \GL(W).
\]
Thus \(G\cong W^*\rtimes \GL(W)\), which proves the stated smoothness, connectedness, and dimension assertions for the group.

Endow the orbit \(O=G\cdot\Lambda\) with its reduced locally closed subscheme structure, and write \(q:G\to O\) for the orbit morphism. The scheme \(O\) is irreducible because it is the image of the irreducible algebraic group \(G\). The smooth locus of \(O\) is nonempty over the perfect field \(k\). It is invariant under the transitive action of \(G(k)\), so it is the entire orbit and \(O\) is smooth.

In characteristic zero, the dominant morphism \(q\) is generically separable. Thus the differential of \(q\) achieves maximal rank \(\dim O\) at some \(k\)-point \(g\in G\). Differentiating the equivariance identity \(q\circ L_g=a_g\circ q\), where \(a_g\) is the action of \(g\) on \(O\), shows that the differential of \(q\) has the same rank at \(g\) and at the identity because the translation automorphisms induce tangent isomorphisms. Consequently, \(d_eq\) surjects onto the tangent space \(T_\Lambda O\). Its image is the normalized coboundary space computed above, which confirms the tangent identification for the reduced orbit even when the ambient law scheme is nonreduced.

Finally, every unital algebra isomorphism preserves the distinguished unit element, so its change of basis belongs to \(G\), matching orbits with unital isomorphism classes.
\end{proof}

\begin{lemma}
\label{lem:family-orbit-closure}
Let \(\mathfrak o\) be a complete equicharacteristic DVR with algebraically closed residue field \(k\), coefficient field \(k\subset\mathfrak o\) and fraction field \(K\).  Let \(\mathscr C\) be a finite free unital \(\mathfrak o\)-algebra with special fibre \(A\).  If \(S\) is a \(d\)-dimensional \(k\)-algebra and
\[
\mathscr C_K\otimes_K\overline K
   \cong S\otimes_k\overline K,
\]
then
\[
A\in\overline{G\cdot S}\subset X_d.
\]
\end{lemma}

\begin{proof}
The unit of \(\mathscr C\) is unimodular because it reduces to the nonzero unit of \(A\).  Extend \(1_{\mathscr C}\) to an \(\mathfrak o\)-basis lifting the prescribed basis of \(A\).  The structure constants define a morphism \(\Spec\mathfrak o\to X_d\) whose closed point is \(A\).

Let \(Y=\overline{G\cdot S}\) with its reduced closed structure. The assumed geometric generic-fibre isomorphism preserves the identity, so the geometric generic point maps to a \(G(\overline K)\)-translate of \(S\).  Every defining equation of \(Y\) therefore maps to zero in \(\overline K\).  Since \(\mathfrak o\hookrightarrow\overline K\) is injective, it already maps to zero in \(\mathfrak o\).  The entire morphism factors through \(Y\), and its closed point \(A\) belongs to \(Y\).
\end{proof}

The next lemma compares the Zariski tangent space of the ambient scheme with the dimension of an orbit. When the tangent dimension matches the orbit dimension, smoothness and local openness follow at points of the orbit closure.

\begin{lemma}
\label{lem:open-orbit-smoothness}
Let a smooth irreducible algebraic group \(H\) act on a finite-type scheme \(X\) over an algebraically closed field \(k\) of characteristic zero.  Let \(y\in X(k)\), and write \(O=H\cdot y\).
\begin{enumerate}[label=\textup{(\roman*)}]
\item If \(\dim T_yX=\dim O\), then \(y\) is smooth and \(O\) is open in \(X\).
\item Suppose \(O\) is open.  If a closed point \(x\in\overline O\) satisfies \(\dim T_xX=\dim O\), then \(x\) is smooth and \(\overline O\) is its unique irreducible component.  Every open orbit whose closure contains \(x\) equals \(O\).
\end{enumerate}
\end{lemma}

\begin{proof}
Consider first the underlying reduced scheme \(X_{\mathrm{red}}\). The image of the orbit morphism is constructible, and therefore contains a nonempty open subset of its irreducible reduced closure. Translating this open subset shows that \(O\) is locally closed. Because \(k\) is perfect, this reduced orbit has a nonempty smooth locus, and transitivity under the group action makes the whole orbit smooth. This is the orbit argument of \cite[Proposition~1.11]{Brion}, applied over the stated algebraically closed field of characteristic zero. The local rings in the dimension comparison that follows are those of the original scheme \(X\).

For~(i), let \(R=\mathcal O_{X,y}\) and write the orbit local ring as \(R/I\). With \(\delta=\dim O\), the assumed tangent dimension gives
\[
 \delta=\dim(R/I)\le\dim R\le\operatorname{edim}R=\delta.
\]
Thus \(R\) is regular, hence an integral domain. The prime ideal \(I\) has height zero, so it is the zero ideal. The locally closed orbit therefore contains an open neighbourhood of \(y\). Its translates show that \(O\) is open in the original scheme \(X\).

For~(ii), the orbit is irreducible since \(H\) is irreducible. As the orbit is open, its closure \(\overline O\) is an irreducible component of dimension \(\delta\). Since \(x\) is closed, the local dimension of \(\overline O\) at \(x\) is \(\delta\). The same dimension inequalities therefore force \(\mathcal O_{X,x}\) to be regular, giving smoothness. A regular local ring has one minimal prime, so this is the unique component through \(x\).

Any other open orbit whose closure contains \(x\) has the same component as its closure. The two nonempty open orbits must intersect, and hence coincide.
\end{proof}

\begin{proposition}
\label{prop:semisimple-orbit-open}
If \(S\in X_d(k)\) is split semisimple, then \(S\) is a smooth point of the scheme \(X_d\) and \(G\cdot S\) is open in \(X_d\).
\end{proposition}

\begin{proof}
The algebra \(S\) is separable, hence projective as an \(S^e\)-module, so \(HH^i(S,S)=0\) for \(i>0\).  In degree two, vanishing of \(HH^2(S,S)\) says that every normalized cocycle is a coboundary.  By Lemma~\ref{lem:fixed-unit-tangents}, the tangent dimension at \(S\) is the orbit dimension.  Lemma~\ref{lem:open-orbit-smoothness}(i) then gives smoothness of \(X_d\) at \(S\) and openness of \(G\cdot S\).
\end{proof}

\Needspace{14\baselineskip}

\begin{theorem}
\label{thm:orbit-smooth-unique}
Let \(A,S\in X_d(k)\), with \(S\) split semisimple.  Assume
\begin{equation}
\label{eq:orbit-criterion-hypotheses}
A\in\overline{G\cdot S},
\qquad
\dim T_AX_d=\dim(G\cdot S).
\end{equation}
Then:
\begin{enumerate}[label=\textup{(\roman*)}]
\item \(A\) is a smooth point of \(X_d\);
\item \(\overline{G\cdot S}\) is the unique irreducible component of \(X_d\) through \(A\);
\item if \(S'\) is split semisimple and \(A\in\overline{G\cdot S'}\), then \(S'\cong S\).
\end{enumerate}
\end{theorem}

\begin{proof}
Proposition~\ref{prop:semisimple-orbit-open} makes the orbits of both semisimple algebras open, and \(G\) is smooth and irreducible by Lemma~\ref{lem:fixed-unit-tangents}.  The assumed equality of tangent and orbit dimensions is the remaining hypothesis of Lemma~\ref{lem:open-orbit-smoothness}(ii), which gives smoothness and the unique component. 

The equality of the two orbits yields a unital algebra isomorphism, again by Lemma~\ref{lem:fixed-unit-tangents}.
\end{proof}

\begin{corollary}
\label{cor:hochschild-orbit-criterion}
Let \(A,S\in X_d(k)\), with \(S\) semisimple and \(A\in\overline{G\cdot S}\).  The dimension equality in Theorem~\ref{thm:orbit-smooth-unique} is equivalent to
\begin{equation}
\label{eq:hochschild-numerical-criterion}
\dim\Der(A)-\dim HH^2(A)=d-\dim Z(S).
\end{equation}
In particular, it holds if
\[
\dim Z(A)=\dim Z(S),
\qquad
\dim HH^1(A)=\dim HH^2(A).
\]
\end{corollary}

\begin{proof}
Lemma~\ref{lem:fixed-unit-tangents} gives
\[
\dim T_AX_d=d(d-1)-\dim\Der(A)+\dim HH^2(A).
\]
For split semisimple \(S\), all derivations are inner and \(\dim\Der(S)=d-\dim Z(S)\).  Therefore
\[
\dim(G\cdot S)=d(d-1)-(d-\dim Z(S)),
\]
where characteristic zero ensures that the stabilizer group scheme is smooth, so its dimension equals the dimension of its Lie algebra.  Equating the two displayed dimensions gives the claimed equivalence.  Finally, \(\dim\Der(A)=d-\dim Z(A)+\dim HH^1(A)\).
\end{proof}

\subsection{Symmetry, periodicity, and Tate duality}
\label{subsec:structural-reduction}

By the deformation constructed in Section~\ref{sec:general-existence}, the contraction algebra lies in the closure of the required semisimple orbit. The next results establish symmetry of the contraction algebra and four-periodicity, hence \(\dim HH^1(A)=\dim HH^2(A)\), and identify the centre with the dual cocenter. The tangent-dimension calculation below uses these two equalities.

\begin{proposition}
\label{prop:contraction-algebra-symmetric}
The contraction algebra is symmetric: \(\Acon\cong D\Acon\) as \(\Acon\)-bimodules.
\end{proposition}

\begin{proof}
Let \(N\) be the nonfree maximal Cohen--Macaulay Van den Bergh summand on the completed base \(R\). The base \(R\) is a three-dimensional isolated compound Du Val (cDV) singularity, and the contraction algebra is the opposite of its stable endomorphism ring, \(A=\underline{\End}_R(N)^{\mathrm{op}}\). The stable endomorphism algebra is symmetric by \cite[Theorem~2.6]{August}, and the opposite of a symmetric algebra is symmetric.
\end{proof}

\begin{proposition}
\label{prop:bimodule-periodicity}
In the stable category of \(A\)-bimodules,
\[
 \Omega^4_{A^e}(A)\simeq A.
\]
Here \(\Omega_{A^e}\) denotes the syzygy functor.
\end{proposition}

\begin{proof}
Put \(\Lambda=\underline{\End}_R(N)\), so that \(A=\Lambda^{\mathrm{op}}\) in our convention. Jasso--Keller--Muro's Theorem~3.2.5 applies because \(R\) is an isolated cDV singularity with a crepant resolution.  Their equation~(3.2.4) identifies the coefficient bimodule with the diagonal \(\Lambda\)-bimodule, and their theorem gives an exact sequence with projective middle terms and \(\Lambda\) at both ends.

To return to the convention for the contraction algebra, exchange the left and right actions. This gives the exact sequence of \(A\)-bimodules
\[
 0\longrightarrow A\longrightarrow P_3\longrightarrow P_2\longrightarrow
 P_1\longrightarrow P_0\longrightarrow A\longrightarrow0
\]
with projective middle terms \cite[Theorem~3.2.5]{JKM}.  The leftmost term of this exact sequence is a fourth syzygy of \(A\), which gives the stated stable bimodule isomorphism.

\end{proof}

\begin{lemma}
\label{lem:functorial-hypersurface-periodicity}
Let \(Q\) be a Noetherian regular local ring of finite dimension, let \(f\in Q\) be a nonunit nonzerodivisor, and put \(R=Q/(f)\). On the ungraded stable category of finitely generated maximal Cohen--Macaulay \(R\)-modules there is a natural isomorphism
\[
 \Omega_R^2\simeq\mathrm{Id}.
\]
\end{lemma}

\begin{proof}
For a nonzero maximal Cohen--Macaulay \(R\)-module \(M\), the depth over the ambient regular local ring is \(\operatorname{depth}_Q M=\dim Q-1\). The resolution bound over the regular local ring \(Q\) \cite[Proposition~10.110.1, Tag~\href{https://stacks.math.columbia.edu/tag/00O7}{00O7}]{Stacks} yields an exact sequence
\[
 0\longrightarrow F_{1,M}\xrightarrow{\varphi_M}F_{0,M}
 \longrightarrow M\longrightarrow0
\]
with finite free \(Q\)-modules and injective \(\varphi_M\).

Because \(fM=0\), multiplication by \(f\) factors uniquely through \(\varphi_M\), giving a homomorphism \(\psi_M:F_{0,M}\to F_{1,M}\) with
\[
 \varphi_M\psi_M=fI,\qquad \psi_M\varphi_M=fI.
\]
The second identity follows by applying the injective \(\varphi_M\). The zero module is handled by zero free modules.

Write a bar for reduction modulo \(f\). These reduced maps assemble into an exact periodic resolution over \(R\). Indeed, if \(\varphi_M(x)=fy=\varphi_M\psi_M(y)\), injectivity gives \(x=\psi_M(y)\), so \(\ker\bar\varphi_M=\operatorname{im}\bar\psi_M\). The other kernel identity \(\ker\bar\psi_M=\operatorname{im}\bar\varphi_M\) follows by the symmetric argument using injectivity of \(\psi_M\), which follows from \(\varphi_M\psi_M=fI\). Consequently, the assignment
\[
 \vartheta_M:M\xrightarrow{\sim}
 K_2(M):=\ker\bar\varphi_M,
 \qquad [x]\longmapsto\bar\psi_M(x),
\]
identifies \(M\) with a second syzygy.

For an arbitrary homomorphism \(a:M\to N\), lift it along the presentations to \(a_0:F_{0,M}\to F_{0,N}\) and then choose \(a_1\) satisfying \(\varphi_Na_1=a_0\varphi_M\). Injectivity of \(\varphi_N\) and the common equation \(f\) give
\[
 a_1\psi_M=\psi_Na_0.
\]
Thus \(K_2(a)\vartheta_M=\vartheta_Na\), where \(K_2(a)\) is the restriction of \(\bar a_1\) to the second kernel. Modifying the initial lift by \(\varphi_Nb\) alters the lifted second map by \(b\varphi_M\), whose reduction vanishes identically on \(K_2(M)\). This proves independence of the chosen lifts, as well as functorial compatibility with addition, identity maps, and composition.

Any homomorphism factoring through a free module continues to factor through a free module, since \(K_2(P)\cong P\) for free \(P\), so these assignments and comparison isomorphisms descend to the stable category.

It remains to compare the functor \(K_2\) with the syzygy functor \(\Omega_R^2\), formed using arbitrary free presentations, in the stable category. Write \(d_{1,M}\) and \(d_{1,N}\) for their degree-one differentials. Lift the identity morphism in both directions between the two presentations. If two chain lifts induce the same module map, their differences satisfy \(\delta_0=d_{1,N}h_0\). Then the difference \(\delta_1-h_0d_{1,M}\) takes values in the target second kernel. Its restriction to the source second kernel factors through the source degree-one free module. Thus the induced differences vanish stably.

The comparison maps are consequently inverse stable isomorphisms and are natural in every module homomorphism. Composing them with \(\vartheta_M\) proves the asserted natural isomorphism, including its restriction to any stable additive subcategory.
\end{proof}

\paragraph{Alternative cluster-tilting argument.}
Dugas's theorem provides a second route to four-periodicity from a cluster-tilting object in a Frobenius category.  Let \(\mathcal B=\operatorname{CM}(R)\), whose projective-injective objects are \(\operatorname{add}R\).  To identify the cluster-tilting object in this category, write the Van den Bergh tilting bundle on the smooth source \(X\) as \(\mathcal O_X\oplus\mathcal N\), and put \(N=f_*\mathcal N\).  The nonfree summand in \cite[Theorem~1.3.1]{JKM} is a \(2\mathbb Z\)-cluster-tilting object of \(\underline{\operatorname{CM}}(R)\).  This property also holds for its \(R\)-dual: duality is an exact contravariant equivalence on maximal Cohen--Macaulay modules, interchanges the two \(\operatorname{Ext}^1\)-orthogonality conditions, and interchanges left and right approximations.  Thus it applies to our perverse convention, with \(\underline{\End}_R(N)=A^{\mathrm{op}}\). Adding the projective-injective summand \(R\) gives the \(2\)-cluster-tilting object \(R\oplus N\) in the Frobenius category \(\mathcal B\).

Write the completed cDV hypersurface as \(R=\C[[x_1,x_2,x_3,x_4]]/(f)\). Lemma~\ref{lem:functorial-hypersurface-periodicity} gives \(\Omega_R^2\simeq\mathrm{Id}\) naturally on its stable category, and hence on \(\operatorname{add}(R\oplus N)\) modulo projectives. In \cite[Theorem~1.1(3)]{Dugas}, we may therefore take \(d_{\mathrm D}=2\) and \(r_{\mathrm D}=1\).  If \(A\ne\C\), the one-curve contraction algebra is local and nonsemisimple, so it has no semisimple block, and Dugas's theorem gives a periodic projective bimodule resolution of \(A^{\mathrm{op}}\), hence of \(A\), with period dividing \((d_{\mathrm D}+2)r_{\mathrm D}=4\).  If \(A=\C\), its stable bimodule class is zero and its positive Hochschild groups vanish.

\begin{lemma}[Tate--Hochschild duality]
\label{lem:tate-hochschild-duality}
For every symmetric finite-dimensional algebra \(A\) over \(\C\),
\[
 \bigl(\widehat{HH}^{-n}(A)\bigr)^\vee
 \cong\widehat{HH}^{n-1}(A)
 \qquad(n\in\mathbb Z).
\]
If \(A\) is four-periodic in the stable bimodule category, then
\begin{equation}
\label{eq:hh1-hh2-duality}
 DHH^1(A)\cong HH^2(A),
 \qquad
 \dim HH^1(A)=\dim HH^2(A).
\end{equation}
\end{lemma}

\begin{proof}
For a symmetric algebra \(A\), the enveloping algebra \(A^e\) is symmetric, and \(\widehat{HH}^*(A)=\widehat{\operatorname{Ext}}^*_{A^e}(A,A)\). The displayed duality is the Tate--Hochschild duality of \cite[equation~(2.11)]{Linckelmann}.  Taking \(n=3\) and using the period-four bimodule isomorphism gives
\[
 DHH^1(A)\cong D\widehat{HH}^{-3}(A)
 \cong\widehat{HH}^{2}(A)=HH^2(A),
\]
where period four identifies Tate degrees \(-3\) and \(1\), and positive Tate degrees agree with ordinary Hochschild cohomology.
\end{proof}

Etingof--Latour--Rains identify the centre of their central extensions with the dual of the cocenter \cite[Corollary~2.2]{EtingofLatourRains}.  The following formulation for symmetric algebras is the one needed here and later for relative traces.

\begin{lemma}
\label{lem:centre-cocenter-duality}
If \(A\) is finite-dimensional symmetric with symmetrizing trace \(\tau\colon A\to\C\), then
\[
 Z(A)\xrightarrow{\sim}HH_0(A)^\vee,
 \qquad z\mapsto\bigl(\bar a\mapsto\tau(za)\bigr).
\]
Consequently,
\begin{equation}
\label{eq:centre-cocenter-dimension}
 \dim Z(A)=\dim HH_0(A).
\end{equation}
\end{lemma}

\begin{proof}
For a central element \(z\),
\[
 \tau(z[a,b])=\tau(zab-zba)=\tau(zab-azb)=0,
\]
so the functional \(a\mapsto\tau(za)\) factors through \(A/[A,A]\). 

Conversely, identify \(A\cong DA\) by the symmetrizing form.  If the functional associated to \(c\in A\) annihilates every commutator, then for all \(a,b\),
\[
 0=\tau(c[a,b])=\tau((ca-ac)b).
\]
Nondegeneracy forces \(ca=ac\) for all \(a\), so \(c\in Z(A)\).  The annihilator of the commutator subspace is therefore exactly the centre, which proves the claimed isomorphism.
\end{proof}

\Needspace{25\baselineskip}
\subsection{Reduction to the cocenter dimension}
\label{subsec:tangent-orbit-dimensions}

For \(G=\Stab_{\GL_d}(1)\), Lemma~\ref{lem:fixed-unit-tangents} and Corollary~\ref{cor:hochschild-orbit-criterion} give the tangent and orbit formulas below.  The second equality uses \(\dim\Der(A)=d-\dim Z(A)+\dim HH^1(A)\). By symmetry and four-periodicity, Lemma~\ref{lem:tate-hochschild-duality} gives equal dimensions for \(HH^1(A)\) and \(HH^2(A)\), so these terms cancel.

The last equality uses \(\dim Z(S_{\mathrm{GV}})=N_{\mathrm{GV}}\) and the fact that all derivations of a semisimple complex algebra are inner:
\begin{align}
 \dim T_AX_d
 &=d(d-1)-\dim\Der(A)+\dim HH^2(A)
 \label{eq:tangent-dimension-raw}\\
 &=d(d-1)-d+\dim Z(A),
 \label{eq:tangent-dimension-final}\\
 \dim(G\cdot S_{\mathrm{GV}})
 &=d(d-1)-d+N_{\mathrm{GV}}.
 \label{eq:semisimple-orbit-dimension}
\end{align}

Subtracting and applying Lemma~\ref{lem:centre-cocenter-duality} gives
\begin{equation}
\label{eq:tangent-orbit-difference}
 \dim T_AX_d-\dim(G\cdot S_{\mathrm{GV}})=\dim HH_0(A)-N_{\mathrm{GV}}.
\end{equation}

Theorem~\ref{thm:general-semisimple-existence} and Lemma~\ref{lem:family-orbit-closure} already give \(A\in\overline{G\cdot S_{\mathrm{GV}}}\).  Thus the only remaining hypothesis of Theorem~\ref{thm:orbit-smooth-unique} is that the right-hand side of \eqref{eq:tangent-orbit-difference} vanish.  This reduces Theorem~\ref{thm:general-uniqueness} to the cocenter count \(\dim HH_0(A)=N_{\mathrm{GV}}\).

\subsection{The cocenter formula: statements and reduction}
\label{subsec:cocenter-statements-reduction}

\begin{corollary}
\label{cor:contraction-cocenter-formula}
Let \(A\) be the contraction algebra of a completed simple one-curve flop of a smooth complex threefold, of length \(\ell\), and let \(n_j\) be the genus-zero GV invariants of that flop. Then
\[
 \dim_\C HH_0(A)=\dim_\C Z(A)=N_{\mathrm{GV}}=\sum_{j=1}^{\ell}n_j.
\]
\end{corollary}

We compute the cocenter of \(A\) by pulling back a universal square presentation from the marked base to the original classifying arc. Once the restricted quotient is identified with \(A\) along that arc, Smith normal form expresses the cocenter dimension directly as the order of vanishing of the determinant of the restricted matrix.

\subsubsection{The universal presentation and the arc comparison}
\label{subsubsec:marked-quotient-downstairs}

Use the primitive marked pair \((\Phi,p)\) attached to the flop in Section~\ref{sec:general-existence}. Here \(p\) is its retained vertex. The marked base is \(H=\mathfrak h/W_0\), where \(W_0=\langle s_i:i\ne p\rangle\). Write \(\mathbb H_\ell=\mathcal O(H)\) for the ring of regular functions on \(H\), and \(R_H=\widehat{\mathcal O}_{H,0}=\widehat{\mathbb H}_\ell\). Let \(\widehat{\mathcal A}_\ell\) be the universal endomorphism algebra of Subsection~\ref{subsec:full-universal-tilting}, completed at the central graded point of its affine target, with its free projector \(e_0\) retained. Put
\[
 \mathcal C_H=
 \widehat{\mathcal A}_\ell/
 \widehat{\mathcal A}_\ell e_0\widehat{\mathcal A}_\ell,
 \qquad \mathsf M_H=HH_0(\mathcal C_H/R_H).
\]
For \(1\le j\le\ell\), let \(P_j\in R_H\) be a reduced equation of the curve divisor \(D_j\), corresponding to curves of class \(j[C]\).

\Needspace{13\baselineskip}
\begin{theorem}
\label{thm:marked-cocenter}
The algebra \(\mathcal C_H\) is finite over \(R_H\). Its cocenter has a free resolution
\begin{equation}
\label{eq:downstairs-presentation}
 0\longrightarrow R_H^q\xrightarrow{\Theta}R_H^q
 \longrightarrow\mathsf M_H\longrightarrow0
\end{equation}
with
\[
 \operatorname{Fitt}^{R_H}_0(\mathsf M_H)
 =(\det\Theta)=(P_1\cdots P_\ell),
\]
where \(P_j=0\) is a reduced equation of the curve divisor \(D_j\). Thus \(\det\Theta=\eta P_1\cdots P_\ell\) for some \(\eta\in R_H^\times\). Each \(P_j\) occurs with exponent \(1\).
\end{theorem}

We prove this presentation in Subsection~\ref{subsec:cocenter-descent-proof}. To see why its determinant matters, first compute its pullback along the original classifying arc.

Let \(R_1=\C[[w]]\), and let \(f_0:\Spec R_1\to H\) be the original classifying arc. Set
\[
 \mathcal C_0=\mathcal C_H\otimes_{R_H,f_0}R_1,\qquad
 \mathsf M_0=HH_0(\mathcal C_0/R_1).
\]
We also need to identify \(\mathcal C_0\) with \(A\):
\[
 \mathcal C_H\otimes_{R_H,f_0}R_1=\mathcal C_0\cong A.
\]
More precisely, Lemma~\ref{lem:arc-quotient-supported} shows that the ordinary pullback \(\mathcal C_0\) is already complete at the contraction point, so the natural map from \(\mathcal C_0\) to its completion is an isomorphism. The comparison of the fixed basic generators identifies the pulled-back endomorphism algebra with the original basic NCCR and preserves the projector \(e_0\) onto the free summand. Taking the quotient by the ideal generated by \(e_0\) gives \(\mathcal C_0\cong A\). These two steps are proved in Subsection~\ref{subsubsec:arc-vs-contraction}.

Lemma~\ref{lem:arc-not-in-discriminant} shows that the original arc lies in no curve divisor, so every \(P_j(f_0(w))\) is nonzero. Proposition~\ref{prop:conserved-multiplicities-gv} identifies its finite order of vanishing with \(n_j\).

\subsubsection{The determinant--length reduction}
\label{subsubsec:pullback-cocenter}

\begin{proof}[Proof of Corollary~\ref{cor:contraction-cocenter-formula} from the stated inputs]
Assume Theorem~\ref{thm:marked-cocenter} and Lemma~\ref{lem:arc-quotient-supported}, whose proofs follow below. By the lemma, the algebra on the original arc is \(\mathcal C_0\cong A\). The cocenter is the cokernel of the commutator map, so it commutes with ordinary base change, without any flatness assumption. Hence \(\mathsf M_0=\mathsf M_H\otimes_{R_H}R_1\). The pulled-back square matrix presents \(\mathsf M_0\), which has finite length over the DVR \(R_1\). Write \(\Theta_0=f_0^*\Theta\). Its determinant is
\begin{equation}
\label{eq:arc-determinant}
 \det\Theta_0=\eta(w)\prod_{j=1}^{\ell}P_j(f_0(w))\ne0.
\end{equation}
Here \(\eta(w)\in R_1^\times\) is the pullback of the unit in the universal determinant factorization. Right exactness of tensor products ensures that the pulled-back presentation computes the cokernel over the arc. Because its determinant is nonzero in the DVR \(R_1\), the square matrix is injective. Even though the classifying map need not be flat, we obtain the exact sequence
\begin{equation}
\label{eq:arc-presentation-exact}
 0\longrightarrow R_1^q\xrightarrow{\Theta_0}R_1^q
 \longrightarrow\mathsf M_0\longrightarrow0.
\end{equation}

Applying the Smith normal form over the discrete valuation ring \(R_1\) yields
\begin{equation}
\label{eq:smith-length}
 \dim_\C\mathsf M_0
 =\operatorname{length}_{R_1}\mathsf M_0
 =\operatorname{ord}_w\det\Theta_0.
\end{equation}

The relative and absolute cocenters share the same underlying \(\C\)-vector space, since \(r[a,b]=[ra,b]\) for central \(r\). Lemma~\ref{lem:arc-quotient-supported} identifies this cocenter with \(HH_0(A)\). The factorization \eqref{eq:arc-determinant} expresses its length as the sum of the vanishing orders along the curve divisors, which Proposition~\ref{prop:conserved-multiplicities-gv} identifies with the GV multiplicities. Thus
\begin{equation}
\label{eq:cocenter-dimension-final}
 \dim_\C HH_0(A)=\sum_j\operatorname{ord}_wP_j(f_0(w))=\sum_j n_j=N_{\mathrm{GV}}.
\end{equation}
Because \(A\) is a finite-dimensional symmetric algebra, centre--cocenter duality identifies \(Z(A)\) with the vector space dual of \(HH_0(A)\), which yields
\begin{equation}
\label{eq:centre-dimension-final}
 \dim_\C Z(A)=N_{\mathrm{GV}}.
\end{equation}
\end{proof}

Subsection~\ref{subsec:universal-marked-cocenter} constructs the cocenter presentation over \(\widehat{\mathcal O}_{\mathfrak h,0}\). Subsection~\ref{subsec:descent-cartan-cover} compares the algebras with the geometric contraction quotients. Finally, Subsection~\ref{subsec:cocenter-descent-proof} descends the presentation to \(R_H\).

\subsection{A cocenter presentation over \texorpdfstring{\(\widehat{\mathcal O}_{\mathfrak h,0}\)}{the completed local ring of h}}
\label{subsec:universal-marked-cocenter}

Keep the primitive marked ADE pair \((\Phi,p)\) attached to the flop.  Here \(p\) is the retained vertex \(\alpha_*\) of Section~\ref{sec:general-existence}, so \(m_p(\alpha)=c_*(\alpha)\). Recall that \(\mathfrak h\) is the Cartan space of rank \(r\), and put
\[
 W_0=\langle s_i:i\ne p\rangle,
 \qquad H=\mathfrak h/W_0.
\]
The completed coefficient ring at the origin is
\[
 \widetilde R=\widehat{\mathcal O}_{\mathfrak h,0}
 \cong\C[[\lambda_1,\ldots,\lambda_r]].
\]
Under \(\mathfrak h\to H\), each reduced curve-divisor equation pulls back to a product of linear root equations, up to a unit. We construct the cocenter presentation over \(\widetilde R\), then descend it to \(R_H=\widehat{\mathcal O}_{H,0}\).

For a positive root \(\alpha=\sum_i m_i(\alpha)\alpha_i\), put
\begin{equation}
\label{eq:N_p-definition}
 N_p:=\#\{\alpha\in\Phi^+:m_p(\alpha)>0\}.
\end{equation}

Choose an orientation \(Q\) of the finite Dynkin diagram, and let \(\overline Q\) be its double.  Over \(\widetilde R\), define the deformed preprojective algebra
\[
 \Pi_\lambda=
 \widetilde R\,\overline Q/
 \left(\sum_{a\in Q_1}[a,a^*]-\sum_i\lambda_i e_i\right).
\]
We use the convention of \cite[Proposition~4.15]{Karmazyn} and write \((\lambda,\alpha)=\sum_i m_i(\alpha)\lambda_i\). Let
\[
 \mathcal D_p=e_p\Pi_\lambda e_p.
\]

\phantomsection
\label{subsubsec:etingof-rains-parameters}
To construct the cocenter presentation, we use the central extension of \citet{EtingofRains}. Take \(\mu=\rho\), normalized by \(\rho_i=(\rho,\alpha_i)=1\). Let \(\Pi\) be the universal deformed preprojective algebra with central parameters \(x_i\) and relation \(\sum_a[a,a^*]=\sum_i x_i e_i\). To form the central extension, adjoin a central polynomial variable \(z\). The coefficient ring for the regular weight \(\rho\) is
\[
 K_\rho=\C[x_1-\rho_1z,\ldots,x_r-\rho_rz].
\]
The fibre at a parameter \(\lambda\) is the algebra \(\Pi^\rho_\lambda\) with relations
\[
 \sum_a[a,a^*]=\sum_i(\rho_i z+\lambda_i)e_i.
\]
Via the coordinates \(\lambda_i:=x_i-\rho_iz\), completing \(K_\rho\) at its augmentation gives the formal power-series ring \(\C[[\lambda_1,\ldots,\lambda_r]]\), namely \(\widetilde R=\widehat{\mathcal O}_{\mathfrak h,0}\).  The \(\lambda_i\) are base coordinates, whereas \(z\) remains a central element of the algebra.  Etingof--Rains prove that \(\Pi[z]\) is finite free over \(K_\rho\) \cite[Theorem~3.2(i)]{EtingofRains}. We now complete this coefficient ring and take the corner at the retained vertex.

Define the marked corner of the completed universal central extension by
\begin{equation}
\label{eq:central-extension-corner}
 E_p:=e_p\bigl(\Pi[z]\widehat\otimes_{K_\rho}\widetilde R\bigr)e_p,
 \qquad \mathsf M_p:=HH_0(E_p/\widetilde R).
\end{equation}
Since the Etingof--Rains central extension is finite free over \(K_\rho\), the corner \(E_p\) is finite over \(\widetilde R\).  Its cocenter \(\mathsf M_p\), being a quotient of \(E_p\), is finite as well.

At the origin of the coefficient space, the closed fibre of \(E_p\) is \(e_p\Pi^\rho_0e_p\).  Since \(HH_0\) is the cokernel of the commutator map, it commutes with arbitrary base change.  In particular,
\[
 \mathsf M_p\otimes_{\widetilde R}\C
 \cong HH_0(e_p\Pi^\rho_0e_p).
\]

\begin{theorem}
\label{thm:universal-marked-cocenter-upstairs}
For every primitive marked ADE pair \((\Phi,p)\), there is an exact sequence
\begin{equation}
\label{eq:universal-cocenter-presentation}
 0\longrightarrow\widetilde R^{N_p}\xrightarrow{\Psi_p}\widetilde R^{N_p}
 \longrightarrow HH_0(\mathcal D_p/\widetilde R)\longrightarrow0
\end{equation}
such that
\begin{equation}
\label{eq:universal-cocenter-determinant}
 \det\Psi_p=\eta\Delta_p,
 \qquad
 \Delta_p:=\prod_{\substack{\alpha\in\Phi^+\\m_p(\alpha)>0}}(\lambda,\alpha),
 \qquad\eta\in\widetilde R^\times.
\end{equation}
\end{theorem}

To prove Theorem~\ref{thm:universal-marked-cocenter-upstairs}, we show that the cocenter of the central extension corner, \(\mathsf M_p=HH_0(E_p/\widetilde R)\), is free of rank \(N_p\). In a basis of this free module, multiplication by \(z\) has matrix \(\Psi_p\), with cokernel \(HH_0(\mathcal D_p/\widetilde R)\). The next lemma describes the deformed preprojective algebras at single-root parameters. These are the factors used to compute the generic rank.

\begin{lemma}
\label{lem:single-root-preprojective-fibre}
Let \(k\) be algebraically closed of characteristic zero and let \(\eta\in k^r\) annihilate exactly one positive root \(\alpha\). Then the ordinary deformed preprojective algebra \(\Pi_\eta\) of the finite Dynkin diagram is \(M_{\operatorname{ht}(\alpha)}(k)\), and its vertex projector \(e_i\) has rank \(m_i(\alpha)\).
\end{lemma}

\begin{proof}
The path-length filtration makes \(\operatorname{gr}\Pi_\eta\) a quotient of the finite-dimensional ordinary preprojective algebra, so \(\Pi_\eta\) is finite-dimensional \cite[Section~2.1]{EtingofRains}. Over the chosen parameter, the simple-module criterion \cite[Theorem~1.2]{CrawleyBoeveyMomentMap} shows that every simple module has dimension vector \(\alpha\), and that such a simple exists: there is no decomposition of \(\alpha\) into two or more positive roots annihilated by \(\eta\).

For two representations \(U,V\), splitting an extension on the vertex spaces gives the complex
\[
 \bigoplus_i\operatorname{Hom}_k(U_i,V_i)
 \xrightarrow{d^0}
 \bigoplus_{a\in\overline Q_1}
       \operatorname{Hom}_k(U_{\mathrm{tail}(a)},V_{\mathrm{head}(a)})
 \xrightarrow{d^1}
 \bigoplus_i\operatorname{Hom}_k(U_i,V_i).
\]
Here \(d^0\) changes the splitting and \(d^1\) is the linearized preprojective relation; the middle cohomology of this complex is \(\operatorname{Ext}^1_{\Pi_\eta}(U,V)\). The trace pairings identify \(d^1(U,V)\) with the transpose, up to sign, of \(d^0(V,U)\). Taking dimensions therefore gives
\[
 \dim\operatorname{Ext}^1(U,V)=\dim\operatorname{Hom}(U,V)+\dim\operatorname{Hom}(V,U)
                 -(\dim U,\dim V),
\]
where the last term is the symmetric Cartan pairing.

Since \((\alpha,\alpha)=2\), two nonisomorphic simple modules of dimension vector \(\alpha\) would have zero homomorphisms in both directions, yielding a negative extension dimension.

Thus the simple \(S\) is unique, and the same formula gives \(\operatorname{Ext}^1(S,S)=0\). Its indecomposable projective cover cannot have nonzero radical. Indeed, if that radical were nonzero, it would admit a quotient isomorphic to \(S\), giving a nonzero self-extension of \(S\). Hence every finite-length module is semisimple. It follows that \(\Pi_\eta=\End_k(S)\), with the stated matrix size and projector ranks.
\end{proof}

\phantomsection
\label{subsubsec:generic-rank}

We now compute the generic rank of \(\mathsf M_p\). Etingof--Rains prove generic semisimplicity in \cite[Theorem~3.2(v) and the display immediately following it]{EtingofRains}. In their decomposition a positive root labels each factor, and that factor is an ordinary deformed preprojective algebra.  The Crawley--Boevey criterion supplies a simple module of dimension vector \(\alpha\) \cite[Theorem~1.2]{CrawleyBoeveyMomentMap}. Combined with the Etingof--Rains dimension calculation, this criterion identifies these factors as matrix algebras.

After passing to a splitting field \(L/\operatorname{Frac}(\widetilde R)\),
\begin{equation}
\label{eq:generic-extension-decomposition}
 (E_p)_L\cong\prod_{m_p(\alpha)>0}M_{m_p(\alpha)}(L).
\end{equation}
Taking the cocenter counts the surviving factors: each matrix algebra contributes one dimension. Thus
\begin{equation}
\label{eq:generic-cocenter-rank}
 \dim_L(\mathsf M_p)_L=N_p.
\end{equation}

\begin{lemma}
\label{lem:closed-fibre-cocenter-count}
For every primitive marked ADE pair \((\Phi,p)\),
\begin{equation}
\label{eq:special-cocenter-dimension}
 \dim_{\C}HH_0(e_p\Pi^\rho_0e_p)=N_p.
\end{equation}
\end{lemma}

\begin{proof}
At \(\lambda=0\), the normalization \(\rho_i=1\) for all simple roots gives the central-extension relation used by Etingof--Latour--Rains (ELR) \cite{EtingofLatourRains}.  They compute the cocenters at nodal vertices for \(\mu=\rho\) \cite[Section~3.2]{EtingofLatourRains}.  To apply that calculation to the six primitive flop markings, we identify their retained vertices and deleted Levi diagrams using Karmazyn's presentations below.

ELR call the branching vertex of a Dynkin diagram of type \(D\) or \(E\) the \emph{nodal vertex}.  Karmazyn's presentations \cite[Examples~4.22--4.27]{Karmazyn} give the following retained vertices and Levi subgroups, using his vertex numbering:

\begin{center}
\begin{tabular}{r|c|c|c|c|l}
\toprule
length & Dynkin & retained & \(W_0\) in & deleted & ELR \\
 & type & vertex & Karmazyn & diagram & status \\
\midrule
1 & \(A_1\) & unique & trivial & empty & immediate \\
2 & \(D_4\) & \(p=4\) & \(S_2^3\) & \(A_1^3\) & nodal \(D\) \\
3 & \(E_6\) & \(p=6\) & \(S_2\times S_3^2\) & \(A_1\sqcup A_2^2\) & nodal \(E_6\) \\
4 & \(E_7\) & \(p=7\) & \(S_3\times S_2\times S_4\) & \(A_2\sqcup A_1\sqcup A_3\) & nodal \(E_7\) \\
5 & \(E_8\) & \(p=4\) & \(S_4\times S_5\) & \(A_3\sqcup A_4\) & non-nodal \\
6 & \(E_8\) & \(p=8\) & \(S_5\times S_2\times S_3\) & \(A_4\sqcup A_1\sqcup A_2\) & nodal \(E_8\) \\
\bottomrule
\end{tabular}
\end{center}

For a positive root \(\alpha\), ELR's phrase ``the root contains the nodal simple root'' is equivalent to \(m_p(\alpha)>0\).  Hence their root count is \(N_p\).

For type \(D\), ELR \cite[Section~3, Lemma~3.4]{EtingofLatourRains} give degreewise upper bounds for the cocenter of the nodal corner.  Summing them gives
\[
 \dim_{\C}HH_0(e_p\Pi^\rho_0e_p)\le N_p.
\]
The finite \(\widetilde R\)-module \(\mathsf M_p\) has generic rank \(N_p\) by \eqref{eq:generic-cocenter-rank}.  For a finite module over the local domain \(\widetilde R\), the closed-fibre dimension is at least the generic rank.  This supplies the reverse inequality and hence equality.

For the nodal \(E_6,E_7,E_8\) corners, ELR \cite[Section~3, Case~2b]{EtingofLatourRains} compute the same graded cocenter dimensions by exact Magma calculations and obtain equality with the positive-root count.  These are the primitive flop lengths \(3,4,6\) in the table.

For length one (type \(A_1\)), the closed corner \(e_p\Pi^\rho_0e_p=\C[z]/(z)\) has one-dimensional cocenter. 

For the remaining length-five case, the retained vertex of \(E_8\) is not nodal, so ELR's nodal calculation does not apply. Appendix~\ref{app:e8-executable-certificate} proves \eqref{eq:e8-equality} by combining the explicit relations, modular-rank upper bounds, a bound on the remaining degrees and the generic-rank lower bound. This gives the required dimension \(104=N_p\) and completes the six primitive cases.
\end{proof}

There is also an independent structural proof of this equality: Theorem~\ref{thm:all-corner-trace-excision}, applied to the closed \(\rho\)-central extension with retained set \(\{p\}\), gives \(\dim HH_0(e_p\Pi^\rho_0e_p)=N_p\), since the roots of the deleted subdiagram are exactly those with marked coefficient zero.  That argument takes place in \(\Pi^\rho_0\) and does not use the universal presentation constructed here.  In particular it includes the non-nodal \(E_8\) marking of length five.

\begin{proof}[Proof of Theorem~\ref{thm:universal-marked-cocenter-upstairs}]
\leavevmode
\begin{romanenu}
\item \emph{Freeness of the auxiliary cocenter module.}
\label{subsubsec:freeness-cp}

Lemma~\ref{lem:closed-fibre-cocenter-count} establishes that the closed-fibre dimension is \(N_p\). Lifting a basis of the closed fibre to the auxiliary cocenter module \(\mathsf M_p\) and invoking Nakayama's lemma yields a surjection
\[
 \widetilde R^{N_p}\twoheadrightarrow \mathsf M_p.
\]
After tensoring with the fraction field \(\operatorname{Frac}(\widetilde R)\), the generic rank identity \eqref{eq:generic-cocenter-rank} shows that this surjection is an isomorphism. Its kernel is a submodule of a free module over the integral domain \(\widetilde R\), hence is torsion-free. Because its generic fibre vanishes, the kernel is zero. Therefore
\begin{equation}
\label{eq:mp-free}
 \mathsf M_p\cong\widetilde R^{N_p}.
\end{equation}

\item \emph{Multiplication by the central element.}
\label{subsubsec:multiplication-by-z}

On the free \(\widetilde R\)-module in \eqref{eq:mp-free}, multiplication by the central element \(z\) defines an endomorphism
\[
 \mu_z\colon\mathsf M_p\to\mathsf M_p.
\]
The quotient algebra of \(E_p\) modulo the ideal generated by \(z\) is precisely the ordinary deformed preprojective corner \(\mathcal D_p\). Since \(z\) is central, passage to its quotient commutes with forming the cocenter, giving
\begin{equation}
\label{eq:cocenter-mod-z}
 HH_0(\mathcal D_p/\widetilde R)\cong\mathsf M_p/z\mathsf M_p
 =\operatorname{coker}\mu_z.
\end{equation}

To evaluate the determinant, let \(V_\alpha\) denote the simple module of the centrally extended algebra indexed by the positive root \(\alpha\). The decomposition displayed immediately following \cite[Theorem~3.2(v)]{EtingofRains} gives the scalar action
\[
 z|_{V_\alpha}=-\frac{(\lambda,\alpha)}{(\rho,\alpha)}.
\]
The same scalar acts on the retained corner of this factor and on its one-dimensional cocenter. Therefore \(\mu_z\) has these eigenvalues, and over the fraction field its determinant factors as
\[
 \det\mu_z=\prod_{m_p(\alpha)>0}
 \left(-\frac{(\lambda,\alpha)}{(\rho,\alpha)}\right).
\]
The denominators are nonzero complex constants, which yields
\[
 \det\mu_z=\eta\Delta_p,
 \qquad\eta\in\widetilde R^\times.
\]

Because this determinant is nonzero, the endomorphism \(\mu_z\) is injective over the integral domain \(\widetilde R\). Together with \eqref{eq:mp-free}--\eqref{eq:cocenter-mod-z}, this gives the required presentation and determinant.\qedhere
\end{romanenu}
\end{proof}

\subsection{Comparison with the geometric contraction algebras}
\label{subsec:descent-cartan-cover}

The algebra comparisons concern two base changes. Extending \(\mathcal C_H\) from \(R_H\) to \(\widetilde R\) gives the corner \(\mathcal D_p\), while pulling \(\mathcal C_H\) back along the original arc gives \(A\). We first prove that \(R_H\to\widetilde R\) is finite and faithfully flat, as needed to compare the completed generators.

Write
\[
 \mathbb H_\Gamma=\mathcal O(\mathfrak h),\qquad
 \mathbb H_\ell=\mathbb H_\Gamma^{W_0},\qquad
 R_H=\widehat{\mathbb H}_\ell,\qquad
 \widetilde R=\widehat{\mathbb H}_\Gamma,
\]
with completions at the graded origins.  The extension \(R_H\to\widetilde R\) is finite free and faithfully flat.

For the primitive markings in the table above, every component of the deleted Levi diagram is of type \(A\). We prove freeness on each component.  The representation of \(W_0\) on \(\mathfrak h\) is a fixed line together with sum-zero permutation representations of \(S_q\), one for each component \(A_{q-1}\).  On such a block the invariant ring is \(\C[e_2,\ldots,e_q]\), where \(e_i\) are the elementary symmetric polynomials in \(y_1,\ldots,y_q\) with \(\sum_i y_i=0\).  The coordinate ring \(\C[y_1,\ldots,y_q]/(\sum_i y_i)\) is free over this invariant ring, with basis
\[
 \{y_1^{a_1}\cdots y_q^{a_q}:0\le a_i<q-i+1\}.
\]
To obtain the displayed basis, start with a generic monic polynomial of degree \(q\), adjoin its first root, divide by the resulting linear factor, and adjoin the next root.  The successive monic quotients have degrees \(q,q-1,\ldots,1\), giving this basis over the coefficient ring. The resulting splitting algebra is the polynomial ring in the ordered roots, and setting their sum equal to zero preserves its free basis. The elementary symmetric polynomials generate the invariants by successive cancellation of the leading monomial of a symmetric polynomial.

Tensoring the block bases and adjoining the fixed coordinate proves freeness of \(\mathbb H_\Gamma\) over \(\mathbb H_\ell\) with rank \(|W_0|>0\).  In the fibre over \(0\in H\), each root coordinate is nilpotent. Thus the topology defined by the extended origin ideal of \(H\) agrees with that defined by the origin ideal of \(\mathfrak h\).  Completion preserves the displayed finite free basis, which proves the assertion for \(R_H\to\widetilde R\).

\subsubsection{The generators after extension to \texorpdfstring{\(\widetilde R\)}{R tilde}}
\label{subsubsec:cartan-comparison}

Let \(\mathcal A_\Gamma\) be Karmazyn's universal extended-Dynkin algebra and put \(e_{0p}=e_0+e_p\). Write \(\mathcal T_{\Gamma,0p}\) for the two-summand generator on the partial resolution over \(\mathfrak h\) associated to this corner \cite[Proposition~4.21]{Karmazyn}, and \(\mathcal T_\ell\) for the generator of Subsection~\ref{subsec:full-universal-tilting}. Completions below are taken at the central graded origins.

\begin{proposition}
\label{prop:formal-marked-comparison}
On the partial resolution obtained by extending to \(\widetilde R\) and completing at the central target point, there is an isomorphism
\begin{equation}
\label{eq:tilting-bundle-comparison}
 \theta:\widehat{\mathcal T}_\ell
 \widehat\otimes_{R_H}\widetilde R
 \xrightarrow{\sim}\widehat{\mathcal T}_{\Gamma,0p}
\end{equation}
preserving the structure-sheaf and retained-vertex summands.  It induces an isomorphism
\begin{equation}
\label{eq:algebra-comparison}
 \widehat{\mathcal A}_\ell\widehat\otimes_{R_H}\widetilde R
 \xrightarrow{\sim}\widehat{e_{0p}\mathcal A_\Gamma e_{0p}}
\end{equation}
sending \(e_0\) to \(e_0\) and \(e_p\) to \(e_p\).
\end{proposition}

\begin{proof}
The geometric comparison in \cite[Section~2.3 and the proof of Proposition~4.21]{Karmazyn} places both generators on the same partial-resolution family over \(\mathfrak h\), obtained from the marked family by pullback along \(\mathfrak h\to H\), or from the simultaneous resolution by contracting the unretained curves. The corresponding coefficient map is \(\mathbb H_\ell\to\mathbb H_\Gamma\).

For the completion comparison, write \(U\) for the affine target ring completed at the maximal ideal of its central point and put \(S=\widetilde R\), \(U_S=U\otimes_{R_H}S\). The finite free parameter extension makes \(U_S\) finite over \(U\), and the quotient
\[
 U_S/\mathfrak m_UU_S\cong S/\mathfrak m_HS
\]
is local Artinian. Hence \(U_S\) is a Noetherian local ring, with maximal ideal \(\mathfrak n_S\), and for some \(c\ge1\) the ideals satisfy
\[
 \mathfrak n_S^c\subseteq\mathfrak m_UU_S\subseteq\mathfrak n_S.
\]
These containments show that the ring is complete at the maximal ideal of the new target point. As a complete Noetherian local ring it is also complete for the parameter ideal \(\mathfrak m_SU_S\): compatible Cauchy representatives for the parameter topology form an \(\mathfrak n_S\)-adic Cauchy sequence, and each fixed power of the parameter ideal is closed in the maximal-ideal topology. These two ideal topologies need not coincide. The source is projective over \(U_S\), and both source and target remain flat over \(S\).

Fix the contraction morphism from the full resolution to the partial resolution
\[
 q:\mathcal Z_S\longrightarrow\mathcal Y_{\ell,S}
\]
for the same chosen simultaneous resolution and the same side of the partial resolution \cite[Section~2.4.2 and the proof of Theorem~2.6]{NabijouWemyss}. The normalized bundle \(\mathcal T_{\Gamma,0p}\) lives on this \emph{partial} resolution and is obtained by lifting its labelled central basic generator. The contraction quotient is formed using the projector onto the structure sheaf. We therefore identify and lift the two labelled summands separately.

On the completed central surface the tilting-bundle construction of \cite[Theorem~4.6 and its proof]{IKWY}, dualized to the zero-perverse convention, gives the labelled basic generator
\[
 \mathcal T_0=\mathcal O\oplus\mathcal T_{0,p}.
\]
The restriction of the fixed \(\mathcal T_\ell\) is this same basic generator. Indeed, its completed summand endomorphism rings are local by the basic normalization in \cite[Remark~3.5]{Karmazyn}, and Proposition~\ref{prop:ordinary-endomorphism-base-change} identifies the specialized endomorphism rings with their quotients by the parameter ideal. The summands therefore remain indecomposable. The positive nonfree summand \((\mathcal T'_\ell)^\vee\) has determinant degree \(1\) on the retained curve, so \(\mathcal T'_\ell\) has determinant degree \(-1\). This nonzero degree persists on the central curve and distinguishes it from \(\mathcal O\). The specialized bundle is a projective generator of the same zero-perverse heart, so its two summands are the indicated labelled basic projectives.

Let \(q_0\) be the restriction of \(q\) to the central surface. The same construction gives the labelled pullback
\[
 q_0^*\mathcal T_0\cong
 \mathcal T_{\mathrm{full},0,0}\oplus
 \mathcal T_{\mathrm{full},0,p}
\]
on the full central resolution, where \(\mathcal T_{\mathrm{full},0,i}\) denotes the summand labelled by \(i\) of its basic generator. The singularities contracted on the central surface are rational, so \(Rq_{0*}\mathcal O=\mathcal O\).

The projection formula below requires this identity over the parameter base. It persists over successive Artin thickenings of \(S\). Indeed, for a small extension, parameter flatness identifies the kernels of structure-sheaf reduction on source and target with their central structure sheaves tensored with the extension kernel. Pushing forward these exact sequences proves the identity and higher-direct-image vanishing on the next thickening. Induction and formal functions on the completed target model \cite[Tag~02OC]{Stacks} give \(Rq_*\mathcal O=\mathcal O\).

With the central labelled summands identified, lift their isomorphisms separately through the parameter thickenings. Each summand of \(\mathcal T_0\) has vanishing higher self-Ext. Choose an isomorphism for each pair of corresponding central summands. Across the parameter thickening from order \(n\) to order \(n+1\), the obstruction to lifting such an isomorphism lies in \(\operatorname{Ext}^1\) of that central summand with itself, tensored with \(\mathfrak m_{\widetilde R}^n/ \mathfrak m_{\widetilde R}^{n+1}\), and therefore vanishes. A lift remains invertible modulo the nilpotent thickening. Induction therefore gives compatible summand-wise isomorphisms on all parameter thickenings.

The completed affine target is also complete for the parameter ideal. Since the source is projective over that target, coherent formal existence and full faithfulness algebraize these compatible isomorphisms and their inverses \cite[Lemma~3.2]{KarmazynDeformations}, \cite[Tag~088F]{Stacks}. Their direct sum is \(\theta\), which preserves both projectors \(e_0\) and \(e_p\) because it was constructed separately on each labelled summand. The choices differ only by automorphisms of the two summands.

Write \(\mathcal T_{\ell,S}\) for the pullback of \(\widehat{\mathcal T}_\ell\) to the partial resolution over \(S\). On the full family, apply the same compatible lifting argument to \(q^*\mathcal T_{\ell,S}\) and the retained direct sum of the full generator. Their common central bundle is a direct summand of a tilting bundle, so its higher self-Ext vanishes. This gives an isomorphism of labelled bundles on \(\mathcal Z_S\). Adjunction and the projection formula, using \(Rq_*\mathcal O=\mathcal O\), now give
\[
\begin{aligned}
 \widehat{\mathcal A}_\ell\otimes_{R_H}S
 &\cong\operatorname{End}_{\mathcal Y_{\ell,S}}
            (\mathcal T_{\ell,S})^{\mathrm{op}}\\
 &\cong\operatorname{End}_{\mathcal Z_S}
            (q^*\mathcal T_{\ell,S})^{\mathrm{op}}\\
 &\cong e_{0p}\widehat{\mathcal A}_\Gamma e_{0p}.
\end{aligned}
\]
The first arrow comes from endomorphism base change in Proposition~\ref{prop:ordinary-endomorphism-base-change}, while the last is \cite[Proposition~4.15]{Karmazyn} applied to the full labelled generator. The same comparison applies to \(\mathcal T_{\Gamma,0p}\), while \(\theta\) was constructed on the partial family itself.

Every arrow is \(S\)-linear and preserves both projectors. Retaining the free projector \(e_0\) is essential because the subsequent contraction quotient divides out the two-sided ideal generated by \(e_0\). The cofinal target ideals checked above identify the ordinary and completed tensors, proving \eqref{eq:algebra-comparison} with the stated completion convention.
\end{proof}

\phantomsection\label{subsubsec:quotient-free-idempotent}
For an algebra \(\Lambda\) with idempotent \(e_0\), scalar extension commutes with the quotient \(\Lambda/\Lambda e_0\Lambda\) by right exactness.  The corner identity
\begin{equation}
\label{eq:corner-quotient}
 \frac{e_{0p}\mathcal A_\Gamma e_{0p}}
 {(e_{0p}\mathcal A_\Gamma e_{0p})e_0(e_{0p}\mathcal A_\Gamma e_{0p})}
 \cong e_p(\mathcal A_\Gamma/\mathcal A_\Gamma e_0\mathcal A_\Gamma)e_p
\end{equation}
then identifies the quotient after extension to \(\widetilde R\). By \cite[Proposition~4.15]{Karmazyn}, \(\mathcal A_\Gamma\) is the universal extended-Dynkin deformed preprojective algebra. Write \((1,d_1,\ldots,d_r)\) for the minimal positive imaginary root of the extended diagram.  The affine parameter relation is \(t_0+\sum_{i=1}^r d_it_i=0\). Killing the affine vertex of \(\mathcal A_\Gamma\) leaves the finite-Dynkin algebra \(\Pi_\lambda\), with the fixed coordinates
\[
 \lambda_i=t_i\quad(1\le i\le r),\qquad
 \lambda_0=-\sum_{i=1}^r d_i\lambda_i.
\]
At this stage the quotient is the completion of \(\mathcal D_p\) at the target origin. We show that this completion is \(\mathcal D_p\) itself. Indeed, \(\mathcal D_p\) is a quotient of the finite \(\widetilde R\)-module \(E_p\), so it is already complete for the parameter ideal. Its closed parameter fibre is the finite-dimensional graded ordinary preprojective corner, so a power of the positive-degree target-origin ideal acts trivially on that fibre. The parameter ideal is contained in the target-origin ideal, and a power of the latter sends \(\mathcal D_p\) into its parameter ideal. Their induced filtrations on \(\mathcal D_p\) therefore have cofinal powers, so the residual target completion introduces no change.

In particular this comparison introduces no parameter automorphism. Thus
\begin{equation}
\label{eq:literal-marked-quotient}
 \left(\widehat{\mathcal A}_\ell/
 \widehat{\mathcal A}_\ell e_0\widehat{\mathcal A}_\ell\right)
 \widehat\otimes_{R_H}\widetilde R\cong e_p\Pi_\lambda e_p=\mathcal D_p.
\end{equation}

\paragraph{The ordinary tensor product.}
The completion is taken at the graded central origin. Because \(\widetilde R\) is a finite free module over \(R_H\), tensoring a complete module with \(\widetilde R\) produces a finite direct sum of copies of that module, which is already complete. The isomorphism \eqref{eq:literal-marked-quotient} therefore descends to the ordinary tensor product
\begin{equation}
\label{eq:downstairs-pullback}
 \mathcal C_H\otimes_{R_H}\widetilde R\cong\mathcal D_p.
\end{equation}

\subsubsection{The quotient on the original classifying arc}
\label{subsubsec:arc-vs-contraction}

Recall \(R_1=\C[[w]]\), the original classifying arc \(f_0\colon\Spec R_1\to H\), and the pulled-back quotient \(\mathcal C_0=\mathcal C_H\otimes_{R_H,f_0}R_1\) from Subsection~\ref{subsec:cocenter-statements-reduction}. We now prove that \(\mathcal C_0\cong A\), as required by the length calculation in that subsection.

Let \(\widehat{\mathcal R}_\ell\) denote the formal completion, at the graded origin, of the affine target coordinate ring \(\mathcal R_\ell\) of the universal family, and define
\[
 Z_0=\widehat{\mathcal R}_\ell\otimes_{R_H}R_1,\qquad
 \Lambda_0=\widehat{\mathcal A}_\ell\otimes_{R_H}R_1.
\]
The relative contraction quotient is then \(\mathcal C_0=\Lambda_0/\Lambda_0e_0\Lambda_0\). We denote by \(\mathfrak m\) the maximal ideal at the contraction point in \(Z_0\).

\begin{lemma}
\label{lem:arc-quotient-supported}
The natural completion map \(\mathcal C_0\to\widehat{\mathcal C_0}_{\mathfrak m}\) is an isomorphism, and
\begin{equation}
\label{eq:completed-contraction}
 \mathcal C_0\cong\widehat{\mathcal C_0}_{\mathfrak m}\cong A.
\end{equation}
In particular,
\begin{equation}
\label{eq:cocenter-identification}
 HH_0(\mathcal C_0)\cong HH_0(A).
\end{equation}
\end{lemma}

\begin{proof}
The finiteness assertion along the arc follows from its nonconstancy. Indeed, the power series \(g(w)=P_1(f_0(w))\) has positive finite order.  Weierstrass division makes \(\C[[w]]\) finite over \(\C[[g(w)]]\).  The latter subring lies in the image of \(R_H\), so \(R_1\) is finite over \(R_H\).

Consequently \(Z_0\) is finite over the complete Noetherian local ring \(\widehat{\mathcal R}_\ell\).  Modulo the maximal ideal of \(\widehat{\mathcal R}_\ell\), the fibre is \(R_1/\mathfrak m_HR_1\), where \(\mathfrak m_H\) is the maximal ideal of \(R_H\).  This fibre is local Artinian with residue field \(\C\). Thus \(Z_0\) is local and complete at \(\mathfrak m\).  The finite \(Z_0\)-algebra \(\Lambda_0\) and its quotient \(\mathcal C_0\) are also complete.

The ordinary tensor products above therefore already compute the completed pullback along the central arc.  Apply Proposition~\ref{prop:literal-special-basic-nccr} with the classifying map independent of its external parameter \(s\). It identifies \(\Lambda_0\) with the original basic NCCR, preserving \(e_0\). Taking the idempotent quotient gives \(\mathcal C_0\cong A\), and \eqref{eq:cocenter-identification} follows on passing to cocenters.
\end{proof}

\subsection{Descent from \texorpdfstring{\(\widetilde R\)}{R tilde} to \texorpdfstring{\(R_H\)}{R H}}
\label{subsec:cocenter-descent-proof}

Theorem~\ref{thm:universal-marked-cocenter-upstairs} establishes a length-one free resolution of the cocenter after base extension to \(\widetilde R\). Using the descent isomorphism \eqref{eq:downstairs-pullback}, we construct a length-one free resolution of \(\mathsf M_H\) over \(R_H\) and determine its zeroth Fitting ideal.

\begin{proof}[Proof of Theorem~\ref{thm:marked-cocenter}]
\phantomsection\label{subsubsec:pd-one-descends}
The algebra \(\mathcal D_p\) is finite over \(\widetilde R\), so finiteness descends along the faithfully flat extension \(R_H\to\widetilde R\) to \(\mathcal C_H\) and \(\mathsf M_H\) over \(R_H\). Because the relative cocenter is the cokernel of the commutator map, it commutes with arbitrary base change, giving
\begin{equation}
\label{eq:cocenter-base-change}
 \mathsf M_H\otimes_{R_H}\widetilde R
 \cong HH_0(\mathcal D_p/\widetilde R).
\end{equation}

Choose a surjection \(R_H^q\twoheadrightarrow\mathsf M_H\) with kernel \(K\). After the faithfully flat base change to \(\widetilde R\), the kernel is \(K\otimes_{R_H}\widetilde R\) and the target is \(\mathsf M_H\otimes_{R_H}\widetilde R\cong HH_0(\mathcal D_p/\widetilde R)\). Theorem~\ref{thm:universal-marked-cocenter-upstairs} supplies a length-one free resolution of this base-changed cocenter. Comparing these two presentations of the same cokernel, Schanuel's lemma shows that \(K\otimes_{R_H}\widetilde R\) is a direct summand of a finite free \(\widetilde R\)-module, and hence is projective.

Finite projectivity descends along the faithfully flat extension \(R_H\to\widetilde R\) \cite[Tag~058S]{Stacks}, so \(K\) is a finite projective \(R_H\)-module. Because \(R_H\) is a local ring, the projective module \(K\) is free.

The base change of \(\mathsf M_H\) to \(\widetilde R\) has a square presentation with nonzero determinant, so \(\mathsf M_H\) has generic rank zero. The two free ranks therefore agree, giving \eqref{eq:downstairs-presentation}.

\phantomsection\label{subsubsec:fitting-descends}
To descend the determinant ideal, use the fact that Fitting ideals commute with base change. Thus
\begin{equation}
\label{eq:fitting-upstairs}
 \operatorname{Fitt}^{R_H}_0(\mathsf M_H)\widetilde R=(\Delta_p).
\end{equation}
For \(i\ne p\), the simple reflection \(s_i\) preserves \(m_p\) and sends every positive root other than \(\alpha_i\) to a positive root. Thus \(W_0\) permutes the factors of \(\Delta_p\) without signs, and \(\Delta_p\in\widetilde R^{W_0}=R_H\). Contracting the ideal in \eqref{eq:fitting-upstairs} along the faithfully flat extension gives
\begin{equation}
\label{eq:fitting-downstairs}
 \operatorname{Fitt}^{R_H}_0(\mathsf M_H)=(\Delta_p).
\end{equation}

By Lemma~\ref{lem:coefficient-orbits}, the roots of each marked coefficient form one \(W_0\)-orbit. Their product descends to the corresponding divisor equation, and Lemma~\ref{lem:root-divisors-and-bad-locus} shows that it has multiplicity one. Grouping the factors of \(\Delta_p\) by these orbits therefore gives
\begin{equation}
\label{eq:fitting-factorization}
 \Delta_p=\eta P_1\cdots P_\ell,\qquad\eta\in R_H^\times.
\end{equation}
\end{proof}

This completes the proof of Theorem~\ref{thm:marked-cocenter}. Combining this presentation with the identification along the original arc in Lemma~\ref{lem:arc-quotient-supported} establishes Corollary~\ref{cor:contraction-cocenter-formula} by the determinant calculation carried out in Subsection~\ref{subsec:cocenter-statements-reduction}.

\subsection{Completion of the uniqueness proof}
\label{subsec:completion-uniqueness}

The cocenter dimension formula establishes the equality between the Zariski tangent space dimension and the orbit dimension. The orbit-closure relation was already established geometrically in Theorem~\ref{thm:general-semisimple-existence}. We now present a second derivation of that relation directly from the square presentation and reduced annihilator of the universal quotient. The \hyperref[par:existence-dependency-map]{dependency discussion} records the inputs used in this second construction.

\begin{proof}[Proof of Theorem~\ref{thm:general-uniqueness}]
The dimension formulas \eqref{eq:tangent-orbit-difference} and \eqref{eq:cocenter-dimension-final} give
\begin{equation}
\label{eq:tangent-orbit-equality}
 \dim T_AX_d=\dim(G\cdot S_{\mathrm{GV}}).
\end{equation}

To apply the orbit criterion, we also need the orbit-closure relation. We construct a finite free family from the universal algebra. Choose a constant translation \(F(s,w)=f_0(w)+sv\) as in Lemmas~\ref{lem:weierstrass-root-polynomials} and \ref{lem:good-to-formal-generic}. The distinguished polynomials \(W_j(s,w)\) have degrees \(n_j\), are separable and pairwise coprime over \(\C((s))\), and their roots map under \(F\) to \(H\setminus\Sigma^{\mathrm{bad}}\). Put
\[
 \mathscr A=\mathcal C_H\otimes_{R_H,F}\C[[s,w]].
\]
Proposition~\ref{prop:weighted-contraction-presentation} gives a square presentation of \(\mathscr A\) over \(\C[[s,w]]\). Its reduction modulo \(s\) has cokernel \(A\) by Lemma~\ref{lem:arc-quotient-supported}. Since \(A\) is finite-dimensional, this reduced matrix is invertible over \(\C((w))\). Lemma~\ref{lem:finite-slice-flatness}, with coefficient ring \(\C[[s]]\), therefore makes \(\mathscr A\) finite free over that ring.

The reduced annihilator in Proposition~\ref{prop:multicurve-reduced-annihilator}, specialized to \(J=\{p\}\), shows that \(\delta=\prod_jP_j\) annihilates \(\mathcal C_H\). Thus \(\mathscr A\) is an algebra over the following finite \(\C[[s]]\)-algebra
\[
 \C[[s,w]]/(\delta(F))
 \cong \C[[s]][w]/\left(\prod_jW_j(s,w)\right).
\]
This finite presentation allows us to retain all roots when inverting \(s\).

Over \(L=\overline{\C((s))}\), the finite algebra in the preceding display splits into one copy of \(L\) for each root of the \(W_j\). Since this algebra acts centrally, its primitive idempotents decompose \(\mathscr A_L\). Each resulting factor is the fibre of \(\mathcal C_H\) at its parameter point.

We identify each factor by lifting its parameter point to the root space. Since \(R_H\to\widetilde R\) is finite and faithfully flat and \(L\) is algebraically closed, this point lifts to an \(L\)-point of \(\Spec\widetilde R\). For a root of \(W_j\), the lifted parameter point lies on a non-Levi root hyperplane of retained coefficient \(j\), since the pullback of \(P_j\) is the product of those root equations. Avoidance of \(\Sigma^{\mathrm{bad}}\) excludes every other root hyperplane. The marked comparison \eqref{eq:downstairs-pullback} identifies the factor with the retained corner of the deformed preprojective algebra at this point, where exactly one positive root vanishes. Lemma~\ref{lem:single-root-preprojective-fibre} shows that the preprojective algebra is a full matrix algebra and that its retained projector has rank \(j\). Its corner is therefore \(M_j(L)\). Consequently
\[
 \mathscr A_L\cong\prod_{j=1}^{\ell}M_j(L)^{\times n_j}.
\]

Lemma~\ref{lem:family-orbit-closure} now gives
\begin{equation}
\label{eq:orbit-closure}
 A\in\overline{G\cdot S_{\mathrm{GV}}}.
\end{equation}
Together with the tangent-space equality proved at the start, this is the hypothesis needed for Theorem~\ref{thm:orbit-smooth-unique}. The theorem shows that \(A\) is smooth, that \(\overline{G\cdot S_{\mathrm{GV}}}\) is its unique component, and that any split semisimple algebra whose orbit closure contains \(A\) is isomorphic to \(S_{\mathrm{GV}}\).

It remains to apply this conclusion to an arbitrary semisimple deformation. Let \(\mathscr D\) be any finite free formal deformation of \(A\) over \(\C[[u]]\) with geometrically semisimple generic fibre. Choose a product \(S'\) of complex matrix algebras with the block sizes and multiplicities of the geometric generic fibre, so that \(\mathscr D\otimes_{\C[[u]]}\overline{\C((u))}\cong S'\otimes_{\C}\overline{\C((u))}\). Lemma~\ref{lem:family-orbit-closure} gives \(A\in\overline{G\cdot S'}\). The orbit criterion therefore gives
\begin{equation}
\label{eq:uniqueness-conclusion}
 S'\cong S_{\mathrm{GV}}
 =\prod_{j=1}^{\ell}M_j(\C)^{\times n_j}.
\end{equation}
Tensoring this isomorphism with \(\overline{\C((u))}\) gives the asserted description of the geometric generic fibre.
\end{proof}

\phantomsection\label{subsec:relative-hh0-discussion}

\begin{corollary}
\label{cor:semisimple-family-cocenter}
Let \(\mathscr B\) be a finite free \(\C[[s]]\)-algebra with special fibre \(A\) and geometrically semisimple generic fibre.  Then \(HH_0(\mathscr B/\C[[s]])\) is free of rank \(N_{\mathrm{GV}}\).
\end{corollary}

\begin{proof}
The relative cocenter is a finite \(\C[[s]]\)-module, and forming it commutes with arbitrary base change. Theorem~\ref{thm:general-uniqueness} gives dimension \(N_{\mathrm{GV}}\) on the closed fibre and, after a splitting extension, on the generic fibre. The equality of these two fibre dimensions excludes its torsion summands by Lemma~\ref{lem:dvr-torsion-exclusion-expanded}, so it is free of rank \(N_{\mathrm{GV}}\).
\end{proof}

\phantomsection\label{subsec:final-status}
For several retained curves, Section~\ref{sec:hua-toda-flatness} uses the exact cocenter sequence associated with deleting vertices to construct the corresponding universal presentations. Lemma~\ref{lem:finite-slice-flatness} makes the modules presented by their pullbacks along classifying maps \(R\to B[[w]]\) finite free over \(B\). The resulting rank formulas are used in Sections~\ref{sec:marked-semisimple} and~\ref{sec:centre-morphism}.
\section{Universal presentations and flat families}
\label{sec:hua-toda-flatness}

Let a smooth threefold flop vary over a complete Noetherian local complex base.  Using a lift of its fixed basic generator, we form the relative contraction quotient.  We prove that this quotient, its cocenter, and its centre are finite free over the base, with ranks determined by the central flop.  This is the content of the three flatness assertions in \cite[Conjecture~4.3]{HuaToda}, in the completed smooth-threefold setting.  Kawamata established algebra flatness over a disk for configurations of exceptional curves \cite[Theorem~3.1(2)]{KawamataNC}.

The simultaneous partial resolutions recalled in \cite[Section~2.3]{Karmazyn} provide the geometric models. The lifting theorem for tilting bundles \cite[Theorem~3.4]{KarmazynDeformations} provides the generators. Below we verify the hypotheses and carry out the classifying and endomorphism comparisons for the families at hand. The quotient by the free projector and its cocenter have square presentations derived from the Etingof--Rains central extension \cite[Theorem~3.2]{EtingofRains}.  The determinants of these two square presentations generate the zeroth Fitting ideals, with factors given by the curve-divisor equations.

Pullback gives matrices over \(B[[w]]\), where \(B\) is the external base and \(w\) is the surface parameter.  Lemma~\ref{lem:finite-slice-flatness} shows that a finite-dimensional cokernel after reduction modulo the maximal ideal of \(B\) forces freeness over \(B\), even when \(B\) has nilpotents. The rank of the cokernel as a \(B\)-module is the order in \(w\) of the determinant after reduction modulo the maximal ideal of \(B\).  We apply the lemma separately to the algebra and its cocenter, then use a symmetrizing trace to identify the centre with the dual cocenter.

The one-curve cocenter presentation is available from Section~\ref{sec:general-uniqueness}.  The pullback and trace arguments below are independent of the number of retained curves.  For several curves, the corresponding presentations are constructed in Subsection~\ref{sec:multicurve-flatness}; the Peirce and cocenter ranks of the multi-curve contraction quotient enter Section~\ref{sec:marked-semisimple}.

\paragraph{The two existence arguments and their dependencies.}
\label{par:existence-dependency-map}
The geometric proof of Theorem~\ref{thm:general-semisimple-existence} identifies its local matrix factors in Theorem~\ref{thm:entire-completed-j-factor} and assembles them in Proposition~\ref{prop:exact-general-generic-product}. It also uses Proposition~\ref{prop:weighted-contraction-presentation} in Lemma~\ref{lem:fullness-away-from-roots} to retain every moving support point. The proof of Theorem~\ref{thm:general-uniqueness} gives another construction of the orbit-closure input directly from \(\mathcal C_H\): Proposition~\ref{prop:weighted-contraction-presentation} and Lemma~\ref{lem:finite-slice-flatness} give a finite free family; Proposition~\ref{prop:multicurve-reduced-annihilator}, with \(J=\{p\}\), and Lemma~\ref{lem:single-root-preprojective-fibre} identify its split factors.  This second construction uses neither the geometric local-factor calculation nor the cocenter count. Both constructions use the comparison of the labelled generators after extension from \(R_H\) to \(\widetilde R\), the identification of the quotient along the original classifying arc, and the geometry of the translated arc and its GV multiplicities from Proposition~\ref{prop:conserved-multiplicities-gv}. The algebra presentation and reduced-annihilator proofs below use the Etingof--Rains central extension independently of semisimple existence and uniqueness.  We use the cocenter presentation to compute the tangent dimension and to prove flatness of the centre and cocenter.

\begin{remark}
\label{rem:hua-toda-support}
A geometric deformation over \((B,\mathfrak m,\C)\) means a projective contraction with affine Noetherian target, source and target flat over \(B\), and the prescribed central contraction as closed fibre.  We complete along the central contraction point and its inverse image. Lift the fixed normalized basic Van den Bergh generator, retaining its structure-sheaf summand and all curve labels. Form the opposite endomorphism algebra of this lift, then quotient by the ideal generated by the structure-sheaf projector.

The finite relative algebra carries all support points that move away from the central point.  Base change is performed on this algebra. For example, \(B[[w]]/(w^2-t)\), with \(B=\C[[t]]\), is free of rank two over \(B\), whereas \(\C((t))[[w]]/(w^2-t)=0\). Thus the generic fibre of the finite relative algebra is obtained by localizing that algebra; completing again at \(w=0\) would remove its moving points. Fixing the generator is equally significant: its original summand multiplicities give the matrix blocks in the split fibres.  A block \(M_j\) contributes \(j^2\) to the algebra dimension, while its basic Morita representative \(\C\) has dimension one.
\end{remark}

\subsection{Freeness and traces over local bases}
\label{subsec:finite-slice-flatness}

A square matrix over \(B[[w]]\) with finite-dimensional cokernel after reduction modulo the maximal ideal of \(B\) gives a finite free \(B\)-module.  Once both an algebra and its cocenter are finite free, a symmetrizing trace on the closed algebra lifts and identifies the centre with the dual cocenter.

\subsubsection{Freeness over a possibly nonreduced base}
\label{subsubsec:finite-slice-matrices}

After reduction modulo the maximal ideal of \(B\), Smith normal form over \(k[[w]]\), where \(k\) is the residue field, computes the cokernel length from the determinant. To prove that the cokernel is free over \(B\), we first lift generators of its closed fibre to show that it is a finite \(B\)-module. Injectivity of the square matrix is then checked through the powers of the maximal ideal, giving the required vanishing of the first Tor group.

\begin{lemma}
\label{lem:finite-slice-flatness}
Let \((B,\mathfrak m,k)\) be a complete Noetherian commutative local ring, put \(D_B=B[[w]]\), and let \(\Phi:D_B^a\to D_B^a\) be a square matrix. Write \(\overline\Phi\) for its reduction modulo \(\mathfrak m\). If \(\det\overline\Phi\ne0\) in \(k[[w]]\), then \(\Phi\) is injective and its cokernel is finite free over \(B\), of rank \(\operatorname{ord}_w\det\overline\Phi\).
\end{lemma}

\begin{proof}
Put \(M=\operatorname{coker}\Phi\).  Smith normal form over \(k[[w]]\) shows that \(M/\mathfrak mM\) is finite-dimensional over \(k\).  Lift a basis to \(m_1,\ldots,m_b\in M\).  Successive reduction modulo \(\mathfrak mM\) writes every \(m\in M\) as
\[
 m=\sum_i b_i^{(n)}m_i+r_n,\qquad
 r_n\in\mathfrak m^nM,\quad
 b_i^{(n+1)}-b_i^{(n)}\in\mathfrak m^n.
\]
Completeness of \(B\) gives limits for the coefficients.  Since \(M\) is finite over the Noetherian local ring \(D_B\), Krull intersection \cite[Tag~\href{https://stacks.math.columbia.edu/tag/00IP}{00IP}]{Stacks} gives \(\bigcap_n\mathfrak m^nM=0\).  The limiting coefficients express \(m\) as a \(B\)-linear combination of the \(m_i\).  Thus \(M\) is finite over \(B\).

Let \(v\in D_B^a\) satisfy \(\Phi v=0\). Since \(\det\overline\Phi\ne0\) in \(k[[w]]\), the reduced map \(\overline\Phi\) is injective on \(k[[w]]^a\), which forces \(v\in\mathfrak mD_B^a\). Because \(D_B=B[[w]]\) is the \((w)\)-adic completion of the polynomial ring \(B[w]\) and completion of a Noetherian ring is flat \cite[Tag~\href{https://stacks.math.columbia.edu/tag/00MB}{00MB}]{Stacks}, \(D_B\) is flat over \(B\). This flatness identifies the \(n\)-th graded piece \(\operatorname{gr}_{\mathfrak m}^n(D_B^a)\) canonically with \((\mathfrak m^n/\mathfrak m^{n+1})\otimes_k k[[w]]^a\). On each graded piece, the map induced by \(\Phi\) is
\[
 1\otimes\overline\Phi:
 (\mathfrak m^n/\mathfrak m^{n+1})\otimes_k k[[w]]^a
 \longrightarrow
 (\mathfrak m^n/\mathfrak m^{n+1})\otimes_k k[[w]]^a,
\]
which is injective because \(\overline\Phi\) is injective. Induction on \(n\) places \(v\) in \(\mathfrak m^nD_B^a\) for every \(n\ge1\). Krull's intersection theorem on the finite \(D_B\)-module \(D_B^a\) then gives \(v=0\), establishing injectivity.

Because \(D_B\) is flat over \(B\), the short exact sequence \(0\to D_B^a\xrightarrow{\Phi}D_B^a\to M\to0\) provides a resolution of \(M\) by flat \(B\)-modules, which gives
\[
 \operatorname{Tor}_1^B(M,k)=\ker\overline\Phi=0.
\]
Take a minimal free surjection \(B^b\to M\), with finite kernel \(L\). Since \(\operatorname{Tor}_1^B(M,k)=0\), the reduction of \(L\) injects into the kernel of \(k^b\to M/\mathfrak m M\). Minimality makes this last map an isomorphism.  Thus \(L/\mathfrak mL=0\), and Nakayama gives \(L=0\).  The rank is the closed-fibre dimension. Smith normal form gives \(b=\dim_k(M/\mathfrak mM)=\operatorname{ord}_w\det\overline\Phi\).
\end{proof}

The preceding freeness criterion, applied to a free resolution over the parameter ring, also identifies the derived pullback with the ordinary pullback.  In the following corollary, \((B,\mathfrak m,k)\) remains a complete Noetherian local ring as in Lemma~\ref{lem:finite-slice-flatness}.

\begin{corollary}
\label{cor:perfect-finite-slice}
Suppose a finite module \(C\) over a Noetherian local ring \(R\) has an exact sequence \(0\to R^a\xrightarrow{\Theta}R^a\to C\to0\). For a local map \(F:R\to B[[w]]\), if \(C\otimes_R k[[w]]\) is finite-dimensional, then
\[
 C\otimes_R^L B[[w]]\simeq
 \bigl(C\otimes_R B[[w]]\bigr)[0].
\]
The ordinary pullback is finite free over \(B\), with rank equal to the dimension of that closed fibre.
\end{corollary}

\begin{proof}
We compute the derived pullback \(C\otimes_R^L B[[w]]\) using the finite free resolution \(0\to R^a\xrightarrow{\Theta}R^a\to C\to0\). Base change along \(F\colon R\to B[[w]]\) yields the two-term complex \(D_B^a\xrightarrow{\Theta(F)}D_B^a\) in degrees \(-1\) and \(0\). Under reduction modulo \(\mathfrak m\), the closed-fibre matrix \(\overline{\Theta(F)}=\Theta(F)\otimes_B k\) has cokernel \(C\otimes_R k[[w]]\), which is finite-dimensional over \(k\) by hypothesis. Inverting \(w\) makes \(\overline{\Theta(F)}\) an invertible matrix over the Laurent field \(k((w))\), so \(\det\overline{\Theta(F)}\ne0\) in \(k[[w]]\).

Applying Lemma~\ref{lem:finite-slice-flatness} to \(\Phi=\Theta(F)\) shows that \(\Theta(F)\) is injective and that its cokernel \(C\otimes_R B[[w]]\) is a finite free \(B\)-module of rank \(\dim_k(C\otimes_R k[[w]])\). The injectivity of \(\Theta(F)\) ensures that the pulled-back resolution has cohomology only in degree zero, establishing the quasi-isomorphism \(C\otimes_R^LB[[w]]\simeq(C\otimes_R B[[w]])[0]\).
\end{proof}

\subsubsection{Symmetrizing traces and base change}
\label{subsubsec:relative-symmetrizing-traces}

The trace pairing identifies the centre with the dual cocenter of a symmetric algebra; see \cite[Corollary~2.2]{EtingofLatourRains} for the central extension.  The relative assertion requires lifting a given closed-fibre trace and preserving its perfect pairing under scalar extension.  Here the base need only be local.

\begin{lemma}
\label{lem:relative-symmetrizing-trace}
Let \((B,\mathfrak m,k)\) be a commutative local ring, and let \(\mathscr A\) be a finite free \(B\)-algebra.  Suppose \(M=HH_0(\mathscr A/B)\) is finite free and the closed fibre \(A_0=\mathscr A\otimes_B k\) of \(\mathscr A\) is symmetric.  Every symmetrizing trace of \(A_0\) lifts to a perfect \(B\)-linear symmetrizing trace \(\tau\) on \(\mathscr A\).  It induces an isomorphism
\[
 Z(\mathscr A)\xrightarrow{\sim}M^\vee,
 \qquad z\longmapsto\bigl(\bar a\longmapsto\tau(za)\bigr).
\]
Here \(M^\vee=\operatorname{Hom}_B(M,B)\). The centre and cocenter commute with every commutative scalar extension \(B\to B'\).  The lifted trace determines the displayed duality, while the base-change maps are canonical.
\end{lemma}

\begin{proof}
As a commutator cokernel, \(M\) commutes with arbitrary scalar extension.  The given closed trace is a functional on \(M\otimes_B k\). Freeness of \(M\) lifts it to a functional \(M\to B\), hence to a trace \(\tau:\mathscr A\to B\).  The determinant of the pairing \((a,b)\mapsto\tau(ab)\) has nonzero reduction, so it is a unit in the local ring \(B\).  The lifted trace is therefore symmetrizing.

To identify the centre, take the annihilator of the commutator submodule under this perfect pairing.  We have
\[
 \tau(z[a,b])=\tau((za-az)b).
\]
This identity shows that a central element annihilates every commutator under the pairing.  Conversely, if \(\tau(z[a,b])=0\) for all \(a,b\), perfection forces \(za-az=0\) for every \(a\).  The annihilator is thus the centre, giving the stated duality and centre freeness.

After any scalar extension to \(B'\), the trace pairing remains perfect. Under centre--cocenter duality, the canonical map \(Z(\mathscr A)\otimes_BB'\to Z(\mathscr A\otimes_BB')\) is the base-change map for the dual of the finite free module \(M\). It is therefore an isomorphism, including when \(B\to B'\) is not flat.
\end{proof}

\subsection{One curve: presentations and pullback along classifying maps}
\label{subsec:one-curve-finite-slices}

At a transverse point of \(D_j\), the matrix algebra \(M_j\) contributes \(j^2\) to the algebra dimension and \(1\) to the cocenter dimension.  These are the exponents in the two determinant formulas.

\subsubsection{A square presentation of the universal quotient}
\label{subsec:weighted-discriminant}

Recall the notation of Subsection~\ref{subsec:descent-cartan-cover}. Here \(R_H=\widehat{\mathcal O}_{H,0}\), \(\widetilde R=\widehat{\mathcal O}_{\mathfrak h,0}\), and \(\mathcal C_H\) is the universal contraction quotient. The extension \(R_H\to\widetilde R\) is finite and faithfully flat. The equation \(P_j=0\) defines the reduced divisor \(D_j\) for curves of class \(j[C]\).  The algebra presentation is derived below from the same central extension that gives the cocenter presentation.

\Needspace{14\baselineskip}
\begin{proposition}
\label{prop:weighted-contraction-presentation}
The finite \(R_H\)-module \(\mathcal C_H\) has a free resolution of the form
\begin{equation}
\label{eq:weighted-contraction-resolution}
 0\longrightarrow R_H^a\xrightarrow{\Theta_A}R_H^a
 \longrightarrow\mathcal C_H\longrightarrow0
\end{equation}
with
\begin{equation}
\label{eq:weighted-contraction-determinant}
 \operatorname{Fitt}^{R_H}_0(\mathcal C_H)
 =(\det\Theta_A)=\left(\prod_{j=1}^{\ell}P_j^{j^2}\right).
\end{equation}
In particular, \(\mathcal C_H\) is a Cohen--Macaulay module of codimension one over \(R_H\).
\end{proposition}

\begin{proof}
Recall from \eqref{eq:central-extension-corner} that \(E_p=e_p(\Pi[z]\widehat\otimes_{K_\rho}\widetilde R)e_p\), the corner at the retained vertex \(p\) of the completed central extension. Theorem~3.2(i) of \cite{EtingofRains} makes \(\Pi[z]\) finite free over \(K_\rho=\C[x_i-\rho_i z]\).  After flat completion, its retained corner is a direct summand of a finite free \(\widetilde R\)-module. It is therefore free, since \(\widetilde R\) is local.  

Before the change of coefficient ring, \(z\) is a polynomial variable, so multiplication by \(z\) is injective.  Flat completion and the corner projection preserve this injection.  Its cokernel is obtained by setting \(z=0\):
\begin{equation}
\label{eq:corner-z-free-resolution}
 0\longrightarrow E_p\xrightarrow{z}E_p
 \longrightarrow\mathcal D_p\longrightarrow0.
\end{equation}

To compute the determinant, use the generic decomposition \eqref{eq:generic-extension-decomposition}.  The factor of \((E_p)_L\) indexed by a positive root with \(m_p(\alpha)>0\) is \(M_{m_p(\alpha)}(L)\), on which \(z\) acts as \(- (\lambda,\alpha)/(\rho,\alpha)\) \cite[Theorem~3.2(v) and the following display]{EtingofRains}. This scalar occurs with multiplicity \(m_p(\alpha)^2\), the dimension of the matrix algebra.  Hence, in the domain \(\widetilde R\),
\[
 \det(z\mid E_p)=\eta
 \prod_{\substack{\alpha\in\Phi^+\\m_p(\alpha)>0}}
 (\lambda,\alpha)^{m_p(\alpha)^2},
 \qquad\eta\in\C^\times.
\]
The determinant is nonzero, so \(\mathcal D_p\) is torsion and has projective dimension one over \(\widetilde R\).

By \eqref{eq:downstairs-pullback}, \(\mathcal C_H\otimes_{R_H}\widetilde R\cong\mathcal D_p\). Choose a finite free surjection \(R_H^a\twoheadrightarrow\mathcal C_H\) with kernel \(K\).  The length-one free resolution \eqref{eq:corner-z-free-resolution} and the base-changed surjection have the same cokernel \(\mathcal D_p\). Schanuel's lemma therefore makes \(K\otimes_{R_H}\widetilde R\) projective.  Finite projectivity descends along the faithfully flat map \(R_H\to\widetilde R\), and projective modules over the local ring \(R_H\) are free.  Thus \(K\) is free.  Since \(\mathcal C_H\) is torsion, \(\operatorname{rk}K=a\), proving \eqref{eq:weighted-contraction-resolution}.

To identify the determinant ideal, recall that Fitting ideals commute with scalar extension. Lemma~\ref{lem:coefficient-orbits} groups the roots of coefficient \(j\) into a single orbit, whose product is the pullback of \(P_j\) with multiplicity one by Lemma~\ref{lem:root-divisors-and-bad-locus}. Every factor in that orbit has exponent \(j^2\), so the displayed determinant ideal is the extension of \((\prod_jP_j^{j^2})\).  Faithful flatness then gives \eqref{eq:weighted-contraction-determinant} over \(R_H\). The Fitting ideal shows that \(\mathcal C_H\) has dimension \(\dim R_H-1\). Since \(R_H\) is regular local and the nonzero module \(\mathcal C_H\) has projective dimension one, Auslander--Buchsbaum \cite[Proposition~10.111.1, Tag~\href{https://stacks.math.columbia.edu/tag/090V}{090V}]{Stacks} gives the same value for its depth. Thus \(\mathcal C_H\) is Cohen--Macaulay.
\end{proof}

The universal contraction quotient and its cocenter have the same support, namely the union of the curve divisors.  Up to units, their determinants are \(\prod_jP_j^{j^2}\) for the algebra and \(\prod_jP_j\) for the cocenter, the latter by Theorem~\ref{thm:marked-cocenter}.  The exponent is the dimension of the generic matrix block for the algebra presentation, and the dimension of the matrix factor's cocenter for the cocenter presentation.

\subsubsection{Classifying an arbitrary geometric deformation}
\label{subsec:arbitrary-base-classification}

Fix a completed simple flop \(f_0:X_0\to Y_0\) of a smooth complex threefold.  Let \(\mathcal T_0=\mathcal O\oplus\mathcal T'_0\) be its selected basic tilting generator, as in Section~\ref{sec:general-existence}, and \(A_0\) its contraction algebra.

The families of Remark~\ref{rem:hua-toda-support} are projective contractions with affine Noetherian target and with source and target flat over the base. Their relative quotients use the fixed generator with its labelled summands.

\Needspace{20\baselineskip}
\begin{proposition}
\label{prop:arbitrary-base-classifying-arcs}
Let \(B\) be a complete Noetherian local \(\C\)-algebra with residue field \(\C\), with the residue identification compatible with its \(\C\)-algebra structure, and let \(f_B\) be such a completed deformation.  A lift of a general surface-section parameter supplies a continuous local map
\[
 F:R_H\longrightarrow D_B=B[[w]]
\]
and identifies \(f_B\) with the completed pullback of the universal marked partial resolution.  There is a relative tilting generator \(\mathcal T_B\) lifting \(\mathcal T_0\).  For any such lift with the fixed summands, its relative contraction quotient satisfies
\begin{equation}
\label{eq:arbitrary-base-quotient-identification}
 \mathscr A_B:=%
 \frac{\End_{X_B}(\mathcal T_B)^{\mathrm{op}}}
 {\End_{X_B}(\mathcal T_B)^{\mathrm{op}}e_0
  \End_{X_B}(\mathcal T_B)^{\mathrm{op}}}
 \cong\mathcal C_H\otimes_{R_H,F}D_B.
\end{equation}
Its reduction modulo \(\mathfrak m\) is \(A_0\).
\end{proposition}

\begin{proof}
\emph{Constructing the classifying map.} We first view the given family as a deformation of its central partial surface resolution. Write the completed target as \(\operatorname{Spec}R\), with maximal ideal \(\mathfrak n\), and retain the projective \(R\)-scheme whose formal completion is the given source. Choose a general central section parameter \(w_0\) and a lift \(w\in\mathfrak n\). Evaluation of power series defines \(D_B=B[[w]]\to R\). Because \(w_0\) is a nonzerodivisor on the central target and source, this regularity persists after their Noetherian completions.

At a source point over the closed target point, let \(Q\) be its local ring. The maps \(B\to D_B\to Q\) are local, \(Q\) is \(B\)-flat, and \(Q/\mathfrak mQ\) is torsion-free, hence flat, over \(\C[[w]]\). The fibrewise criterion with module \(Q\), finite over itself, gives \(D_B\)-flatness \cite[Tag~\href{https://stacks.math.columbia.edu/tag/00MP}{00MP}]{Stacks}. Every source point specializes to a point over the closed target point by properness over \(R\). Localizing at its generalizations then gives \(D_B\)-flatness at every source point. The same criterion applies to \(R\).

Let \(p:S\to V=\operatorname{Spec}A\) be the algebraic marked partial ADE resolution modeling the closed surface, and let \(\widehat S\) be its completion over the vertex ideal \(\mathfrak q\subset A\). Reid's section lemma gives this partial resolution \cite[Section~1, preprint p.~5]{KatzMorrison}. The strict transform of the central surface section agrees with its scheme-theoretic inverse image, because smallness excludes an exceptional divisorial component, while the Cartier divisor on the smooth central source has no embedded component.

Fix the central marked identification with \(\widehat S\). The standard simultaneous partial-resolution source family over \(H=\mathfrak h/W_0\) is semiuniversal for ordinary deformations of \(S\), its Kodaira--Spencer map is an isomorphism, and \(H\) is smooth \cite[Section~6, preprint p.~38]{KatzMorrison} \cite[Section~2.3, preprint pp.~9--10]{Karmazyn}. The equivariant terminology in these references includes ordinary deformations, and no action on the given family is assumed.

We compare the deformation tangent space of the algebraic surface \(S\) with that of its formal completion \(\widehat S\). Since \(S\) has only ADE singularities, \(L_{S/\C}\) is perfect of amplitude \([-1,0]\). Put \(K=R\mathcal Hom(L_{S/\C},\mathcal O_S)\). The finite \(A\)-module \(H^1(S,K)\) is supported at the vertex. Away from it, \(p\) is an isomorphism and \(K\) is a tangent bundle in degree zero. Thus completion leaves \(H^1(S,K)\) unchanged.

The continuous cotangent complex of \(\widehat S\) is the completion of \(L_{S/\C}\) \cite[Proposition~3.9, preprint p.~16]{PerezRodriguezFormalDeformations}. Dualization commutes with this completion because \(L_{S/\C}\) is perfect. Proper formal functions therefore gives
\[
 T^1_{\mathrm{cont}}(\widehat S)
 =H^1(\widehat S,\widehat K)
 =H^1(S,K)^{\wedge}_{\mathfrak q}=T^1(S)
\]
\cite[Tags~\href{https://stacks.math.columbia.edu/tag/0A0G}{0A0G} and~\href{https://stacks.math.columbia.edu/tag/0A0H}{0A0H}]{Stacks}. The comparison respects Kodaira--Spencer maps.

The tangent-space comparison allows the classifying map and its source isomorphism to be lifted together. Put \(\mathfrak a=(\mathfrak m,w)\) and \(D_n=D_B/\mathfrak a^n\). At each stage of the induction we choose both \(F_n:R_H\to D_n\) and an isomorphism \(\phi_n\) of the completed source pullback with the given completed source over \(D_n\). Suppose the pair is chosen at stage \(n\). The power-series ring \(R_H\) allows an arbitrary lift of \(F_n\) to \(D_{n+1}\). Compare the source of this universal pullback with the given source over \(D_{n+1}\), using \(\phi_n\) to fix the identification over \(D_n\). For the square-zero extension with kernel \(J_n=\mathfrak a^n/\mathfrak a^{n+1}\), annihilated by the maximal ideal of \(D_{n+1}\), the classes of such liftings with this reduction identification fixed form a torsor under \(T^1_{\mathrm{cont}}(\widehat S)\otimes_\C J_n\) \cite[Propositions~4.2--4.3 and Corollary~4.4, preprint pp.~19--20]{PerezRodriguezFormalDeformations}.

Here the formal sources are locally Noetherian, and their reductions with respect to the ideals defining completion at the central target points are of finite type over \(\C\). Their maps to the Artin bases are therefore pseudo-finite-type in the sense of \cite[Section~1.1.1]{PerezRodriguezFormalDeformations}. No properness or adicness over \(D_n\) is required. To apply Corollary~4.4 of that paper, prescribe the extension kernel \(\mathcal O_{\widehat S}\otimes_\C J_n\) and its canonical coefficient map. For a flat lifting this map is an isomorphism. Conversely, an extension with this prescribed map is flat by \cite[Lemma~91.11.2, Tag~\href{https://stacks.math.columbia.edu/tag/08MQ}{08MQ}]{Stacks}.

Thus the classification is of flat liftings with the specified reduction isomorphism, rather than arbitrary coherent extensions.

Write \(\mathfrak X_n\) for the identified formal source over \(D_n\) and \(\widehat L_{\mathfrak X_n/D_n}\) for its completed cotangent complex. Ordinary flat cotangent base change on the already identified algebraic standard pullback, followed by completion, gives
\[
 \operatorname{Ext}^1_{\mathfrak X_n}
 \bigl(\widehat L_{\mathfrak X_n/D_n},
       \mathcal O_{\widehat S}\otimes_\C J_n\bigr)
 \cong T^1_{\mathrm{cont}}(\widehat S)\otimes_\C J_n
\]
\cite[Tag~\href{https://stacks.math.columbia.edu/tag/08QQ}{08QQ}]{Stacks}. This uses the standard pullback at stage \(n\), not an algebraization of the unknown lifting. The torsor classifies liftings with the reduction \(\phi_n\) fixed. Equal classes therefore give an isomorphism that extends \(\phi_n\).

Changing the lifted base coordinates changes the torsor class by the Kodaira--Spencer map. Since this map is surjective, choose the lifted coordinates so that the difference class vanishes. This gives \(F_{n+1}\) and an isomorphism \(\phi_{n+1}\) extending the specified \(\phi_n\). Hence \(F=\varprojlim F_n\) is a continuous local map \(R_H\to D_B\), with compatible formal source isomorphisms.

\emph{Recovering the contraction.} The compatible source identifications determine the targets through global functions on the sources. This requires both parameter-adic completion and completion at the maximal ideal of the target, which we compare separately. Set \(I=\mathfrak aR\), \(R_n=R/I^n\), and let \(X_n\) be the base change of the chosen projective source to \(R_n\). The closed \(D_B\)-fibre is the partial resolution of a rational double point, so \(H^0(X_1,\mathcal O)=R_1\) and \(H^{>0}(X_1,\mathcal O)=0\). For the reduction from \(D_{n+1}\) to \(D_n\), flatness identifies the square-zero extension kernels with \(\mathcal O_{X_1}\otimes_\C J_n\) and \(R_1\otimes_\C J_n\), respectively. Induction on these exact sequences gives \(H^0(X_n,\mathcal O)=R_n\) and \(H^{>0}(X_n,\mathcal O)=0\), also for the standard pullback.

Since \(R_n\) is complete for its maximal ideal, formal functions identifies the global functions on the completed source at this stage with \(R_n\) \cite[Tag~\href{https://stacks.math.columbia.edu/tag/02OC}{02OC}]{Stacks}. Thus \(\phi_n\) induces an isomorphism of the complete affine targets, and respects their canonical contraction maps. After identifying those targets, Grothendieck's full faithfulness for morphisms, applied over \(R_n\) with its maximal ideal, algebraizes \(\phi_n\) and its inverse to scheme isomorphisms \cite[Tag~\href{https://stacks.math.columbia.edu/tag/0A4Z}{0A4Z}]{Stacks}. Uniqueness makes these algebraizations compatible in \(n\).

For passage from the thickenings to the full target, we use the following fact. A complete Noetherian local ring \((R,\mathfrak n)\) is also complete for \(I\subset\mathfrak n\), because compatible representatives are \(\mathfrak n\)-adically Cauchy, and each \(I^r\) is \(\mathfrak n\)-adically closed. Separatedness follows from Krull intersection. This does not assert cofinality of the two topologies.

The target isomorphisms therefore pass to the actual complete target rings. Applying the same full-faithfulness theorem again, now with ideal \(I\), identifies the projective sources and their contraction morphisms, hence their formal completions. Properness here is over the affine target throughout.

\emph{The marking and the contraction quotient.} Having identified the contraction over the base, we compare its chosen generator before taking the quotient. The central identification preserves the curve labels. For a marking given by degree line bundles, the lifted line-bundle classes are unique since \(H^1(\widehat S,\mathcal O)= H^2(\widehat S,\mathcal O)=0\). The square-zero units sequence also makes their isomorphisms compatible \cite[Tag~\href{https://stacks.math.columbia.edu/tag/0C6R}{0C6R}]{Stacks}.

Indeed, for the square-zero coefficient ideal \(\mathcal J\), the units sequence has kernel \(1+\mathcal J\cong\mathcal J\). It is exact on the same underlying formal topological space, just as for schemes. The vanishing of \(H^2(\mathcal J)\) lifts line bundles, and that of \(H^1(\mathcal J)\) gives uniqueness and surjectivity on the units needed to lift their isomorphisms. This retains the fixed marked source family and its positive class, without choosing a lift of the classifying map to \(\mathfrak h\).

The flat universal source and target make the classifying square Tor-independent by Lemma~\ref{lem:arbitrary-pullback-tor-independent}. Lemma~\ref{lem:perverse-generator-survives} then ensures that the pulled-back bundle remains a perverse projective generator. Its central restriction is \(\mathcal T_0\), with its specified summands, by Proposition~\ref{prop:literal-special-basic-nccr} and the same rank-and-determinant argument on the central arc.

The vanishing of the second self-Ext group allows each summand of \(\mathcal T_0\) to lift through a square-zero thickening. The first self-Ext vanishing allows the isomorphisms between corresponding lifts to extend as well. As in Proposition~\ref{prop:formal-marked-comparison}, projectivity over the complete target allows coherent formal existence and full faithfulness to identify the compatible bundle lifts. The isomorphisms are constructed summand by summand, so they preserve \(e_0\). Thus any other lift of the fixed generator gives an isomorphic marked quotient.

Ordinary endomorphism base change
(Proposition~\ref{prop:ordinary-endomorphism-base-change}) and right exactness
of the idempotent quotient identify the marked quotient with the completion,
at the central target point, of the ordinary polynomial pullback. It remains to show that the ordinary polynomial quotient is already complete at this target point.

The reduction \(F_0\) modulo \(\mathfrak m\) classifies the original flop, and the argument of Lemma~\ref{lem:arc-quotient-supported} identifies the quotient along \(F_0\) with \(A_0\).  This closed fibre is finite-dimensional, so Proposition~\ref{prop:weighted-contraction-presentation} and Corollary~\ref{cor:perfect-finite-slice} make the right side of \eqref{eq:arbitrary-base-quotient-identification} finite over \(B\). On this finite \(B\)-module, \(w\) acts nilpotently modulo \(\mathfrak m\). Consequently the \((\mathfrak m,w)\)-adic and \(\mathfrak m\)-adic topologies agree.  The ordinary and completed quotient expressions therefore coincide. This proves \eqref{eq:arbitrary-base-quotient-identification} for the ordinary pullback, retaining all moving support points.
\end{proof}

\subsubsection{Flatness and base change for one curve}
\label{subsec:hua-toda-theorem}

\begin{proposition}
\label{thm:hua-toda-arbitrary-base}
Let \(f_0\) and \(f_B\) be as in Proposition~\ref{prop:arbitrary-base-classifying-arcs}.  Put
\[
 d=\dim_\C A_0=\sum_jj^2n_j,\qquad N_{\mathrm{GV}}=\sum_jn_j.
\]
Then \(\mathscr A_B\) is finite free over \(B\) of rank \(d\), and \(HH_0(\mathscr A_B/B)\) and \(Z(\mathscr A_B)\) are finite free of rank \(N_{\mathrm{GV}}\).  For every commutative \(B\)-algebra \(B'\), the canonical maps
\begin{align*}
 HH_0(\mathscr A_B/B)\otimes_BB'
 &\xrightarrow{\sim}HH_0(\mathscr A_B\otimes_BB'/B'),\\
 Z(\mathscr A_B)\otimes_BB'
 &\xrightarrow{\sim}Z(\mathscr A_B\otimes_BB')
\end{align*}
are isomorphisms.  Moreover
\[
 \mathcal C_H\otimes_{R_H,F}^L B[[w]]\simeq\mathscr A_B[0].
\]
The conclusion holds for every complete Noetherian local base \(B\) and every geometric family specified in the hypotheses.
\end{proposition}

\begin{proof}
For the algebra, pull back \eqref{eq:weighted-contraction-resolution} along \(F\).  The closed-fibre cokernel is \(A_0\), so Corollary~\ref{cor:perfect-finite-slice} gives algebra freeness, rank \(d\), and the derived statement.  In terms of the divisor equations,
\[
 d=\operatorname{ord}_w\det\Theta_A(F_0)
   =\sum_jj^2\operatorname{ord}_wP_j(F_0)
   =\sum_jj^2n_j,
\]
where Corollary~\ref{cor:multicurve-cartan-gram}, applied with \(J=\{p\}\), identifies the orders along the possibly nonconvergent central classifying arc \(F_0\) with the same GV multiplicities. The comparison uses the intrinsic local Hilbert length and applies to both formal and analytic classifying arcs of the same contraction.

The cocenter is a commutator cokernel, so it commutes with pullback of \eqref{eq:downstairs-presentation}. This gives
\[
 HH_0(\mathscr A_B/B)\cong
 \operatorname{coker}\bigl(\Theta(F):D_B^q\to D_B^q\bigr).
\]
Here the relative cocenter can be computed over either \(B\) or \(D_B\): the \(B\)-span of commutators is already a \(D_B\)-module, since \(r[a,b]=[ra,b]\) for a central scalar \(r\). Lemma~\ref{lem:finite-slice-flatness} and \(\operatorname{ord}_w\prod_jP_j(F_0)=N_{\mathrm{GV}}\) make the cocenter free of rank \(N_{\mathrm{GV}}\).

To obtain the centre, recall that \(A_0\) is symmetric by Proposition~\ref{prop:contraction-algebra-symmetric}. Lemma~\ref{lem:relative-symmetrizing-trace} proves centre freeness, constructs the lifted symmetrizing trace, and proves centre base change.  Cocenter base change already follows from its cokernel presentation.
\end{proof}

The same freeness criterion applies to each individual root-divisor equation. The resulting rank is the intersection multiplicity along the classifying arc of the closed fibre.

\begin{corollary}
\label{cor:arbitrary-base-root-schemes}
For each \(j\), the scheme
\[
 \operatorname{Spec}\bigl(B[[w]]/(P_j(F))\bigr)
 \longrightarrow\Spec B
\]
is finite flat of degree \(n_j\).
\end{corollary}

\begin{proof}
Apply Lemma~\ref{lem:finite-slice-flatness} to the one-by-one matrix \(P_j(F)\), whose closed specialization has order \(n_j\).
\end{proof}

For several curves, the remaining task is to construct the universal presentations with every retained projector kept in the comparison.

\subsection{Idempotent subalgebras associated to several curves}
\label{sec:multicurve-flatness}

Let \(f_0:X_0\to Y_0\) be the completion at an exceptional image point of a projective flopping contraction of a smooth complex threefold, with any finite exceptional configuration.  We use the geometric families and contraction quotients of Remark~\ref{rem:hua-toda-support}. For one curve, each root determines a matrix size. For several curves, we also need the rank of each retained projector on that matrix block. On each generic matrix block, the Peirce dimensions are products of these projector ranks. The cocenter dimension of the generic semisimple algebra is the number of its simple factors.

\subsubsection{Restricted roots and the marked universal quotient}
\label{subsec:multicurve-marked-quotient}

Let \(\Delta\) be the ADE diagram of a general surface section, and let \(J\subseteq\Delta\) be the nonempty set of retained vertices. Write \(\mathfrak h\) for the corresponding Cartan space. We use the Katz--Morrison discriminant description in the restricted-root form of \citet[Theorem~2.6]{NabijouWemyss}. Here a retained vertex labels a curve that remains on the partial resolution.  The complementary set \(I\) labels the curves contracted from the minimal resolution.  Put
\[
 I=\Delta\setminus J,\qquad W_I=\langle s_i:i\in I\rangle,
 \qquad S=\widehat{\mathcal O}_{\mathfrak h,0},\qquad
 R_J=S^{W_I}=\widehat{\mathcal O}_{\mathfrak h/W_I,0}.
\]
Write \(S=\C[[\lambda_i:i\in\Delta]]\), and let \(q:\mathfrak h\to\mathfrak h/W_I\) denote the quotient map.  In the notation of \cite{NabijouWemyss}, the contracted subset is \(I\) and the retained subset is \(J\).  Each generator \(s_i\), \(i\in I\), is an order-two reflection on \(\mathfrak h\), and \(W_I\) is finite.  For this invariant-ring calculation, put \(H=\C[\mathfrak h]^{W_I}\). Chevalley's theorem makes \(H\) a polynomial ring and \(\C[\mathfrak h]\) a finite free \(H\)-module of rank \(|W_I|\) \cite[Theorem~(A), Lemma~3 and the proof of Theorem~(B)]{Chevalley1955}. This applies to disconnected complements, their fixed subspaces and components of types \(D\) and \(E\) as well.

Let \(\mathfrak m=H_{>0}\) and \(\mathfrak n=\C[\mathfrak h]_{>0}\). The quotient fibre over \(\mathfrak m\) is supported only at the origin: for a nonzero \(v\in\mathfrak h\), choose a linear form nonzero at all points of its finite \(W_I\)-orbit and take the product of its translates.  This is a positive-degree invariant nonzero at \(v\). Thus \(\mathfrak n^c\subseteq\mathfrak m\C[\mathfrak h]\subseteq \mathfrak n\) for some \(c\ge1\).  Completing the finite free module therefore gives \(S=\C[\mathfrak h]\otimes_H\widehat H_{\mathfrak m}\). The Reynolds projector \(|W_I|^{-1}\sum_{g\in W_I}g\) commutes with this scalar extension, so \(S^{W_I}=\widehat H_{\mathfrak m}=R_J\). Consequently both completed rings are regular and \(R_J\to S\) is finite free of rank \(|W_I|\), hence faithfully flat.

For \(\alpha=\sum_i m_i(\alpha)\alpha_i\in\Phi^+(\Delta)\), set
\[
 \pi_I(\alpha)=(m_j(\alpha))_{j\in J},\qquad
 \mathcal B_J=\pi_I(\Phi^+(\Delta))\setminus\{0\},\qquad
 s_\beta=\sum_{j\in J}\beta_j.
\]
Thus \(\pi_I(\alpha)\) retains exactly the coefficients relevant to the partial resolution, while \(s_\beta\) will be the size of the corresponding generic matrix block.  Proportional but distinct restricted roots remain distinct elements of \(\mathcal B_J\).

\Needspace{12\baselineskip}
\begin{lemma}
\label{lem:multicurve-root-products}
For \(\beta\in\mathcal B_J\), the product
\[
 p_\beta=\prod_{\pi_I(\alpha)=\beta}\ell_\alpha,
 \qquad \ell_\alpha=\sum_i m_i(\alpha)\lambda_i,
\]
is invariant under \(W_I\).  It descends to a reduced irreducible element \(P_\beta\in R_J\), whose zero locus is the curve divisor for the retained class \(\beta\) of the marked partial-resolution family.  Distinct \(\beta\)'s give distinct divisors, even when they are proportional.
\end{lemma}

\begin{proof}
The fibres of \(\pi_I\) on non-Levi positive roots are single \(W_I\)-orbits \cite[Lemma~2.3]{NabijouWemyss}.  For \(i\in I\), \(s_i\) preserves the retained coefficients and sends every positive root other than \(\alpha_i\) to a positive root.  It therefore permutes the factors of \(p_\beta\) without signs.

Those factors are pairwise nonassociated linear primes in the factorial ring \(S\).  A proper factorization in \(R_J\) would select a proper nonempty invariant subset of them, contrary to transitivity.  This proves irreducibility in the completed invariant ring \(R_J\).  The geometric interpretation of the descended divisor, including the distinction between proportional restricted roots, is \cite[Theorem~2.6 and Remark~2.7]{NabijouWemyss}.  Finally, at a general point of a non-Levi root hyperplane the \(W_I\)-stabilizer is trivial.  The quotient is therefore generically \'etale there, so the descended equation has multiplicity one.
\end{proof}

To identify the algebra supported on these divisors, let \(\mathcal X_J\to\mathcal Y_J\) be the completed universal family over \(R_J\).  Write the normalized generator and contraction quotient of this completed universal family as
\[
 T_J=\mathcal O\oplus\bigoplus_{j\in J}T_j,\qquad
 \Lambda_J=\End(T_J)^{\mathrm{op}},\qquad
 C_J=\Lambda_J/\Lambda_Je_0\Lambda_J.
\]
Let \(e_0\) and \(e_j\) be the projectors onto \(\mathcal O\) and \(T_j\), respectively. These projectors remain fixed throughout the comparison.  Put \(e_J=\sum_{j\in J}e_j\) and \(e_{0J}=e_0+e_J\).

On the closed surface, Iyama--Kalck--Wemyss--Yang identify the pullback of the tilting bundle from a partial resolution with the retained summands of the generator on the full resolution \cite[Theorem~4.6]{IKWY}.  Karmazyn's uniqueness theorem for bundle lifts extends the central comparison to the formal family.

For ordinary classifying pullbacks, also retain the algebraic model before completion.  Put \(H_J^{\mathrm{alg}}=\C[\mathfrak h]^{W_I}\), and write
\[
 \pi_J^{\mathrm{alg}}:\mathcal X_J^{\mathrm{alg}}
 \longrightarrow\Spec\mathcal R_J^{\mathrm{alg}}
\]
for the standard simultaneous partial resolution over this polynomial marked base.  Its labelled algebraic generator, constructed in the proof, is denoted \(T_J^{\mathrm{alg}}\), with \(\Lambda_J^{\mathrm{alg}}=\End(T_J^{\mathrm{alg}})^{\mathrm{op}}\). After completing the parameter ring, while leaving the target uncompleted at its central point, the model has
\[
 \mathcal R_J^{\mathrm{par}}
 =\mathcal R_J^{\mathrm{alg}}\otimes_{H_J^{\mathrm{alg}}}R_J,
 \qquad
 \Lambda_J^{\mathrm{par}}
 =\Lambda_J^{\mathrm{alg}}\otimes_{H_J^{\mathrm{alg}}}R_J.
\]
Completion at the maximal ideal of its central target point gives the previously defined \(\mathcal X_J,T_J,\Lambda_J\).

\begin{proposition}
\label{prop:multicurve-cartan-comparison}
There is an isomorphism preserving every labelled projector
\[
 \Lambda_J\widehat\otimes_{R_J}S
 \cong e_{0J}\widehat\Pi_{\widetilde\Delta,\lambda}e_{0J}.
\]
Killing the free projector gives
\begin{equation}
\label{eq:multicurve-cartan-quotient}
 C_J\otimes_{R_J}S\cong D_J:=e_J\Pi_{\Delta,\lambda}e_J.
\end{equation}
For every commutative Noetherian \(R_J\)-algebra \(B\), form the ordinary scheme pullback of the model over \(R_J\) whose target has not been completed at its central point, and write \(T_{J,B}^{\mathrm{ord}}\) for its pulled-back labelled generator. It remains a tilting generator, and the canonical marked algebra map gives
\[
 \End(T_{J,B}^{\mathrm{ord}})^{\mathrm{op}}
 \cong\Lambda_J^{\mathrm{par}}\otimes_{R_J}B.
\]
These assertions persist under central localization and ideal-adic completion of the resulting Noetherian affine target.  On a completed formal classifying pullback the full endomorphism algebra is the completed tensor product over \(R_J\). The contraction quotient is finite over \(R_J\), so its pullback is still given by the ordinary tensor product.
\end{proposition}

\begin{proof}
The algebraic model in \cite[Section~2.3]{Karmazyn} is defined for every coloured ADE diagram.  Its affine target is finite type over \(H_J^{\mathrm{alg}}\); source and target are flat over that base, the morphism is projective, its fibres have dimension at most one, and \(R\pi_*\mathcal O=\mathcal O\), as specified at the start of \cite[Section~3.1]{Karmazyn}.  The weighted ADE construction gives \(\C^*\)-actions on the affine target and parameter base whose coordinate rings are positively graded, with degree-zero part \(\C\).  The source is normal, since it is flat over the regular parameter base with geometrically normal partial-resolution surface fibres.  Locally this follows from Serre's criterion: base parameters provide the required depth, while a height-one point lies either in a normal generic fibre or over a height-one base point with a zero-dimensional normal fibre.

Since \(\mathcal X_J^{\mathrm{alg}}\) is projective over an affine variety, it is quasi-projective.  Choose an auxiliary ample line bundle \(\mathcal L\) on it.  By normality, a positive tensor power of \(\mathcal L\) admits a \(\C^*\)-linearization \cite[Lemma~4.10 and Remark~4.11]{KarmazynDeformations}. Replace \(\mathcal L\) by this power.  Its pullback to the completion at the central target point is relatively ample and linearized, as required by \cite[Appendix, Propositions~(A.5)--(A.6)]{NamikawaFlops}.

The normalized bundles on the central partial resolution have equivariant structures.  In the positive convention they are the torsion-free pullbacks of the retained graded special McKay modules, and the zero-perverse bundles are their duals; see \cite[Theorem~6.2 and Corollary~6.6]{KarmazynDeformations} and \cite[Theorem~4.6 and its proof]{IKWY}.

Apply \cite[Lemma~5.3 and Theorem~5.4]{KarmazynDeformations} to the algebraic equivariant family.  Here ``graded complete'' has the graded meaning of \cite[Definitions~4.1 and Remark~4.2]{KarmazynDeformations}, which includes the positively graded polynomial ring \(H_J^{\mathrm{alg}}\).  The positive gradings, normality, projectivity and flat affine target required by that theorem have just been checked. It gives an algebraic tilting bundle and an \(H_J^{\mathrm{alg}}\)-flat endomorphism algebra.

The summands and their inclusion and retraction maps may be lifted and algebraized equivariantly, preserving all labels: self-Ext groups in positive degrees vanish for each summand, taking \(\C^*\)-invariants is exact, and equivariant coherent algebraization preserves the specified maps \cite[Appendix, Lemma~(A.4), Propositions~(A.5)--(A.6)]{NamikawaFlops}. Thus the resulting bundle is \(T_J^{\mathrm{alg}}\), and uniqueness of the compatible lifts of each labelled summand identifies its completion with the fixed \(T_J\). The structure-sheaf summand is kept throughout.  Properness makes \(\Lambda_J^{\mathrm{alg}}\) finite over \(\mathcal R_J^{\mathrm{alg}}\).

Before taking a classifying pullback, we prove that the positive dual of this algebraic lift is globally generated and has relatively ample determinant.  For this calculation abbreviate
\[
 H=H_J^{\mathrm{alg}},\quad U=\mathcal R_J^{\mathrm{alg}},\quad
 X=\mathcal X_J^{\mathrm{alg}},\quad T=T_J^{\mathrm{alg}},\quad
 P=T^\vee,\quad E=\End_X(T).
\]
\phantomsection\label{par:multicurve-algebraic-positivity}
Let \(\mathfrak m=H_{>0}\), set \(U_c=U/\mathfrak mU\) and \(X_c=X\times_H\C\), and let \(o\in\Spec U\) be the homogeneous origin.  Write \(T_c=T|_{X_c}\) and \(P_c=T_c^\vee\).

The specialization in \cite[Theorem~5.4 and its proof]{KarmazynDeformations} is induced by restriction of endomorphisms.  It therefore preserves the fixed projector \(e_0\).  With composition multiplication on \(E\), its direct summand
\[
 M=\operatorname{Hom}_X(T,\mathcal O_X)=\Gamma(X,P)\cong e_0E
\]
satisfies
\[
 M\otimes_H\C\xrightarrow{\sim}\Gamma(X_c,P_c)
\]
by restriction of sections.  Properness makes \(M\) a finite \(U\)-module.  Consider the equivariant evaluation and its coherent cokernel:
\[
 \mathcal O_X\otimes_U M\longrightarrow P\longrightarrow
 \mathcal Q\longrightarrow0.
\]
Its restriction to \(X_c\) is precisely the evaluation of \(P_c\).

Let \(\mathfrak n_c\) be the homogeneous maximal ideal of \(U_c\), and let \(\widehat U_c=\widehat{(U_c)_{\mathfrak n_c}}\). Proper flat base change identifies
\[
 \Gamma(X_c,P_c)\otimes_{U_c}\widehat U_c
 \cong\Gamma(X_c\times_{U_c}\widehat U_c,\widehat P_c),
\]
compatibly with evaluation.  The specified labelled normalization on this completed \emph{central surface} identifies \(\widehat P_c\) with its positive generator, which is globally generated. Consequently \(\mathcal Q|_F=0\) for the scheme-theoretic fibre \(F=\pi^{-1}(o)\); Nakayama then makes \(\operatorname{Supp}\mathcal Q\) disjoint from \(F\). The image \(\pi(\operatorname{Supp}\mathcal Q)\) is a closed invariant subset of the affine target, because \(\pi\) is proper. Every nonempty such subset contains \(o\): the strictly positive grading extends the orbit of any closed point to \(\mathbb A^1\), with value \(o\) at zero.  Thus this image is empty and \(\mathcal Q=0\).  Choosing \(N\) generators of \(M\) over \(U\) gives a surjection \(\mathcal O_X^{\oplus N}\twoheadrightarrow P\).

For the total determinant \(L=\det P\), the central normalization gives
\[
 \deg_{C_i}L=\sum_{j\in J}\deg_{C_i}\det(T_j^\vee)
           =\sum_{j\in J}\delta_{ij}=1.
\]
The disjoint union of the projective-line normalizations of the retained curves maps finitely and surjectively to \(F\), through its reduction. The pullback of \(L|_F\) is ample, so \(L|_F\) is ample even if \(F\) is nonreduced \cite[Tag~\href{https://stacks.math.columbia.edu/tag/0B5V}{0B5V}]{Stacks}. For the proper morphism \(\pi\), fibrewise ampleness extends to relative ampleness over a target neighbourhood \cite[Tag~\href{https://stacks.math.columbia.edu/tag/0D2S}{0D2S}]{Stacks}. Thus the locus \(\{y\in\Spec U:L|_{X_y}\text{ is ample}\}\) is open. It is invariant and contains \(o\), so its closed invariant complement is empty.  Consequently \(L\) is \(\pi\)-ample.  This concerns the total determinant; the individual summand determinants can have degree zero on other retained curves. Finally,
\[
 H^1(X,T)=\operatorname{Ext}_X^1(\mathcal O_X,T)=0,
\]
as a direct summand of the higher self-Ext groups of \(T\). Global generation, relative ampleness and this vanishing, together with the structure-sheaf summand, verify the hypotheses of \cite[Propositions~3.2.6--3.2.7 and Corollary~3.2.8]{VanDenBergh}. Thus \(P\) and \(T\) are the \((-1)\)- and zero-perverse projective generators, respectively, with their bounded coherent derived equivalences.

We next compare the generators after extension from \(R_J\) to \(S=\widehat{\mathcal O}_{\mathfrak h,0}\). Let \((U,\mathfrak n)\) be the affine target ring completed at the maximal ideal of its central point. Since \(S\) is finite free over \(R_J\), the ring \(U_S=U\otimes_{R_J}S\) is finite over \(U\), with closed fibre \(S/\mathfrak m_{R_J}S\).  This fibre is local Artinian, so \(U_S\) has one maximal ideal \(\mathfrak n_S\), and for some \(c\ge1\),
\[
 \mathfrak n_S^c\subseteq\mathfrak nU_S\subseteq\mathfrak n_S.
\]
It is therefore complete for the maximal ideal \(\mathfrak n_S\).  The finite endomorphism modules are complete as well.  Thus \(\Lambda_J\otimes_{R_J}S\) is already complete, giving the completed tensor product in the proposition.

For parameter thickenings we use the ideal \(\mathfrak m_SU_S\) on this already completed target.  A complete Noetherian local ring is complete for every ideal contained in its maximal ideal: a compatible sequence of coset representatives is maximal-ideal Cauchy, and each ideal power is closed because its finite quotient is complete \cite[Tag~\href{https://stacks.math.columbia.edu/tag/00MA}{00MA}]{Stacks}. The same holds for finite modules.  This justifies parameter-adic formal existence here without identifying the parameter and target topologies.

On the central ADE surface, let \(q_0:Z_0\to X_{J,0}\) be the map from the minimal to the partial resolution.  The tilting-bundle construction in \cite[Theorem~4.6 and its proof]{IKWY} gives
\[
 q_0^*T_{J,0}\cong T_{\mathrm{full},0,0}
              \oplus\bigoplus_{j\in J}T_{\mathrm{full},0,j}.
\]
For the positive dual bundles, pull back the line bundles with degrees \(\delta_{ij}\) and their maximal extensions. The first Ext groups from the middle terms of these pulled-back extensions to \(\mathcal O\) vanish, so the extensions are maximal.  Dualizing gives the displayed zero-perverse convention. The nonfree summands therefore have determinant degrees \(-\delta_{ij}\).

Fix the full simultaneous resolution compatible with the selected standard partial resolution. In this comparison, write \(q\) for the morphism from the full resolution to the partial resolution after base change from \(R_J\) to \(S\). Rationality of the central partial resolution gives \(Rq_{0*}\mathcal O=\mathcal O\).  On a small parameter thickening, flatness identifies the structure-sheaf kernel with the central structure sheaf tensored by the square-zero kernel of the base.  Proper pushforward and the central higher-image vanishing therefore extend this identity through every Artin thickening.  Formal functions \cite[Tag~\href{https://stacks.math.columbia.edu/tag/02OC}{02OC}]{Stacks} give \(Rq_*\mathcal O=\mathcal O\) on the completed family.

The pullback of \(T_J\) and the retained full-resolution summands are now bundles on the same \(S\)-flat model with identified central restrictions. Vanishing of the first self-Ext extends each specified labelled isomorphism compatibly to the next thickening; its inverse lifts with it. Grothendieck existence and full faithfulness \cite[Tag~\href{https://stacks.math.columbia.edu/tag/088F}{088F}]{Stacks} produce the summand-preserving isomorphism on the completed model. These are chosen isomorphisms, determined up to automorphisms of the labelled summands.  Projection formula and \(Rq_*\mathcal O=\mathcal O\) then identify the endomorphism algebras of these bundles. The full universal McKay algebra is the extended-Dynkin deformed preprojective algebra \cite[Proposition~4.15]{Karmazyn}. Use its non-affine parameters \(\lambda_i=t_i\) and \(\lambda_0=-\sum_i d_i\lambda_i\), where \((1,d_i)\) is the minimal positive imaginary root.  Retaining the summands labelled by \(0\) and \(J\) gives the isomorphism between \(\Lambda_J\widehat\otimes_{R_J}S\) and the corner \(e_{0J}\widehat\Pi_{\widetilde\Delta,\lambda}e_{0J}\), with these coordinates.

For the unrestricted ordinary pullback, return to the finite-type model. The ring \(\mathcal R_J^{\mathrm{par}}\) is finite type over \(R_J\), so \(\mathcal R_J^{\mathrm{par}}\otimes_{R_J}B\) is Noetherian for \emph{every} Noetherian \(B\) in the statement.  Its projective source is Noetherian too.  For affine source and target rings \(R_J\to A\to Q\) of this model, flatness over \(R_J\) gives
\[
 Q\otimes_A^L(A\otimes_{R_J}B)
 \simeq Q\otimes_{R_J}^LB\simeq Q\otimes_{R_J}B.
\]
Derived base change \cite[Tag~\href{https://stacks.math.columbia.edu/tag/08IB}{08IB}]{Stacks}, applied to the endomorphism vector bundle, together with tilting and \(R_J\)-flatness of \(\Lambda_J^{\mathrm{par}}\), concentrates the comparison in degree zero.  The degree-zero map is the ordinary endomorphism pullback, and preserves multiplication and every projector. The same comparison for \(\mathcal O\) gives \(R\pi_{B*}\mathcal O=\mathcal O\).

The finite evaluation surjection for \((T_J^{\mathrm{alg}})^\vee\), proved on the algebraic model above, pulls back to a surjection onto \((T_{J,B}^{\mathrm{ord}})^\vee\). The determinant of \((T_{J,B}^{\mathrm{ord}})^\vee\) remains relatively ample, and the pulled-back bundle retains its fixed structure-sheaf summand.  Moreover, \(H^1(T_{J,B}^{\mathrm{ord}})=0\), since \(\operatorname{Ext}^1(\mathcal O,T_{J,B}^{\mathrm{ord}})\) is a direct summand of the vanishing group \(\operatorname{Ext}^1(T_{J,B}^{\mathrm{ord}},T_{J,B}^{\mathrm{ord}})\). Projectivity, the fibre-dimension bound and the structure-sheaf pushforward identity have all been checked on this new Noetherian target. The criterion and derived equivalence of \cite[Propositions~3.2.6--3.2.7 and Corollary~3.2.8]{VanDenBergh} therefore apply on this target and give both specified perverse projective generators and their bounded coherent derived equivalences.  No flatness, locality, reducedness or equivariance hypothesis on \(B\) is used.

Central localization and ideal-adic completion of the new Noetherian affine target are flat ring maps \cite[Tag~\href{https://stacks.math.columbia.edu/tag/00MB}{00MB}]{Stacks}. Proper flat base change identifies the endomorphism algebras of the pulled-back bundles with ordinary scalar extensions over the target ring.  For a completed formal classifying pullback, the same calculation on finite Artin parameter quotients and formal functions gives the completed tensor product of the full endomorphism algebra over \(R_J\).

Finally, scalar extension commutes with the idempotent quotient. For an algebra \(\Lambda\) carrying these projectors,
\[
 \frac{e_{0J}\Lambda e_{0J}}
      {e_{0J}\Lambda e_0\Lambda e_{0J}}
 \cong e_J(\Lambda/\Lambda e_0\Lambda)e_J.
\]
Removing the affine vertex leaves the finite-Dynkin corner. That quotient is finite over \(S\) by the central extension freeness in \cite[Theorem~3.2(i)]{EtingofRains}, followed by setting its central variable to zero.  Its closed parameter fibre has nilpotent positive-degree ideal, so the module topology defined by the parameter ideal agrees with that defined by the central target maximal ideal.  Exactness of Noetherian completion therefore leaves this already complete finite module unchanged, proving the ordinary tensor formula \eqref{eq:multicurve-cartan-quotient}. Faithful flatness of \(R_J\to S\) makes \(C_J\) finite over \(R_J\).  On the algebraic model this also follows by graded Nakayama: a homogeneous basis of the finite central quotient lifts to generators over \(H_J^{\mathrm{alg}}\), since each remaining term is a positive-degree parameter times an element of smaller weight, and the graded quotient is bounded below. Completing this algebraic quotient in the parameter ideal therefore gives \(C_J\). For subsequent pullbacks whose support is finite over the deformation base, the quotient is still given by an ordinary tensor product.
\end{proof}

\subsubsection{The presentation of the algebra and its Fitting ideal}
\label{subsec:multicurve-algebra-presentation}

Let \(\Pi_x\) be the universal deformed preprojective algebra of \(\Delta\), with central parameters \(x_i\). Adjoin a central polynomial variable \(z\).  Multiplication by \(z\) will give the square presentation of the retained corner; its determinant will give the divisor multiplicities.  Set \(\rho_i=1\), so \(\langle\rho,\alpha\rangle=\sum_i m_i(\alpha)\), and put \(\lambda_i=x_i-\rho_i z\) and
\[
 \mathscr E=\Pi_x[z]\otimes_{\C[\lambda_i]}S,\qquad \mathscr E_J=e_J\mathscr E e_J.
\]
Etingof--Rains prove freeness over the shifted parameter ring and semisimplicity over its generic point \cite[Theorem~3.2(i),(v)]{EtingofRains}.  The simple factors of the generic central extension are indexed by positive roots. The scalar by which \(z\) acts on a factor vanishes on its root hyperplane. After splitting the generic algebra, write \(V_\alpha\) for the simple module indexed by \(\alpha\). Taking the retained corner replaces \(\End(V_\alpha)\) by \(\End(e_JV_\alpha)\).

\Needspace{15\baselineskip}
\begin{lemma}
\label{lem:multicurve-root-blocks}
The module \(\mathscr E_J\) is finite free over \(S\), and
\begin{equation}
\label{eq:multicurve-z-resolution}
 0\longrightarrow \mathscr E_J\xrightarrow{z}\mathscr E_J\longrightarrow D_J
 \longrightarrow0
\end{equation}
is exact.  If \(K=\operatorname{Frac}S\), then
\begin{equation}
\label{eq:multicurve-generic-root-blocks}
 (\mathscr E_J)_{\overline K}
 \cong\prod_{\substack{\alpha\in\Phi^+\\\pi_I(\alpha)\ne0}}
 M_{s_{\pi_I(\alpha)}}(\overline K).
\end{equation}
On the factor indexed by \(\alpha\), the central element \(z\) acts as \(-\ell_\alpha/\langle\rho,\alpha\rangle\).
\end{lemma}

\begin{proof}
The whole extension is finite free over the shifted parameter ring \cite[Theorem~3.2(i)]{EtingofRains}. The corner projection \(a\mapsto e_Jae_J\) is \(S\)-linear and idempotent, so its image is projective and hence free over the local ring \(S\).  Multiplication by the polynomial variable \(z\) is injective before completion and remains injective after flat completion and on the corner.  Its cokernel is obtained by setting \(z=0\), which proves \eqref{eq:multicurve-z-resolution}.

The generic simple factors of the whole extension are \(\End(V_\alpha)\), with \(\dim e_iV_\alpha=m_i(\alpha)\) \cite[Theorem~3.2(v) and its proof]{EtingofRains}.  Orthogonality of the projectors gives
\[
 e_JV_\alpha=\bigoplus_{j\in J}e_jV_\alpha,
 \qquad\dim e_JV_\alpha=s_{\pi_I(\alpha)}.
\]
Taking the corner of \(\End(V_\alpha)\) gives \(\End(e_JV_\alpha)\).  Taking traces in the preprojective relation gives \(0=\ell_\alpha+z\langle\rho,\alpha\rangle\), determining the scalar action of \(z\).
\end{proof}

\begin{lemma}
\label{lem:multicurve-single-root-block}
Let \(S\to\Omega\) define a geometric point of \(\Spec S\), where \(\Omega\) is an algebraically closed extension of \(\C\).  Suppose exactly one positive root \(\alpha\) has \(\ell_\alpha=0\) at this point, and put \(\beta=\pi_I(\alpha)\ne0\).  Then
\[
 D_J\otimes_S\Omega
 \cong\End_\Omega\left(\bigoplus_{j\in J}\Omega^{\beta_j}\right),
\]
with each labelled projector \(e_j\) of rank \(\beta_j\). This is also the geometric fibre of \(C_J\) at the image point in \(\Spec R_J\), with the same labelled projectors.
\end{lemma}

\begin{proof}
We first work with the full central extension before taking the retained corner.  Write \(h_\gamma=\langle\rho,\gamma\rangle\) and \(c_\gamma=-\ell_\gamma/h_\gamma\) for each positive root \(\gamma\). Freeness and the generic decomposition in Lemma~\ref{lem:multicurve-root-blocks}, applied to the full vertex set, show that
\[
 \prod_{\gamma\in\Phi^+}(z-c_\gamma)=0
\]
in the central extension: the identity holds over \(\operatorname{Frac}S\) and therefore over \(S\), by freeness.

Work on \(H_\alpha=(\ell_\alpha=0)\) inside \(\Spec S\), and invert the equation of every other positive root.  On this integral open set put
\[
 q_\alpha(Z)=\prod_{\gamma\ne\alpha}(Z-c_\gamma),\qquad
 \varepsilon_\alpha=\frac{q_\alpha(z)}{q_\alpha(0)}.
\]
The denominator is a unit.  Since \(zq_\alpha(z)=0\) and \(q_\alpha(Z)-q_\alpha(0)\) is divisible by \(Z\), \(\varepsilon_\alpha\) is a central idempotent and \(z\varepsilon_\alpha=0\).  B\'ezout for \(Z\) and \(q_\alpha(Z)\) shows that \(z\) is invertible on the complementary summand. Thus setting \(z=0\) leaves precisely the \(\varepsilon_\alpha\)-summand, a finite projective module on this integral open set.

The \(\varepsilon_\alpha\)-summand has rank \(h_\alpha^2\).  To compute it at the generic point of \(H_\alpha\), write \(\lambda=v+t\rho\), with \(v\in H_\alpha\) and \(t=\ell_\alpha(\lambda)/h_\alpha\). The change \(z_\alpha=z+t\) identifies the full central extension with the extension restricted to \(H_\alpha\), followed by adjoining this parameter.  In its generic decomposition, \(z_\alpha\) acts on the \(\gamma\)-block by \(-\ell_\gamma(v)/h_\gamma\). These functions are distinct on \(H_\alpha\): if \(\gamma/h_\gamma-\eta/h_\eta\) were a multiple of \(\alpha\), pairing with \(\rho\) would make that multiple zero, forcing \(\gamma=\eta\).  Only the \(\alpha\)-block has scalar zero, and its dimension is \(h_\alpha^2\).  The rank of the projective summand is constant on the integral open set, so every fibre of the finite-Dynkin algebra \(\Pi_\lambda\) at which only \(\pm\alpha\) vanish has this dimension.

Crawley--Boevey's simple-module criterion \cite[Theorem~1.2]{CrawleyBoeveyMomentMap} supplies a simple \(\Pi_\lambda\)-module of dimension vector \(\alpha\). Indeed \(\alpha\) is positive, \(\lambda\cdot\alpha=0\), and there is no decomposition of \(\alpha\) into two or more positive roots all annihilated by \(\lambda\): each such root would have to equal \(\alpha\).  This verifies the criterion's strict inequality condition vacuously; in Dynkin type its function \(p(\gamma)=1-(\gamma,\gamma)/2\) is zero on every root. Since the ground field is algebraically closed, the action on this simple module maps \(\Pi_\lambda\) onto its full endomorphism algebra.  Both algebras have dimension \(h_\alpha^2\), so the map is an isomorphism.  Its vertex projectors have ranks \(m_i(\alpha)\) by the dimension-vector statement. Taking the retained corner gives the displayed marked matrix algebra. Finally, Proposition~\ref{prop:multicurve-cartan-comparison} identifies this corner with the corresponding fibre of \(C_J\), preserving every projector.
\end{proof}

\begin{theorem}
\label{thm:multicurve-weighted-presentation}
There is a free resolution of the form
\begin{equation}
\label{eq:multicurve-weighted-resolution}
 0\longrightarrow R_J^a\xrightarrow{\Theta_J}R_J^a
 \longrightarrow C_J\longrightarrow0
\end{equation}
with
\begin{equation}
\label{eq:multicurve-weighted-fitting}
 \operatorname{Fitt}^{R_J}_0(C_J)
 =\left(\prod_{\beta\in\mathcal B_J}P_\beta^{s_\beta^2}\right).
\end{equation}
For \(i,j\in J\), the module \(e_iC_Je_j\) also has a free resolution with two free modules of equal rank, and
\[
 \operatorname{Fitt}^{R_J}_0(e_iC_Je_j)
 =\left(\prod_{\beta\in\mathcal B_J}P_\beta^{\beta_i\beta_j}\right).
\]
For the zero module, the zeroth Fitting ideal is the unit ideal. The nonzero module \(C_J\) is Cohen--Macaulay of codimension one.
\end{theorem}

\begin{proof}
Over the algebraic closure of the fraction field, the factor belonging to \(\alpha\) in \eqref{eq:multicurve-generic-root-blocks} has algebra dimension \(s_{\pi_I(\alpha)}^2\).  Therefore
\[
 \det(z\mid \mathscr E_J)
 =u\prod_{\pi_I(\alpha)\ne0}\ell_\alpha^{s_{\pi_I(\alpha)}^2}
 =u\prod_\beta(q^*P_\beta)^{s_\beta^2},\qquad u\in\C^\times.
\]
The determinant is nonzero, so \(D_J\) is a finite torsion \(S\)-module with a free resolution consisting of two free modules of equal rank.  Choose a finite free surjection onto \(C_J\). After base change, its target is \(D_J\) by \eqref{eq:multicurve-cartan-quotient}. Comparing this surjection with \eqref{eq:multicurve-z-resolution} by Schanuel's lemma makes its kernel projective, as in the proof of Proposition~\ref{prop:weighted-contraction-presentation}.  Faithfully flat descent \cite[Tag~\href{https://stacks.math.columbia.edu/tag/058S}{058S}]{Stacks} makes the original kernel projective over \(R_J\), hence free because \(R_J\) is local; torsion forces the two ranks to agree.  This proves \eqref{eq:multicurve-weighted-resolution}.  Base change and faithful contraction of Fitting ideals give \eqref{eq:multicurve-weighted-fitting}.

For the Peirce entries, restrict multiplication by \(z\) to the free direct summand \(e_i\mathscr E e_j\).  The cokernel of this restricted map is \(e_iD_Je_j\), whose square presentation descends by the same syzygy construction.  The \(\alpha\)-summand of the generic fibre of \(e_i\mathscr E e_j\) is \(\operatorname{Hom}(e_jV_\alpha,e_iV_\alpha)\), of dimension \(m_i(\alpha)m_j(\alpha)=\beta_i\beta_j\), giving the entrywise formula. The Fitting ideal shows that the support has codimension one. Regularity of \(R_J\) and Auslander--Buchsbaum \cite[Proposition~10.111.1, Tag~\href{https://stacks.math.columbia.edu/tag/090V}{090V}]{Stacks} give \(\operatorname{depth}_{R_J}C_J=\dim R_J-1=\dim C_J\).
\end{proof}

\begin{proposition}
\label{prop:multicurve-reduced-annihilator}
Put \(\delta_J=\prod_{\beta\in\mathcal B_J}P_\beta\).  Then
\[
 \operatorname{Ann}_{R_J}(C_J)=(\delta_J).
\]
There is a matrix \(\Xi_J\) such that \(\Theta_J\Xi_J=\Xi_J\Theta_J=\delta_J I_a\).
\end{proposition}

\begin{proof}
Consider the polynomial
\[
 Q(Z)=\prod_{\pi_I(\alpha)\ne0}
       \bigl(\langle\rho,\alpha\rangle Z+\ell_\alpha\bigr)
\]
Each generic matrix factor is annihilated by its corresponding linear factor in this product, by Lemma~\ref{lem:multicurve-root-blocks}. Thus \(Q\) annihilates \(\mathscr E_J\) over the fraction field.  Freeness over the domain \(S\) therefore implies \(Q(z)\mathscr E_J=0\).  Setting \(z=0\) gives \((q^*\delta_J)D_J=0\).  The Fitting ideal shows that the support is exactly the union of the indicated hyperplanes, whose reduced equation is \(q^*\delta_J\).  Hence
\[
 (q^*\delta_J)\subseteq\operatorname{Ann}_S(D_J)
 \subseteq\sqrt{\operatorname{Ann}_S(D_J)}=(q^*\delta_J).
\]
Faithful flatness gives the annihilator equality over \(R_J\).  Since \(\delta_J\) annihilates the cokernel of \(\Theta_J\), lift \(\delta_J I_a\) through \(\Theta_J\) to obtain \(\Xi_J\). For the reverse product, apply \(\Theta_J\) to \(\Xi_J\Theta_J-\delta_J I_a\). The first identity makes the result zero, and injectivity of \(\Theta_J\) gives the second identity.
\end{proof}

\subsubsection{An exact sequence for cocenters after deleting vertices}
\label{subsec:multicurve-trace-excision}

Each matrix factor contributes one dimension to the cocenter.  To obtain the cocenter presentation over \(S=\widehat{\mathcal O}_{\mathfrak h,0}\), we compute the dimension of the cocenter of the retained corner at \(\lambda=0\).  Write \(P_\Delta=\Pi^\rho_0\) for the ordinary central extension specialized at \(\lambda=0\), with \(\rho_i=1\); the element \(z\) remains part of this algebra.  Put
\[
 Q_I=P_\Delta/P_\Delta e_JP_\Delta.
\]
The quotient deletes the retained vertices, leaving the components of \(I\).  Etingof--Latour--Rains compute the cocenter dimension of the full central extension for each of these diagrams \cite[Theorem~3.1]{EtingofLatourRains}. The following exact sequence identifies the cocenter of the retained corner with the kernel of \(HH_0(P_\Delta)\to HH_0(Q_I)\).

\Needspace{15\baselineskip}
\begin{theorem}
\label{thm:all-corner-trace-excision}
For every nonempty \(J\subseteq\Delta\), there is an exact sequence
\begin{equation}
\label{eq:all-corner-trace-excision}
 0\longrightarrow HH_0(e_JP_\Delta e_J)
 \longrightarrow HH_0(P_\Delta)
 \longrightarrow HH_0(Q_I)\longrightarrow0.
\end{equation}
In particular, if \(N_J=|\Phi^+(\Delta)|-|\Phi^+(I)|\), then
\begin{equation}
\label{eq:all-corner-trace-dimension}
 \dim_\C HH_0(e_JP_\Delta e_J)=N_J.
\end{equation}
With arrows of internal degree one and \(z\) of degree two, the odd pieces of \(HH_0(e_JP_\Delta e_J)\) vanish and
\[
 \dim HH_0(e_JP_\Delta e_J)_{2p}
 =\#\{\alpha\in\Phi^+(\Delta):\pi_I(\alpha)\ne0,
                         \operatorname{ht}(\alpha)=p+1\}.
\]
\end{theorem}

The main point is injectivity of the corner map.  In the exact sequence of Lemma~\ref{lem:vertex-deletion-traces}, the kernel of that corner map is the image of the boundary from \(HH_1(Q_I)\).  This boundary vanishes once \(HH_1(P_\Delta)\to HH_1(Q_I)\) is surjective. We prove surjectivity by lifting Hochschild cycles involving the image of the central element \(z\).

Keller's localization theorem gives a sequence for dg Hochschild homology.  The following comparison relates the low-degree dg Hochschild groups in that sequence to those of the ordinary Jacobi algebras.  For these ordinary Jacobi algebras we use polynomial Ginzburg dg algebras, with cohomological grading for dg algebras and homological grading for Hochschild homology.

\begin{lemma}
\label{lem:ginzburg-low-degree-comparison}
Let \(\Gamma\) be a nonpositive, homologically smooth, bimodule \(3\)-Calabi--Yau dg algebra over \(\C\), with \(A=H^0(\Gamma)\) finite-dimensional, basic, and self-injective.  Then \(H^{-1}(\Gamma)=0\), and the canonical projection gives natural isomorphisms
\[
 HH_i(\Gamma)\cong HH_i(A),\qquad i=0,1.
\]
For any nonpositive dg algebra, the degree-zero comparison holds and its negative Hochschild homology vanishes.
\end{lemma}

\begin{proof}
Write \(\mathcal D(\Gamma)\) for the derived category of right dg modules. Its full subcategories of perfect modules and of modules with finite-dimensional total cohomology are denoted \(\operatorname{per}\Gamma\) and \(\mathcal D_{\mathrm{fd}}\Gamma\), respectively. Amiot's theorem makes \(\mathcal C=\operatorname{per}\Gamma/\mathcal D_{\mathrm{fd}}\Gamma\) a Hom-finite \(2\)-Calabi--Yau category.  The image \(T\) of \(\Gamma\) is cluster tilting, with endomorphism algebra \(A\), \cite[Theorem~2.1]{Amiot}. The projection to \(\mathcal C\) is fully faithful on the fundamental domain
\[
 \mathcal F=\operatorname{per}\Gamma\cap\mathcal D^{\le0}
                \cap{}^\perp\mathcal D^{\le-2}
 \qquad\text{\cite[Proposition~2.9]{Amiot}}.
\]

Both \(\Gamma\) and \(\Gamma[1]\) lie in this domain.  Indeed, they are perfect and have cohomology in nonpositive degrees.  For \(M\in\mathcal D^{\le-2}\), the two groups \(\operatorname{Hom}_{\mathcal D(\Gamma)}(\Gamma,M)=H^0(M)\) and \(\operatorname{Hom}_{\mathcal D(\Gamma)}(\Gamma[1],M)=H^{-1}(M)\) vanish, giving the required left orthogonality.

Since \(A\) is basic and self-injective, the criterion of \citet[Proposition~3.6]{IyamaOppermann}, in the formulation \cite[Proposition~4.4]{HerschendIyama}, gives \(T\cong T[2]\).  Hence
\[
 H^{-1}(\Gamma)
 =\operatorname{Hom}_{\operatorname{per}\Gamma}(\Gamma[1],\Gamma)
 \cong\operatorname{Hom}_{\mathcal C}(T[1],T)
 \cong\operatorname{Hom}_{\mathcal C}(T,T[1])=0.
\]
The last equality holds because \(T\) is cluster tilting.

For the Hochschild comparison, use the reversed internal grading \(H_i\Gamma=H^{-i}(\Gamma)\) and filter the Hochschild bicomplex by chain length.  Over \(\C\), this produces the first-quadrant spectral sequence
\[
 E^2_{p,q}=HH_p(H_*\Gamma)_q\Longrightarrow HH_{p+q}(\Gamma).
\]
Under this reversed convention, the vanishing just proved is \(H_1\Gamma=0\), so the row of internal degree one is zero.  In total degrees zero and one, the only possible terms are \(HH_0(H_0\Gamma)\) and \(HH_1(H_0\Gamma)\), at \((0,0)\) and \((1,0)\), and no higher differential enters or leaves them.  The edge maps are induced by the canonical projection, which proves naturality and both comparisons.  The degree-zero statement does not use the \(H_1\)-vanishing.  Finally, there are no chains of negative total degree, so negative Hochschild homology vanishes.
\end{proof}

\begin{lemma}
\label{lem:vertex-deletion-traces}
Let \(P=\operatorname{Jac}(Q,W)\) be the ordinary Jacobi algebra of a finite quiver with polynomial potential.  Delete a set of vertices, with sum of projectors \(e\), and denote the resulting Jacobi algebra by \(Q'=P/PeP\).  Suppose \(P\) and \(Q'\) are finite-dimensional, basic, and self-injective.  There is a natural exact sequence
\begin{equation}
\label{eq:vertex-deletion-five-term}
 HH_1(P)\longrightarrow HH_1(Q')\longrightarrow HH_0(ePe)
 \longrightarrow HH_0(P)\longrightarrow HH_0(Q')\longrightarrow0.
\end{equation}
The map from the corner sends \([a]\) to \([a]\).
\end{lemma}

\begin{proof}
Let \(\Gamma\) and \(\Gamma'\) denote the polynomial Ginzburg dg algebras before and after the deletion.  Both are nonpositive, homologically smooth, and bimodule \(3\)-Calabi--Yau \cite[Theorem~6.3]{KellerCY}.  The path dg categories of these finite quivers are finitely cellular, and the potential supplies the negative cyclic lift of the deformation class required by the corrected theorem \cite[Theorem~2.1]{KellerCYErratum}. Vertex deletion is a dg localization \cite[Corollary~7.3]{KellerCY}. Its derived kernel is the localizing subcategory \(\operatorname{Loc}(e\Gamma)\), by \cite[Theorem~5.8(b)]{KellerCY}.  For several deleted vertices, compose the single-vertex localizations: each next deleted representable is the image of its original representable, so the composite kernel is generated by their union.

Put \(E=e\Gamma e\).  The bimodule \(e\Gamma\) gives an adjunction
\[
 -\otimes_E^L e\Gamma:
 \mathcal D(E)\rightleftarrows\mathcal D(\Gamma):
 R\operatorname{Hom}_\Gamma(e\Gamma,-).
\]
The right adjoint is \((-)e\) and preserves coproducts, since \(e\Gamma\) is a summand of the free right \(\Gamma\)-module. The objects on which the adjunction unit is an isomorphism form a localizing subcategory. This subcategory contains the generator \(E\), so it is all of \(\mathcal D(E)\).  Thus the left adjoint is fully faithful with image \(\operatorname{Loc}(e\Gamma)\). Consequently \(E\to\Gamma\to\Gamma'\) is an exact sequence in the Morita category by \cite[Theorem~4.11(iii)]{KellerDG}.

Morita invariance and localization of Hochschild mixed complexes \cite[Theorem~5.2]{KellerDG} give a triangle
\[
 C_*(e\Gamma e)\longrightarrow C_*(\Gamma)
 \longrightarrow C_*(\Gamma')\longrightarrow C_*(e\Gamma e)[1].
\]
The dg hypotheses of Lemma~\ref{lem:ginzburg-low-degree-comparison} have been verified for \(\Gamma\) and \(\Gamma'\), and their degree-zero cohomology algebras \(P\) and \(Q'\) are finite-dimensional, basic and self-injective by assumption. The lemma therefore identifies the two degree-one terms with \(HH_1(P)\) and \(HH_1(Q')\), and their degree-zero terms with the ordinary cocenters. For \(e\Gamma e\), only the general nonpositive case of that lemma is needed: \(H^0(e\Gamma e)=ePe\) and \(HH_{-1}(e\Gamma e)=0\).  Thus the corresponding segment of the long exact sequence is \eqref{eq:vertex-deletion-five-term}. Under Morita invariance, the trace of an endomorphism of the projective \(e\Gamma\) is its cyclic class in \(\Gamma\).  The corner map is therefore the stated inclusion on cyclic classes, and is induced by the bimodule \(e\Gamma\), even though the corner inclusion need not be unital.
\end{proof}

\begin{lemma}
\label{lem:central-cycle-generation}
Let \(Q\) be a finite-dimensional symmetric complex algebra and \(c\in Z(Q)\).  If every derivation of \(Q\) annihilating \(c\) is inner, then the map
\[
 Q\longrightarrow HH_1(Q),\qquad q\longmapsto[q\otimes c]
\]
is surjective.  If \(P\twoheadrightarrow Q\) lifts \(c\) to a central element \(z\in P\), then \(HH_1(P)\to HH_1(Q)\) is surjective.
\end{lemma}

\begin{proof}
The cycles in the statement are available because the Hochschild boundary of \(q\otimes c\) is \(qc-cq=0\).  A symmetrizing trace \(\tau\) identifies the dual chain complex with the Hochschild cochain complex and gives a perfect pairing
\[
 HH_1(Q)\times HH^1(Q)\longrightarrow\C,
 \qquad\langle[q\otimes c],[D]\rangle=\tau(qD(c)).
\]
Inner derivations kill \(c\), so this formula depends only on the cohomology class of \(D\).  If a class pairs trivially with every displayed cycle, then \(\tau(qD(c))=0\) for every \(q\in Q\). Nondegeneracy of the symmetrizing form gives \(D(c)=0\).  The hypothesis then makes \(D\) inner, so its cohomology class is zero.  Finite-dimensional duality proves the asserted surjectivity. For the second assertion, a lift \(a\) of \(q\) gives the cycle \(a\otimes z\), which maps to \(q\otimes c\).
\end{proof}

\begin{proof}[Proof of Theorem~\ref{thm:all-corner-trace-excision}]
For the vertex-deletion sequence, use the following polynomial potential.  Fix an orientation \(Q\) of \(\Delta\), use left-to-right path composition, and put
\[
 r_i=e_i\left(\sum_a[a,a^*]\right)e_i,
 \qquad W_\Delta=\frac12\sum_i r_i^2.
\]
For an arrow \(a:i\to j\), its two cyclic derivatives are
\[
 \partial_{a^*}W_\Delta=r_i a-a r_j,
 \qquad\partial_aW_\Delta=a^*r_i-r_j a^*.
\]
It follows that in the Jacobi algebra \(z=\sum_i r_i\) commutes with every arrow and that \(r_i=ze_i\).  Conversely, these central-extension relations recover both displayed cyclic derivatives. Fixing the arrows gives homomorphisms in both directions, with \(z\mapsto\sum_i r_i\) from the central extension. These maps of polynomial algebras are inverse, and hence
\[
 \operatorname{Jac}(\overline Q,W_\Delta)\cong P_\Delta.
\]
Deleting \(J\) removes the incident arrows.  The potential restricts to the sum of the analogous potentials on the connected components \(I_\nu\) of \(I\), and consequently
\[
 Q_I\cong\prod_{I_\nu}P_{I_\nu},
\]
compatibly with the common central element \(z\).

If \(I\) is empty, the asserted sequence is the identity sequence on the full algebra. Assume henceforth that \(I\) is nonempty. For every connected ADE diagram, Etingof--Latour--Rains prove symmetry \cite[Proposition~2.1]{EtingofLatourRains}, the cocenter dimension formula \cite[Theorem~3.1]{EtingofLatourRains}, and that derivations annihilating \(z\) are inner \cite[Corollary~3.5 and the following remark]{EtingofLatourRains}. The normalization \(\rho_i=1\) induces the same normalization on every component.  The one-vertex algebra is \(\C\), where the corresponding claims are verified directly.  Hence \(P_\Delta\) and \(Q_I\) are finite, basic, symmetric Jacobi algebras. In particular they are self-injective, so Lemma~\ref{lem:vertex-deletion-traces} applies to this deletion. A derivation of a finite product preserves its central idempotents, so the condition on derivations holds componentwise for \(Q_I\).  Apply Lemma~\ref{lem:central-cycle-generation} to \(P_\Delta\twoheadrightarrow Q_I\), with \(c\) the image of \(z\).  It gives \(HH_1(P_\Delta)\twoheadrightarrow HH_1(Q_I)\). In \eqref{eq:vertex-deletion-five-term}, this makes the boundary from \(HH_1(Q_I)\) vanish.  The remaining sequence is \eqref{eq:all-corner-trace-excision}.

The cocenter dimensions of the full algebras in \cite[Theorem~3.1]{EtingofLatourRains} are the numbers of positive roots, and degree \(2p\) counts the roots of height \(p+1\). Subtract the formula for the components of \(I\) from the formula for \(\Delta\).  The roots removed in this subtraction are precisely those with zero retained coefficients, equivalently those supported on \(I\). This gives both asserted dimensions.  Since the ordinary corner and quotient maps preserve the internal grading, the same subtraction gives the degreewise statement.
\end{proof}

The dimension of the corner cocenter is therefore the number of positive roots with at least one nonzero retained coefficient.  This also equals the generic cocenter rank of the retained central extension.

\subsubsection{The universal cocenter and its reduced divisor}
\label{subsec:multicurve-universal-cocenter}

\begin{theorem}
\label{thm:multicurve-universal-cocenter}
The module \(M_J=HH_0(\mathscr E_J/S)\) is free of rank \(N_J\). For \(\mathcal H_J=HH_0(C_J/R_J)\), there is a free resolution of the form
\begin{equation}
\label{eq:multicurve-cocenter-resolution}
 0\longrightarrow R_J^b\xrightarrow{\Theta_{\mathrm{tr},J}}R_J^b
 \longrightarrow\mathcal H_J\longrightarrow0
\end{equation}
with
\begin{equation}
\label{eq:multicurve-cocenter-fitting}
 \operatorname{Fitt}^{R_J}_0(\mathcal H_J)=(\delta_J)
 =\left(\prod_{\beta\in\mathcal B_J}P_\beta\right).
\end{equation}
\end{theorem}

\begin{proof}
Because \(M_J\) is a commutator cokernel, it commutes with scalar extension.  Its closed fibre is \(HH_0(e_JP_\Delta e_J)\), whose dimension is \(N_J\) by Theorem~\ref{thm:all-corner-trace-excision}. On the other hand, \eqref{eq:multicurve-generic-root-blocks} gives generic rank \(N_J\), since each matrix algebra has a one-dimensional cocenter.  Lift a basis of the closed fibre and use Nakayama to obtain \(S^{N_J}\twoheadrightarrow M_J\).  Its kernel is a submodule of a free module over the domain \(S\); its generic fibre is zero, so the kernel itself is zero.  Thus \(M_J\) is free over \(S\).

Since \(z\) is central,
\[
 HH_0(D_J/S)=M_J/zM_J.
\]
On the factor indexed by \(\alpha\), the scalar action of \(z\) is \(-\ell_\alpha/\langle\rho,\alpha\rangle\), by Lemma~\ref{lem:multicurve-root-blocks}. The cocenter of that factor is one-dimensional, so this eigenvalue occurs once. Therefore
\[
 \det(z\mid M_J)
 =u\prod_{\pi_I(\alpha)\ne0}\ell_\alpha
 =u\,q^*\delta_J,\qquad u\in\C^\times.
\]
The determinant is nonzero. Multiplication by \(z\) is therefore injective on the free module \(M_J\), giving a square presentation of \(HH_0(D_J/S)\). By \eqref{eq:multicurve-cartan-quotient}, this module is \(\mathcal H_J\otimes_{R_J}S\).  As in the proof of Theorem~\ref{thm:multicurve-weighted-presentation}, the kernel of a finite free surjection onto \(\mathcal H_J\) becomes projective after faithfully flat base change.  Faithfully flat descent \cite[Tag~\href{https://stacks.math.columbia.edu/tag/058S}{058S}]{Stacks} makes the original kernel projective over \(R_J\), hence free because \(R_J\) is local.  Since \(\mathcal H_J\) is torsion, the two free ranks agree, proving \eqref{eq:multicurve-cocenter-resolution}.  Fitting ideals commute with base change, and their equality over \(S\) descends because \(R_J\to S\) is faithfully flat. This gives \eqref{eq:multicurve-cocenter-fitting}.
\end{proof}

The product defining \(P_\beta\) contains one linear factor for each positive root restricting to \(\beta\).  The Fitting exponents of the three presentations are
\[
\begin{array}{c|c}
\text{universal module}&\text{exponent of }P_\beta\\ \hline
 e_iC_Je_j&\beta_i\beta_j\\
 C_J&s_\beta^2\\
 HH_0(C_J/R_J)&1
\end{array}
\]
For a family of classifying arcs, these exponents determine the ranks over the deformation base: restrict the determinants to the arc of the closed fibre and take their orders.

\subsubsection{Expressing a geometric family as a pullback}
\label{subsec:multicurve-geometric-flatness}

Let \((B,\mathfrak m,\C)\) be complete Noetherian local, and let \(f_B:X_B\to Y_B\) be a deformation of \(f_0\) as a projective contraction, with source and affine target flat over \(B\).  We complete the target at the central contraction point and the source along its inverse image, retaining every support point of the finite relative quotient. Write \(T_0=\mathcal O\oplus\bigoplus_{j\in J}T_{0,j}\) for the fixed normalized basic generator.

The variable \(w\) is the surface-section parameter; \(\Spec B\) is the base of the external deformation.

\begin{proposition}
\label{prop:multicurve-classifying-pullback}
There is a continuous local map \(F:R_J\to D_B=B[[w]]\) and a relative tilting generator \(T_B\) lifting \(T_0\), such that the relative contraction quotient satisfies
\begin{equation}
\label{eq:multicurve-family-identification}
 \mathscr A_B:=\End(T_B)^{\mathrm{op}}/
 \End(T_B)^{\mathrm{op}}e_0\End(T_B)^{\mathrm{op}}
 \cong C_J\otimes_{R_J,F}B[[w]].
\end{equation}
All labelled projectors are preserved.  Any other lift of the same marked generator gives an isomorphic marked quotient, and \(\mathscr A_B/\mathfrak m\mathscr A_B\cong A_0\).
\end{proposition}

\begin{proof}
Lift a general surface-section parameter \(w_0\) to \(w\). The general section \(w_0\) is a nonzerodivisor on the smooth central source and on the cDV target. The central source \(X_0\) and target \(Y_0\) are therefore flat over \(\C[[w]]\). As in Proposition~\ref{prop:arbitrary-base-classifying-arcs}, this closed-fibre flatness and the assumed \(B\)-flatness give \(D_B\)-flat source and target by the fibrewise criterion \cite[Tag~\href{https://stacks.math.columbia.edu/tag/00MP}{00MP}]{Stacks}. The closed surface is the partial resolution with retained set \(J\). The source-deformation argument in Proposition~\ref{prop:arbitrary-base-classifying-arcs} applies to every nonempty retained subgraph, with base \(\mathfrak h/W_I\). It constructs compatible classifying maps and source isomorphisms on successive Artin quotients, then recovers the target from global functions on the source. Passing to the limit gives \(F\) and the marked contraction isomorphism.

Proposition~\ref{prop:multicurve-cartan-comparison} fixes the normalized labelled universal generator.  We first identify its central restriction with \(T_0\). The restriction to the central threefold is a basic generator: modulo \(w\), these summands are the distinct indecomposable projectives on the partial surface.  A nontrivial decomposition of one such summand would restrict to a nontrivial decomposition on that partial surface, because a summand whose reduction is zero vanishes by Nakayama in the completed neighbourhood.  The determinant degrees distinguish the curve labels. Van den Bergh's classification of normalized perverse projectives \cite[Proposition~3.5.4 and Theorem~3.5.5]{VanDenBergh} consequently identifies the central generator with \(T_0\).

The higher self-Ext groups of these summands vanish.  Karmazyn's lifting lemma \cite[Lemma~3.2]{KarmazynDeformations} identifies their lifts through every external-base thickening and in the formal limit while preserving the chosen projectors. Put \(D_n=D_B/(\mathfrak m,w)^n\).  The universal source, target, and endomorphism algebra are flat over \(R_J\).  The Tor-independence calculation in Proposition~\ref{prop:multicurve-cartan-comparison} and derived base change \cite[Tag~\href{https://stacks.math.columbia.edu/tag/08IB}{08IB}]{Stacks} therefore identify the endomorphism algebra on the \(n\)-th pullback with \(\Lambda_J\otimes_{R_J}D_n\), preserving multiplication and every projector.  Taking the idempotent quotient gives \(C_J\otimes_{R_J}D_n\).

The endomorphism algebra and its projector ideal are finite modules over the Noetherian completed affine target.  Formal functions and exactness of completion on finite modules \cite[Tags~\href{https://stacks.math.columbia.edu/tag/02OC}{02OC} and~\href{https://stacks.math.columbia.edu/tag/00MA}{00MA}]{Stacks} identify the geometric quotient with \(\varprojlim_n(C_J\otimes_{R_J}D_n)\). Since \(C_J\) is finite over \(R_J\), its ordinary pullback is a finite module over the complete ring \(D_B\), and this limit is \(C_J\otimes_{R_J}D_B\).  This proves \eqref{eq:multicurve-family-identification}.

The ordinary tensor expression must still be compared with the completed support.  The central classifying arc lies in no curve divisor: otherwise the generic fibre of the pulled-back partial resolution would contain an exceptional curve, whose closure in \(X_0\) would be an exceptional surface.  Choose a divisor equation with pullback \(g(w)\ne0\). Since the classifying arc is local, \(g(w)\) has positive order. Then \(\C[[w]]\) is finite over \(\C[[g(w)]]\), and hence over \(R_J\).  Let \(U_J\) be the completed universal affine target ring, with maximal ideal \(\mathfrak n_J\), and set
\[
 Z_0=U_J\otimes_{R_J,f}\C[[w]],\qquad f=F\bmod\mathfrak m.
\]
The finiteness of \(\C[[w]]\) over \(R_J\) makes \(Z_0\) finite over \(U_J\).  Its closed fibre over \(U_J\) is
\[
 Z_0/\mathfrak n_JZ_0
 \cong\C[[w]]/\mathfrak m_{R_J}\C[[w]],
\]
which is local Artinian: the ideal on the right contains \(g(w)\ne0\) and is contained in \((w)\).  Since a finite algebra is integral, every maximal ideal of \(Z_0\) lies over \(\mathfrak n_J\).  The local closed fibre therefore gives a unique maximal ideal \(\mathfrak n_0=\mathfrak n_JZ_0+wZ_0\). For some \(c\ge1\) one has \(\mathfrak n_0^c\subseteq\mathfrak n_JZ_0\subseteq\mathfrak n_0\), so the two adic topologies coincide. Every finite module over the complete Noetherian local ring \(U_J\) is complete \cite[Tag~\href{https://stacks.math.columbia.edu/tag/00MA}{00MA}]{Stacks}. Consequently \(Z_0\), the central pullback of the completed universal endomorphism algebra, and its idempotent quotient are already \(\mathfrak n_0\)-adically complete. Ordinary scalar extension commutes with the idempotent quotient, since the image of the extended ideal is exactly the ideal generated by the extended projector.  Thus the marked identification with the completed central generator established above identifies the ordinary quotient along \(f\) with \(A_0\), as in Lemma~\ref{lem:arc-quotient-supported}.

The ordinary pullback on the right-hand side of \eqref{eq:multicurve-family-identification} is finite over \(D_B\) by Theorem~\ref{thm:multicurve-weighted-presentation}, and its closed \(B\)-fibre is finite-dimensional.  Corollary~\ref{cor:perfect-finite-slice} applied to the pulled-back square presentation makes this ordinary pullback finite over \(B\).  On the resulting finite \(B\)-module \(w\) is nilpotent modulo \(\mathfrak m\), so the \((\mathfrak m,w)\)-adic and \(\mathfrak m\)-adic topologies agree. The ordinary pullback is already complete, which identifies the ordinary tensor expression with the quotient carrying the full relative support.
\end{proof}

We first express the ranks in terms of the orders \(m_\beta\) of the divisor equations on the closed classifying arc. We then use the perturbation count of Nabijou--Wemyss to identify these orders with the GV invariants.

\Needspace{20\baselineskip}
\begin{proposition}
\label{thm:multicurve-algebra-flatness}
In the preceding setting, \(\mathscr A_B\) is finite free over \(B\) of rank \(\dim_\C A_0\), and
\[
 C_J\otimes_{R_J,F}^LB[[w]]\simeq\mathscr A_B[0].
\]
Write \(f=F\bmod\mathfrak m\) and \(m_\beta=\operatorname{ord}_wP_\beta(f(w))\).  These integers are positive and finite, and
\begin{align}
 \operatorname{rank}_B\mathscr A_B
 &=\sum_{\beta\in\mathcal B_J}m_\beta s_\beta^2,
 \label{eq:multicurve-algebra-rank}\\
 \operatorname{rank}_B(e_i\mathscr A_Be_j)
 &=\sum_{\beta\in\mathcal B_J}m_\beta\beta_i\beta_j.
 \label{eq:multicurve-entry-rank}
\end{align}
Moreover \(\Spec B[[w]]/(P_\beta(F))\) is finite flat over \(B\) of degree \(m_\beta\).
\end{proposition}

\begin{proof}
Pull back the algebra presentation of Theorem~\ref{thm:multicurve-weighted-presentation} along \(F\).  Its reduction modulo \(\mathfrak m\) has cokernel \(A_0\), which is finite-dimensional, so the reduced matrix is invertible over \(\C((w))\). Lemma~\ref{lem:finite-slice-flatness} gives finite freeness over \(B\) and injectivity of the pulled-back presentation.  The resulting exact sequence proves the derived statement.

Smith normal form identifies the closed-fibre length with the order of the determinant. The Fitting formula \eqref{eq:multicurve-weighted-fitting} supplies its factors and exponents, giving \eqref{eq:multicurve-algebra-rank}.  Each \(P_\beta(f)\) is nonzero because their product, with positive exponents, is nonzero.  Since \(f\) is local, \(P_\beta(f)\) vanishes at the origin, and hence \(m_\beta>0\).  Apply Lemma~\ref{lem:finite-slice-flatness} to the presentations of the Peirce spaces and to the one-by-one matrices \(P_\beta(F)\).  This gives the remaining rank formula and the stated finite flat schemes cut out by the divisor equations.
\end{proof}

For the numerical comparison, we need the exceptional locus to lie in a scheme proper over the parameter base. We also need the stable-sheaf point of a split curve to be reduced. The next two lemmas establish these facts.

\begin{lemma}
\label{lem:multicurve-critical-support}
For the classifying family above, the relative critical scheme of the affine ADE surface family, defined by the surface equation and its partial derivatives in the surface coordinates, is finite over \(B\).  Its projective inverse image contains every geometric exceptional curve and is projective over \(B\). For a holomorphic classifying bidisc, one may choose a fixed closed surface-parameter disc and a small external disc so that this critical locus has no points on the surface-parameter boundary, is proper with finite fibres over the external disc, and has proper inverse image under the contraction.
\end{lemma}

\begin{proof}
\emph{Finiteness over the polynomial base.} Work over the polynomial parameter ring \(H_J^{\mathrm{alg}}\), with parameters \(a\), and write the universal affine surface equation as \(f_{\mathrm{ADE}}(x_1,x_2,x_3;a)=0\).  The algebra
\[
 \mathcal J_{\mathrm{crit}}
 =\frac{H_J^{\mathrm{alg}}[x_1,x_2,x_3]}
 {(f_{\mathrm{ADE}},\partial_{x_1}f_{\mathrm{ADE}},
                     \partial_{x_2}f_{\mathrm{ADE}},
                     \partial_{x_3}f_{\mathrm{ADE}})}
\]
is finite over \(H_J^{\mathrm{alg}}\).  Indeed the central Jacobian algebra is finite-dimensional.  For the standard ADE equations its relations are, up to nonzero constants,
\[
\begin{array}{c|l}
 A_n&x_1=x_2=x_3^n=0\\
 D_n&x_1=x_2x_3=x_2^2+(n-1)x_3^{n-2}=0\\
 E_6&x_1=x_2^2=x_3^3=0\\
 E_7&x_1=3x_2^2+x_3^3=x_2x_3^2=0\\
 E_8&x_1=x_2^2=x_3^4=0.
\end{array}
\]
In type \(D_n\) these imply \(x_3^{n-1}=x_2^3=0\); in type \(E_7\) they imply \(x_3^5=x_2^3=0\). Choose homogeneous lifts of a basis of the central critical algebra. To express a homogeneous element in terms of the lifts, first reduce it modulo the positive parameter ideal. The remainder is a sum of positive-weight parameters times elements of lower weight. Induction on weight shows that the lifts generate \(\mathcal J_{\mathrm{crit}}\) over the parameter ring.

\emph{Restriction to the original arc.} Let \(\Delta_{\mathrm{tot}}\) cut out the full surface discriminant, including both curve and residual-singularity components \cite[Remark~2.8 and the following equality]{NabijouWemyss}. One may use the invariant square of the product of all positive-root forms.  Since the support of the finite critical module is this hypersurface, some power \(\Delta_{\mathrm{tot}}^q\) annihilates it. We claim that \(\Delta_{\mathrm{tot}}\) does not vanish identically on the central classifying arc.  Otherwise the critical algebra over \(\C[[w]]\) would have a nonzero generic fibre, hence a field quotient finite over \(\C((w))\).  At the corresponding critical point the central threefold equation \(G_0\) and its three surface partial derivatives vanish.  Extend \(d/dw\) to this finite separable field extension and differentiate \(G_0=0\); it follows that \(\partial_wG_0=0\) as well.  This gives a branch of the singular locus dominating the \(w\)-disc, contradicting the isolated central threefold singularity.  Equivalently, its finite-length four-variable Jacobian quotient cannot map to a field in which \(w\) is nonzero.

\emph{Finiteness over the complete base.} Put \(g=\Delta_{\mathrm{tot}}(F)\in B[[w]]\). Since \(g^q\) has nonzero reduction modulo \(\mathfrak m\), Lemma~\ref{lem:finite-slice-flatness} applied to this scalar matrix makes \(B[[w]]/(g^q)\) finite over \(B\).  The pulled-back critical algebra is finite over this quotient ring and therefore finite over \(B\).

Each partial resolution is an isomorphism over the smooth locus of its surface fibre. Every geometric exceptional curve therefore maps into the relative critical subscheme. The inverse image of this finite scheme under the projective contraction is projective over \(B\).  This finite support is already complete at the central target point.

\emph{Proper support over an analytic disc.} For a holomorphic bidisc, choose a closed \(w\)-disc on whose boundary the surface-discriminant equation is nonzero after pullback by the central classifying arc. Shrink the external disc to preserve that nonvanishing.  The pulled-back discriminant has finitely many roots in each \(w\)-disc, with finitely many critical points over each root.

To prove properness, we also need a bound on the surface coordinates. The affine coordinates on the critical scheme satisfy monic equations over the parameter ring. On a compact set of parameter values these give uniform bounds on the coordinates.  The inverse image of a compact external subdisc in the critical scheme is therefore closed in a compact set.  The absence of \(w\)-boundary points proves properness.  Its points specialize only to the central critical point, so the neighbourhood can be chosen inside the original target germ.  Projectivity of the contraction then gives properness of the inverse image.
\end{proof}

\begin{lemma}
\label{lem:multicurve-odp-sheaf-contribution}
Let \(F\) be a field of characteristic zero.  On the split local threefold \(Z_{\mathrm{loc}}=\operatorname{Tot}_{\mathbb{P}^{1}_{F}} (\mathcal O(-1)^{\oplus2})\), or its formal completion along the zero section \(C_{\mathrm{loc}}\), the only pure one-dimensional semistable sheaf with \(\chi=1\) and support on \(C_{\mathrm{loc}}\) is \(\mathcal O_{C_{\mathrm{loc}}}\).  Its first and second self-Ext groups vanish.  It therefore gives one reduced isolated stable-sheaf point, rational over \(F\), with local virtual contribution \(+1\) under the usual stable-sheaf obstruction theory. For finitely many disjoint such local curves, there are no additional stable objects supported on higher multiples or on more than one curve.
\end{lemma}

\begin{proof}
A pure one-dimensional sheaf with Hilbert polynomial \(dn+1\) cannot be strictly semistable. Indeed, equality of reduced Hilbert polynomials for a proper subsheaf of degree \(d'\le d\) would give \(d\chi'=d'\), hence \(d'=d\) and \(\chi'=1\). The quotient would then have dimension zero and Euler characteristic zero, forcing it to be zero, a contradiction.

\emph{Support on a finite thickening.} To classify these sheaves, we first show that they are annihilated by a power of the zero-section ideal. Let \(p:Z_{\mathrm{loc}}\to\mathbb{P}^{1}_{F}\) be the projection, and let \(\mathcal F\) be such a sheaf of generic multiplicity \(m\). In the algebraic total space, Noetherian coherence and the stated set-theoretic support imply that a power of the zero-section ideal annihilates \(\mathcal F\). We justify this implication separately in the formal case. The split contraction
\[
 \operatorname{Tot}_{\mathbb P^1_F}
       (\mathcal O(-1)^{\oplus2})
 \longrightarrow\operatorname{Spec}F[a,b,c,d]/(ad-bc)
\]
is projective and is an isomorphism away from the origin. In homogeneous coordinates \([s_0:s_1]\), it sends the two fibre coordinates \(x,y\) to \((s_0x,s_1x,s_0y,s_1y)\), and a nonzero column recovers \([s_0:s_1]\). On each standard chart the pulled-back origin ideal is the zero-section ideal.

Complete the target at the origin, put \(\widehat Q=F[[a,b,c,d]]/(ad-bc)\), and let \(X_{\widehat Q}\) be the corresponding projective source. The completion of \(X_{\widehat Q}\) along its closed fibre is the formal neighbourhood in the statement.

Coherent formal existence algebraizes \(\mathcal F\) to a coherent sheaf \(\mathcal F^{\mathrm{alg}}\) on \(X_{\widehat Q}\) \cite[Remark~30.27.2, Tag~\texttt{088F}]{Stacks}. Formal functions \cite[Theorem~30.20.5, Tag~\texttt{02OC}]{Stacks} and completeness of finite \(\widehat Q\)-modules identify
\[
 H^0(X_{\widehat Q},\mathcal F^{\mathrm{alg}})
 \cong H^0(\widehat Z_{\mathrm{loc}},\mathcal F).
\]
Here \(\widehat Z_{\mathrm{loc}}\) denotes the zero-section formal completion. The space on the right is finite-dimensional over \(F\), as required for the stated Euler characteristic to be defined. The coherent pushforward is therefore a finite-length \(\widehat Q\)-module, so it vanishes away from the closed point.

The contraction is an isomorphism away from the closed target point, so \(\mathcal F^{\mathrm{alg}}\) is supported on the closed fibre of the contraction. Noetherian coherence now implies annihilation by a power of the closed-fibre ideal, which is the zero-section ideal. The completed sheaf is annihilated by the same power of the zero-section ideal.

\emph{Identification of the sheaf.} In both cases, \(\mathcal F\) is a sheaf on a finite zero-section thickening, and \(V=p_*\mathcal F\) is coherent on \(\mathbb P^1_F\). Its torsion subsheaf is invariant under the normal-coordinate action, and would give a zero-dimensional subsheaf of \(\mathcal F\). Purity therefore makes \(V\) a vector bundle, of rank \(m\), with \(\chi(V)=1\).

We use a line subbundle of maximal degree to test semistability. Choose the largest integer \(a\) with \(\operatorname{Hom}(\mathcal O(a),V)\ne0\). It exists because a rational frame embeds \(V\) into \(\mathcal O(D)^{\oplus m}\) for an effective divisor \(D\), and \(a\ge0\) because \(h^0(V)-h^1(V)=1\). Since \(\operatorname{Pic}(\mathbb{P}^{1}_{F})=\mathbb Z[\mathcal O(1)]\), a nonzero map \(\mathcal O(a)\to V\) is saturated: otherwise its saturation is a line bundle of larger degree, contradicting maximality.

The normal-coordinate algebra acts by \(\mathcal O(1)^{\oplus2}\otimes V\to V\). The normal-coordinate action on \(\mathcal O(a)\) is zero, since \(\operatorname{Hom}(\mathcal O(a+1),V)=0\). Thus \(\mathcal O(a)\) gives a subsheaf on the reduced zero section.

For a polarization of degree \(h>0\) there, its reduced constant Hilbert term is \((a+1)/h\), whereas that of \(\mathcal F\) is \(1/(mh)\). If \(m>1\), this contradicts semistability. Hence \(m=1\), \(V=\mathcal O\), and the normal action is zero because \(H^0(\mathcal O(-1))=0\). This proves \(\mathcal F=\mathcal O_{C_{\mathrm{loc}}}\), which is stable by the same Hilbert-polynomial criterion.

\emph{The local contribution.} Having identified the sheaf, we compute its deformations and obstructions. The zero-section Koszul resolution computes the self-Ext groups of \(\mathcal O_{C_{\mathrm{loc}}}\) as
\[
 \operatorname{Ext}^k(\mathcal O_{C_{\mathrm{loc}}},
                      \mathcal O_{C_{\mathrm{loc}}})
 \cong\bigoplus_{q=0}^2
 H^{k-q}\left(\mathbb{P}^{1}_{F},\bigwedge^q(\mathcal O(-1)^{\oplus2})\right).
\]
After applying sheaf Hom into the zero-section sheaf, the Koszul differentials vanish because the normal coordinates act by zero. The only nonzero groups are \(F\) in degrees zero and three, from \(H^0(\mathbb{P}^{1}_{F},\mathcal O)=F\) and \(H^1(\mathbb{P}^{1}_{F},\mathcal O(-2))=F\), respectively. Thus the tangent and obstruction spaces are both zero, and the local virtual point has degree \(+1\). This calculation also holds in the formal neighbourhood, since the same finite Koszul resolution computes the proper-support Ext groups. The same argument applies after every field extension of \(F\), so the point is geometrically reduced and isolated.

For a disjoint union of curves, it remains to exclude stable sheaves supported on more than one component. The comaximal support ideals split a sheaf into the direct sum of its pieces. If two pieces are nonzero, one has Euler characteristic at least one and degree strictly smaller than the total degree, and hence destabilizes. Only one curve can support a stable object, and the calculation above applies to it.
\end{proof}

\begin{lemma}
\label{lem:central-cycle-hilbert-length}
Choose a normal partial ADE surface \(S_0\) in the central smooth threefold, and let \(p\) be the exceptional image.  For a nonzero effective cycle \(\beta=\sum_{i\in J}\beta_i[C_i]\), put
\[
 \mathcal I_{Z_\beta,S_0}
 =\bigcap_{\beta_i>0}\mathcal I_{C_i,S_0}^{(\beta_i)}.
\]
Every pure one-dimensional stable sheaf of Euler characteristic one and cycle \(\beta\), supported on the exceptional locus, is \(\mathcal O_{Z_\beta}\).  Such a sheaf exists precisely when \(Z_\beta\subset f_0^{-1}(p)\) scheme-theoretically and \(\chi(\mathcal O_{Z_\beta})=1\).  It then has \(H^0=\C\) and \(H^1=0\), is stable for every ample polarization, and determines an isolated point of the threefold Hilbert space.  Denote its local length by \(k_\beta\).

Consider a completed classifying deformation over \(R_s=\C[[s]]\), with flat source and target, and with the retained determinant line bundles fixed.  Suppose its geometric generic complete curves are finitely many disjoint split \((-1,-1)\)-curves.  If the central point above exists, its relative local Hilbert algebra is finite free of rank \(k_\beta\) over \(R_s\).  Its geometric generic points are reduced and correspond exactly to the curves of retained degree vector \(\beta\).  Conversely, a generic curve of that vector specializes to \(Z_\beta\), so its existence implies the stated central conditions.
\end{lemma}

\begin{proof}
A stable sheaf \(F\) with proper support has only scalar endomorphisms. Indeed a nonzero endomorphism with nonzero kernel would have an image of slope greater than that of \(F\), whereas stability forces the slope of that proper subsheaf to be smaller.  An injective endomorphism is surjective, since its cokernel has zero Hilbert polynomial. Thus \(\End(F)\) is a finite-dimensional division algebra over \(\C\), hence is \(\C\). Every function in the maximal ideal of \(p\) acts nilpotently on \(F\); being scalar, that action is zero.  Consequently \(F\) lies on the scheme-theoretic contraction fibre and on \(S_0\).

Rationality of the partial ADE surface gives \(H^1(S_0,\mathcal O)=0\). For any closed subscheme \(Z\) supported on its exceptional fibre, formal functions and the one-dimensional fibres give \(H^2(S_0,\mathcal I_Z)=0\) \cite[Tag~\href{https://stacks.math.columbia.edu/tag/02OC}{02OC}]{Stacks}. The cohomology sequence therefore gives \(H^1(\mathcal O_Z)=0\), and \(\chi(\mathcal O_Z)\ge1\) when \(Z\ne\varnothing\). Since \(\chi(F)=1\), \(F\) has a nonzero section.  Its image is a pure cyclic subsheaf \(\mathcal O_Z\subset F\).  The preceding inequality, together with its degree being at most the degree of \(F\), makes its reduced Hilbert polynomial at least that of \(F\).  Stability forces \(\mathcal O_Z=F\), so \(H^0(F)=\C\) and \(H^1(F)=0\).

On the normal surface \(S_0\), purity of \(\mathcal O_Z\) excludes embedded height-two primes from \(\mathcal O_{S_0}/\mathcal I_Z\).  At a height-one prime the local ring is a DVR, and the primary ideal is determined by its length.  Intersecting these primary ideals gives exactly the displayed symbolic-power ideal. Thus \(Z=Z_\beta\). Conversely, a proper quotient of \(\mathcal O_{Z_\beta}\) is the structure sheaf of a nonempty proper subscheme and has Euler characteristic at least one.  Every proper one-dimensional subsheaf therefore has Euler characteristic at most zero, which proves stability for every ample polarization.

For the infinitesimal statement, an Artin flat deformation \(F_A\) has \(H^0(F_A)\cong A\) and \(H^1(F_A)=0\), by induction over small extensions.  A generator of \(H^0(F_A)\) gives a surjective evaluation map \(\mathcal O\to F_A\); changing the generator by a unit preserves its kernel.  The endomorphism algebra of \(F_A\) is \(A\), so the affine target acts through a unique moving target section.  Sending \(F_A\) to the kernel of its evaluation map identifies the pointed sheaf and embedded Hilbert functors.  This identification also applies to prescribed ambient deformations and to completion at the central target point.

We now construct the open-and-closed locus of the relative Hilbert space on which the cycle is \(\beta\).  Write \(L_i\) for the retained determinant line bundles, so that \(\deg_{C_j}L_i=\delta_{ij}\).  For either the central contraction or the prescribed \(R_s\)-family, start with \(M=\operatorname{Hilb}_{X/Y}\) and restrict to its one-dimensional Hilbert-polynomial strata.  On an \'etale scheme chart \(U\to M\), let \(\mathcal Z_U\to U\) be the universal proper flat family.  Define
\[
 d_i(Z)=\chi(\mathcal O_Z\otimes L_i)-\chi(\mathcal O_Z).
\]
The functions \(\chi(\mathcal O_Z)\) and \(d_i(Z)\) are locally constant and commute with base change: apply \cite[Lemma~36.32.2, Tag~\href{https://stacks.math.columbia.edu/tag/0B9T}{0B9T}]{Stacks} to the proper finitely presented family \(\mathcal Z_U\to U\), with its flat finitely presented structure sheaf and its twists by \(L_i\). The conditions
\[
 \chi(\mathcal O_Z)=1,\qquad d_i(Z)=\beta_i\quad(i\in J)
\]
therefore define an open-and-closed subspace \(M_\beta\subset M\). The proper morphism used here is the universal Hilbert family \(\mathcal Z_U\to U\); the ambient contraction need not be proper over the parameter base. For a one-dimensional subscheme with cycle \(\sum_j a_j[C_j]\), additivity modulo zero-dimensional sheaves gives
\[
 d_i(Z)=\sum_j a_j\deg_{C_j}L_i=a_i.
\]
Indeed a zero-dimensional sheaf has the same Euler characteristic after an invertible twist. Modulo zero-dimensional sheaves, the class of \(\mathcal O_Z\) is the sum of its generic lengths times the integral component classes. Thus these numerical conditions specify the full retained cycle even when the Hilbert subscheme is not initially assumed pure.

Assume first that the central point \(Z_\beta\) exists.  Every point of the central \(M_\beta\) lies over \(p\): its Hilbert subscheme has positive one-dimensional degree, whereas the central contraction has zero-dimensional fibres away from \(p\).  Every closed subscheme of \(f_0^{-1}(p)\) with the prescribed invariants equals \(Z_\beta\). To see this without assuming purity, let \(T\) be the maximal zero-dimensional subsheaf of its structure sheaf.  This is an ideal, and its quotient is pure and has cycle \(\beta\).  The symbolic-power argument above identifies that quotient with \(\mathcal O_{Z_\beta}\).  Additivity and \(\chi(\mathcal O_{Z_\beta})=1\) give \(1=\operatorname{length}T+1\), so \(T=0\).

We show that the local Hilbert ring at \([Z_\beta]\) is Artinian.  By \cite[Proposition~99.9.4, Tag~0D01]{Stacks}, the relative Hilbert space is locally of finite presentation over the Noetherian target.  Choose an affine \'etale chart near \([Z_\beta]\).  Its central restriction is Noetherian and of finite type over the central target.  The preceding set-theoretic support assertion says that the target maximal ideal \(\mathfrak m_p\) is contained in the nilradical of this chart. Consequently some power of \(\mathfrak m_p\) annihilates its structure sheaf.  The chart is therefore of finite type over an Artinian local \(\C\)-algebra.  Its fibre over \(p\) has only the Hilbert point \([Z_\beta]\); on an \'etale chart the points above it form a discrete zero-dimensional fibre.  The finite-type \(\C\)-scheme chart is Jacobson, so it has dimension zero near a chosen lift of this point. Its local ring is Noetherian of dimension zero, hence Artinian. The evaluation-map construction above identifies the pointed Artin functor of this local ring with the embedded threefold Hilbert functor.  It follows that the latter has the same Artinian local ring and intrinsic length \(k_\beta\).  In particular the threefold Hilbert point is isolated.

For the \(R_s\)-family use the analogous \(M_\beta\), choose an \'etale chart and a lift of the central point, and let \(\mathscr H_\beta\) be its completed local ring.  The preceding argument identifies \(\mathscr H_\beta/(s)\) with the central threefold Hilbert ring.  The local Artin functor is unchanged by taking this open-and-closed locus.  The tangent space and a complete obstruction space for this quotient deformation problem are
\[
 \begin{aligned}
 T_\beta&=\operatorname{Hom}(I_\beta,\mathcal O_{Z_\beta})
          \cong\operatorname{Ext}^1(\mathcal O_{Z_\beta},\mathcal O_{Z_\beta}),\\
 V_\beta&=\operatorname{Ext}^1(I_\beta,\mathcal O_{Z_\beta})
          \cong\operatorname{Ext}^2(\mathcal O_{Z_\beta},\mathcal O_{Z_\beta}).
 \end{aligned}
\]
These comparisons follow from the ideal sequence using \(H^0=\C\), \(H^1=H^2=0\), and evaluation at \(1\). For a small extension of \(R_s\)-algebras, ambient flatness gives the canonical quotient extension of Tag~0CZU; its splitting is equivalent to a flat quotient lift by Tag~08MQ.  This construction commutes with linear pushouts of the small kernel and keeps the prescribed ambient family and \(s\) fixed.  The lifted quotient remains proper: the central target ideal acts trivially modulo the Artin maximal ideal, so a power of it annihilates the lift.  The support consequently lies in the projective inverse image of an Artinian target thickening.  Thus \(V_\beta\) is a complete obstruction space for the stated local Hilbert functor.

Proper-support duality on the central smooth crepant threefold gives \(V_\beta\cong T_\beta^\vee\).  Here \(\mathcal O_{Z_\beta}\) is perfect and has proper exceptional support, and the canonical bundle is trivial near that support.  The self-Ext comparison in the proof of Proposition~\ref{prop:conserved-multiplicities-gv} transfers these groups to the chosen model; this is the only part of that proof used here. The closed Hilbert algebra has finite length \(k_\beta\), and the complete obstruction space \(V_\beta\) has the same dimension as the tangent space \(T_\beta\). These are the hypotheses of Lemma~\ref{lem:complete-obstruction-flatness}, which gives
\[
 \mathscr H_\beta\text{ finite free over }R_s,\qquad
 \operatorname{rank}_{R_s}\mathscr H_\beta=k_\beta.
\]
Finally, the local ring of the chosen chart injects into its completion by Krull intersection.  Its image is an \(R_s\)-submodule of the finite module \(\mathscr H_\beta\), hence is finite and closed in the \(s\)-adic topology.  It is also dense in that topology, since the closed quotient is Artinian and the maximal-ideal and \(s\)-adic topologies of the completion are cofinal.  The local ring therefore equals \(\mathscr H_\beta\).  The universal Hilbert family pulls back to \(\Spec\mathscr H_\beta\). The relative tangent spaces at its generic-fibre parameter points are consequently those of the relative Hilbert functor.

It remains to identify the generic points.  Lemma~\ref{lem:multicurve-critical-support} places the exceptional image in a target scheme finite over \(\C[[s]]\). Its coordinates are therefore integral, so after a finite extension of \(\C((s))\) the image extends to a section over the valuation ring. The projective inverse image of that section contains the unique saturated schematic closure of the curve.  The closure is DVR-flat, with central cycle \(\beta\) and Euler characteristic one.  The central fibre of the closure lies scheme-theoretically in \(f_0^{-1}(p)\subset S_0\).  Its pure quotient is \(\mathcal O_{Z_\beta}\); the rational-surface cohomology calculation above gives its Euler characteristic at least one.  Additivity then forces that value to be one and the zero-dimensional kernel to vanish. Thus the central fibre is exactly \(Z_\beta\).  This proves both the existence of the central point and the specialization assertion.

Conversely, a generic Hilbert object has positive one-dimensional degree and lies in one contraction fibre.  Every such generic fibre is a reduced integral \(\mathbb P^1\); its only one-dimensional closed subscheme is itself.  The locally constant Euler differences for the retained determinant line bundles preserve the vector \(\beta\). Uniqueness of saturated closure and of the pointed \'etale lift prevents duplicate Hilbert points.  For each resulting split \((-1,-1)\)-curve \(E\), Lemma~\ref{lem:multicurve-odp-sheaf-contribution} gives \(\operatorname{Ext}^1(\mathcal O_E,\mathcal O_E)=0\). Together with \(H^0(\mathcal O_E)\) equal to the ground field and \(H^1(\mathcal O_E)=0\), the ideal sequence makes the embedded Hilbert tangent space zero.  The finite geometric generic local rings are therefore fields.  Their number is the rank \(k_\beta\).
\end{proof}

\begin{corollary}
\label{cor:multicurve-cartan-gram}
The Cartan matrix with respect to the idempotents \((e_i)_{i\in J}\) is
\[
 \bigl(\dim_\C e_iA_0e_j\bigr)_{i,j\in J}
 =\sum_{\beta\in\mathcal B_J}m_\beta\,\beta\beta^{\mathsf T}.
\]
It is symmetric positive definite.  For the genus-zero invariants of the central flopping contraction, \(m_\beta=n_\beta\).
\end{corollary}

\begin{proof}
Restricting \eqref{eq:multicurve-entry-rank} to the closed fibre gives the displayed matrix formula.  Every retained simple root projects to a standard basis vector and every \(m_\beta\) is positive, so the sum is positive definite.

We compare the formal classifying arc \(f\) with an analytic marked arc \(g_{\mathrm{an}}\) for the same original contraction.  Fix the same standard partial resolution, its chosen side, and the retained curve basis.  Let \(g\) be the completion of \(g_{\mathrm{an}}\). The convention for \(n_\beta\) is the one fixed in the Introduction: the curve count of \cite[Section~3.1, equation~(3.A)]{NabijouWemyss}. In that reference, the contracted subset is \(I\), and the retained subset is our \(J\). Lemma~\ref{lem:multicurve-root-products} identifies their class-\(\beta\) component with our reduced divisor \(P_\beta=0\).

Choose a small disk for \(g_{\mathrm{an}}\) on whose boundary the full surface discriminant is nonzero.  A sufficiently small transverse perturbation avoids the multiple-root strata and satisfies \cite[Section~2.6, condition~(2.G)]{NabijouWemyss}. The argument principle and their curve--intersection correspondence \cite[Lemma~2.9 and Corollary~2.10]{NabijouWemyss} give
\[
 n_\beta=\operatorname{ord}_w P_\beta(g(w))=:m_\beta(g).
\]
These are intersections with the reduced divisor on \(\mathfrak h/W_I\).  Distinct proportional restricted roots retain their separate components.

To compare this order with the order along \(f\), form a transverse formal translation of each of the two arcs, using Lemma~\ref{lem:formal-transverse-translation} for the full surface discriminant and its root-intersection strata.  For each translation, take the ordinary pullback of the fixed algebraic model over \(\C[[s,w]]\), followed by full central target completion. Universal parameter flatness and the flatness of this completion make both source and target flat over \(\C[[s]]\); their marked central contraction is the prescribed \(f_0\).

For an arc \(a\in\{f,g\}\), the finite algebra defined by the translated equation \(P_\beta=0\) has rank \(m_\beta(a)\) over \(\C[[s]]\), by Lemma~\ref{lem:finite-slice-flatness}. At a point where exactly one positive root \(\alpha\) vanishes and \(\pi_I(\alpha)\ne0\), the chosen simultaneous partial resolution has the transverse \(A_1\) model: the full Weyl quotient has normal coordinate \(t=\zeta^2\), while \(\mathfrak h\to\mathfrak h/W_I\) is \'etale at the chosen point of that non-Levi root hyperplane.  The transverse pullback is the split ODP \(xy-z^2=\zeta^2\) after a finite splitting extension, with its selected small resolution and one \((-1,-1)\)-curve.  If exactly one positive root \(\alpha\) vanishes and \(\pi_I(\alpha)=0\), no curve is retained.  The same local equations underlie \cite[Lemma~2.9]{NabijouWemyss} and apply over the formal splitting field. Distinct roots of the translated discriminant polynomials have different \(w\)-coordinates, and the union of the curve divisors is the full locus of complete curves.  Thus there are exactly \(m_\beta(a)\) complete curves with the indicated label and no further complete curves.  Lemma~\ref{lem:multicurve-critical-support} ensures that their full target coordinates, and hence their projective inverse images, survive central target completion.

The line bundles used on these two families are the pullbacks of the specified labelled bundles in Proposition~\ref{prop:multicurve-cartan-comparison}.  Choose a local lift to the root hyperplane \(H_\alpha\subset\mathfrak h\). Specialize the projective flat family of divisors over \(H_\alpha\) provided by Lemma~\ref{lem:chosen-root-divisor-family}. Constancy of the determinant-line-bundle degrees in this flat family gives, for each retained \(i\),
\[
 \deg_E\det(P_i)
   =\sum_{j\in J}\beta_j\deg_{C_j}\det(P_i)=\beta_i.
\]
Here the central degrees are \(\delta_{ij}\); the root specialization uses Lemma~\ref{lem:chosen-root-divisor-family} and the same morphism from the full resolution to the partial resolution, and the same labelled bundles, as in Remark~\ref{rem:root-family-full-resolution}. Consequently every retained coordinate, and not only their sum, agrees with the chosen class basis.

Apply Lemma~\ref{lem:central-cycle-hilbert-length} to the two completed families just constructed, whose sources and targets are flat and whose generic curves are disjoint split \((-1,-1)\)-curves.  Each contains a class-\(\beta\) generic curve, since its root degree is positive, and hence has central point \(Z_\beta\). The marked identification of their central completed contractions identifies this unique cyclic sheaf and its embedded threefold Hilbert functor.  Its local length \(k_\beta\) is therefore the same for both families.  The lemma gives the individual equalities
\[
 m_\beta(f)=k_\beta=m_\beta(g)=n_\beta.
\]
This comparison does not require the two classifying maps to be equal.

Returning to the original formal translation of \(f\), the algebra defined by the translated equation \(P_\beta=0\) has rank \(m_\beta(f)=n_\beta\) over \(\C[[s]]\) and the determinant degrees just computed retain its original marking.  We form its generic fibre by localizing the full relative family before completing at moving points. For a single retained curve, Proposition~\ref{prop:conserved-multiplicities-gv} identifies the same count with the original Katz Hilbert invariant.
\end{proof}

\subsubsection{Hua--Toda flatness over general bases}
\label{subsec:multicurve-hua-toda}

Kawamata proved algebra flatness over a disk \cite[Theorem~3.1(2)]{KawamataNC}.  To extend this to general complete local bases, we apply the rank calculation above to the algebra and its Peirce spaces, and Lemma~\ref{lem:finite-slice-flatness} to the cocenter presentation. A relative symmetrizing trace then identifies the centre with the dual cocenter.  These steps give the following freeness and base-change statements.

\Needspace{25\baselineskip}
\begin{theorem}[Hua--Toda flatness and base change]
\label{thm:multicurve-hua-toda}
Let \(f_B\), \(T_B\), and \(\mathscr A_B\) be as in Proposition~\ref{prop:multicurve-classifying-pullback}, with arbitrary complete Noetherian local complex base \((B,\mathfrak m,\C)\). Let \(n_\beta=m_\beta\) be the central curve multiplicities of Corollary~\ref{cor:multicurve-cartan-gram}, and put \(N_{\mathrm{GV}}=\sum_{\beta\in\mathcal B_J}n_\beta\).  Then
\[
 \mathscr A_B\cong B^{\sum_\beta n_\beta s_\beta^2},\qquad
 HH_0(\mathscr A_B/B)\cong B^{N_{\mathrm{GV}}},\qquad Z(\mathscr A_B)\cong B^{N_{\mathrm{GV}}}
\]
as \(B\)-modules.  For every commutative \(B\)-algebra \(B'\), the canonical maps
\begin{align*}
 HH_0(\mathscr A_B/B)\otimes_BB'
 &\xrightarrow{\sim}HH_0(\mathscr A_B\otimes_BB'/B'),\\
 Z(\mathscr A_B)\otimes_BB'
 &\xrightarrow{\sim}Z(\mathscr A_B\otimes_BB')
\end{align*}
are isomorphisms.  Every symmetrizing trace of \(A_0\) lifts to a perfect \(B\)-linear symmetrizing trace of \(\mathscr A_B\). Moreover,
\begin{align*}
 C_J\otimes_{R_J,F}^L B[[w]]&\simeq\mathscr A_B[0],\\
 \mathcal H_J\otimes_{R_J,F}^L B[[w]]
 &\simeq HH_0(\mathscr A_B/B)[0].
\end{align*}
The free bases and lifted trace can be chosen, and the base-change maps are canonical.  The auxiliary \(B\)-algebra \(B'\) need not be Noetherian, local, or complete.  These are ordinary tensors of the finite relative algebra; no replacement by \(B'[[w]]\) is made.
\end{theorem}

\begin{proof}
Proposition~\ref{thm:multicurve-algebra-flatness} gives the asserted rank and freeness of \(\mathscr A_B\), and shows that the derived pullback of \(C_J\) is concentrated in degree zero.  Put \(H_B=HH_0(\mathscr A_B/B)\).  The \(B\)-span of commutators is already a \(D_B\)-module, since \(d[a,b]=[da,b]\) for a central scalar \(d\in D_B\).  Hence
\[
 H_B\cong\mathcal H_J\otimes_{R_J,F}D_B
 \cong\operatorname{coker}\Theta_{\mathrm{tr},J}(F).
\]
On the closed classifying arc, the determinant of this presentation is a unit times \(\prod_\beta P_\beta(f(w))\).  It is nonzero and has order \(N_{\mathrm{GV}}\).  Lemma~\ref{lem:finite-slice-flatness} gives freeness of rank \(N_{\mathrm{GV}}\) and injectivity of the pulled-back presentation.  This also proves the second derived statement over \(B\), without assuming that \(B\) is reduced.

Finally, the central contraction algebra is symmetric: the stable endomorphism algebra of the entire maximal Cohen--Macaulay generator is covered by \cite[Theorem~2.6]{August}.  For any symmetrizing trace of \(A_0\), apply Lemma~\ref{lem:relative-symmetrizing-trace} to the finite free algebra and cocenter just obtained.  The lemma lifts the closed trace to a perfect relative symmetrizing trace and identifies the centre with the dual cocenter. The trace depends on the lift chosen; the two base-change maps are canonical.
\end{proof}

For the local consequences below, write \(d=\dim_\C A_0\).

\begin{corollary}
\label{cor:hua-toda-local-bases}
The commutator submodule of \(\mathscr A_B\) is a \(B\)-direct summand. For every ideal \(L\subseteq B\),
\[
 [\mathscr A_B,\mathscr A_B]\cap L\mathscr A_B
 =L[\mathscr A_B,\mathscr A_B].
\]
For a continuous local map \(B\to B'\) of complete Noetherian local complex bases, \(\mathscr A_B\otimes_BB'\) is the geometric relative quotient formed from the pulled-back fixed generator and full finite support.  Over a Noetherian algebraic or analytic base whose relative quotient is finite over the base, the algebra, centre, and cocenter are locally free near a fibre of the stated kind, with ranks \(d,N_{\mathrm{GV}},N_{\mathrm{GV}}\), and the centre and cocenter commute with arbitrary base change.  Several exceptional image points are treated by their completed pieces.
\end{corollary}

\begin{proof}
The cokernel of the commutator submodule is finite free, so the commutator sequence splits.  Intersecting the resulting direct-sum decomposition with \(L\mathscr A_B\) proves the saturation formula.

By Proposition~\ref{prop:multicurve-classifying-pullback}, geometric base change reduces to comparing the power-series pullbacks of \(C_J\). Let \(W(w)\) be the characteristic polynomial of multiplication by \(w\) on the finite free \(B\)-module \(\mathscr A_B\).  Cayley--Hamilton makes \(W\) annihilate this module.  Since \(w\) is nilpotent on its closed fibre, \(W\) is monic of degree \(d\) with reduction \(w^d\). Both \(B[w]/(W)\) and \(B[[w]]/(W)\) are free of rank \(d\), the latter by Lemma~\ref{lem:finite-slice-flatness}.  The inclusion of polynomials into power series induces a map between these two free \(B\)-modules.  It is an isomorphism modulo \(\mathfrak m\), hence an isomorphism over \(B\). The same argument after a continuous local map \(B\to B'\) gives
\[
 B'[[w]]/(W')\cong B'[w]/(W')
 \cong(B[w]/(W))\otimes_BB'.
\]
Tensoring the module killed by \(W\) over these finite polynomial quotients therefore identifies its ordinary \(B\)-base change with its power-series base change.  This retains every point of its finite support.

For the noncomplete local assertion, retain the standing projective, source-flat, target-flat and fixed-generator hypotheses, and write \(A\) for the relative quotient, assumed finite over the parameter ring \((B,\mathfrak m,\C)\). First complete at a central target point \(y\), with local ring \((R,\mathfrak n)\).  Proper coherence makes \(E=(f_*\mathcal End(T)^{\mathrm{op}})_y\), its projector ideal, and its quotient finite \(R\)-modules. Formal functions is used for this target-point completion: in the algebraic case use \cite[Tag~\href{https://stacks.math.columbia.edu/tag/02OC}{02OC}]{Stacks}, and in the analytic case use \cite[Theorems~19.57 and~19.59]{RiemenschneiderComplexSpaces}. It identifies \(\widehat E\) with the endomorphisms of the completed marked bundle.  Exact completion on finite target modules then identifies the completed quotient with \(\widehat E/\widehat Ee_0\widehat E\) \cite[Tag~\href{https://stacks.math.columbia.edu/tag/00MA}{00MA}]{Stacks}.

Here the analytic completed contraction is a projective scheme over \(\widehat R\).  Choose the given projective embedding and put \(D_k=R/\mathfrak n^k\).  At each level the coherent structure algebra, its defining ideal and the pushed-forward marked bundle are coherent modules on the reduced projective space \(\mathbb P^r_\C\), with finite \(D_k\)-actions.  Apply projective GAGA on that reduced space: its exact, fully faithful coherent-sheaf equivalence algebraizes both the modules and all their morphisms \cite[no.~6, Proposition~3, p.~9; no.~9, Proposition~10, p.~17; no.~12, Theorems~1--3, pp.~19--20]{SerreGAGA}.  Choose a \(\C\)-basis of \(D_k\); its action is finitely many endomorphisms whose multiplication and unit identities remain valid by faithfulness.  The same argument transfers the quotient-algebra action, ideal inclusion, marked inclusions and retractions, and their labelled projectors.  Exactness preserves the defining quotient.

Compatibility with tensor products and cokernels identifies reduction to \(D_{k-1}\) for \(k\ge2\), so the algebraized structure algebras and bundles form one compatible infinitesimal system; no flatness over \(R\) is used here.  Finally, algebraic and analytic local rings on \(\mathbb P^r_\C\) have the same completions. This remains true after the finite \(D_k\)-extension and the defining quotient, and for the completed bundle modules.  Faithfully flat local completion therefore shows that each algebraized marked module is locally free. Grothendieck algebraization and existence then give the projective \(\widehat R\)-scheme and marked bundle \cite[Tags~\href{https://stacks.math.columbia.edu/tag/089A}{089A} and~\href{https://stacks.math.columbia.edu/tag/088F}{088F}]{Stacks}. Their completed local rings and bundle modules along the closed target fibre are those of the original analytic family.

Flatness of local maps is preserved by completion \cite[Tag~\href{https://stacks.math.columbia.edu/tag/0C4G}{0C4G}]{Stacks}, so the target and source are flat over \(\widehat B\): descend flatness from these completed local rings, and use that every point of the projective source specializes to the closed target fibre. The bundle is locally free there by faithfully flat completion, hence everywhere by projectivity and the closed nonfree locus. Formal full faithfulness identifies the central algebraization with the prescribed completed contraction \cite[Tag~\href{https://stacks.math.columbia.edu/tag/0A4Z}{0A4Z}]{Stacks}; all marked summands are retained, so the complete-base theorem applies.  Algebraically the corresponding model is the projective base change to \(\operatorname{Spec}\widehat R\).

The target and parameter completions agree only after taking the finite quotient.  When the central support is \(y\), its target maximal ideal acts nilpotently on \(A/\mathfrak mA\), so for some \(c\ge1\)
\[
 \mathfrak n^{ck}A\subseteq\mathfrak m^kA
 \subseteq\mathfrak n^kA\qquad(k\ge1).
\]
Thus \(\widehat A_{\mathfrak n}=A\otimes_B\widehat B\). For several central support points, let \(C\) be the commutative image of the target acting on \(A\).  It is finite over \(B\), being a \(B\)-submodule of \(\End_B(A)\). The finite-algebra completion decomposition \cite[Tag~\href{https://stacks.math.columbia.edu/tag/07N9}{07N9}]{Stacks} identifies \(A\otimes_B\widehat B\) with the product of the quotients completed at all central support points.  No splitting or selection of a local branch over the original base is required.

Apply the complete-base theorem to these pieces.  The cocenter commutes with scalar extension as the cokernel of the commutator map. Faithfully flat descent of finite projectivity \cite[Tag~\href{https://stacks.math.columbia.edu/tag/058S}{058S}]{Stacks} therefore makes \(A\) and its cocenter finite free over \(B\). Lift a central symmetrizing trace through this free cocenter; its pairing is perfect since its determinant is a unit on the closed fibre. The perfect trace identifies the centre with the dual of the finite free cocenter, so the centre is finite free and commutes with arbitrary scalar extension.  These bases and the perfect pairing extend over a neighbourhood of the chosen base point, proving the local assertion.
\end{proof}

\subsubsection{The exponents in two examples}
\label{subsec:multicurve-examples}

For \(\Delta=A_3\), with chain \(1-2-3\), take \(J=\{1,3\}\). The restricted roots are \((1,0),(0,1),(1,1)\), and
\[
 P_{10}=\lambda_1(\lambda_1+\lambda_2),\qquad
 P_{01}=\lambda_3(\lambda_2+\lambda_3),\qquad
 P_{11}=\lambda_1+\lambda_2+\lambda_3.
\]
Thus the algebra and cocenter Fitting ideals are respectively
\[
 (P_{10}P_{01}P_{11}^4),\qquad (P_{10}P_{01}P_{11}).
\]
If the three intersection orders are \(a,b,c\), the algebra rank is \(a+b+4c\), its centre and cocenter have rank \(a+b+c\), and its Cartan matrix with respect to the two curve idempotents is
\[
 \begin{pmatrix}a+c&c\\c&b+c\end{pmatrix}.
\]
This example distinguishes the number of positive roots with a given restriction from the exponent contributed by its matrix block.

For the length-five marking of \(E_8\), deleting the retained vertex leaves \(A_3\sqcup A_4\).  Theorem~\ref{thm:all-corner-trace-excision} gives
\[
 \dim HH_0(e_p\Pi^\rho_0e_p)=120-6-10=104.
\]
This recovers \eqref{eq:e8-equality} from the cocenter formula of Etingof--Latour--Rains for the central extension of the full preprojective algebra.  The explicit rank certificate in Subsection~\ref{subsec:e8-certificate} gives the direct corner calculation used in the single-curve proof.

The Peirce dimensions and cocenter count computed here are the numerical quantities used in the orbit argument with fixed idempotents in Section~\ref{sec:marked-semisimple}.

\section{Semisimple generic fibres with the curve projectors retained}
\label{sec:marked-semisimple}

The universal presentations determine the dimensions of the Peirce spaces and the number of simple factors for each retained subset. We now construct a semisimple deformation with those dimensions and with the curve projectors fixed.

For an analytic disk deformation satisfying the hypotheses of Kawamata's Theorem~3.1, the relative contraction algebra is flat and its generic fibre is a product of matrix algebras; the multiplicities of the original summands are recorded by the integers \(m_{j,i}\) \cite[Theorem~3.1(2)--(3)]{KawamataNC}.  We keep the separable algebra of labelled projectors fixed in the formal family. After a finite extension, the lattice construction produces a polynomially split algebra.  We then use the orbit argument to prove that every geometrically semisimple generic fibre has the same matrix factors and the same ranks for the labelled projectors.

Throughout this section, \(X_0\to Y_0\) is the completion at an exceptional image point of a projective flopping contraction of a smooth complex threefold.  We use the marked generator, retained set \(J\), restricted roots \(\mathcal B_J\), and quotient \(C_J\) of Subsection~\ref{sec:multicurve-flatness}.  The original curves are labelled \((C_i)_{i\in J}\).  In the positive-dual convention the generator is
\[
 P_0=\mathcal O\oplus\bigoplus_{i\in J}P_{0,i},\qquad
 \deg_{C_j}\det P_{0,i}=\delta_{ij},\qquad
 A=\End(P_0)/\End(P_0)e_0\End(P_0).
\]
Dualizing identifies this quotient with the convention \(\End(T_0)^{\mathrm{op}}/(e_0)\) used earlier, preserving all projectors.  Put
\begin{equation}
\label{eq:marked-target-definition}
 E=\bigoplus_{i\in J}\C e_i,\qquad
 V_\beta=\bigoplus_{i\in J}\C^{\beta_i},\qquad
 S_E=\prod_{\beta\in\mathcal B_J}\End_\C(V_\beta)^{\times n_\beta}.
\end{equation}
On each \(\beta\)-factor, \(e_i\) acts as the projection onto \(\C^{\beta_i}\).  These labelled idempotents are generally noncentral.

A marking consists of the fixed embedding \(E\hookrightarrow A\), which prescribes the projectors and their left and right actions. Equivalently, the simple module of a \(\beta\)-factor has dimension vector \(\bigl(\dim_\C e_iV_\beta\bigr)_{i\in J}=\beta\) relative to the idempotents \((e_i)_{i\in J}\), and this dimension vector occurs \(n_\beta\) times. Write \(s_\beta=\sum_i\beta_i\), keeping distinct proportional restricted roots separate.

\subsection{Constructing a semisimple generic fibre with fixed projectors}
\label{subsec:marked-existence}

\begin{theorem}
\label{thm:marked-semisimple-existence}
There is a finite free unital \(\C[[s]]\)-algebra \(\mathscr A\), obtained from the universal marked geometric family, with a marked identification \(\mathscr A/(s)\cong A\).  Its labelled orthogonal idempotents \(e_i\) sum to its identity.  For a finite extension \(K'/\C((s))\), there is an isomorphism preserving these idempotents
\[
 \mathscr A\otimes_{\C[[s]]}K'\cong S_E\otimes_\C K'.
\]
In particular, \(e_i\) has rank \(\beta_i\) on each of the \(n_\beta\) factors of size \(s_\beta\). After writing \(K'=\C((u))\), \(s=u^e\), there is also a finite free marked \(\C[u]\)-algebra \(\mathscr P\) such that
\[
 \mathscr P/(u)\cong A,\qquad
 \mathscr P[u^{-1}]\cong S_E\otimes_\C\C[u,u^{-1}].
\]
Thus the deformation can be chosen polynomially split, with all projector ranks preserved.
\end{theorem}

\begin{lemma}
\label{lem:formal-transverse-translation}
Let \(R=\C[[t_1,\ldots,t_r]]\), and let \(f_0(w)=(f_1(w),\ldots,f_r(w))\) be any formal arc through the origin, with \(f_\nu(w)\in w\C[[w]]\). Let \(P_1,\ldots,P_q\in R\) be pairwise nonassociated prime elements such that \(d_b=\operatorname{ord}_wP_b(f_0(w))\) is positive and finite.  Let \(\Sigma\) be a reduced closed formal subscheme of \(\bigcup_bV(P_b)\), of codimension at least two in \(\Spec R\). There is a nonempty Zariski open set of translation vectors \(v=(v_1,\ldots,v_r)\in\C^r\) with the following property. Fix any \(v\) in this set and put \(F(s,w)=f_0(w)+sv\). Thus the \(\nu\)-th coordinate is \(f_\nu(w)+sv_\nu\), and \(F(0,w)=f_0(w)\).  Each \(P_b(F)\) is a unit times a distinguished polynomial in \(w\) of degree \(d_b\); over \(\C((s))\) these polynomials are separable, pairwise coprime, and have no root mapping into \(\Sigma\). After a finite ramified extension their roots lie in \(u\C[[u]]\).
\end{lemma}

\begin{proof}
Introduce independent parameters \(a_1,\ldots,a_r\), and put \(B_{\mathrm{tr}}=\C[[a_1,\ldots,a_r]]\), \(D=B_{\mathrm{tr}}[[w]]\), and \(G_b=P_b(f_0(w)+a)\). At \(a=0\), the specialization of \(G_b\) is \(P_b(f_0(w))\), which has finite order \(d_b\) and is therefore nonzero. Lemma~\ref{lem:finite-slice-flatness} applies to the one-by-one matrix \(G_b\), making its cokernel free of rank \(d_b\) over \(B_{\mathrm{tr}}\), with basis \(1,w,\ldots,w^{d_b-1}\), since these classes are a basis on the closed fibre.  Expressing \(w^{d_b}\) in this basis gives a monic polynomial \(W_b(a,w)\), whose lower coefficients vanish at \(a=0\).  It belongs to \((G_b)\). The same freeness criterion applied to \(W_b\) shows that the surjection \(D/(W_b)\to D/(G_b)\) is an isomorphism between free modules with the same displayed basis.  Hence \((W_b)=(G_b)\); since \(D\) is a domain, \(G_b\) is a unit times this monic distinguished polynomial of degree \(d_b\).

The coordinate change \(t_\nu=a_\nu+f_\nu(w)\) gives
\[
 D/(G_b)\cong (R/(P_b))[[w]],
\]
which is a domain.  Its generic fibre over \(B_{\mathrm{tr}}\) is therefore a finite field extension.  In characteristic zero that extension is separable, so the discriminant of \(W_b\) is nonzero.

For \(b\ne c\), the incidence module \(D/(G_b,G_c)\) is finite over \(B_{\mathrm{tr}}\).  The coordinate change identifies it with \((R/(P_b,P_c))[[w]]\), of dimension at most \(r-1\), whereas \(\dim B_{\mathrm{tr}}=r\).  It is consequently torsion, so the resultant of \(W_b,W_c\) is nonzero.  To impose avoidance of \(\Sigma\), let \(K_\Sigma\) be its radical ideal.  Since \(\prod_bP_b\in K_\Sigma\), the module
\[
 D/K_\Sigma(f_0(w)+a)
\]
is finite over \(B_{\mathrm{tr}}\).  After the coordinate change its dimension is at most \(r-1\), whereas \(\dim B_{\mathrm{tr}}=r\); it is therefore torsion.  Choose a nonzero annihilator, taking \(1\) when the module is zero.  Multiply it by the discriminants and resultants already obtained, and call the resulting nonzero element \(b(a)\).

Let \(b_h\) be the lowest nonzero homogeneous term of \(b\). For every \(v\) with \(b_h(v)\ne0\), \(b(sv)=s^h(b_h(v)+O(s))\ne0\).  This element is invertible over \(\C((s))\), proving separability, pairwise coprimality, and avoidance of \(\Sigma\).  Division by the monic polynomial \(\prod_bW_b\) shows that the power-series quotients used to test avoidance of \(\Sigma\) can be computed in \(B_{\mathrm{tr}}[w]/(\prod_bW_b)\).  Thus ordinary substitution \(a=sv\) agrees with completed substitution and preserves the annihilation relation for \(D/K_\Sigma(f_0(w)+a)\), with the annihilating element \(b(a)\) specialized to \(b(sv)\).

The roots of a distinguished monic polynomial are integral over \(\C[[s]]\) and reduce to zero.  A finite splitting extension is a complete discretely valued field with residue field \(\C\), so its valuation ring is \(\C[[u]]\).  Absorbing a unit into \(u\) gives \(s=u^e\), and every root lies in \(u\C[[u]]\).
\end{proof}

The local equations of Katz--Morrison's simultaneous partial resolution describe the contraction at each root.  Nabijou--Wemyss identify its discriminant components with restricted roots and the unique curve class at a general point of each component \cite[Theorem~2.6]{NabijouWemyss}.  Apply this description to the translated slice of Lemma~\ref{lem:formal-transverse-translation}.

\begin{proposition}
\label{prop:formal-marked-morsification}
Let \(f_0(w)\) be a formal arc classifying the marked contraction \(X_0\to Y_0\). There is a translated arc \(F(s,w)\), corresponding to a homomorphism \(R_J\to\C[[s,w]]\), with \(F(0,w)=f_0(w)\). Its geometric generic contraction has exactly \(n_\beta\) disjoint split \((-1,-1)\)-curves of each class \(\beta\).  At each such curve the completed target is a split ordinary double point and the source is its chosen small resolution.  The generic fibre is formed from the full finite relative family.
\end{proposition}

\begin{proof}
The marked base \(R_J\) is regular, and its curve-divisor equations \(P_\beta\) are distinct primes by Lemma~\ref{lem:multicurve-root-products}.  Their orders on the central arc are
\[
 \operatorname{ord}_wP_\beta(f_0(w))=n_\beta
\]
by Proposition~\ref{thm:multicurve-algebra-flatness} and Corollary~\ref{cor:multicurve-cartan-gram}.  We impose transversality on every component of the surface discriminant, including the components for roots with zero retained coefficients. At a point corresponding to a root whose retained coefficients are all zero, this ensures that the threefold slice is smooth.

The central arc is contained in none of these components.  Write the central isolated cDV hypersurface as \(G(x_1,x_2,x_3,w)=0\).  The relative critical scheme defined by \(G,G_{x_1},G_{x_2},G_{x_3}\) is finite over \(\C[[w]]\), because the central surface of the cDV hypersurface is ADE; this is the topological finiteness argument of Lemma~\ref{lem:topological-finiteness-expanded}.  If the relative critical scheme had a component dominating the \(w\)-axis, the normalization of that component would have a parameter \(u\) with \(w=u^e\).  Differentiating \(G(x(u),u^e)=0\), while using the three vanishing partial derivatives, would give \(G_w=0\) along that curve.  This contradicts the isolated singularity of the threefold.

Take \(\Sigma\) to be the reduced image, under \(\mathfrak h\to\mathfrak h/W_I\), of intersections of two distinct root hyperplanes in \(\mathfrak h\). The finite quotient preserves the codimension of this locus, so \(\Sigma\) has codimension at least two; \(\Sigma\) contains the singular and intersection loci of the full discriminant components.  Apply Lemma~\ref{lem:formal-transverse-translation} to the reduced equations of the full discriminant components.  For an intersection of the resulting slice with a curve divisor, choose a lift of its parameter point to \(\mathfrak h\). Exactly one positive root \(\alpha\) vanishes at this lift, with \(\pi_I(\alpha)=\beta\ne0\).

The stabilizer in the full Weyl group is \(\{1,s_\alpha\}\): a real vector's stabilizer is generated by reflections in its vanishing roots, as follows by moving it to the closed dominant chamber and expanding \(v-wv\) along a reduced reflection word into positive roots with nonnegative coefficients.  For a complex vector \(a+ib\), choose a generic real \(c\) so that \(a+cb\) has the same stabilizer and vanishing roots.  These are finite linear fixed-subspace conditions, so the assertion also holds over algebraically closed extension fields. Since \(W_I\) fixes the retained fundamental weights and \(\pi_I(\alpha)\ne0\), it does not contain \(s_\alpha\). Thus \(\mathfrak h\to\mathfrak h/W_I\) is \'etale here.  Completing a finite quotient at its orbit gives invariants under the stabilizer of one lift.  A coordinate \(\zeta\) normal to \(H_\alpha\subset\mathfrak h\) is negated by \(s_\alpha\), while the tangential coordinates are fixed; hence the unmarked normal parameter is \(\tau=\zeta^2\).

The fibre of the semiuniversal ADE family at the corresponding point of \(\mathfrak h/W\) has one \(A_1\) singularity, because exactly one positive root vanishes at the chosen lift.  This is the single-coincidence case in the \(A,D\) constructions and the single effective-root configuration in the \(E\) construction \cite[Theorem~1, Sections~4--5 and Table~4]{KatzMorrison}. The same paper's preferred versal equations contain a constant-term base coefficient, so the total unmarked hypersurface is smooth \cite[Section~3, Tables~2--3]{KatzMorrison}.

To obtain the local equation, use the nondegenerate Hessian in the fibre variables to solve for a formal critical section.  After translating by that section, eliminate a fibre-homogeneous term \(r_k\) of degree \(k\ge3\) by the substitution \(\mathbf{x}\mapsto\mathbf{x}-k^{-1}H^{-1}\nabla r_k\), where \(\mathbf{x}\) is the column of three fibre variables and \(H\) is the quadratic Hessian.  Euler's identity removes \(r_k\); the successive substitutions converge.  Diagonalizing the quadratic form gives \(xy-\xi^2=g\).  The critical-value function \(g\) vanishes exactly on the smooth discriminant \(\tau=0\), and total-space smoothness gives \(dg\ne0\).  Thus \(g=\tau\) times a unit.  Absorbing that unit in \(\tau\) and its square root in \(\zeta\) gives the marked local model
\[
 xy-\xi^2=\zeta^2.
\]

No Levi root vanishes at this point, so the morphism from the full resolution to the partial resolution contracts nothing near its sole exceptional curve. Formal versality for the resolved \(A_1\) surface identifies its resolved parameter \(\sigma\), whose affine smoothing parameter is \(\sigma^2\), with one of \(\zeta,-\zeta\) \cite[Section~2.3]{Karmazyn}.  These are the two explicit simultaneous \(A_1\) resolutions in \cite[Section~4]{KatzMorrison}. The simple root of \(P_\beta(F)\) makes \(\zeta\) a parameter along the slice.  After a finite splitting extension the equation is \(xy-p_+p_-=0\), and the chosen source is its small resolution. The charts of Lemma~\ref{lem:split-odp-projectives} give the exceptional curve its normal bundle \(\mathcal O(-1)^{\oplus2}\).  Its retained curve class is \(\beta\) \cite[Theorem~2.6 and its proof]{NabijouWemyss}. Only finitely many intersection points with the discriminant and lifts to \(\mathfrak h\) occur; their residue field extensions and the constants needed to split the quadratic forms lie in a single finite extension.  The formal coordinate recursions then remain over that field.

If the vanishing root \(\alpha\) has \(\pi_I(\alpha)=0\), then \(s_\alpha\in W_I\).  Here the marked normal parameter is already \(\tau\), so the transverse model \(xy-\xi^2=\tau\) has smooth threefold slice and no contracted curve. Away from the full discriminant the surface family is smooth.  Hence the roots with nonzero retained coefficients account for every exceptional curve.  The distinguished polynomial for \(P_\beta(F)\) has degree \(n_\beta\) and simple roots, which gives the stated multiplicity.  Forming the finite relative family before inverting \(s\) retains every moving root, as in Remark~\ref{rem:generic-fibre-ring}.
\end{proof}

\begin{proposition}
\label{prop:marked-odp-projector-ranks}
At a split ordinary double point of class \(\beta\) in Proposition~\ref{prop:formal-marked-morsification}, the completed relative contraction quotient is
\[
 \End_{K'}\left(\bigoplus_{i\in J}(K')^{\beta_i}\right),
\]
and \(e_i\) is the projection onto \((K')^{\beta_i}\). The maximal ideal of the completed target acts trivially on this quotient.
\end{proposition}

\begin{proof}
Fix the normalized positive generator and its labelled summands. The comparison in Proposition~\ref{prop:multicurve-cartan-comparison} and the completed pullback in Proposition~\ref{prop:multicurve-classifying-pullback} preserve those summands separately.  All full endomorphisms in this paragraph are formed on the specified completed target model.

After a finite ramified extension, lift the moving intersection point under \(\mathfrak h\to\mathfrak h/W_I\). The lift is a formal arc through \(0\in\mathfrak h\) on one root hyperplane \(H_\alpha\), and its generic point lies on no other root hyperplane.  Write \(\alpha=\sum_k c_k\alpha_k\), so that \(\beta_i=c_i\) for \(i\in J\). By Lemma~\ref{lem:chosen-root-divisor-family}, on the chosen standard full simultaneous resolution the effective Cartier divisor \(\sum_k c_kC_k\) is the scheme-theoretic central fibre of a projective flat family over \(H_\alpha\). Pulling this family over \(H_\alpha\) back to the formal arc gives a proper flat curve family over its complete discrete valuation ring.  Its generic fibre is the unique root curve, and its central fundamental cycle is \(\sum_k c_k[C_k]\).  The morphism from the full resolution to the partial resolution contracts the curves indexed by \(I\) and is an isomorphism on this generic non-Levi root curve \cite[Theorem~2.6 and its proof]{NabijouWemyss}. The family is defined over \(H_\alpha\), where the root \(\alpha\) remains specified.

Pull back \(\det P_i\) to this proper family of curves. The compatible formal restrictions of \(\det P_i\) algebraize to a line bundle by projective formal existence \cite[Lemma~30.24.3, Tag~\texttt{0885}]{Stacks}; the nonfree locus is empty since its closed image in the local base would otherwise meet the closed fibre.  The degree is constant in the flat proper family.  More explicitly, the derived pushforwards of this line bundle and of the structure sheaf are perfect and commute with base change \cite[Lemmas~36.35.5 and~36.35.10, Tags~\texttt{0DI4} and~\texttt{0DJT}]{Stacks}; over the local base, alternating sums of the ranks of finite free complexes give their fibrewise Euler characteristics.  Their difference is the degree.  On the central fundamental cycle it is the sum of the component degrees with multiplicities \(c_k\).  The pulled-back line bundle has degree zero on deleted curves and degree \(\delta_{ij}\) on each retained curve. Consequently, separately for every label,
\[
 \deg_{C_\beta}\det P_i
 =\sum_{k\in J}c_k\delta_{ik}=\beta_i.
\]
For one retained vertex, this specializes to the determinant-degree calculation in Lemma~\ref{lem:marked-degree-specialization}.

The positivity proved for the algebraic lift in Proposition~\ref{prop:multicurve-cartan-comparison} applies to \(P=(T_J^{\mathrm{alg}})^\vee\), with the prescribed labelled summands.  Its finite evaluation surjection pulls back to a surjection onto the positive generator here, and its total determinant remains relatively ample under this pullback and the specified completion at the central target point.  The structure-sheaf summand is preserved, while \(H^1(P^\vee)=0\) is the \(\mathcal O,P^\vee\)-component of the Ext vanishing in the tilting comparison.  Together with \(R\pi_*\mathcal O=\mathcal O\) and the one-dimensional fibre bound, these are the hypotheses of the positive perverse-generator criterion \cite[Propositions~3.2.6--3.2.7]{VanDenBergh}.  Thus the individual summands of \(P\) are projectives in the fixed positive perverse heart on the completed split ODP resolution.

The split Atiyah calculation of Lemma~\ref{lem:split-odp-projectives} gives exactly two indecomposable positive perverse projectives, \(\mathcal O\) and \(\mathcal Q\), with determinant degrees zero and one.  Consequently, if \(q_i=\operatorname{rank}P_i\), Krull--Schmidt and the displayed degree give, separately for each label,
\[
 \widehat P_i\cong
 \mathcal O^{\oplus(q_i-\beta_i)}\oplus\mathcal Q^{\oplus\beta_i}.
\]
After pushforward, the generator is a direct sum of free modules and \(\bigoplus_i I^{\oplus\beta_i}\), where \(I\) is the rank-one Atiyah module.  The original \(e_0\) projects onto one free summand. The ideal it generates consists of all maps factoring through finite free modules, so it also kills every additional free summand in this decomposition.  Lemma~\ref{lem:atiyah-stable-endomorphisms} identifies \(\underline{\End}(I)\) with \(K'\), with the target maximal ideal acting by zero.  The surviving matrix entries are therefore copies of \(K'\), with their usual composition.  Grouping them by the original summands gives precisely the asserted projectors and multiplication.
\end{proof}

\begin{lemma}
\label{lem:marked-polynomialization}
Let \(L\) be a finite free \(\C[[u]]\)-algebra with marked special fibre \(A\), and suppose a marked generic isomorphism identifies \(L[u^{-1}]\) with \(S_E\otimes\C((u))\).  There is a finite free marked \(\C[u]\)-algebra \(\mathscr P\) whose completion is \(L\) and whose localization at \(u\) is \(S_E\otimes\C[u,u^{-1}]\).
\end{lemma}

\begin{proof}
The generic marked isomorphism lets us regard \(L\) as a lattice in \(S_E((u))\).  For some \(m\ge0\),
\[
 u^mS_E[[u]]\subseteq L\subseteq u^{-m}S_E[[u]].
\]
Define \(\mathscr P=L\cap S_E[u,u^{-1}]\).  It is a unital subalgebra containing the constant labelled projectors, and
\[
 u^mS_E[u]\subseteq\mathscr P\subseteq u^{-m}S_E[u].
\]
The two inclusions show that \(\mathscr P\) is finite and torsion-free over the principal ideal domain \(\C[u]\), hence free, and its localization is the stated split algebra.

To compare completions, take an element of \(L\) and truncate its Laurent expansion far enough.  The truncation lies in \(\mathscr P\), because the discarded Laurent tail lies in \(u^mS_E[[u]]\subseteq L\), and the truncation approximates the original element to arbitrarily high order.  The same inclusions show that the ambient \(u\)-adic topology agrees with the lattice topology.  Completion therefore identifies \(\mathscr P\otimes_{\C[u]}\C[[u]]\) with \(L\), preserving products and projectors.  Reduction modulo \(u\) gives the prescribed marked algebra \(A\).
\end{proof}

For finite-dimensional algebras without distinguished idempotents, Jelisiejew--Smoktunowicz relate formal, polynomially split, and orbit-closure formulations of deformation \cite[Theorems~2.19--2.20]{JS}. The lattice argument retains the chosen embedding of \(E\) and the constant projectors while passing from the formal family to the polynomial one.

\begin{proof}[Proof of Theorem~\ref{thm:marked-semisimple-existence}]
Let \(f_0(w)\) be a formal arc classifying the given marked contraction. Apply Lemma~\ref{lem:formal-transverse-translation} to this arc and the curve-divisor equations \(P_\beta\). Fix a translation vector \(v\) in the open set supplied by the lemma. In the chosen coordinates the translated arc is \(F(s,w)=f_0(w)+sv\). Its pullback sends \(P\in R_J\) to \(P(F(s,w))\in\C[[s,w]]\).  For its excluded locus, intersect the curve divisors with the reduced image of all intersections of two distinct root hyperplanes in \(\mathfrak h\). Each such intersection has codimension two in \(\mathfrak h\), and the finite quotient \(\mathfrak h\to\mathfrak h/W_I\) preserves its dimension, so the excluded locus has codimension at least two.  The orders of the central root functions are \(n_\beta\), by Corollary~\ref{cor:multicurve-cartan-gram}.

Pull back the finite-type standard simultaneous partial resolution and then complete at the contraction point of the target and its inverse image. Parameter flatness is preserved by this base change and by the flat Noetherian target completion.  The prescribed central contraction and all its labelled generator summands are retained. Proposition~\ref{prop:multicurve-classifying-pullback} identifies the relative contraction quotient with the ordinary tensor
\[
 \mathscr A=C_J\otimes_{R_J,F}\C[[s,w]].
\]
The square presentation in Theorem~\ref{thm:multicurve-weighted-presentation} has closed determinant \(\prod_\beta P_\beta(f_0(w))^{s_\beta^2}\), up to a unit. Every factor is nonzero and has order \(n_\beta\), as checked above. Lemma~\ref{lem:finite-slice-flatness} therefore makes \(\mathscr A\) finite free over \(\C[[s]]\), of rank \(\sum_\beta n_\beta s_\beta^2\), with marked special fibre \(A\).  The tensor statement here concerns the contraction quotient, which is finite over \(R_J\). The full endomorphism algebras are compared on the specified completed geometric model.

By Proposition~\ref{prop:multicurve-reduced-annihilator}, \(\delta_J=\prod_\beta P_\beta\) annihilates \(C_J\). Its pullback is a unit times a product \(W=\prod_\beta W_\beta\) of distinguished polynomials, where \(\deg_w W_\beta=n_\beta\). The finite algebra
\[
 \C[[s,w]]/(W)=\C[[s]][w]/(W)
\]
is formed before inverting \(s\).  After a finite splitting extension \(K_1/\C((s))\), its distinct roots \(w_{\beta,a}\) therefore give, by the Chinese remainder theorem, a product decomposition
\[
 \mathscr A\otimes K_1
 =\prod_{\beta,a}
 C_J\otimes_{R_J,\,F(s,w_{\beta,a})}K_1,
 \qquad 1\le a\le n_\beta.
\]
Each root has positive valuation, so evaluation of the formal classifying map is defined in the valuation ring before localization. This construction retains every support point of the finite relative algebra.

Pass temporarily to an algebraic closure of \(K_1\). Lift the parameter point of one factor to \(\Spec S\) along the finite faithfully flat morphism \(\Spec S\to\Spec R_J\).  Avoidance of the excluded locus means that exactly one positive root \(\alpha\) vanishes there.  Its retained vector is \(\beta\), because the image of \(P_\beta\) in \(S\) is the product of the root forms whose retained coefficients are \(\beta\).  Lemma~\ref{lem:multicurve-single-root-block}, together with the comparison of labelled generator summands in Proposition~\ref{prop:multicurve-cartan-comparison}, identifies this factor with
\[
 \End\left(\bigoplus_{i\in J}\overline{K_1}^{\,\beta_i}\right),
\]
with each \(e_i\) its indicated projector.  Thus the matrix multiplication and every individual label are identified.  A matrix for each marked isomorphism and its inverse has finitely many coefficients algebraic over \(K_1\).  A single finite extension \(K'/K_1\) contains all these coefficients and splits every factor, giving the stated marked product \(S_E\otimes K'\).

The finite-extension argument in Lemma~\ref{lem:formal-transverse-translation} writes the valuation ring of \(K'\) as \(\C[[u]]\), with \(s=u^e\). Apply Lemma~\ref{lem:marked-polynomialization} to \(L=\mathscr A\otimes_{\C[[s]]}\C[[u]]\). It gives the \(\C[u]\)-algebra family with its prescribed special fibre, all constant projectors, and completion equal to this entire marked formal family.
\end{proof}

\begin{remark}
\label{rem:marked-existence-kawamata}
Kawamata begins with an analytic disk deformation and proves flatness of the relative contraction algebra; when the exceptional loci in the noncentral fibres consist of disjoint \((-1,-1)\)-curves, he obtains the generic matrix product \cite[Theorem~3.1(2)--(3)]{KawamataNC}.  His proof also decomposes the restriction of each original summand with multiplicities \(m_{j,i}\), the geometric source of the projector ranks.

The formal classifying-arc statement needed here is established through the universal quotient with the embedding of \(E\) fixed, and is then replaced by a polynomial family with the same completion.  The comparison after forgetting the idempotents with other deformation formulations is supplied by \cite[Theorems~2.19--2.20]{JS}; the marked lattice step is written out here.
\end{remark}

\subsection{Fixed-idempotent laws and the orbit criterion}
\label{subsec:marked-orbit-criterion}

The preceding theorem constructs a family whose generic fibre is a product of matrix algebras with fixed idempotents.  The unmarked structure-constant and orbit method used here is reviewed by Jelisiejew--Smoktunowicz: semisimple orbits are open, and the formal, polynomially split, and orbit-closure deformation notions agree after the appropriate finite extensions \cite[Section~2.1, Theorems~2.14 and~2.20]{JS}.  With a fixed separable algebra \(E\), the method shows that finite free associative formal deformations of \(A\) with geometrically semisimple generic fibre have the same simple matrix factors and the same dimension vectors relative to the lifted idempotents.

The relative tangent formula uses the centralizer of \(E\), which is generally larger than the centre.  Peirce dimensions and the number of simple factors are necessary numerical inputs, but they do not determine the marked factorization by themselves.

Let \(S\) be a semisimple algebra whose simple modules have dimension vectors \(v_a\) relative to \((e_i)_{i\in J}\), and put \(c_{ij}=\dim e_iSe_j\).  The Gram matrix \(C=(c_{ij})=\sum_a v_av_a^{\mathsf T}\) and the number of factors \(\#\{a\}\) do not, by themselves, determine the vectors.  For example, the two multisets
\[
 \{(2,1),(2,1),(2,1),(0,3),(1,0)\},\qquad
 \{(3,0),(1,2),(1,2),(1,2),(1,0)\}
\]
both have five vectors, Gram matrix \(\left(\begin{smallmatrix}13&6\\6&12\end{smallmatrix}\right)\), and unmarked algebra \(M_3(\C)^4\times\C\).  They are distinct even after exchanging the two labels.  Thus the Gram matrix and the number of simple modules admit several multisets of dimension vectors.  The orbit argument in Theorem~\ref{thm:marked-general-uniqueness} starts from the product constructed geometrically in Theorem~\ref{thm:marked-semisimple-existence} and proves that every geometrically semisimple generic fibre, together with the lifted idempotents, is isomorphic to it.

Let \(V=\bigoplus_{i,j\in J}V_{ij}\) be a fixed \(E\)-bimodule, with \(c_{ij}=\dim_\C V_{ij}\) and distinguished nonzero vectors \(e_i\in V_{ii}\).  Let \(X_E(V)\) be the scheme of associative laws with these fixed idempotents.  It parametrizes associative \(E\)-balanced products \(V\otimes_E V\to V\) with these prescribed left and right \(E\)-actions.  Thus \(1=\sum_i e_i\) is fixed before the multiplication varies, and the associativity equations define the scheme structure. Its change-of-basis group is
\[
 G_E=\{g\in\Aut_{E\text{-}E}(V):g(e_i)=e_i\text{ for all }i\}.
\]
This group is smooth and irreducible: its factors are \(\GL(V_{ij})\) for \(i\ne j\) and the stabilizers of the nonzero vectors \(e_i\) in \(\GL(V_{ii})\).  Set
\begin{equation}
\label{eq:marked-group-dimensions}
 b=\sum_i c_{ii},\qquad
 g=\dim G_E=\sum_{i,j}c_{ij}^2-b.
\end{equation}
For a law \(B\in X_E(V)\), its centralizer \(B^E=\{x:xe=ex\text{ for all }e\in E\}\) is \(\bigoplus_i e_iBe_i\), of dimension \(b\).

To compute tangent vectors with the prescribed \(E\)-actions, we use relative Hochschild cochains. Their normalization keeps both the left and the right action fixed.

\begin{lemma}
\label{lem:marked-tangent-dimensions}
For \(B\in X_E(V)\), normalized relative Hochschild cochains compute ordinary Hochschild cohomology.  If \(Z_E^2(B)\) and \(B_E^2(B)\) are their degree-two cocycles and coboundaries, then
\[
 T_BX_E(V)=Z_E^2(B),\qquad
 T_B(G_E\cdot B)=B_E^2(B).
\]
In particular,
\begin{equation}
\label{eq:marked-tangent-dimension}
 \dim T_BX_E(V)
 =g-b+\dim Z(B)-\dim HH^1(B)+\dim HH^2(B).
\end{equation}
\end{lemma}

\begin{proof}
The normalized relative complex has
\[
 C_E^0(B)=B^E,\qquad
 C_E^n(B)=\operatorname{Hom}_{E^e}
       ((B/E)^{\otimes_E n},B)\quad(n\ge1).
\]
Insertion of the identity gives a contracting homotopy for the augmentation of the relative bar complex, so that complex is exact.  Normalization removes its split degenerate subcomplex and leaves the terms \(B\otimes_E(B/E)^{\otimes_E n}\otimes_E B\). Each is induced from an \(E^e\)-module.  Since \(E^e\) is a product of copies of \(\C\), every such module is projective, and induction makes each term \(B^e\)-projective.  Applying \(\operatorname{Hom}_{B^e}(-,B)\) to this projective resolution of \(B\) gives the normalized relative cochain complex. Its cohomology is therefore ordinary \(HH^*(B)\).

Linearizing associativity gives the relative two-cocycle equation: fixed \(E\)-actions require the cochain to vanish when an input belongs to \(E\).  Linearizing a change of \(E\)-bimodule basis gives the coboundaries.  The group scheme of algebra automorphisms fixing \(E\) is a closed subgroup scheme of \(G_E\), hence is smooth in characteristic zero \cite[Lemma~39.8.2, Tag~\texttt{047N}]{Stacks}.  Consequently the orbit map \(G_E\to G_E\cdot B\) is smooth \cite[Proposition~1.11 and Theorem~1.16]{Brion}.  The orbit has the asserted tangent space, and the stabilizer has Lie algebra \(\Der_E(B)\).

The degree-zero part of the relative Hochschild complex gives the exact sequence
\[
 0\longrightarrow Z(B)\longrightarrow B^E
 \longrightarrow\Der_E(B)\longrightarrow HH^1(B)\longrightarrow0.
\]
Moreover, \(\dim C_E^1(B)=g\), because these cochains are precisely the \(E\)-bimodule endomorphisms vanishing on the distinguished vectors.  It follows that \(\dim B_E^2(B)=g-b+\dim Z(B)-\dim HH^1(B)\). Adding \(\dim HH^2(B)\) proves the formula.
\end{proof}

\begin{lemma}
\label{lem:marked-semisimple-orbits}
A finite-dimensional semisimple complex algebra with a fixed complete labelled system of orthogonal idempotents is classified, up to marked isomorphism, by a multiset of nonzero vectors \(v_a=(v_{a,i})_{i\in J}\) of nonnegative integers.  Its algebra and Peirce dimensions are
\[
 S(\{v_a\})=\prod_a\End_\C\left(\bigoplus_i\C^{v_{a,i}}\right),
 \qquad c_{ij}=\sum_a v_{a,i}v_{a,j}.
\]
Every such orbit in \(X_E(V)\) is open, has dimension \(g-b+\dim Z(S)\), and is dense in an irreducible component.
\end{lemma}

\begin{proof}
On a simple matrix factor, simultaneous orthogonal projections are classified by the dimensions of their images.  Inner automorphisms change bases in those images, and an algebra automorphism may permute factors carrying the same vector.  This proves the classification and the Peirce formula; in particular \(\dim Z(S)\) is the number of vectors.  A marked isomorphism between laws on \(V\) is exactly an element of \(G_E\).

Semisimplicity gives \(HH^n(S)=0\) for \(n>0\), so Lemma~\ref{lem:marked-tangent-dimensions} makes the tangent and orbit dimensions equal to \(g-b+\dim Z(S)\).  The stabilizer is a product of groups \((\prod_i\GL_{v_{a,i}})/\mathbb G_m\), followed by finite permutations of identical factors, and is smooth.  Apply Lemma~\ref{lem:open-orbit-smoothness}(i) to the smooth irreducible group \(G_E\).  The orbit is open, so its closure is an irreducible component.
\end{proof}

\begin{theorem}
\label{thm:marked-orbit-criterion}
Let \(B\in X_E(V)\), and suppose a marked semisimple algebra \(S\in X_E(V)\) satisfies
\[
 B\in\overline{G_E\cdot S},\qquad
 \dim Z(B)=\dim Z(S),\qquad
 \dim HH^1(B)=\dim HH^2(B).
\]
Then \(B\) is a smooth point of \(X_E(V)\), and \(\overline{G_E\cdot S}\) is its unique irreducible component. If another marked semisimple law \(S'\) has \(B\in\overline{G_E\cdot S'}\), then \(S'\cong S\) as marked algebras.
\end{theorem}

\begin{proof}
In the tangent formula of Lemma~\ref{lem:marked-tangent-dimensions}, the assumed equality of the two Hochschild dimensions cancels the cohomology terms. The equality of centre dimensions then matches the semisimple orbit dimension in Lemma~\ref{lem:marked-semisimple-orbits}:
\[
 \dim T_BX_E(V)=g-b+\dim Z(S)=\dim(G_E\cdot S).
\]
The semisimple orbits are open by Lemma~\ref{lem:marked-semisimple-orbits}.  Applying Lemma~\ref{lem:open-orbit-smoothness}(ii) proves smoothness, the unique-component assertion, and equality of the two marked orbits.
\end{proof}

\subsection{Periodicity and general uniqueness}
\label{subsec:marked-general-uniqueness}

Return to the multi-curve contraction algebra \(A\), and write
\begin{equation}
\label{eq:marked-numerical-inputs}
 c_{ij}=\dim e_iAe_j=\sum_\beta n_\beta\beta_i\beta_j,
 \qquad N_{\mathrm{GV}}=\sum_\beta n_\beta=\dim Z(A)=\dim HH_0(A).
\end{equation}
Corollary~\ref{cor:multicurve-cartan-gram} and Theorem~\ref{thm:multicurve-hua-toda}, restricted to the central fibre, place \(A\) and the geometrically constructed product on the same \(E\)-bimodule.  The Peirce dimensions determine the fixed-idempotent law scheme. Together with the Hochschild equality proved below, the cocenter dimension gives the required tangent dimension.  The formal family proves that \(A\) lies in the closure of the constructed semisimple orbit.

\begin{proposition}
\label{prop:multicurve-hochschild-equality}
The algebra \(A\) is symmetric and satisfies
\[
 \Omega^4_{A^e}(A)\simeq A
 \quad\text{in the stable bimodule category},\qquad
 DHH^1(A)\cong HH^2(A).
\]
\end{proposition}

\begin{proof}
The target is an isolated cDV hypersurface with a smooth crepant resolution, so August's symmetry theorem applies to the stable endomorphism algebra of the entire nonfree maximal Cohen--Macaulay generator \cite[Theorem~2.6]{August}.  The four-term projective bimodule extension in \cite[Setting~3.2.1 and Theorem~3.2.5]{JKM} also applies to this generator and has the diagonal bimodule at both ends.  Thus the untwisted four-periodicity used in Proposition~\ref{prop:bimodule-periodicity} applies here.  For a symmetric algebra with this four-periodicity, Lemma~\ref{lem:tate-hochschild-duality} identifies the dual of \(HH^1(A)\) with \(HH^2(A)\).  Each semisimple direct factor has zero positive Hochschild cohomology.
\end{proof}

\begin{lemma}
\label{lem:marked-idempotent-lifting}
Let \((B,\mathfrak m,\C)\) be complete Noetherian local, and let \(\mathscr D\) be a finite free \(B\)-algebra with identified special fibre \(A\).  The system \((e_i)\) lifts to orthogonal idempotents \((\widetilde e_i)\) summing to one.  Any two such lifts are simultaneously conjugate by a unit congruent to one modulo \(\mathfrak m\mathscr D\).  Each \(\widetilde e_i\mathscr D\widetilde e_j\) is free of rank \(c_{ij}\), and an adapted trivialization defines a morphism \(\Spec B\to X_E(V)\) whose closed point is \(A\).
\end{lemma}

\begin{proof}
Across a square-zero ideal, if \(x^2-x=q\) lies in that ideal, then \(x+(1-2x)q\) is idempotent.  Lift the first projector and work successively in its complementary corner to preserve orthogonality; take the last projector to be the remaining complement.  Applying this construction to successive powers of \(\mathfrak m\) and using completeness gives the system.

For two lifts \((p_i)\) and \((q_i)\), the element
\[
 h=\sum_iq_ip_i\equiv1\pmod{\mathfrak m\mathscr D}
\]
is invertible by completeness and satisfies \(hp_i=q_ih\). This proves simultaneous conjugacy.

The maps \(x\mapsto p_ixp_j\) are mutually orthogonal \(B\)-linear projectors.  Their images are direct summands of the free module \(\mathscr D\), hence free over the local ring, with ranks determined by the special fibre.  Choose bases lifting fixed bases of \(e_iAe_j\); in each diagonal block include \(p_i\), which is unimodular because its reduction is nonzero.  These choices give the claimed \(E\)-bimodule trivialization and classifying morphism.
\end{proof}

\Needspace{27\baselineskip}
\begin{theorem}
\label{thm:marked-general-uniqueness}
Let \(\mathscr D\) be any finite free unital \(\C[[t]]\)-algebra with special fibre identified with \(A\), and suppose its geometric generic fibre is semisimple.  For any lifts of the labelled original idempotents, that fibre is marked-isomorphic to \(S_E\otimes_\C\overline{\C((t))}\).  Equivalently, the multiset of dimension vectors of its simple modules, relative to the lifted idempotents, is exactly
\begin{equation}
\label{eq:marked-unique-multiset}
 \{\,\beta\text{ repeated }n_\beta\text{ times}:
                    \beta\in\mathcal B_J\,\}.
\end{equation}
Thus there is exactly one isomorphism class of semisimple geometric generic fibre together with the lifted idempotents.  After forgetting the idempotents, this fibre is obtained by scalar extension from the complex semisimple algebra
\[
 \prod_{m\ge1}M_m(\C)^{\times N_m},\qquad
 N_m=\sum_{s_\beta=m}n_\beta.
\]
Moreover \(A\) is smooth in \(X_E(V)\), and the unique irreducible component through \(A\) is \(\overline{G_E\cdot S_E}\).
\end{theorem}

\begin{proof}
The Peirce formula in \eqref{eq:marked-numerical-inputs} places \(A\) and \(S_E\) on the same \(E\)-bimodule \(V\).  Theorem~\ref{thm:marked-semisimple-existence} places \(A\) in \(\overline{G_E\cdot S_E}\): choose bases in every Peirce summand of the free family, including each projector, and use the generic marked isomorphism.  Every equation of the orbit closure vanishes at the generic point after a finite field extension, hence on the whole parameter DVR, and therefore at \(A\).  Equation \eqref{eq:marked-numerical-inputs} gives \(\dim Z(A)=N_{\mathrm{GV}}=\dim Z(S_E)\), while Proposition~\ref{prop:multicurve-hochschild-equality} gives \(\dim HH^1(A)=\dim HH^2(A)\).  Theorem~\ref{thm:marked-orbit-criterion} now proves smoothness and the unique-component assertion.

Let \(\mathscr D\) be an arbitrary family as in the statement.  Lift the projectors and choose the compatible bases supplied by Lemma~\ref{lem:marked-idempotent-lifting}.  The geometric generic fibre of \(\mathscr D\) is represented by a constant semisimple complex algebra \(S'\) with lifted idempotents, classified by the multiset of dimension vectors of its simple \(S'\)-modules.  The finitely many coefficients of a marked isomorphism lie in a finite extension of \(\C((t))\).  The same orbit-closure argument gives \(A\in\overline{G_E\cdot S'}\), so the orbit criterion forces \(S'\cong S_E\).  This proves \eqref{eq:marked-unique-multiset}.  Changing the lifted idempotents conjugates them simultaneously, so the answer does not depend on their choice.  Forgetting the projectors combines precisely the factors with the same total size \(s_\beta\), giving the stated unmarked algebra.
\end{proof}

\begin{corollary}
\label{cor:marked-smooth-laws}
Put \(d=\dim A\), and let \(b,g\) be as in \eqref{eq:marked-group-dimensions}.  Then
\begin{align*}
 \widehat{\mathcal O}_{X_E(V),A}&\cong
     \C[[x_1,\ldots,x_{g-b+N_{\mathrm{GV}}}]],\\
 \widehat{\mathcal O}_{X_d,A}&\cong
     \C[[y_1,\ldots,y_{d(d-1)-d+N_{\mathrm{GV}}}]].
\end{align*}
The unique component of \(X_d\) through \(A\) is the closure of the unmarked orbit of \(S_E\).  In every finite free associative deformation of \(A\) over a complete Noetherian local complex base \((B,\mathfrak m,\C)\), every semisimple geometric fibre has the multiset of dimension vectors of simple modules given in \eqref{eq:marked-unique-multiset}, with respect to any lifts of the original projectors.
\end{corollary}

\begin{proof}
Theorem~\ref{thm:marked-general-uniqueness} gives smoothness at \(A\) in \(X_E(V)\). For smoothness in \(X_d\), forgetting the fixed idempotents in the existence theorem places \(A\) in the unmarked orbit closure of \(S_E\). Proposition~\ref{prop:multicurve-hochschild-equality} gives equal first and second Hochschild dimensions, and \(\dim Z(A)=N_{\mathrm{GV}}=\dim Z(S_E)\). The unmarked tangent calculation therefore gives equality with the orbit dimension, so both hypotheses of Theorem~\ref{thm:orbit-smooth-unique} hold.  The dimension formulas then follow from the corresponding orbit dimensions and the power-series description of complete regular local \(\C\)-algebras with residue field \(\C\) \cite[Lemma~10.160.10(2), Tag~\texttt{0C0S}]{Stacks}.  This description uses the given scalar embedding and any lifts of a cotangent-space basis.

For an arbitrary complete local associative family, choose a smooth irreducible open neighbourhood \(U_E\) of \(A\) avoiding all other components.  The orbit \(G_E\cdot S_E\) is dense in \(U_E\).  Every semisimple complex orbit meeting \(U_E\) is open, hence meets that dense orbit and equals it.

This conclusion persists over an algebraically closed extension \(\Omega/\C\).  In fixed dimension \(d\), there are only finitely many multisets of nonzero integral vectors \(v_a\) with \(\sum_a(\sum_i v_{a,i})^2=d\); the Peirce equations select a subset. Each gives a constant complex semisimple marked law.  By Lemma~\ref{lem:marked-semisimple-orbits}, every such law has an open orbit.  All these orbits except the distinguished one are disjoint from \(U_E\), and their open subschemes remain disjoint after base change.  Every semisimple \(\Omega\)-law is marked-isomorphic to the scalar extension of one of these complex laws: its simple modules give exactly the corresponding integral vectors.  The orbit map is smooth and surjective, and its nonempty fibres over the algebraically closed field \(\Omega\) have \(\Omega\)-points.  Thus those base-changed orbits exhaust the semisimple geometric laws, and only the distinguished orbit meets \((U_E)_\Omega\).

By Lemma~\ref{lem:marked-idempotent-lifting}, the given family defines \(\Spec B\to X_E(V)\).  This map factors through \(U_E\), since an open subset of a local spectrum containing its closed point is the whole spectrum.  Every semisimple geometric fibre therefore has the asserted dimension-vector multiset, for any choice of the lifts.
\end{proof}

The smooth open neighbourhood of \(A\) in \(X_d\) is used in Section~\ref{sec:centre-morphism} for arbitrary finite free associative families.  Section~\ref{sec:hua-toda-flatness} concerns the relative quotients arising from geometric deformations.

\begin{remark}
\label{rem:marked-uniqueness-scope}
Theorem~\ref{thm:marked-general-uniqueness} determines the multiset of dimension vectors of the simple modules after a splitting extension.  Over \(\C((t))\), Galois action may permute simple factors with the same full dimension vector, since every \(e_i\) is defined over the ground field.  Different integral families can therefore have isomorphic geometric generic fibres. For the basic algebra \(A\), complete primitive idempotent systems are conjugate after a permutation of their labels.  Hence the multiset of dimension vectors is intrinsic to \(A\) up to such a permutation.  The underlying semisimple algebra retains only the aggregated numbers \(N_m\). The example in the preceding subsection shows why the individual dimension vectors cannot be recovered from those matrix sizes.
\end{remark}

\Needspace{9\baselineskip}
\subsection{Projector ranks in examples}
\label{subsec:marked-checks}

Let \(A\) be the normalized basic contraction algebra of a completed smooth-threefold flop whose surface section is the \(A_3\), \(J=\{1,3\}\) configuration in Subsection~\ref{subsec:multicurve-examples}.  Let \(a,b,c\) be the finite positive orders of \(P_{10},P_{01},P_{11}\) on its central marked classifying arc.  For any complete local finite free associative family with identified special fibre \(A\) and any complete labelled orthogonal lifts of \(e_1,e_3\), Corollary~\ref{cor:marked-smooth-laws} identifies every semisimple geometric fibre with the scalar extension, to its geometric fibre field, of the complex marked algebra
\[
 \C^a\times\C^b\times M_2(\C)^c,
 \qquad e_1=(1^a,0^b,E_{11}^c),\quad
        e_3=(0^a,1^b,E_{22}^c).
\]
The three vectors \((1,0),(0,1),(1,1)\) recover both the matrix sizes and the Cartan matrix computed from the same universal presentations.

Consider the truncated two-cycle algebra with basis \(e_1,e_2,a,b,ab,ba\), in which every path of length three is zero.  It has
\[
 C=\begin{pmatrix}2&1\\1&2\end{pmatrix},\qquad
 b=4,\quad g=6,\quad \dim HH^0=3,\quad
 \dim HH^1=\dim HH^2=1.
\]
The relative cochain spaces in degrees zero through three have dimensions \(4,6,12,24\), and their differentials have ranks \(1,4,7\).  Thus the tangent dimension in \(X_E(V)\) is \(5\), equal to the orbit dimension for the algebra with its fixed idempotents \(\C\times\C\times M_2(\C)\), with
\[
 e_1=(1,0,E_{11}),\qquad e_2=(0,1,E_{22}).
\]
The accompanying exact rational calculation records these ranks and the two Gram decompositions.  The example shows that the relative tangent formula agrees with the orbit dimension for the dimension vectors supplied by the construction.

\section{Consequences for centres and their deformations}
\label{sec:centre-morphism}

Section~\ref{sec:marked-semisimple} proves that the contraction algebra \(A\) defines a smooth point of \(X_d\).  On a smooth open neighbourhood of that point, distinguish two questions.  For a given finite free associative family, how do the centre and cocenter of that family vary?  For a prescribed commutative deformation of the centre, does an associative deformation of the whole algebra realize it?  For the first question, we prove that centres and cocenters form vector bundles compatible with base change. For the second, we compute the centre tangent map and construct a type \(D\) first-order centre deformation whose class lies outside its image.

\begin{conjecture}[{\citet[Conjecture~1.4]{JS}}]
\label{conj:js-centre}
Let \(A\) be a contraction algebra and \(d=\dim A\).
\begin{enumerate}[label=\textup{(\arabic*)}]
\item The algebra \(A\) defines a smooth point of the moduli space of \(d\)-dimensional unital associative algebras; in other words, \(A\) is unobstructed.
\item The dimension of the centre is locally constant near this point, and the map sending an algebra to its centre, regarded as a commutative algebra, is smooth.\footnote{In personal correspondence with the author, Joachim Jelisiejew confirmed that Conjecture~1.4(2) is formulated for the map to commutative algebra structures on the centre (October~2026).}
\end{enumerate}
\end{conjecture}

For the completed smooth-threefold contractions considered here, the smoothness already proved establishes Conjecture~\ref{conj:js-centre} \textup{(1)}.  The bundle theorem gives local constancy and arbitrary base change for centres and cocenters in given associative families.  The centre tangent map gives the first-order lifting criterion for prescribed commutative centre deformations.  Corollary~\ref{cor:centre-morphism-not-smooth} gives counterexamples \(D_{n,2}\), \(n\ge2\), to the smoothness assertion in Conjecture~\ref{conj:js-centre}\textup{(2)}.

The split semisimple fibre has \(N_{\mathrm{GV}}=\sum_\beta n_\beta\) matrix factors. This is also the common rank of the centre and cocenter near \(A\). In the scheme \(\Com_c\) of commutative laws, \(c\) denotes the dimension of the underlying vector space. For the centres studied here, \(c=N_{\mathrm{GV}}\).  In this section, a marking of a centre deformation means an identification of its special fibre with \(Z(A)\); it is distinct from the labelled projectors used in Section~\ref{sec:marked-semisimple}.

Recall the fixed-unit law scheme \(X_d\) and its change-of-basis group \(G\) from Subsection~\ref{sec:orbit-criterion}.  The points of \(X_d\) are unital associative products on a fixed \(d\)-dimensional vector space.  Write \(\mathcal V=V\otimes\mathcal O_{X_d}\) for the tautological algebra and \(\mathcal V_U\) for its restriction to an open subscheme \(U\).

The structure-constant space and its Hochschild tangent description come from associative deformation theory \cite{Gerstenhaber}; the fixed-unit form used here is reviewed in \cite[Section~2.1 and Proposition~2.8]{JS}.  The preceding description of semisimple fibres gives the local bundle statement, from which we construct the centre morphism and its induced map \(\chi_B\) on first-order deformation classes.

\subsection{Unobstructedness}
\label{subsec:unobstructedness}

The smooth local lifting argument for Conjecture~\ref{conj:js-centre} is phrased in the fixed-unit deformation framework of \cite{JS}.

\begin{corollary}
\label{thm:unobstructed}
Let \(A\) be the basic contraction algebra of a completed projective flopping contraction of a smooth complex threefold.  Put \(d=\dim A\) and \(N_{\mathrm{GV}}=\sum_\beta n_\beta\); for one curve this is \(N_{\mathrm{GV}}=\sum_jn_j\).  Then
\begin{equation}
\label{eq:unobstructed-power-series}
 \widehat{\mathcal O}_{X_d,A}\cong\C[[t_1,\ldots,t_{d(d-1)-d+N_{\mathrm{GV}}}]].
\end{equation}
Consequently \(A\) is unobstructed: if \(R'\to R\) is a surjection of local Artinian \(\C\)-algebras with residue field \(\C\), every finite free unital \(R\)-algebra \(\mathscr B\) with \(\mathscr B\otimes_R\C\cong A\) lifts to a finite free unital \(R'\)-algebra \(\mathscr B'\) with \(\mathscr B'\otimes_{R'}R\cong\mathscr B\).  In particular, Conjecture~\ref{conj:js-centre}\textup{(1)} holds for these algebras.
\end{corollary}

\begin{proof}
Corollary~\ref{cor:marked-smooth-laws} identifies the completed local ring at \(A\) with a power-series ring in the stated number of variables; in the one-curve case it already follows from Theorem~\ref{thm:general-uniqueness}.

For the lifting statement, interpret the multiplication on \(\mathscr B\) as a point of \(X_d\) specializing to the smooth point \(A\).  Choose an \(R\)-basis of \(\mathscr B\) whose first vector is \(1\) and which lifts the fixed basis of \(A\); the unit is unimodular because it reduces to the unit of \(A\). The structure constants define a homomorphism \(\mathcal O(X_d)\to R\). It sends the maximal ideal of the point \(A\) into the nilpotent maximal ideal of \(R\), and hence factors through \(\widehat{\mathcal O}_{X_d,A}\).  Lift the images of the variables \(t_i\) arbitrarily to \(R'\).  The resulting homomorphism \(\widehat{\mathcal O}_{X_d,A}\to R'\) gives structure constants over \(R'\) satisfying the associativity and unit equations, hence an algebra \(\mathscr B'\) reducing to \(\mathscr B\).
\end{proof}

\begin{remark}
\label{rem:unobstructed-framed}
Smoothness of \(X_d\) at \(A\) gives the lifting property even when \(HH^3(A)\) is nonzero \cite[Section~1]{JS}.  The number \(d(d-1)-d+N_{\mathrm{GV}}\) is the local dimension at \(A\) of the scheme of algebra laws on a fixed vector space with fixed identity.  The deformation functor of \(A\) up to isomorphism has tangent space \(HH^2(A)\), of dimension \(\dim HH^1(A)\) by Lemma~\ref{lem:tate-hochschild-duality}; the second assertion of Corollary~\ref{thm:unobstructed} says that this functor is unobstructed.
\end{remark}

\subsection{Centres and cocenters in associative families}
\label{subsec:local-centre-structure}

Section~\ref{sec:hua-toda-flatness} proves freeness for relative contraction quotients arising from geometric families.  The result below concerns an arbitrary finite free associative family with special fibre \(A\).  In a smooth open neighbourhood of \(A\) in \(X_d\), the semisimple orbit is dense. The commutator maps used to define the centre and cocenter have the same ranks on that orbit and at \(A\).  Reducedness and density then give constant rank after shrinking the neighbourhood, as the proof below shows.  Hua--Toda's generic matrix decomposition and flatness question concern geometric families \cite[Theorem~4.2 and Conjecture~4.3]{HuaToda}; the smooth neighbourhood in \(X_d\) extends the centre and cocenter conclusions to associative families.

\begin{theorem}
\label{thm:centre-cocenter-bundles}
Let \(A\), \(d\), and \(N_{\mathrm{GV}}\) be as in Corollary~\ref{thm:unobstructed}. There is a smooth irreducible open neighbourhood \(U\) of \(A\) in \(X_d\) with the following properties.
\begin{enumerate}[label=\textup{(\roman*)}]
\item The centre \(\mathsf Z_U=Z(\mathcal V_U)\) and the cocenter \(\mathsf M_U=HH_0(\mathcal V_U/\mathcal O_U)\) are locally free of rank \(N_{\mathrm{GV}}\), and their formation commutes with arbitrary base change: for every morphism \(g\colon U'\to U\), the natural maps \(g^*\mathsf Z_U\to Z(g^*\mathcal V_U)\) and \(g^*\mathsf M_U\to HH_0(g^*\mathcal V_U/\mathcal O_{U'})\) are isomorphisms.
\item Let \(\tau_0\) be a symmetrizing trace of \(A\).  After shrinking \(U\), there is an \(\mathcal O_U\)-linear map \(\tau\colon\mathcal V_U\to\mathcal O_U\) that vanishes on commutators, restricts to \(\tau_0\) at \(A\), and for which \((a,b)\mapsto\tau(ab)\) is a perfect pairing.  Every fibre of \(\mathcal V_U\) is then a symmetric algebra, and
\[
 \mathsf Z_U\xrightarrow{\ \sim\ }\mathsf M_U^\vee,
 \qquad z\longmapsto\bigl(\bar a\mapsto\tau(za)\bigr).
\]
\end{enumerate}
In particular, the dimension of the centre is constant on \(U\), as the first half of Conjecture~\ref{conj:js-centre}\textup{(2)} predicts.
\end{theorem}

\begin{proof}
By Corollary~\ref{cor:marked-smooth-laws}, \(A\) is smooth in \(X_d\), and the unique component of \(X_d\) through \(A\) is the closure of the unmarked orbit of \(S_E\).  Choose a smooth irreducible open neighbourhood \(U\) meeting no other component.  Its semisimple locus is dense, and each semisimple geometric fibre has \(N_{\mathrm{GV}}\) simple factors.

For the centre, consider the commutator operator
\[
 \kappa\colon\mathcal V_U\longrightarrow
 \mathcal E nd_{\mathcal O_U}(\mathcal V_U),
 \qquad a\longmapsto(b\mapsto ab-ba),
\]
a map of free modules of ranks \(d\) and \(d^2\) with kernel \(\mathsf Z_U\).  At a point \(B\), the rank of \(\kappa\) is \(d-\dim Z(B)\). This rank is \(d-N_{\mathrm{GV}}\) both on the semisimple locus, because \(S_E\) has \(N_{\mathrm{GV}}\) simple factors, and at \(A\), by \eqref{eq:marked-numerical-inputs}.  All minors of size \(d-N_{\mathrm{GV}}+1\) consequently vanish on the dense semisimple locus and hence on the reduced scheme \(U\).  Some minor of size \(d-N_{\mathrm{GV}}\) is nonzero at \(A\); shrink \(U\) so that it is a unit.  A matrix with an invertible \(r\times r\) minor and vanishing \((r+1)\)-minors is equivalent to \(\left(\begin{smallmatrix}1_r&0\\0&0\end{smallmatrix}\right)\).  Thus \(\ker\kappa\) is a direct summand of rank \(N_{\mathrm{GV}}\), and \(\ker(g^*\kappa)=g^*\ker\kappa\) for every \(g\).  Since the centre of \(g^*\mathcal V_U\) is \(\ker(g^*\kappa)\), this proves the assertion for \(\mathsf Z_U\).

For the cocenter, use the parallel cokernel calculation with
\[
 b\colon\mathcal V_U\otimes_{\mathcal O_U}\mathcal V_U\longrightarrow
 \mathcal V_U,\qquad a\otimes a'\longmapsto aa'-a'a,
\]
whose cokernel is \(\mathsf M_U\).  The rank of \(b\) at \(B\) is \(d-\dim HH_0(B)\), which equals \(d-N_{\mathrm{GV}}\) on the semisimple locus and, by \eqref{eq:marked-numerical-inputs}, at \(A\).  The same minor calculation as for \(\kappa\) gives a local matrix with an identity block of size \(d-N_{\mathrm{GV}}\) and zero elsewhere.  Its cokernel \(\mathsf M_U\) is therefore locally free of rank \(N_{\mathrm{GV}}\).  Formation of a cokernel commutes with base change by right exactness.

For~(ii), choose an extension of the functional \(\tau_0\) to a local section of \(\mathsf M_U^\vee\), and compose it with \(\mathcal V_U\to\mathsf M_U\).  At the local ring of \(A\), Lemma~\ref{lem:relative-symmetrizing-trace} makes the pairing defined by the extended trace perfect. Shrink \(U\) so that the determinant remains invertible.  The same trace formula gives \(\mathsf Z_U\cong\mathsf M_U^\vee\) throughout \(U\) and after base change.  Every fibre is therefore symmetric.
\end{proof}

\begin{corollary}
\label{cor:multicurve-associative-families}
Let \((B,\mathfrak m,\C)\) be a Noetherian local complex algebra, and let \(\mathscr D\) be a finite free unital \(B\)-algebra with special fibre \(A\).  Then \(Z(\mathscr D)\) and \(HH_0(\mathscr D/B)\) are free of rank \(N_{\mathrm{GV}}\), commute with every commutative scalar extension \(B\to B'\), and every symmetrizing trace of \(A\) lifts to a perfect trace on \(\mathscr D\). The assertion applies to every local associative family specified here, including all complete local bases in the preceding sections.
\end{corollary}

\begin{proof}
Choose a basis beginning with the unit and lifting a fixed basis of \(A\). The structure constants of \(\mathscr D\) in this basis define \(\Spec B\to X_d\).  This map factors through \(U\), because an open subset of a local spectrum containing its closed point is the whole spectrum.  Theorem~\ref{thm:centre-cocenter-bundles} identifies the pulled-back bundles with the centre and cocenter of \(\mathscr D\), including after arbitrary scalar extension. They are finite free because \(B\) is local. The algebra \(\mathscr D\) and its cocenter are now free, and the special fibre \(A\) is symmetric, so Lemma~\ref{lem:relative-symmetrizing-trace} lifts any prescribed symmetrizing trace.
\end{proof}

In particular, these conclusions apply when \(R\) is a local Artinian \(\C\)-algebra with residue field \(\C\) and \(\mathscr B\) is a finite free unital \(R\)-algebra with \(\mathscr B\otimes_R\C\cong A\), including every local Artinian scalar extension \(R\to R''\). The base-change isomorphism \(Z(\mathscr B)\otimes_R\C\cong Z(A)\) shows that every central element of \(A\) lifts to a central element of \(\mathscr B\).

The same neighbourhood also gives equality of the first and second Hochschild cohomology dimensions at every point, as follows.

\begin{corollary}
\label{cor:hh1-hh2-near-A}
For every \(B\in U(\C)\), \(\dim HH^1(B)=\dim HH^2(B)\).
\end{corollary}

\begin{proof}
The scheme \(U\) is smooth and irreducible of dimension \(d(d-1)-d+N_{\mathrm{GV}}\), so \(\dim T_BX_d=d(d-1)-d+N_{\mathrm{GV}}\).  As in the proof of Corollary~\ref{cor:hochschild-orbit-criterion}, \(\dim T_BX_d=d(d-1)-\dim\Der(B)+\dim HH^2(B)\) and \(\dim\Der(B)=d-\dim Z(B)+\dim HH^1(B)\).  Since \(\dim Z(B)=N_{\mathrm{GV}}\) by Theorem~\ref{thm:centre-cocenter-bundles}, the claim follows.
\end{proof}

\begin{remark}
\label{rem:hua-toda-flatness}
The universal presentations in Theorem~\ref{thm:multicurve-hua-toda} give cocenter freeness on geometric slices and the numerical input for marked uniqueness. Marked uniqueness then gives the smooth neighbourhood of the algebra law used here.  The bundle theorem on this neighbourhood extends the centre and cocenter properties to arbitrary associative families.
\end{remark}

\subsection{The centre morphism}
\label{subsec:centre-morphism-setup}

Lifting central elements inside a fixed associative family is weaker than lifting an arbitrarily prescribed multiplication on the centre. To formulate the latter problem, we need a morphism to a scheme of commutative algebra laws.

The centre bundle has constant rank but no preferred basis.  A basis whose first vector is the identity, or a unit-preserving frame, expresses the centre multiplication by structure constants, and hence as a point of a fixed scheme of commutative algebra laws.  The fixed-unit scheme \(X_d\) used above is the structure-constant scheme of \cite[Section~2.1 and Proposition~2.8]{JS}; see also \cite{Gerstenhaber,Gabriel}.  The commutative analogue of \(X_d\) is as follows.

For \(c\ge1\), let \(\Com_c\) be the affine scheme of commutative associative multiplications on \(\C^c\) with identity \(e_1\), and let \(G_c=\Stab_{\GL_c}(e_1)\), a smooth connected group of dimension \(c(c-1)\).  Poonen constructs the scheme of based finite-rank commutative algebras and its fixed-identity locus \cite[Definition~1.5 and Proposition~3.1]{Poonen}; \(\Com_c\) is the corresponding complex fixed-unit law scheme used here.  For a commutative algebra \(Z\), write \(Z^2_{\mathrm{com}}(Z)\) for the symmetric normalized Hochschild \(2\)-cocycles of \(Z\) with values in \(Z\), and \(B^2_{\mathrm{com}}(Z)=\delta\bigl(C^1_{\mathrm{norm}}(Z,Z)\bigr)\), where \((\delta f)(a,b)=af(b)-f(ab)+f(a)b\).  Put
\[
 T^1_{\mathrm{com}}(Z)=Z^2_{\mathrm{com}}(Z)/B^2_{\mathrm{com}}(Z).
\]
This is Harrison's cohomology in cohomological degree two, which governs first-order commutative deformations of \(Z\) \cite[Section~1, p.~192]{Harrison}.  To compare with Harrison's original unnormalized definition, a symmetric cocycle \(f\) satisfies \(f(1,a)=a f(1,1)\); subtracting the coboundary of \(a\mapsto a f(1,1)\) normalizes it.  If a coboundary \(\delta h\) is normalized, then \(h(1)=\delta h(1,1)=0\), so the normalized quotient is unchanged.  The following tangent calculation identifies this commutative deformation space with the quotient of the tangent space to the framed law scheme by infinitesimal changes of frame.

\Needspace{10\baselineskip}
\begin{lemma}
\label{lem:commutative-tangents}
Let \(Z\in\Com_c(\C)\).  Then
\[
 T_Z\Com_c=Z^2_{\mathrm{com}}(Z),
 \qquad
 T_Z(G_c\cdot Z)=B^2_{\mathrm{com}}(Z),
\]
and \(T^1_{\mathrm{com}}(Z)\) is a subspace of \(HH^2(Z)\).  The \(G_c\)-orbits are the isomorphism classes.  The orbit \(G_c\cdot\C^c\) of the split \'etale algebra is open in \(\Com_c\) and has dimension \(c(c-1)\).
\end{lemma}

\begin{proof}
The fixed-unit calculation is the commutative specialization of Lemma~\ref{lem:fixed-unit-tangents}.  A dual-number multiplication \(m+\varepsilon\psi\) remains commutative exactly when \(\psi\) is symmetric.  Since \(Z\) is commutative, every normalized coboundary is symmetric, so \(B^2_{\mathrm{com}}(Z)\) is the space of all normalized coboundaries and \(T^1_{\mathrm{com}}(Z)\subset HH^2(Z)\).  The algebra \(\C^c\) is separable, so \(HH^2(\C^c)=0\) and its tangent space equals the orbit tangent space.  Its stabilizer is the finite group permuting the primitive idempotents.  The proof of Proposition~\ref{prop:semisimple-orbit-open} then shows that the orbit is open, of dimension \(\dim G_c=c(c-1)\).
\end{proof}

We call the closure of \(G_c\cdot\C^c\), with its reduced structure, the \emph{smoothable component} \(\Sm_c\subset\Com_c\).  It is the component reached by degenerating the split \(c\)-point algebra, and it is irreducible of dimension \(c(c-1)\).

A unit-preserving frame places a commutative deformation of the centre in \(\Com_c\).  The following construction gives a local model of the map in Conjecture~\ref{conj:js-centre}\textup{(2)}, relating algebra deformations to commutative deformations of \(Z(A)\) \cite[Section~1]{JS}.  Smoothness of the centre morphism would lift a prescribed commutative centre deformation.  We test smoothness on \(\Com_c\) and on its smoothable component. Corollary~\ref{cor:centre-morphism-not-smooth} gives counterexamples even after restriction to the scheme-theoretic image.

\Needspace{12\baselineskip}
\begin{definition}[The centre morphism]
\label{def:centre-morphism}
Let \(U\subset X_d\) be an open subscheme on which the centre \(\mathsf Z_U=Z(\mathcal V_U)\) is a locally free \(\mathcal O_U\)-module of rank \(c\) whose formation commutes with arbitrary base change.  Let \(\mathscr F_U\to U\) be the scheme of unit-preserving frames: for a \(U\)-scheme \(U'\), the \(U'\)-points of \(\mathscr F_U\) are the isomorphisms \(\mathcal O_{U'}^{\,c}\to\mathsf Z_U\otimes\mathcal O_{U'}\) sending \(e_1\) to \(1\).  The unit section spans a subbundle of \(\mathsf Z_U\), so \(\mathscr F_U\) is a \(G_c\)-torsor that is locally trivial for the Zariski topology.  Transporting the multiplication of \(\mathsf Z_U\) along the tautological frame defines a \(G_c\)-equivariant morphism
\[
 \cent\colon\mathscr F_U\longrightarrow\Com_c .
\]
We say that the \emph{centre morphism is smooth at} \(B\in U(\C)\) if \(\cent\) is smooth at a point of \(\mathscr F_U\) over \(B\).  By equivariance, this does not depend on the choice of the point.
\end{definition}

This smoothness condition is equivalent to smoothness at \(B\) of the morphism \(U\to[\Com_c/G_c]\) classifying \(\mathsf Z_U\): the quotient stack classifies locally free rank-\(c\) commutative algebras \cite[Definition~3.2]{Poonen}, and the frame bundle \(\mathscr F_U\) is its base change along the smooth atlas \(\Com_c\to[\Com_c/G_c]\).  An object of this fibre product is a law on \(\mathcal O^c\) and a unit-preserving algebra isomorphism to the centre; forgetting the law leaves exactly a centre frame. Smoothness and its open locus are detected by this smooth surjective base change \cite[Lemmas~101.33.4 and~101.33.6, Tags~\texttt{0DN7} and~\texttt{0DZR}]{Stacks}.  In a local trivialization of \(\mathsf Z_U\), the unit is unimodular because its image in every residue-field algebra is nonzero.  Its span is therefore a direct summand, which justifies the \(G_c\)-torsor description above.

\subsection{The centre tangent map and the smooth locus}
\label{subsec:centre-tangent-map}

Let \(U\) be as in Theorem~\ref{thm:centre-cocenter-bundles}, so that the centre morphism \(\cent\colon\mathscr F_U\to\Com_c\) of Definition~\ref{def:centre-morphism} is defined.  Fix \(B\in U(\C)\) and put \(Z_B=Z(B)\).  For \(\phi\in Z^2_{\mathrm{norm}}(B,B)\), let \(B_\phi\) be \(B[\varepsilon]\) with multiplication \(aa'+\varepsilon\phi(a,a')\), where \(\varepsilon^2=0\).  By Lemma~\ref{lem:fixed-unit-tangents}, these are the tangent vectors to \(X_d\) at \(B\).  The base-change property of the centre bundle lets us lift elements of \(Z_B\) to central elements of \(B_\phi\); their induced first-order multiplication, modulo changes of the chosen frame, defines \(\chi_B\).

\Needspace{12\baselineskip}
\begin{lemma}
\label{lem:centre-tangent-map}
Let \(\phi\in Z^2_{\mathrm{norm}}(B,B)\).
\begin{enumerate}[label=\textup{(\roman*)}]
\item There is a linear map \(\tau\colon Z_B\to B\) with \(\tau(1)=0\) such that \(z+\varepsilon\tau(z)\) is central in \(B_\phi\) for all \(z\in Z_B\); equivalently,
\begin{equation}
\label{eq:central-lift}
 [\tau(z),a]=\phi(a,z)-\phi(z,a)\qquad(a\in B).
\end{equation}
Any two such maps differ by a linear map \(Z_B\to Z_B\) vanishing at \(1\).
\item For such \(\tau\), the bilinear map
\begin{equation}
\label{eq:induced-centre-cocycle}
 \psi(z,z')=\phi(z,z')+\tau(z)z'+z\tau(z')-\tau(zz')
\end{equation}
takes values in \(Z_B\) and lies in \(Z^2_{\mathrm{com}}(Z_B)\).  The map \(z+\varepsilon z''\mapsto z+\varepsilon(\tau(z)+z'')\) is an isomorphism from \(Z_B[\varepsilon]\), with multiplication \(zz'+\varepsilon\psi(z,z')\), onto \(Z(B_\phi)\).
\item The class of \(\psi\) in \(T^1_{\mathrm{com}}(Z_B)\) depends only on the class of \(\phi\) in \(HH^2(B)\).  This defines a linear map
\[
 \chi_B\colon HH^2(B)\longrightarrow T^1_{\mathrm{com}}(Z_B),
\]
the \emph{centre tangent map} of \(B\).
\end{enumerate}
\end{lemma}

\begin{proof}
(i) The algebra \(B_\phi\) gives a dual-number point of \(X_d\) with closed point in the open subscheme \(U\), so the corresponding morphism factors through \(U\).  By Theorem~\ref{thm:centre-cocenter-bundles}(i), \(Z(B_\phi)\) is free of rank \(c\) over \(\C[\varepsilon]\) and reduces to \(Z_B\).  Every \(z\) thus lifts to a central element \(z+\varepsilon b\); choose such lifts on a basis containing \(1\) to define \(\tau\).  Comparing \((z+\varepsilon\tau(z))a\) with \(a(z+\varepsilon\tau(z))\) in \(B_\phi\) gives~\eqref{eq:central-lift}.  An element \(\varepsilon b\) is central in \(B_\phi\) exactly when \(b\in Z_B\), which proves the last assertion.

(ii) Multiply these chosen central lifts.  The product of \(z+\varepsilon\tau(z)\) and \(z'+\varepsilon\tau(z')\) is
\[
 zz'+\varepsilon\bigl(\phi(z,z')+\tau(z)z'+z\tau(z')\bigr)
 =\bigl(zz'+\varepsilon\tau(zz')\bigr)+\varepsilon\psi(z,z').
\]
It is central, so \(\varepsilon\psi(z,z')\) is central and \(\psi(z,z')\in Z_B\).  Central elements commute, so \(\psi\) is symmetric, and \(\psi(1,z')=0\) because \(\phi\) is normalized and \(\tau(1)=0\).  The map in (ii), from the proposed centre deformation to \(Z(B_\phi)\), is \(\C[\varepsilon]\)-linear and bijective, since both sides are free of rank \(c\) and it reduces to the identity of \(Z_B\).  By the computation above it is multiplicative for the stated product.  That product is therefore associative, so \(\psi\) is a Hochschild cocycle.

(iii) Finally, changing \(\tau\) by \(f\colon Z_B\to Z_B\) with \(f(1)=0\) changes \(\psi\) by \(\delta f\).  If \(\phi=\delta g\) with \(g\in C^1_{\mathrm{norm}}(B,B)\), then \(1+\varepsilon g\) is an algebra isomorphism \(B_\phi\to B[\varepsilon]\).  Hence \(z-\varepsilon g(z)\) is central in \(B_\phi\), and the choice \(\tau=-g|_{Z_B}\) in~\eqref{eq:induced-centre-cocycle} gives \(\psi=0\).

For linearity, choose a fixed vector-space complement to \(Z_B\) in \(B\). Projection to this complement selects the unique complement-valued solution \(\tau\) of~\eqref{eq:central-lift}; since that equation has a right-hand side linear in \(\phi\), this choice and the induced class of \(\psi\) are linear in \(\phi\).
\end{proof}

\begin{proposition}
\label{prop:centre-differential}
Let \((B,\mathfrak f)\) be a \(\C\)-point of \(\mathscr F_U\), and identify \(\C^c\) with \(Z_B\) by the frame \(\mathfrak f\).  The image of \(d\cent\colon T_{(B,\mathfrak f)}\mathscr F_U\to T_{Z_B}\Com_c =Z^2_{\mathrm{com}}(Z_B)\) is spanned by \(B^2_{\mathrm{com}}(Z_B)\) and the cocycles~\eqref{eq:induced-centre-cocycle}.  Consequently
\begin{equation}
\label{eq:centre-differential-rank}
 \operatorname{rank}d\cent_{(B,\mathfrak f)}
 =\dim B^2_{\mathrm{com}}(Z_B)+\operatorname{rank}\chi_B,
\end{equation}
and \(d\cent_{(B,\mathfrak f)}\) is surjective if and only if \(\chi_B\) is surjective.
\end{proposition}

\begin{proof}
A tangent vector to \(\mathscr F_U\) at \((B,\mathfrak f)\) is a first-order deformation \(B_\phi\) together with a unit-preserving frame of \(Z(B_\phi)\) lifting \(\mathfrak f\).  By Lemma~\ref{lem:centre-tangent-map}(i), these frames are \(z\mapsto z+\varepsilon\tau(z)\), where \(\tau\) runs through a torsor under \(C^1_{\mathrm{norm}}(Z_B,Z_B)=\operatorname{Lie}G_c\).  By Lemma~\ref{lem:centre-tangent-map}(ii), \(d\cent\) sends this tangent vector to \(\psi\).  Varying \(\tau\) gives all of \(\psi+B^2_{\mathrm{com}}(Z_B)\), which proves the description of the image.  The image of \(d\cent\) contains \(B^2_{\mathrm{com}}(Z_B)\). Modulo \(B^2_{\mathrm{com}}(Z_B)\), the image of \(d\cent\) is the image of \(\chi_B\) in \(T^1_{\mathrm{com}}(Z_B)\).
\end{proof}

The rank formula separates changes of frame from the induced commutative deformation class. After quotienting the tangent space to \(\Com_c\) at \(Z_B\) by changes of frame, the cokernel of \(\chi_B\) measures which commutative first-order deformation classes cannot be induced by an associative deformation.

To infer smoothness of \(\cent\) from this criterion, its source must be smooth.  The frame bundle is smooth because it is locally a product of the smooth neighbourhood \(U\) and the smooth group \(G_c\).  The following lemma applies even when the target is singular.

\begin{lemma}
\label{lem:smooth-source-criterion}
Let \(g\colon X\to Y\) be a morphism of schemes of finite type over \(\C\), let \(x\in X(\C)\) and \(y=g(x)\), and suppose that \(X\) is smooth at \(x\).  Then \(g\) is smooth at \(x\) if and only if \(dg_x\colon T_xX\to T_yY\) is surjective.  The target \(Y\) may be singular.
\end{lemma}

\begin{proof}
Apply the smooth-source Jacobian criterion of \cite[Section~8.1, Proposition~8.1]{GilleParimalaSuresh} with base \(\Spec\C\).  That statement concerns quasi-separated algebraic stacks of finite presentation and their tangent categories. Here we may restrict to affine neighbourhoods of \(x\) and \(y\), so its finite-presentation and quasi-separatedness hypotheses hold. Both residue fields are \(\C\).  For schemes, the tangent category is the discrete category of maps from \(\Spec\C[\varepsilon]/(\varepsilon^2)\) with the prescribed closed point, so essential surjectivity of the tangent functor is exactly surjectivity of \(dg_x\). The source-smoothness hypothesis is precisely the assumption on \(X\). The cited criterion gives the claimed equivalence, without any reducedness or smoothness assumption on \(Y\); it also implies that \(Y\) is smooth at \(y\) whenever these equivalent conditions hold.
\end{proof}

The smooth-source hypothesis cannot be omitted.  For example, the closed immersion
\[
 \Spec\C[t]/(t^2)\longrightarrow\Spec\C[t]/(t^3)
\]
induces an isomorphism on tangent spaces and is not flat.

\begin{theorem}
\label{thm:centre-smooth-locus}
Let \(U\) be as in Theorem~\ref{thm:centre-cocenter-bundles} and \(B\in U(\C)\).  The centre morphism is smooth at \(B\) if and only if the centre tangent map
\[
 \chi_B\colon HH^2(B)\longrightarrow T^1_{\mathrm{com}}(Z(B))
\]
is surjective.  These points form the \(\C\)-points of a dense open subscheme \(U^{\cent}\subset U\), which contains every semisimple algebra in \(U\).
\end{theorem}

\begin{proof}
The torsor \(\mathscr F_U\) is locally trivial with smooth fibre \(G_c\) over the smooth scheme \(U\), so it is smooth.  Lemma~\ref{lem:smooth-source-criterion} therefore reduces smoothness of \(\cent\) to surjectivity of its differential. Proposition~\ref{prop:centre-differential} shows that changes of frame supply all commutative coboundaries; modulo those directions, the image of the differential is exactly the image of \(\chi_B\). Thus the differential is surjective if and only if \(\chi_B\) is.

The smooth locus of \(\cent\) is open and, by equivariance, \(G_c\)-stable.  The torsor projection is smooth and surjective, hence open, so the image of this locus in \(U\) is open; call it \(U^{\cent}\).  If \(B\) is semisimple, then \(HH^2(B)=0\), and \(Z(B)\cong\C^c\) has \(T^1_{\mathrm{com}}=0\) by Lemma~\ref{lem:commutative-tangents}; so \(\chi_B\) is surjective.  The semisimple locus is dense in \(U\).
\end{proof}

We next construct a type \(D\) centre deformation whose class is nonzero in the cokernel of \(\chi_A\).  The criterion just proved then shows that the centre morphism is not smooth at \(A\).

\subsection{Type \texorpdfstring{\(D\)}{D} counterexamples to Conjecture~1.4(2)}
\label{subsec:type-d-counterexample}

For \(n\ge2\), consider the algebra \(D_{n,2}\) of \cite[Definition~3.2 and Proposition~3.3]{JS},
\begin{equation}
\label{eq:dn2-presentation}
 A_n:=D_{n,2}
 =\frac{\C\langle x,y\rangle}
 {\bigl(xy+yx,\ y^2+x^2+x^{2n-1},\ x^{2n+1}\bigr)}.
\end{equation}
It is the completed Jacobi algebra of the potential \(\tfrac1{2n}x^{2n}+\tfrac13x^3+xy^2\); the ideal displayed in \eqref{eq:dn2-presentation} contains a power of \((x,y)\), so completion changes nothing \cite[Proposition~3.3]{JS}.  By Brown--Wemyss, \(A_n\) is the contraction algebra of a simple flop of a smooth threefold, of type \(D\) and length two \cite[Theorems~1.5 and~8.11]{BrownWemyss}.  It has dimension \(4n\) by Lemma~\ref{lem:dn2-structure}. Lemma~\ref{lem:dn2-split-family} below constructs a split associative family and proves that every geometrically semisimple generic fibre of a finite free formal associative deformation of \(A_n\) has the matrix factors and multiplicities of \(M_2(\C)^{\times(n-1)}\times\C^{4}\). This type is also recorded in \cite[Theorems~3.7 and~3.13]{JS}; our construction and uniqueness argument here are independent.  Theorem~\ref{thm:general-uniqueness} therefore identifies the GV invariants of \(A_n\) as \(n_1=4\) and \(n_2=n-1\), with \(c=n+3\).

We use this presentation to construct a polynomial deformation of \(Z(A_n)\) and show directly that it cannot be the centre of an associative deformation of \(A_n\) with the prescribed special-fibre identification.

Throughout this subsection we write
\[
 \su=x^2,\qquad \sv=x^{2n-1},\qquad \sw=x^{2n-2}y,
\]
using sans-serif letters for elements of the centre.  The symbols \(u\) and \(w\) continue to denote geometric parameters.

The obstruction uses the contrast between \(\sv\) and \(\sw\).  The former is a power of the noncentral element \(x\), whereas the latter lies outside the subalgebra \(\C[x]\).  A commutative centre family makes \(\sv^2\) acquire a \(\sw\)-component.  The calculation below excludes this symbolically for every \(n\ge2\); for \(n=2\), the later socle argument removes the special-fibre marking, and the exact computation describes the entire tangent image.

\begin{lemma}
\label{lem:dn2-structure}
Let \(n\ge2\).
\begin{enumerate}[label=\textup{(\roman*)}]
\item In \(A_n\), \(x^{2n-1}y=0\), and the elements \(1,x,\ldots,x^{2n}\) and \(y,xy,\ldots,x^{2n-2}y\) form a basis.  With \(x^i=0\) for \(i>2n\) and \(x^iy=0\) for \(i>2n-2\),
\[
 x^iy\cdot x^j=(-1)^jx^{i+j}y,\qquad
 x^iy\cdot x^jy=(-1)^{j+1}\bigl(x^{i+j+2}+x^{i+j+2n-1}\bigr).
\]
\item The centralizer of \(x\) is \(\C[x]\oplus\C\sw\), where \(\C[x]=\operatorname{span}(1,x,\ldots,x^{2n})\).
\item The centre is
\[
 Z_n:=Z(A_n)=\operatorname{span}(1,\su,\ldots,\su^n,\sv,\sw)
 \cong\frac{\C[\su,\sv,\sw]}
 {(\su^{n+1},\su\sv,\su\sw,\sv^2,\sv\sw,\sw^2)},
\]
of dimension \(n+3\).
\item \(\sw x=x\sw=0\), and \(\sv A_n+A_n\sv\subseteq\C[x]\).
\end{enumerate}
\end{lemma}

\begin{proof}
(i) Since \(yx^i=(-1)^ix^iy\), multiplying \(y^2=-x^2-x^{2n-1}\) by \(y\) on the left and on the right gives \(-x^2y+x^{2n-1}y=-x^2y-x^{2n-1}y\), so \(x^{2n-1}y=0\).  The listed elements span.  To prove independence directly, take a vector space with basis \(E_0,\ldots,E_{2n},F_0,\ldots,F_{2n-2}\), interpreting larger indices as zero, and define operators
\[
 XE_i=E_{i+1},\quad XF_j=F_{j+1},\quad
 YE_i=(-1)^iF_i,\quad
 YF_j=(-1)^{j+1}(E_{j+2}+E_{j+2n-1}).
\]
On each basis vector, \(XY+YX=0\) and \(X^{2n+1}=0\). Moreover \(Y^2E_i=-E_{i+2}-E_{i+2n-1}\) and \(Y^2F_j=-F_{j+2}\); since \(X^{2n-1}F_j=0\), these identities give \(Y^2+X^2+X^{2n-1}=0\).  Thus \(x\mapsto X\), \(y\mapsto Y\) defines a representation of \(A_n\).  Applying the spanning monomials to \(E_0\) gives the distinct basis vectors \(E_i\) and \(F_j\), so the monomials are independent.  The products follow from \(yx^j=(-1)^jx^jy\) and the relation for \(y^2\).

(ii) Write \(a=f(x)+g(x)y\).  Then \([a,x]=[g(x)y,x]=-2xg(x)y\).  By (i), this vanishes exactly when \(g\) is a multiple of \(x^{2n-2}\).

(iii) By (ii), a central element has the form \(f(x)+\lambda\sw\).  Since \(2n-2\) is even, \(\sw\) commutes with \(y\).  Moreover \([f(x),y]=2\sum_{i\text{ odd}}f_ix^iy\), which vanishes exactly when \(f_i=0\) for all odd \(i\le2n-3\).  This gives the displayed basis.  The relations hold because \(x^{2n+1}=0\), \(x^{2n-1}y=0\) and \(\sw^2=x^{4n-4}y^2=-x^{4n-2}-x^{6n-5}=0\); comparing dimensions gives the presentation.

(iv) This follows from the products in (i).
\end{proof}

\begin{lemma}
\label{lem:dn2-split-family}
For every \(n\ge2\), there is a finite free unital \(\C[[h]]\)-algebra \(\mathscr B_n\), with special fibre \(A_n\), whose generic fibre is
\[
 \mathscr B_n\otimes_{\C[[h]]}\C((h))
 \cong \bigl(M_2(\C)^{\times(n-1)}\times\C^4\bigr)\otimes_\C\C((h)).
\]
Every geometrically semisimple generic fibre of a finite free formal associative deformation of \(A_n\) has these matrix sizes and multiplicities. The scheme \(X_{4n}\) is smooth at \(A_n\), with local dimension \(16n^2-7n+3\).  The unique irreducible component of \(X_{4n}\) through \(A_n\) is the closure of this semisimple orbit.
\end{lemma}

\begin{proof}
Put \(R_h=\C[[h]]\), \(D_h=R_h[[x]]\), and
\[
 \begin{gathered}
 p_h=x\prod_{i=1}^{n-1}(x^2-ih^2),\qquad q_h=x^2-nh^2,\qquad
 g_h=p_h+q_h,\\
 T_h=D_h/(p_hg_h),\qquad M_h=D_h/(p_h).
 \end{gathered}
\]
Let \(\pi:T_h\to M_h\) be the quotient and let \(\sigma\) be the involution of \(M_h\) induced by \(x\mapsto-x\). For \(k,l\in M_h\), the expression \(\mathsf b(k,l)=-\widetilde k\,\sigma(\widetilde l)g_h\) is well defined in \(T_h\), where tildes denote lifts to \(D_h\) and \(\sigma\) acts on \(D_h\) by \(x\mapsto-x\).  Changing either lift changes the expression by a multiple of \(p_hg_h\).  On \(T_h\oplus M_hy\) set
\[
 (f,k)(f',k')
 =\bigl(ff'+\mathsf b(k,k'),\ \pi(f)k'+k\sigma(\pi(f'))\bigr).
\]
The unit is \((1,0)\).

Associativity with at most one \(M_h\)-factor follows from the commuting left and right actions \(\pi,\sigma\pi\). With two such factors it follows from
\[
 \mathsf b(\pi(f)k,l)=f\mathsf b(k,l),\quad
 \mathsf b(k,l\sigma(\pi(f)))=\mathsf b(k,l)f,\quad
 \mathsf b(k\sigma(\pi(f)),l)=\mathsf b(k,\pi(f)l).
\]
With three factors, the two products are \(-k\sigma(l)\pi(g_h)m\) and \(-k\sigma(l)m\sigma(\pi(g_h))\), equal because \(\pi(g_h)=\pi(q_h)\) is fixed by \(\sigma\).

We next compute the ranks over \(R_h\).  The equation \(g_h(hX)=0\) becomes
\[
 X^2-n+h^{2n-3}X\prod_{i=1}^{n-1}(X^2-i)=0.
\]
At \(X=\pm\sqrt n\), each successive coefficient in a formal solution is determined by division by the nonzero number \(\pm2\sqrt n\). The equation therefore has solutions \(\rho_\pm/h\in R_h\), with \(\rho_\pm=\pm\sqrt n\,h+O(h^2)\). Division by \(x-\rho_+\) is valid in \(R_h[[x]]\), since \(\rho_+\in hR_h\) makes evaluation and the quotient converge \(h\)-adically.  The resulting quotient vanishes at \(\rho_-\): its value times \(\rho_--\rho_+\) is zero, and \(R_h\) is a domain. Dividing once more gives \(g_h=(x-\rho_+)(x-\rho_-)v_h\), where \(v_h(0,x)=1+x^{2n-3}\), so \(v_h\) is a unit.

Consequently \(T_h\) is the quotient by the monic distinguished polynomial \(p_h(x-\rho_+)(x-\rho_-)\) of degree \(2n+1\). For a monic distinguished polynomial of degree \(a\), division modulo \(h\) is division by \(x^a\); successive corrections modulo \(h^r\) give a unique remainder of degree below \(a\) over \(R_h\). Thus \(T_h\) and \(M_h\) are free of ranks \(2n+1\) and \(2n-1\). Our algebra \(\mathscr B_n=T_h\oplus M_hy\) is free of rank \(4n\). Its reduction has the basis of Lemma~\ref{lem:dn2-structure} and relations \(xy=-yx\), \(y^2=-x^2-x^{2n-1}\), \(x^{2n+1}=0\); hence its special fibre is \(A_n\).

For each \(1\le i\le n-1\), the matrices
\[
 X_i=h\sqrt i\begin{pmatrix}1&0\\0&-1\end{pmatrix},
 \qquad
 Y_i=h\sqrt{n-i}\begin{pmatrix}0&1\\1&0\end{pmatrix}
\]
and the four scalar pairs
\[
 (0,h\sqrt n),\quad(0,-h\sqrt n),\quad(\rho_+,0),\quad(\rho_-,0)
\]
give representations of \(\mathscr B_n\). Indeed, in every case \(XY=-YX\), \(Y^2=-g_h(X)\), \(p_h(X)Y=0\) and \(p_h(X)g_h(X)=0\); evaluation of the formal series converges since \(X\) is divisible by \(h\).

Over \(\C((h))\), the values of \(x^2\) are \(ih^2,0,\rho_+^2,\rho_-^2\). They are distinct except for the two zero-\(x\) characters. For the last two values, equality would imply \(\rho_-=-\rho_+\), hence \(p_h(\rho_+)=q_h(\rho_+)=0\); but \(q_h(\rho_+)=0\) gives
\[
 p_h(\rho_+)=\rho_+h^{2n-2}\prod_{i=1}^{n-1}(n-i)\ne0.
\]
Interpolation in \(x^2\) separates these factors; \(y\) separates the two zero-\(x\) characters.  On each matrix factor, \(X_i\) gives the two diagonal projectors and \(Y_i\) gives both nonzero off-diagonal entries.  The resulting map onto \(M_2(\C((h)))^{\times(n-1)}\times\C((h))^4\) is surjective. Both sides have dimension \(4n\), so it is an isomorphism.

Finally apply Theorem~\ref{thm:marked-orbit-criterion} with \(E=\C1\).  The constructed family places \(A_n\) in this semisimple orbit closure; the centre of \(A_n\) has dimension \(n+3\) by Lemma~\ref{lem:dn2-structure}, and the semisimple product has \((n-1)+4=n+3\) matrix factors.  The Hochschild equality is Proposition~\ref{prop:multicurve-hochschild-equality}, applicable to the geometric realization of \(A_n\). The criterion proves smoothness, the unique component and uniqueness for any other semisimple orbit closure through \(A_n\). A formal deformation with geometrically semisimple generic fibre places its special point in the orbit closure of the complex semisimple algebra with the same matrix sizes and multiplicities: each defining equation vanishes after the generic field extension, hence in the parameter domain and on its special fibre. This gives the stated uniqueness for every such formal family. The orbit dimension is \((4n)(4n-1)-4n+(n+3)=16n^2-7n+3\).
\end{proof}

\begin{lemma}
\label{lem:dn2-centre-family}
Let
\begin{equation}
\label{eq:dn2-centre-family}
 \mathscr Z_n=\frac{\C[t,\su,\sv,\sw]}
 {\bigl(\su^{n+1},\su\sv,\su\sw,\sv\sw,\sw^2,\ \sv^2-t\sw\bigr)}.
\end{equation}
\begin{enumerate}[label=\textup{(\roman*)}]
\item The algebra \(\mathscr Z_n\) is free over \(\C[t]\) with basis \(1,\su,\ldots,\su^n,\sv,\sw\), and \(\mathscr Z_n/(t)\cong Z_n\).
\item For \(t\ne0\), the fibre is isomorphic to \(\C[\su,\sv]/(\su^{n+1},\su\sv,\sv^3)\), which is smoothable.
\item In the basis of (i), \(\mathscr Z_n\) defines a morphism \(\mathbb A^1\to\Com_{n+3}\) that factors through \(\Sm_{n+3}\).  Its tangent vector at \(t=0\) is the cocycle \(\psi_0\in Z^2_{\mathrm{com}}(Z_n)\) with \(\psi_0(\sv,\sv)=\sw\) and \(\psi_0=0\) on all other pairs of basis vectors.
\end{enumerate}
\end{lemma}

\begin{proof}
(i) Let the powers of \(\su\) multiply as in \(\C[\su]/(\su^{n+1})\), put \(\sv^2=t\sw\), and let all other products of two basis vectors other than \(1\) vanish.  The element \(\sw\) annihilates \(\su,\ldots,\su^n,\sv,\sw\), and \(\sv\) annihilates the powers of \(\su\).  Hence every product of three such basis vectors is zero in either order, unless all three are powers of \(\su\).  The multiplication is therefore associative and commutative.  It satisfies the relations of~\eqref{eq:dn2-centre-family}, and those relations reduce every monomial to the stated basis; this gives inverse isomorphisms.  Setting \(t=0\) recovers Lemma~\ref{lem:dn2-structure}(iii).

(ii) For \(t\ne0\), eliminate \(\sw=t^{-1}\sv^2\).  The family
\[
 \frac{\C[h,\su,\sv]}{\bigl(\su\sv,\ \su^{n+1}-h^n\su,\
 \sv^3-h^2\sv\bigr)}
\]
is free over \(\C[h]\) with basis \(1,\su,\ldots,\su^n,\sv,\sv^2\). Indeed, it maps onto the fibre product over \(\C[h]\) of \(\C[h][\su]/(\su^{n+1}-h^n\su)\) and \(\C[h][\sv]/(\sv^3-h^2\sv)\) along \(\su\mapsto0\) and \(\sv\mapsto0\); that fibre product is free with the images of these elements as a basis, and they span the source.  For \(h\ne0\), its fibre has the \(n+3\) distinct points \((\su,\sv)=(h\xi,0)\) with \(\xi^n=1\), and \((0,0)\), \((0,h)\), \((0,-h)\), so it is isomorphic to \(\C^{n+3}\).  In this basis the family is a morphism \(\mathbb A^1\to\Com_{n+3}\) mapping \(h\ne0\) into the irreducible set \(G_{n+3}\cdot\C^{n+3}\).  Its fibre at \(h=0\), which is the algebra in (ii), therefore lies in the closure \(\Sm_{n+3}\).

(iii) By (ii), the points with \(t\ne0\) map to \(\Sm_{n+3}\), which is closed.  Since \(\mathbb A^1\) is reduced, the whole morphism factors through \(\Sm_{n+3}\).  The structure constants are affine-linear in \(t\), and the only \(t\)-dependent one is that of \(\sv\cdot\sv\).
\end{proof}

\begin{theorem}
\label{thm:dn2-nonlifting}
Let \(n\ge2\), let \(R_\varepsilon=\C[\varepsilon]/(\varepsilon^2)\), and let \(\mathscr Z_{n,\varepsilon}=\mathscr Z_n\otimes_{\C[t]}R_\varepsilon\) with \(t\mapsto\varepsilon\).  There is no finite free \(R_\varepsilon\)-algebra \(\mathscr B\) with \(\mathscr B/\varepsilon\mathscr B\cong A_n\) together with a homomorphism of \(R_\varepsilon\)-algebras \(\iota\colon\mathscr Z_{n,\varepsilon}\to\mathscr B\) whose reduction modulo \(\varepsilon\) is the inclusion \(Z_n\subset A_n\).
\end{theorem}

The theorem does not require the image of \(\iota\) to be central. If \(\mathscr Z_{n,\varepsilon}\) were the centre of an associative deformation with the prescribed marking, its inclusion would supply exactly the homomorphism excluded by the theorem.

\begin{proof}
Suppose that \((\mathscr B,\iota)\) exists.  The argument compares two descriptions of the square of a lift of \(\sv\): the relation on \(x\) forces it into \(\varepsilon\C[x]\), whereas the proposed centre family forces it to equal \(\varepsilon\sw\).

Since \(\mathscr B\) is free over \(R_\varepsilon\), multiplication by \(\varepsilon\) induces an isomorphism of \(A_n\)-bimodules
\begin{equation}
\label{eq:epsilon-bimodule}
 A_n\xrightarrow{\ \sim\ }\varepsilon\mathscr B,\qquad
 a\longmapsto\varepsilon\tilde a,
\end{equation}
where \(\tilde a\) is any lift of \(a\).  We write \(\varepsilon a\) for the image of \(a\).  Choose a lift \(\tilde x\in\mathscr B\) of \(x\), and put \(\tilde{\sv}=\iota(\sv)\) and \(\tilde{\sw}=\iota(\sw)\).

Since \(x^{2n+1}=0\), we have \(\tilde x^{2n+1}=\varepsilon b\) for some \(b\in A_n\).  As \(\tilde x\) commutes with its own powers, \(0=[\tilde x^{2n+1},\tilde x]=\varepsilon[b,x]\), so \([b,x]=0\) by \eqref{eq:epsilon-bimodule}.  Lemma~\ref{lem:dn2-structure}(ii) gives \(b\in\C[x]\oplus\C\sw\).

Now \(\tilde{\sv}\) lifts \(\sv=x^{2n-1}\), so \(\tilde{\sv}=\tilde x^{2n-1}+\varepsilon g\) for some \(g\in A_n\).  Since \(\varepsilon^2=0\) and \(4n-2\ge2n+1\),
\[
 \tilde{\sv}^{\,2}=\tilde x^{2n+1}\tilde x^{2n-3}
   +\varepsilon\bigl(\sv g+g\sv\bigr)
 =\varepsilon\bigl(bx^{2n-3}+\sv g+g\sv\bigr).
\]
Here \(2n-3\ge1\) and \(\sw x=0\), so \(bx^{2n-3}\in\C[x]\), and \(\sv g+g\sv\in\C[x]\) by Lemma~\ref{lem:dn2-structure}(iv).  Hence \(\tilde{\sv}^{\,2}\in\varepsilon\,\C[x]\).  On the other hand, \(\tilde{\sv}^{\,2}=\iota(\sv^2)=\iota(\varepsilon\sw) =\varepsilon\tilde{\sw}=\varepsilon\sw\).  By~\eqref{eq:epsilon-bimodule}, \(\sw\in\C[x]\), contrary to the basis in Lemma~\ref{lem:dn2-structure}(i), in which \(\sw=x^{2n-2}y\) is independent of the powers of \(x\).
\end{proof}

\Needspace{22\baselineskip}
\begin{corollary}
\label{cor:centre-morphism-not-smooth}
Let \(n\ge2\), \(c=n+3\), and let \(U\) be a neighbourhood of \(A_n\) as in Theorem~\ref{thm:centre-cocenter-bundles}.
\begin{enumerate}[label=\textup{(\roman*)}]
\item The class \([\psi_0]\) is not in the image of \(\chi_{A_n}\).  Hence the centre morphism is not smooth at \(A_n\), and \(A_n\notin U^{\cent}\).
\item The centre morphism factors through \(\Sm_c\), which is the scheme-theoretic image of \(\cent\colon\mathscr F_U\to\Com_c\).  The morphism \(\cent\colon\mathscr F_U\to\Sm_c\) is not smooth at the points over \(A_n\) either.
\item The family \(\mathscr Z_n\) does not lift near \(t=0\): there is no finite free algebra over \(\C[t]_{(t)}\), over the convergent power series in \(t\), or over \(\C[[t]]\), with special fibre \(A_n\) and centre identified with \(\mathscr Z_n\) compatibly with \(Z_n\subset A_n\).
\end{enumerate}
Thus \(A_n=D_{n,2}\), \(n\ge2\), gives a counterexample to the smoothness assertion in Conjecture~\ref{conj:js-centre}\textup{(2)}.  The local-constancy assertion holds for \(A_n\) by Theorem~\ref{thm:centre-cocenter-bundles}.
\end{corollary}

\begin{proof}
(i) Suppose that \([\psi_0]=\chi_{A_n}([\phi])\).  By Lemma~\ref{lem:centre-tangent-map}, there is a map \(\tau\) for which \eqref{eq:induced-centre-cocycle} equals \(\psi_0\), after changing \(\tau\) by an element of \(C^1_{\mathrm{norm}}(Z_n,Z_n)\).  The map \(z\mapsto z+\varepsilon\tau(z)\) is then an algebra homomorphism \(\mathscr Z_{n,\varepsilon}\to(A_n)_\phi\) lifting \(Z_n\subset A_n\), contradicting Theorem~\ref{thm:dn2-nonlifting}.  Thus \(\chi_{A_n}\) is not surjective. Theorem~\ref{thm:centre-smooth-locus} equates this surjectivity with smoothness of the centre morphism, giving the conclusion.

(ii) Over the dense semisimple locus of \(U\), the centres are isomorphic to \(\C^c\), so \(\cent\) maps a dense open subscheme of the reduced scheme \(\mathscr F_U\) into \(G_c\cdot\C^c\).  Hence \(\cent\) factors through \(\Sm_c\).  By equivariance its image contains the open orbit.  Moreover, the kernel of \(\mathcal O(\Com_c)\to\Gamma(\mathscr F_U,\mathcal O)\) is radical, because \(\mathscr F_U\) is reduced.  Thus the scheme-theoretic image is exactly the reduced orbit closure \(\Sm_c\).

Let \(\mathfrak f\) be the frame of \(Z_n\) given by the basis of Lemma~\ref{lem:dn2-centre-family}(i).  By part (iii) of that lemma, \(\psi_0\) is tangent to \(\Sm_c\) at \(\cent(A_n,\mathfrak f)\).  By (i) and Proposition~\ref{prop:centre-differential}, it is not in the image of \(d\cent_{(A_n,\mathfrak f)}\).  So the differential of \(\cent\colon\mathscr F_U\to\Sm_c\) at \((A_n,\mathfrak f)\) is not surjective, and \(\cent\) is not smooth there by Lemma~\ref{lem:smooth-source-criterion}; by equivariance, the same holds at every point over \(A_n\).

(iii) Reducing such a lift modulo \(t^2\) would give a pair \((\mathscr B,\iota)\) as excluded by Theorem~\ref{thm:dn2-nonlifting}.
\end{proof}

\begin{remark}
\label{rem:dn2-mechanism}
The socle element \(\sv=x^{2n-1}\) of \(Z_n\) is a power of the noncentral element \(x\).  In every first-order deformation, the square of a lift of \(\sv\) stays in the subalgebra generated by a chosen lift of \(x\), because of the relation \(x^{2n+1}=0\) and the shape of the centralizer of \(x\).  The commutative centre deformation \(\sv^2=t\sw\) instead gives the square a component along \(\sw=x^{2n-2}y\), outside that subalgebra.

The obstruction is first order.  A ramified base change \(t=h^r\) with \(r\ge2\) has zero first-order class, so lifts of that pullback are a higher-order problem.
\end{remark}

\subsection{Exact tangent data for \texorpdfstring{\(D_{2,2}\)}{D2,2}}
\label{subsec:d22-tangent-data}

The preceding symbolic obstruction applies to every \(n\ge2\).  For \(n=2\), exact linear algebra computes the full first-order image of \(\chi_A\).  This computation also supplies an independent certificate for non-lifting without a prescribed identification of the special-fibre centre.  Here \(A=A_2\) has basis \(1,x,x^2,x^3,x^4,y,xy,x^2y\), and \(Z_2\) has basis \(1,\su,\sv,\sw,\sq\) with \(\su=x^2\), \(\sv=x^3\), \(\sw=x^2y\) and \(\sq=\su^2=x^4\).  The only nonzero product of two of \(\su,\sv,\sw,\sq\) is \(\su^2=\sq\).

\Needspace{23\baselineskip}
\begin{proposition}
\label{prop:d22-tangent-data}
For \(A=D_{2,2}\) and \(Z_2=Z(A)\), the following hold.
\begin{enumerate}[label=\textup{(\roman*)}]
\item \(\dim\Der(A)=6\), \(\dim HH^1(A)=\dim HH^2(A)=3\), and \(\dim Z^2_{\mathrm{norm}}(A,A)=53=\dim_AX_8\).
\item \(\dim B^2_{\mathrm{com}}(Z_2)=10\), \(\dim Z^2_{\mathrm{com}}(Z_2)=26\), and \(\dim T^1_{\mathrm{com}}(Z_2)=16\).
\item The centre tangent map \(\chi_A\) has rank \(3\).  Its image is spanned by three explicitly given cocycle classes \([\psi_1],[\psi_2],[\psi_3]\), recorded in Appendix~\ref{app:d22-representatives}.  The representatives satisfy \(\psi_i(s,s)\in\C\sq\) for every \(s\in\Soc(Z_2)\).
\item At every point of \(\mathscr F_U\) over \(A\), \(\operatorname{rank}d\cent=10+3=13\), whereas \(\dim T_{Z_2}\Sm_5\ge\dim\Sm_5=20\) and \(\dim T_{Z_2}\Com_5=26\).
\end{enumerate}
\end{proposition}

\begin{proof}
The cochain dimensions and ranks in (i)--(iii) are exact linear algebra over \(\Q\), described in Appendix~\ref{app:d22-computation}. The equality with the local dimension in (i) uses Lemma~\ref{lem:dn2-split-family} with \(n=2\): its split family places \(A\) in the orbit closure of \(M_2(\C)\times\C^4\), of dimension \(8\cdot7-8+5=53\).  Hence \(53\le\dim_AX_8\le\dim T_AX_8=53\), so \(A\) is smooth and that orbit closure is the unique component through \(A\). For (iv), Proposition~\ref{prop:centre-differential} adds the dimension of the frame directions in (ii) to the rank of \(\chi_A\) in (iii). Lemma~\ref{lem:commutative-tangents} gives the dimension of the open split-\'etale orbit, hence of its closure \(\Sm_5\). The tangent space at a closed point has dimension at least this local dimension; the ambient tangent dimension is the cocycle dimension in (ii).
\end{proof}

At \(D_{2,2}\), the image of the differential has codimension at least seven in the tangent space to the scheme-theoretic image of the centre morphism. Theorem~\ref{thm:dn2-nonlifting} fixes the special-fibre identification \(\mathscr Z_{n,\varepsilon}/(\varepsilon)\cong Z_n\).  For \(n=2\), we remove that identification by proving a socle condition that is invariant under every automorphism of \(Z_2\).

\begin{proposition}
\label{prop:d22-unmarked}
No first-order deformation of \(A=D_{2,2}\) over \(R_\varepsilon\) has centre isomorphic to \(\mathscr Z_{2,\varepsilon}\) as an \(R_\varepsilon\)-algebra.  Equivalently, the \(\Aut(Z_2)\)-orbit of \([\psi_0]\) does not meet the image of \(\chi_A\).
\end{proposition}

\begin{proof}
We first prove the socle constraint directly for any finite free commutative \(R_\varepsilon\)-algebra \(\mathscr D\) with special fibre \(Z_2\) admitting a homomorphism to a finite free associative deformation \(\mathscr B\) of \(A\) that reduces to \(Z_2\subset A\). Write the multiplication of \(\mathscr D\), in a unit-preserving lifted basis, as \(ab+\varepsilon\psi(a,b)\).  Put \(\mathfrak m=(\su,\sv,\sw,\sq)\).  Then \(\mathfrak m^2=\C\sq\) and \(\Soc(Z_2)=\C\sv\oplus\C\sw\oplus\C\sq\). We will show
\begin{equation}
\label{eq:socle-diagonal}
 \psi(s,s)\in\C s+\mathfrak m^2
 \qquad\text{for every }s\in\Soc(Z_2).
\end{equation}

Associativity on \((\su,s,s)\) gives \(\su\psi(s,s)=s\psi(\su,s)\). If the \(\sv\)- or \(\sw\)-coordinate of \(s\) is nonzero, the left side lies in \(\C\su+\C\sq\), while the right side is a scalar multiple of \(s\).  Comparing coefficients forces the \(1\)- and \(\su\)-coordinates of \(\psi(s,s)\) to vanish.  These coordinates are polynomials in \(s\) and vanish on the dense open set just considered, so \(\psi(s,s)\in\Soc(Z_2)\) for every socle element.

For \(a,b,c\in\C\) with \(d=a^2-b^2\ne0\), put \(z=ax+by\) and \(s=a\sv+b\sw+c\sq\). Lemma~\ref{lem:dn2-structure} gives
\[
 z^2=d\su-b^2\sv,\qquad
 z^3=d(a\sv+b\sw)-ab^2\sq,\qquad
 z^4=d^2\sq,\qquad z^5=0.
\]
For \(r\in A\), let \(r_x,r_y,r_{xy}\) be its coefficients on \(x,y,xy\); the other basis elements are central.  Then
\[
 [r,z]=2(br_x-ar_y)xy-2r_{xy}(a\sw+b\sv+b\sq).
\]
It follows that
\[
 C_A(z)=\C[z]\oplus\C(b\sv+a\sw),\qquad
 (b\sv+a\sw)z=0,\qquad
 \C[z]\cap\Soc(Z_2)=\C s+\C\sq.
\]
Indeed, the centralizer is \(Z_2+\C z\), and the displayed powers of \(z\), together with \(b\sv+a\sw\), are independent because \(d\ne0\).  Also
\[
 s=P(z),\qquad
 P(T)=d^{-1}T^3+(c+ab^2/d)d^{-2}T^4.
\]

Choose a lift \(X\) of \(z\) in \(\mathscr B\).  The isomorphism \eqref{eq:epsilon-bimodule} writes \(X^5=\varepsilon r\) with \(r\in C_A(z)\), since \(X\) commutes with its powers.  For every \(j\ge6\), therefore, \(X^j=\varepsilon r z^{j-5}\) lies in \(\varepsilon\C[z]\).  If \(S\) is the image in \(\mathscr B\) of the chosen lift of \(s\) in \(\mathscr D\), then \(S=P(X)+\varepsilon g\) for some \(g\in A\).  The products of \(s\) with \(1,x,y\) are \(s,a\sq,-b\sq\), respectively, and its products with every other basis element are zero.  Since \(s\) is central,
\[
 S^2=P(X)^2+\varepsilon(sg+gs)\in\varepsilon\C[z].
\]
On the other hand \(S^2=\varepsilon\psi(s,s)\), because \(s^2=0\). Thus \(\psi(s,s)\in\C[z]\cap\Soc(Z_2)=\C s+\C\sq\) when \(d\ne0\).

To include the remaining socle directions, write
\[
 \psi(a\sv+b\sw+c\sq,a\sv+b\sw+c\sq)
   =F(a,b,c)\sv+G(a,b,c)\sw+H(a,b,c)\sq .
\]
The quadratic polynomials \(F,G\) satisfy \(bF-aG=0\) on \(a^2-b^2\ne0\), hence identically.  Coprimality of \(a,b\) in \(\C[a,b,c]\) gives \(F=aL\) and \(G=bL\) for a polynomial \(L\). The displayed diagonal is \(Ls+(H-cL)\sq\), proving \eqref{eq:socle-diagonal} also when \(s\in\C\sq\).

This condition depends only on the cohomology class: for a normalized linear map \(f:Z_2\to Z_2\), \((\delta f)(s,s)=2sf(s)\in\C s\). It is invariant under \(\Aut(Z_2)\), which preserves \(\Soc(Z_2)\) and \(\mathfrak m^2\). The class \([\psi_0]\) violates it, since \(\psi_0(\sv,\sv)=\sw\notin\C\sv+\C\sq\). The centre tangent map is constructed from precisely such lifted multiplications, so its image satisfies the condition; no automorphic transport of \([\psi_0]\) belongs to that image.

Finally suppose \(Z(\mathscr B)\cong\mathscr Z_{2,\varepsilon}\) as an \(R_\varepsilon\)-algebra.  Then \(Z(\mathscr B)\) is free of rank five.  The natural map \(Z(\mathscr B)/\varepsilon Z(\mathscr B)\to Z_2\) is injective: if a central element lies in \(\varepsilon\mathscr B\), it is annihilated by \(\varepsilon\), and freeness of \(Z(\mathscr B)\) puts it in \(\varepsilon Z(\mathscr B)\). Dimension five makes the reduced map an isomorphism. The putative centre isomorphism consequently transports \([\psi_0]\) by an automorphism of \(Z_2\), up to a coboundary, to a class satisfying \eqref{eq:socle-diagonal}.  The invariance and failure just proved give a contradiction.
\end{proof}

\paragraph{Centres in families and the scheme of algebra laws.}
\label{subsec:conjecture-summary}
The smoothness proved in Section~\ref{sec:marked-semisimple} gives an open neighbourhood of \(A\) in \(X_d\) on which centre and cocenter have constant rank.  Corollary~\ref{thm:unobstructed} lifts associative deformations, and Theorem~\ref{thm:centre-cocenter-bundles} describes their centres and cocenters.  A prescribed commutative centre deformation has the additional first-order condition that its class lie in the image of \(\chi\). Corollary~\ref{cor:centre-morphism-not-smooth} gives counterexamples \(D_{n,2}\), \(n\ge2\), to the smoothness assertion in Conjecture~\ref{conj:js-centre}\textup{(2)}, even when the target is restricted to the smoothable component or to the scheme-theoretic image.

The type \(A\) and \(D\) algebras studied in \cite{JS} belong to this geometric class \cite[Theorem~1.5]{BrownWemyss}.  The type \(D\) examples show that base-change compatibility of centres in associative families does not guarantee that a prescribed commutative centre deformation lifts.

\section*{Acknowledgements}
\addcontentsline{toc}{section}{Acknowledgements}
I made extensive use of large language models, including ChatGPT, and the mathematical exploration agents Danus and Rethlas throughout this work. These tools assisted with mathematical exploration and exposition, suggested proof strategies and intermediate arguments, and helped check and reorganize calculations. I formulated the research questions, developed the strategies, and selected and revised the arguments presented here. I have made every effort to verify these arguments against the cited literature and accompanying exact certificates, and I take full responsibility for the manuscript and for making it public.

Neither the existence nor the priority date of this thesis should discourage a wholly independent human proof of these results. If a subsequent work obtains them through an argument developed wholly independently by humans, I claim no special priority or credit for the proof ideas beyond standard bibliographic acknowledgement of prior or concurrent results.

If later work makes substantive use of the results, constructions, or proof strategies in this thesis, ordinary scholarly citation is appropriate regardless of the tools used. In particular, I would appreciate explicit attribution when closely related arguments appear in AI-assisted work, since the origin of ideas generated through related language models can otherwise be difficult to trace.

\appendix
\Needspace{12\baselineskip}
\section{Root enumeration for the primitive marked pairs}
\label{app:general-certificates}

The root enumeration verifies the assertion used in Lemma~\ref{lem:coefficient-orbits} that the positive roots with a fixed nonzero coefficient at the retained vertex form one orbit under the Levi Weyl group.  Karmazyn's primitive presentations determine the six marked pairs \cite[Examples~4.22--4.27]{Karmazyn}.  For the enumeration, number the vertices in the diagrams of Figure~\ref{fig:marked-pairs} from left to right along the horizontal chain, and number the lower vertex last.  Thus the retained vertices are \(1,2,3,3,4,3\) for lengths \(1,\ldots,6\), respectively. The table in the proof of Lemma~\ref{lem:closed-fibre-cocenter-count} instead uses Karmazyn's numbering, so the two conventions must not be mixed.

Represent a root by its vector \(v=(v_1,\ldots,v_r)\) in the simple-root basis.  If \(i\sim h\) denotes adjacency in the diagram, the simple reflection acts by
\[
 s_i(v)_h=v_h\quad(h\ne i),\qquad
 s_i(v)_i=-v_i+\sum_{h\sim i}v_h.
\]
To enumerate the full root system, start from the simple roots. Apply every simple reflection \(s_i\) to each new vector, retaining each vector only once.  Since the ADE root system is finite, this terminates and gives the full Weyl orbits of the simple roots, hence the entire root system.  Keeping the vectors with nonnegative coordinates produces \(\Phi^+\).

Let \(p\) be the retained vertex.  For each retained coefficient \(j>0\), choose a positive root with \(v_p=j\) and repeat the procedure using only reflections with \(i\ne p\).  The resulting orbit agrees, in all six primitive cases, with the complete set of positive roots of coefficient \(j\).  The numbers of positive roots are
\[
 1,\quad12,\quad36,\quad63,\quad120,\quad120,
\]
and the non-Levi totals are
\[
 N_p=1,\quad9,\quad29,\quad53,\quad104,\quad106.
\]
The individual orbit cardinalities are those of Lemma~\ref{lem:coefficient-orbits}.  The same calculation checks that the largest marked coefficient is the stated length.

The integer implementation is included in \path{tools/check_mathematics.py} in the source distribution.  At every step it verifies
\[
 2\sum_i v_i^2-2\sum_{\{i,h\}\text{ an edge}}v_iv_h=2.
\]
Integer vector operations make the enumeration exact.  Distinct coefficient orbits give disjoint sets of hyperplanes, because proportional roots in a reduced root system differ only by sign.  The curve divisors descended from these disjoint sets of hyperplanes are therefore distinct.

\Needspace{16\baselineskip}
\section{The \texorpdfstring{\(E_8\)}{E8} calculation for length five}
\label{app:e8-executable-certificate}

\subsection{Relations for the central extension}
\label{app:e8-specialization}

We compute the dimension of the cocenter of the closed corner for the primitive \(E_8\) marking of length five. Subsection~\ref{app:e8-specialization} derives the homogeneous relations, and Subsection~\ref{app:e8-ranks} records the relation matrices and their modular ranks. Subsection~\ref{subsec:e8-certificate} combines these bounds with the tail estimate and the generic-rank lower bound to prove cocenter dimension \(104\).  The central extension is that of \cite[Theorem~3.2]{EtingofRains}, and the form of the relations at a nodal vertex comes from \cite[Section~3]{EtingofLatourRains}.  The specialization at the vertex of coefficient five, relation reduction, and exact row-rank certificate are carried out here.

We use the normalization \(\rho_i=(\rho,\alpha_i)=1\) from Subsection~\ref{subsec:universal-marked-cocenter}.  The algebra \(\Pi^\rho_0\) is the path algebra over \(\C[z]\) of the double of a finite \(E_8\) quiver \(Q\), with \(z\) central, modulo
\[
 \sum_{a\in Q_1}[a,a^*]=z\sum_ie_i .
\]
Different orientations of \(Q\) give isomorphic algebras \cite[Section~2.1]{EtingofLatourRains}.

We compose paths from left to right: if \(a:u\to v\) and \(b:v\to w\), then \(ab:u\to w\).  For an edge joining \(v\) and \(v'\), let \(\ell_{v,v'}\) denote the summand of the left-hand side that is a loop at \(v\) through \(v'\). It is \(aa^*\) or \(-a^*a\), according to the orientation of the arrow \(a\) of \(Q\) joining them.  The relation at \(v\) reads \(\sum_{v'\sim v}\ell_{v,v'}=z\,e_v\).  If \(p\colon v\to v'\) and \(q\colon v'\to v\) are the two arrows of the double quiver, then \(\ell_{v,v'}=\pm pq\) and \(\ell_{v',v}=\mp qp\). Hence
\begin{equation}
\label{eq:arm-commutation}
 \ell_{v,v'}\,p=-p\,\ell_{v',v},\qquad q\,\ell_{v,v'}=-\ell_{v',v}\,q .
\end{equation}

\begin{lemma}
\label{lem:arm-relations}
Let \(v_0,v_1,\ldots,v_r\) be a chain in the Dynkin graph such that \(v_1,\ldots,v_r\) have no neighbours other than those in the chain.  Then \(X=\ell_{v_0,v_1}\) satisfies
\[
 X(X+z)(X+2z)\cdots(X+rz)=0 .
\]
\end{lemma}

\begin{proof}
Put \(Y=\ell_{v_1,v_0}\).  By~\eqref{eq:arm-commutation}, \(qf(X)=f(-Y)q\) for every polynomial \(f\) over \(\C[z]\). Hence \(Xf(X)=\pm pf(-Y)q\).  It is therefore enough to find a degree-\(r\) polynomial \(f\) with \(f(-Y)=0\).  

If \(r=1\), the relation at \(v_1\) is \(Y=z\), and \(f(T)=T+z\) has the required property.  If \(r>1\), the relation at \(v_1\) gives \(\ell_{v_1,v_2}=z-Y\), and the induction hypothesis gives \(\prod_{m=0}^{r-1}(\ell_{v_1,v_2}+mz)=0\).  Substitution of \(z-Y\) for \(\ell_{v_1,v_2}\) gives \(\prod_{m=1}^{r}(mz-Y)=0\), and \(f(T)=\prod_{m=1}^r(T+mz)\) works.
\end{proof}

The nodal corners of Etingof--Latour--Rains are presented by relations of the same form \cite[equation~(3.1)]{EtingofLatourRains}.  We number the vertices of \(E_8\) as in \cite[Example~4.26]{Karmazyn}: the finite diagram is the chain \(1,2,3,4,8,7,6\) with \(5\) attached to the branching vertex \(8\), and the retained vertex is \(p=4\).  Orient \(Q\) so that its arrow \(a\) between \(4\) and \(8\) points from \(4\) to \(8\), and put
\[
 d=-\ell_{4,3},\qquad b=-\ell_{8,5},\qquad c=-\ell_{8,7}.
\]
Apply Lemma~\ref{lem:arm-relations} to the three arms, of lengths \(3\), \(1\), and \(2\), and then use the relations at vertices \(4\) and \(8\). This gives
\begin{equation}
\label{eq:e8-arm-relations}
\begin{gathered}
 d(d-z)(d-2z)(d-3z)=0,\qquad b(b-z)=0,\qquad c(c-z)(c-2z)=0,\\
 aa^*=d+z\,e_4,\qquad a^*a+b+c+z\,e_8=0.
\end{gathered}
\end{equation}
These are the relations of Karmazyn's intermediate presentation \cite[Example~4.26]{Karmazyn} at \(t_i=z\), after removal of the free vertex.

To compute the corner at \(\lambda=0\), put
\[
 B=b-2z\,e_8,\qquad C=c,\qquad H=a^*a=-B-C-3z\,e_8,
\]
so that
\begin{equation}
\label{eq:e8-B-C-relations}
 B^2=-3zB-2z^2,\qquad C^3=3zC^2-2z^2C .
\end{equation}
In the corner \(B_5=e_4\Pi^\rho_0e_4\), define
\[
 \beta=aa^*+3z\,e_4=d+4z\,e_4,\qquad
 \gamma=aBa^*-2z^2e_4 .
\]

\begin{lemma}
\label{lem:e8-generators}
The elements \(z,\beta,\gamma\) generate the algebra \(B_5\).
\end{lemma}

\begin{proof}
To prove that \(z,\beta,\gamma\) generate \(B_5\), we reduce closed paths at the two vertices \(4\) and \(8\). A closed path at \(4\) is a product of paths meeting \(4\) only at their endpoints. Those entering the arm \(\{1,2,3\}\) are polynomials in \(d\) and \(z\), by the arm relations and~\eqref{eq:arm-commutation}. The remaining paths have the form \(aqa^*\) with \(q\) closed at \(8\).

A closed path at \(8\) is a product of polynomials in \(b,z\), polynomials in \(c,z\), and elements \(a^*ra\) with \(r\) closed at \(4\). Since \(da=a(H-z)\) by~\eqref{eq:e8-arm-relations}, we have \(a^*d^ka=H(H-z)^k\). Induction on path length now shows that every closed path at \(8\) is a polynomial in \(B\), \(C\), and \(z\).

We now reduce the remaining paths based at the retained vertex. Substitute \(C=-a^*a-B-3z\) in an expression \(aq(B,C)a^*\). Each occurrence of \(a^*a\) splits the expression into a product of two elements of the same form, while powers of \(B\) reduce by~\eqref{eq:e8-B-C-relations}. Thus \(B_5\) is generated by \(z\), \(aa^*\), and \(aBa^*\), equivalently by \(z\), \(\beta\), and \(\gamma\).
\end{proof}

\begin{lemma}
\label{lem:e8-relations}
In \(B_5\), the elements
\begin{align*}
 \mathsf P&=(\beta-4z)(\beta-5z)(\beta-6z)(\beta-7z),\\
 \mathsf R_4&=(\gamma+\beta^2)^2+\beta\gamma\beta-5z^2(\gamma+\beta^2)+4z^4,\\
 \mathsf R_5&=\gamma\beta\gamma+\gamma^2\beta+\gamma\beta^3-5z^2\gamma\beta
\end{align*}
vanish.
\end{lemma}

\begin{proof}
Since \(\beta=d+4z\), the element \(\mathsf P\) is \(d(d-z)(d-2z)(d-3z)=0\).

For the remaining two relations, multiply elements of the form \(s\,e_4+aPa^*\), with \(s\in\C[z]\) and \(P\) a noncommutative polynomial in \(B\) and \(C\), using the rule
\[
 (s\,e_4+aPa^*)(s'e_4+aP'a^*)=ss'e_4+a(sP'+s'P+PHP')a^* .
\]
Only~\eqref{eq:e8-B-C-relations} and \(H=-B-C-3z\) enter the reduction. They give \(HB=2z^2-CB\), \(BH=2z^2-BC\) and
\[
 \gamma+\beta^2=7z^2+a(3z-C)a^*,\qquad
 \beta\gamma=-6z^3+a(3z-C)Ba^*,\qquad
 \gamma\beta=-6z^3+aB(3z-C)a^* .
\]

For \(\mathsf R_5\), the same reduction rules give \(\beta\gamma+(\gamma+\beta^2)\beta-5z^2\beta=a(C-z)(C-2z)a^*\).  Hence
\[
 \mathsf R_5
 =\gamma\bigl(\beta\gamma+(\gamma+\beta^2)\beta-5z^2\beta\bigr)
 =a(BH-2z^2)(C-z)(C-2z)a^*
 =-aBC(C-z)(C-2z)a^*,
\]
which vanishes by~\eqref{eq:e8-B-C-relations}.

For \(\mathsf R_4\), put \(g=\gamma+\beta^2\), so that \(\mathsf R_4=(g-z^2)(g-4z^2)+\beta\gamma\beta\).  The same rules give
\begin{align*}
 (g-z^2)(g-4z^2)
 &=18z^4+a\bigl(27z^3-9z^2C+(3z-C)H(3z-C)\bigr)a^*,\\
 (3z-C)H(3z-C)
 &=-27z^3-9z^2B+11z^2C+3zBC+3zCB-CBC,\\
 \beta\gamma\beta
 &=-18z^4+a\bigl(9z^2B-2z^2C-3zBC-3zCB+CBC\bigr)a^*,
\end{align*}
where the second line uses \(C^3=3zC^2-2z^2C\).  Adding the reduced terms gives zero.
\end{proof}

The program \path{tools/check_e8_rho_relations_formal.py} independently verifies these reductions over the polynomial ring \(\mathbb Z[z]\), with \(z\) retained as a variable.

\begin{remark}[Sign conventions in Karmazyn's arXiv version]
\label{rem:karmazyn-sign-conventions}
The comparison in this remark concerns \href{https://arxiv.org/abs/1709.02720v1}{arXiv:1709.02720v1} (8 September 2017), Example~4.26, pp.~23--25 of \cite{Karmazyn}. The invariant generators are listed on p.~23, the intermediate presentation is on pp.~24--25, and the final display is on p.~25. The relations used here are derived independently in Lemma~\ref{lem:e8-relations}. The rank matrices below use the relations \(\mathsf P,\mathsf R_4,\mathsf R_5\) of Lemma~\ref{lem:e8-relations}.  These agree with the last two relations of Karmazyn's intermediate presentation, written in the invariants \(R_1,R_2,R_3,S_2\) of the five-arm parameters \(\tau_1,\ldots,\tau_5\).  At the specialization \(\rho_i=1\), \(\lambda=0\), \((\tau_1,\ldots,\tau_5)=(2,1,0,-1,-2)z\) and \((R_1,R_2,R_3,S_2)=(3z,2z^2,0,2z^2)\); here \(\beta=a_4a_4^*+R_1\) and \(\gamma=a_4Ba_4^*-S_2\).  For the quadratic and cubic relations displayed there, the shifts must be \(B=b+\tau_5\) and \(C=c+\tau_3\).

The final display on p.~25 of that version rewrites these relations through \(T_3=R_2+S_2+R_1S_1\), \(T_2=R_3+R_2S_1+R_1S_2\), \(T_1=-R_3R_1+R_2S_2\) and \(T_0=R_3S_2\).  The expressions \(T_3,T_2,T_1,T_0\) are \(+\sigma_2,\ldots,+\sigma_5\) in the parameters \(\tau_i\), whereas Karmazyn's list of invariant generators sets \(T_i=-\sigma(\tau)\).  Substituting the sign convention \(T_i=-\sigma(\tau)\) into Karmazyn's final display reverses the signs of the parameter terms and gives
\[
 \mathsf R_4'=(\gamma+\beta^2)^2+\beta\gamma\beta+5z^2(\gamma+\beta^2)-4z^4,
 \qquad
 \mathsf R_5'=\gamma\beta\gamma+\gamma^2\beta+\gamma\beta^3+5z^2\gamma\beta .
\]

Neither expression vanishes in \(B_5\). We prove nonvanishing in the cocenter by determining the relation spaces in degrees four and five without using the later rank table. Put \(\mathcal F=\C[z]\langle\beta,\gamma\rangle\). The matrix Hilbert polynomial in \cite[Theorem~3.2(iv)]{EtingofRains} gives
\[
 (\dim(B_5)_m)_{m=0}^{5}=(1,2,4,7,10,14),\qquad
 (\dim\mathcal F_m)_{m=0}^{5}=(1,2,4,7,12,20).
\]
For the first calculation, let \(C_Q\) be the adjacency matrix of the numbered \(E_8\) tree. The diagonal coefficients at vertex \(4\) of \((1-C_Qt+t^2)^{-1}\) in arrow degrees \(0,2,4,6,8,10\) are \(1,1,2,3,3,4\), obtained from the recurrence \(v_{d}=C_Qv_{d-1}-v_{d-2}\), with \(v_{-1}=0\) and \(v_0\) the fourth standard basis vector. Multiplying by \(1+t^2+\cdots+t^{58}\) gives the displayed partial sums. The second calculation counts words of weights \(1\) and \(2\), together with powers of the central variable.

Thus the degree-four kernel of \(\mathcal F\to B_5\) has dimension two. The independent relations \(\mathsf P,\mathsf R_4\) span it.

In degree five the kernel has dimension six and is spanned by
\[
 z\mathsf P,\quad\beta\mathsf P,\quad z\mathsf R_4,\quad
 \beta\mathsf R_4,\quad\mathsf R_4\beta,\quad\mathsf R_5.
\]
Their independence follows by comparing successively the coefficients of \(\gamma\beta\gamma\), \(\beta\gamma^2\), \(\gamma^2\beta\), \(z\gamma^2\), \(\beta^5\), and \(z\beta^4\). Write \(Q_\rho=\mathcal F/(\mathsf P,\mathsf R_4,\mathsf R_5)\). Consequently \(Q_\rho\to B_5\) and the induced cocenter maps are isomorphisms through degree five.

We can therefore test the two expressions in these explicit cyclic relation spaces. In degree four, the cyclic relations are exactly the span of \([\mathsf P]\) and \([\mathsf R_4]\). In cyclic words,
\begin{align*}
 [\mathsf P]&=[\beta^4]-22z[\beta^3]+179z^2[\beta^2]-638z^3[\beta]+840z^4,\\
 [\mathsf R_4]&=[\gamma^2]+3[\beta^2\gamma]+[\beta^4]-5z^2[\gamma]
 -5z^2[\beta^2]+4z^4,
\end{align*}
and \([\mathsf R_4']-[\mathsf R_4]=10z^2[\gamma]+10z^2[\beta^2]-8z^4\). If \([\mathsf R_4']\) were in the span, comparing the coefficients of \([\gamma^2]\) and of \(z[\beta^3]\) would force \([\mathsf R_4']=[\mathsf R_4]\), which is false. Hence \(\mathsf R_4'\ne0\) in \(B_5\).

For the second expression, \([\mathsf R_5']-[\mathsf R_5]=10z^2[\gamma\beta]\). The degree-five cyclic relation space is spanned by \([z\mathsf P],[\beta\mathsf P],[z\mathsf R_4], [\beta\mathsf R_4],[\mathsf R_5]\). The coefficient of \(z[\gamma^2]\) first eliminates \([z\mathsf R_4]\) from a possible expression for \(z^2[\gamma\beta]\). The coefficients of \([\gamma^2\beta]\) and \([\gamma\beta^3]\) in \([\beta\mathsf R_4]\) and \([\mathsf R_5]\) are respectively \((1,3)\) and \((2,1)\); their determinant is \(-5\). They therefore eliminate these two classes as well. The remaining classes \([z\mathsf P]\) and \([\beta\mathsf P]\) contain no \(\gamma\), and cannot produce \(z^2[\gamma\beta]\). Thus \([\mathsf R_5']\ne0\), and hence \(\mathsf R_5'\ne0\) in \(B_5\).

The rank table below is computed from \(\mathsf P,\mathsf R_4,\mathsf R_5\).
\end{remark}

\subsection{Relation matrices and modular rank bounds}
\label{app:e8-ranks}

The C++17 program \path{tools/check_e8_corrected_rho.cpp} constructs the finite relation matrices used in Proposition~\ref{prop:e8-exact-modular}. In each degree \(m\le23\), using half the path length as the grading, it enumerates the cyclic-word columns \(z^a[w]_{\mathrm{cyc}}\), with \(a\ge0\) and \(a+\operatorname{wt}(w)=m\), taking one word from each cyclic class. For each \(\mathsf R\in\{\mathsf P,\mathsf R_4,\mathsf R_5\}\), \(a\ge0\), and word \(w\) in \(\beta,\gamma\) satisfying \(a+\deg\mathsf R+\operatorname{wt}(w)=m\), it includes the row \(z^a[\mathsf Rw]_{\mathrm{cyc}}\). The empty word is included and has weight zero; the weights of \(z,\beta,\gamma\) are \(1,1,2\), respectively. The program then performs exact Gaussian elimination over \(\mathbb F_{1000003}\).

The command \texttt{make math} compiles the program, repeats the calculation over \(\mathbb F_{1000033}\) and compares both outputs with the table below.

The second prime checks the rank calculation.  For the characteristic-zero estimate, a nonzero minor modulo either prime remains nonzero over \(\Q\).  Hence \(r_m\) is a lower bound for the characteristic-zero relation rank and \(u_m=c_m-r_m\) is an upper bound for the characteristic-zero cocenter nullity.  The tail and generic-rank arguments in Subsection~\ref{subsec:e8-certificate} establish the exact total dimension.  The table lists \(m\), \(c_m\), \(r_m\), and \(u_m\).
\begin{center}
\small
\renewcommand{\arraystretch}{0.92}
\begin{tabular}{@{}rrrr@{}}
\toprule
Degree \(m\) & Columns \(c_m\) & Rank \(r_m\) over \(\mathbb F_{1000003}\),
\(\mathbb F_{1000033}\) & Nullity bound \(u_m\) \\
\midrule
0  & 1    & 0    & 1 \\
1  & 2    & 0    & 2 \\
2  & 4    & 0    & 4 \\
3  & 6    & 0    & 6 \\
4  & 9    & 2    & 7 \\
5  & 12   & 5    & 7 \\
6  & 17   & 10   & 7 \\
7  & 22   & 16   & 6 \\
8  & 30   & 24   & 6 \\
9  & 40   & 34   & 6 \\
10 & 55   & 49   & 6 \\
11 & 74   & 69   & 5 \\
12 & 105  & 100  & 5 \\
13 & 146  & 142  & 4 \\
14 & 210  & 206  & 4 \\
15 & 304  & 300  & 4 \\
16 & 447  & 443  & 4 \\
17 & 658  & 655  & 3 \\
18 & 987  & 984  & 3 \\
19 & 1480 & 1478 & 2 \\
20 & 2246 & 2244 & 2 \\
21 & 3416 & 3414 & 2 \\
22 & 5227 & 5225 & 2 \\
23 & 8014 & 8013 & 1 \\
\bottomrule
\end{tabular}
\end{center}
The SHA-256 checksums of the two \(E_8\) programs are:
\begin{center}
\scriptsize
\begin{tabular}{@{}l@{\hspace{1em}}l@{}}
\toprule
File & SHA-256 \\
\midrule
\texttt{check\_e8\_corrected\_rho.cpp} &
\texttt{2545d84e337ea8b4ac650c428b90880a36ab47cd168e8935e0d4598ac1dbabd9}\\
\texttt{check\_e8\_rho\_relations\_formal.py} &
\texttt{208305e7e5a7613b75b66959b4ac7022c31644b54c988cbd642baad7382d5f1a}\\
\bottomrule
\end{tabular}
\end{center}

\Needspace{14\baselineskip}
\subsection{The cocenter dimension and the tail}
\label{subsec:e8-certificate}

The relation-matrix bounds and the grading bounds together determine the cocenter for the \(E_8\) marking of length five. The lower bound comes from the generic decomposition indexed by positive roots in \eqref{eq:generic-cocenter-rank}, established before the closed-fibre count in Lemma~\ref{lem:closed-fibre-cocenter-count}.

Deleting the retained vertex for the length-five marking of \(E_8\) leaves \(A_3\times A_4\), so
\begin{equation}
\label{eq:e8-N_p}
 N_p=|\Phi^+(E_8)|-|\Phi^+(A_3)|-|\Phi^+(A_4)|
 =120-6-10=104.
\end{equation}
Let \(B_5=e_p\Pi^\rho_0e_p\) be the closed marked corner. Subsection~\ref{app:e8-specialization} derives its relations directly from the preprojective relations specialized at \(\rho_i=1\), \(\lambda=0\).  It constructs elements \(\beta,\gamma\in B_5\) that, together with \(z\), generate \(B_5\) (Lemma~\ref{lem:e8-generators}). The three relations
\begin{align*}
 \mathsf P&=(\beta-4z)(\beta-5z)(\beta-6z)(\beta-7z),\\
 \mathsf R_4&=(\gamma+\beta^2)^2+\beta\gamma\beta-5z^2(\gamma+\beta^2)+4z^4,\\
 \mathsf R_5&=\gamma\beta\gamma+\gamma^2\beta+\gamma\beta^3-5z^2\gamma\beta
\end{align*}
vanish in \(B_5\) (Lemma~\ref{lem:e8-relations}).  Give \(\C[z]\langle\beta,\gamma\rangle\) and \(B_5\) the grading equal to half the path length:
\[
 \deg z=\deg\beta=1,
 \qquad\deg\gamma=2.
\]
The three relations are homogeneous, so there is a surjection of graded algebras
\begin{equation}
\label{eq:e8-presentation}
 Q_\rho:=\frac{\C[z]\langle\beta,\gamma\rangle}
 {(\mathsf P,\mathsf R_4,\mathsf R_5)}\twoheadrightarrow B_5 .
\end{equation}
The induced surjection \(HH_0(Q_\rho)\twoheadrightarrow HH_0(B_5)\) transfers cocenter upper bounds from \(Q_\rho\) to \(B_5\).  Thus it is enough first to bound the cocenter of the quotient defined by the three displayed relations. Remark~\ref{rem:karmazyn-sign-conventions} explains why the signs of the parameter terms differ from a direct substitution into the displayed universal relations in the arXiv version of \cite{Karmazyn} specified in that remark.

\begin{proposition}
\label{prop:e8-exact-modular}
Let \(h_m=\dim_{\C}HH_0(B_5)_m\).  In degrees \(0,\ldots,23\), the following sequence gives upper bounds for \(h_m\):
\[
 (1,2,4,6,7,7,7,6,6,6,6,5,5,4,4,4,4,3,3,2,2,2,2,1).
\]
The sum of these bounds is \(99\).
\end{proposition}

\begin{proof}
Consider degree \(m\) of the integral matrix construction in Subsection~\ref{app:e8-ranks}: its columns are the cyclic words \(z^a[w]_{\mathrm{cyc}}\), and its rows are all multiples \(z^a[\mathsf Rw]_{\mathrm{cyc}}\) of that degree, for \(\mathsf R\in\{\mathsf P,\mathsf R_4,\mathsf R_5\}\), including the empty word.  These rows span all cocenter relations: a two-sided multiple \(u\mathsf Rv\) has the same cyclic class as \(\mathsf Rvu\), and central powers of \(z\) can be moved to the front. The cokernel in degree \(m\) is therefore \(HH_0(Q_\rho)_m\).

The relation matrices have integer entries.  Exact Gaussian elimination over \(\mathbb F_{1000003}\) gives the ranks recorded in Subsection~\ref{app:e8-ranks}; running the same computation over \(\mathbb F_{1000033}\) gives the same ranks as an additional modular check.

A minor that is nonzero modulo either prime is nonzero over \(\Q\). The modular rank therefore bounds the characteristic-zero rank from below.  Subtracting from the column count gives an upper bound for \(\dim HH_0(Q_\rho)_m\). The graded surjection \eqref{eq:e8-presentation} induces a surjection on cocenters and hence gives the stated upper bound for \(h_m\).
\end{proof}

\subsubsection{Tail bound and exact grading comparison}
\label{subsubsec:e8-tail}

The finite calculation gives upper bounds through degree \(23\). We first bound the remaining degrees, then compare the resulting total with the generic rank of the universal cocenter. For \(E_8\) the Coxeter number is \(h=30\). In the Etingof--Rains grading, arrows have degree \(1\) and \(z\) has degree \(2\). Their Hilbert-series formulas \cite[Section~2.1, equation~(1), and Theorem~3.2(iii)--(iv)]{EtingofRains} show that the ordinary preprojective algebra has top arrow degree \(h-2=28\), while the central extension \(\Pi^\rho_0\) has top degree \(2(h-2)=56\). On the marked corner \((z,\beta,\gamma)\) have degrees \((2,2,4)\). The finite computation uses exactly half this grading, so \(B_5\) vanishes in degrees above \(28\).

For the quotient by \(z\) we have a sharper bound. Setting \(z=0\) recovers the ordinary preprojective corner \(B_5/(z)\). Because the \(E_8\) Dynkin graph is bipartite, all closed paths have even arrow length. The top arrow degree bound \(h-2=28\) in the Etingof--Rains grading therefore forces the corner \(B_5/(z)\) to vanish in degrees above \(14\) in the halved grading.

Because \(z\) is central in \(B_5\), reduction modulo \(z\) commutes with the cocenter functor, giving
\[
 HH_0(B_5)/(z)\cong HH_0(B_5/(z)).
\]
The vanishing of \(B_5/(z)\) in degrees above \(14\) implies that \(HH_0(B_5/(z))_m=0\) for \(m>14\). Consequently, multiplication by \(z\) on the cocenter is surjective
\[
 HH_0(B_5)_{m-1}\twoheadrightarrow HH_0(B_5)_m
 \qquad(m>14).
\]
The surjectivity of multiplication by \(z\) makes the dimensions nonincreasing for \(m>14\). The finite calculation in Proposition~\ref{prop:e8-exact-modular} gives \(\dim_\C HH_0(B_5)_{23}\le1\). Successive surjectivity therefore forces every remaining degree \(24\le m\le 28\) to have dimension at most one, while degrees above \(28\) vanish. Adding this tail bound of \(5\times 1=5\) to the sum of \(99\) through degree \(23\) yields the total cocenter upper bound
\[
 \dim HH_0(B_5)\le 99+5=104.
\]

We now establish the matching lower bound from the universal family. The finite module \(\mathsf M_p\) has generic rank \(N_p=104\) by \eqref{eq:generic-cocenter-rank}, and its closed fibre over the origin is \(HH_0(B_5)\). Because the closed-fibre dimension of a finite module over a local ring is at least its generic rank, we obtain the reverse inequality
\[
 \dim HH_0(B_5)\ge104.
\]
The two inequalities together force the exact dimension
\begin{equation}
\label{eq:e8-equality}
 \dim HH_0(B_5)=104=N_p.
\end{equation}
This proves \eqref{eq:special-cocenter-dimension} in the remaining case. The degreewise inequalities now have equal sums on the two sides. Thus every bound of Proposition~\ref{prop:e8-exact-modular} is attained, each degree \(24,\ldots,28\) contributes exactly one, and \(HH_0(Q_\rho)_m\to HH_0(B_5)_m\) is an isomorphism for \(m\le23\).

\section{Exact linear algebra for \texorpdfstring{\(D_{2,2}\)}{D2,2}}
\label{app:d22-computation}

Theorem~\ref{thm:dn2-nonlifting} proves the centre-lifting obstruction symbolically for every \(D_{n,2}\).  This appendix gives an independent exact-rational calculation for \(D_{2,2}\): it computes the tangent data in Proposition~\ref{prop:d22-tangent-data} and the automorphism-invariant socle condition \eqref{eq:socle-diagonal} used in Proposition~\ref{prop:d22-unmarked}.  The program \path{tools/verify_conjecture_1_4_counterexample.py} computes the cochain ranks, and \path{tools/check_d22_tangent_image.py} verifies the tangent image and socle condition.  Both use exact rational arithmetic (\texttt{fractions.Fraction} in Python) and are run by \texttt{make math}. The computation is specific to \(n=2\).  Theorem~\ref{thm:dn2-nonlifting} and Corollary~\ref{cor:centre-morphism-not-smooth} prove non-lifting symbolically for every \(n\ge2\).

\Needspace{7\baselineskip}
\paragraph{Cochains and ranks.}
The algebra \(A=D_{2,2}\) is given in the basis \((1,x,x^2,x^3,x^4,y,xy,x^2y)\) by the multiplication table of Lemma~\ref{lem:dn2-structure}(i), and the first program checks associativity and the unit on all triples of basis vectors.   The centre \(Z_2\) has the basis \((1,\su,\sv,\sw,\sq)\), embedded as \((1,x^2,x^3,x^2y,x^4)\), and the same program checks that the embedding is multiplicative.

Normalized cochains are multilinear maps on the span of the basis vectors other than \(1\); for \(Z_2\), the second cochains are symmetric.  The coboundary maps are sparse rational matrices.  Exact row reduction gives the following dimensions.
\begin{center}
\small
\begin{tabular}{@{}lrr@{}}
\toprule
 & \(A\), associative & \(Z_2\), commutative \\
\midrule
normalized \(1\)-cochains & \(56\) & \(20\) \\
rank of \(\delta\) on \(1\)-cochains & \(50\) & \(10\) \\
normalized \(2\)-cochains & \(392\) & \(50\) (symmetric) \\
rank of \(\delta\) on \(2\)-cochains & \(339\) & \(24\) \\
\(2\)-cocycles & \(53\) & \(26\) \\
cocycles modulo coboundaries & \(HH^2(A)\): \(3\) & \(T^1_{\mathrm{com}}(Z_2)\): \(16\) \\
\bottomrule
\end{tabular}
\end{center}
The derivations of \(A\) form the kernel on \(1\)-cochains, of dimension \(56-50=6\); inner derivations have dimension \(8-5=3\), so \(\dim HH^1(A)=3\).

\paragraph{The centre tangent map.}
For each of the \(53\) basis cocycles \(\phi\) of \(Z^2_{\mathrm{norm}}(A,A)\), the first program solves the linear system~\eqref{eq:central-lift} for \(\tau\) with values in the chosen complement \(\operatorname{span}(x,y,xy)\) of \(Z_2\).  Requiring \(\tau\) to take values in this complement fixes the representative used in the calculation.  The first program checks that the solution exists and is unique, forms \(\psi\) by \eqref{eq:induced-centre-cocycle}, and checks that \(\psi\) is symmetric, takes values in \(Z_2\), and is a cocycle.

It also verifies that the coboundaries of \(A\) map to commutative coboundaries and that the \(53\) cocycles \(\psi\) have rank \(3\) modulo \(B^2_{\mathrm{com}}(Z_2)\).  Finally, the coefficient of \(\sw\) in \(\psi(\sv,\sv)\) vanishes on all these cocycles and on \(B^2_{\mathrm{com}}(Z_2)\), but equals \(1\) on \(\psi_0\).  This is an independent exact verification of Corollary~\ref{cor:centre-morphism-not-smooth}(i) for \(n=2\).

\paragraph{Representatives and the socle condition.}
\label{app:d22-representatives}
The representatives in Proposition~\ref{prop:d22-tangent-data}(iii) are normalized and vanish on every pair not displayed below, apart from the values imposed by symmetry:
\[
\begin{array}{ll}
 \psi_1(\su,\sw)=-\sw,\quad \psi_1(\su,\sq)=-2\sq,\quad
 \psi_1(\sv,\sv)=\sq;\\[2pt]
 \psi_2(\su,\sv)=-\tfrac12\sw,\quad \psi_2(\su,\sw)=\tfrac12\sv,\quad
 \psi_2(\sv,\sw)=\tfrac12\sq;\\[2pt]
 \psi_3(\su,\sw)=-2\sw,\quad \psi_3(\su,\sq)=-2\sq,\quad
 \psi_3(\sw,\sw)=2\sq.
\end{array}
\]
These formulas also give \(\psi_i(s,s)\in\C\sq\) for every \(s\in\Soc(Z_2)\).

The second program checks that \(\psi_1,\psi_2,\psi_3\) of Proposition~\ref{prop:d22-tangent-data}(iii) are commutative cocycles whose classes lie in the image of \(\chi_A\) and span it.  Reducing them modulo \(B^2_{\mathrm{com}}(Z_2)\) by the normalized row-echelon basis of the coboundaries leaves the coefficients of \(\sw\) in \(\psi(\su,\sv)\), of \(\sw\) in \(\psi(\su,\sw)\), and of \(\sq\) in \(\psi(\su,\sq)\) as non-pivot coordinates.  The coefficient vectors of \(\psi_1,\psi_2,\psi_3\) in these three coordinates form the matrix
\[
 \begin{pmatrix}0&-1&-2\\-\frac12&0&0\\0&-2&-2\end{pmatrix},
 \qquad\det=-1,
\]
whose rows correspond to \(\psi_1,\psi_2,\psi_3\).  The program then verifies~\eqref{eq:socle-diagonal} for all \(53\) image cocycles and for the \(20\) generators of \(B^2_{\mathrm{com}}(Z_2)\), as a polynomial identity in the coordinates of \(s=\alpha\sv+\beta\sw+\gamma\sq\), and checks that \(\psi_0\) violates \eqref{eq:socle-diagonal}.

\paragraph{Checks for small values of \(n\).}
For \(n=2,\ldots,6\), the first program also checks the identities of Lemma~\ref{lem:dn2-structure} used in the proof of Theorem~\ref{thm:dn2-nonlifting}: \(\dim A_n=4n\), \(\dim Z(A_n)=n+3\), the centralizer of \(x\) has dimension \(2n+2\), \(\sw x=x\sw=0\), and \(\sv A_n+A_n\sv\subseteq\C[x]\).  It also checks that the multiplication of~\eqref{eq:dn2-centre-family} for \(n=2\) is associative over \(\Q[t]\). These are finite checks of identities whose general proof is already given in Lemma~\ref{lem:dn2-structure}.

\Needspace{5\baselineskip}
The SHA-256 checksums of the two programs are:
\begin{center}
\scriptsize
\begin{tabular}{@{}l@{\hspace{1em}}l@{}}
\toprule
File & SHA-256 \\
\midrule
\texttt{verify\_conjecture\_1\_4\_counterexample.py} &
\texttt{e51dd6b0816f019f3cd2f38c9fb12312ec4b9f49d18323a18748055bdbb7ccef}\\
\texttt{check\_d22\_tangent\_image.py} &
\texttt{2140afdf251c2fb1f123a02cdb92fac87b4d029e822f336a4a2ca1e9e92dc9ae}\\
\bottomrule
\end{tabular}
\end{center}

\clearpage
\section*{List of notation}
\addcontentsline{toc}{section}{List of notation}
\begingroup
\small
\setlength{\tabcolsep}{0pt}
\renewcommand{\arraystretch}{1.12}
\newenvironment{notationgroup}[1]{%
  \par\medskip\noindent
  \begin{tabular}{@{}>{\raggedright\arraybackslash}p{.20\linewidth}@{\hspace{.025\linewidth}}>{\raggedright\arraybackslash}p{.595\linewidth}@{\hspace{.025\linewidth}}>{\raggedright\arraybackslash}p{.155\linewidth}@{}}
  \multicolumn{3}{@{}l}{\normalsize\bfseries #1}\\[2pt]\midrule
}{\end{tabular}\par}
\begin{notationgroup}{The contraction and the marked surface}
\(f\colon X\to\Spec R,\ C\) & Completed simple flop and its reduced exceptional curve \(C\cong\mathbb P^1\). & \S\,\ref{sec:general-existence}\\[2pt]
\(R\oplus N,\ \Lambda\) & Basic generator and its opposite endomorphism algebra \(\Lambda=\End_R(R\oplus N)^{\mathrm{op}}\). & \S\,\ref{sec:introduction}\\[2pt]
\(e_0,\ A=\Acon\) & Projector onto the free summand and the contraction algebra \(A=\Lambda/\Lambda e_0\Lambda\). & \S\,\ref{sec:introduction}\\[2pt]
\(\Phi,\ \Phi^+,\ \mathfrak h,\ W\) & ADE root system of the surface section, its positive roots, its Cartan space and its Weyl group. & Def.~\ref{def:marked-partial-resolution-data}\\[2pt]
\(*,\ p,\ c_*(\alpha)\) & The retained vertex is written \(*\) or \(p\). Its coefficient in a root \(\alpha\) is \(c_*(\alpha)=m_p(\alpha)\). & Def.~\ref{def:marked-partial-resolution-data}, \S\,\ref{subsec:universal-marked-cocenter}\\[2pt]
\(\ell,\ \theta\) & Length of the flop and highest root, related by \(\ell=c_*(\theta)\). & Def.~\ref{def:marked-partial-resolution-data}\\[2pt]
\(W_0,\ H\) & Weyl subgroup generated by the reflections at the deleted vertices, and the marked deformation base \(H=\mathfrak h/W_0\). & Def.~\ref{def:marked-partial-resolution-data}\\[2pt]
\(q,\ \pi_W\) & Quotient map \(q\colon\mathfrak h\to H\) and induced finite map \(\pi_W\colon H\to\mathfrak h/W\). & Def.~\ref{def:marked-partial-resolution-data}, \S\,\ref{subsec:full-marked-morsification}\\[2pt]
\(H_\alpha,\ D_j,\ P_j\) & Root hyperplane in \(\mathfrak h\), curve divisor in \(H\) for class \(j[C]\), and a reduced equation of that divisor. & \S\,\ref{subsec:full-marked-morsification}\\[2pt]
\(D^{\mathrm{curv}},\ \Sigma^{\mathrm{bad}}\) & Union of the curve divisors, and the closed subset removed so that exactly one root pair vanishes at each lift of a remaining point. & Lem.~\ref{lem:root-divisors-and-bad-locus}\\[2pt]
\end{notationgroup}

\begin{notationgroup}{The classifying arc and the curve counts}
\(f_0,\ \nu_0,\ R_1\) & Original analytic classifying arc, its formal completion at \(w=0\), and the arc ring \(R_1=\C[[w]]\). & Lem.~\ref{lem:marked-analytic-classifying-arc}\\[2pt]
\(s,\ w,\ v,\ F(s,w)\) & In affine coordinates on \(H\), \(v\in\C^{\dim H}\) is a fixed translation vector. For each \(s\), the arc \(w\mapsto F(s,w)=f_0(w)+sv\) is the translate of the original arc by \(sv\). & Lem.~\ref{lem:good-constant-translate-expanded}\\[2pt]
\(R_0,\ B\) & For the translated one-curve family, \(R_0=\C[[s]]\) is the deformation ring and \(B=R_0[[w]]\) also contains the arc coordinate. & \S\,\ref{subsec:full-general-assembly}\\[2pt]
\(C_j,\ C^{(j)}\) & The symbolic thickening cut out by \(\mathcal I_C^{(j)}\) on the central surface, and Toda's pure thickening. They define the same subscheme. & Lem.~\ref{lem:symbolic-pure-thickening}\\[2pt]
\(n_j\) & Genus-zero GV invariant of class \(j[C]\). & Prop.~\ref{prop:general-source-package}\\[2pt]
\(k_j\) & Length of the original threefold's local Hilbert ring at \([C_j]\). & Prop.~\ref{prop:conserved-multiplicities-gv}\\[2pt]
\(m_j,\ W_j\) & \(m_j=\operatorname{ord}_{w=0}P_j(f_0(w))\). The roots of the monic polynomial \(W_j(s,T)\) give the intersections of the translated arc with \(D_j\). Its degree is \(m_j=k_j=n_j\). & Lem.~\ref{lem:weierstrass-root-polynomials}, Prop.~\ref{prop:conserved-multiplicities-gv}\\[2pt]
\(\mathscr H_j\) & Completed local ring of the relative Hilbert space at the point representing \(C_j\), taken on an \'etale chart. & Prop.~\ref{prop:conserved-multiplicities-gv}\\[2pt]
\end{notationgroup}

\clearpage
\begin{notationgroup}{The relative algebras and their cocenters}
\(\mathbb H_\ell,\ R_H,\ \widetilde R\) & Ring of regular functions on \(H\), its completion at the origin, and the completed local ring of \(\mathfrak h\) at the origin, respectively. & \S\,\ref{subsubsec:marked-quotient-downstairs}, \S\,\ref{subsec:universal-marked-cocenter}\\[2pt]
\(\mathcal T,\ \mathcal P,\ \mathcal A_\ell\) & Universal tilting bundle, its positive dual \(\mathcal P=\mathcal T^\vee\), and the universal flopping algebra \(\mathcal A_\ell=\End(\mathcal T)^{\mathrm{op}}\). & \S\,\ref{subsec:full-universal-tilting}\\[2pt]
\(S_\nu,\ \Lambda_\nu,\ \mathscr A\) & Completed target ring after pullback along the translated family, the opposite endomorphism algebra of the pulled-back tilting bundle, and the contraction quotient \(\mathscr A=\Lambda_\nu/\Lambda_\nu e_0\Lambda_\nu\). & \S\,\ref{subsec:full-general-assembly}\\[2pt]
\(D_B\) & For deformations over a complete local base \(B\), the classifying homomorphism takes values in \(D_B=B[[w]]\). & Prop.~\ref{prop:arbitrary-base-classifying-arcs}\\[2pt]
\(K,\ K'\) & Generic field \(K=\C((s))\) and finite splitting extension \(K'=\C((u))\), with \(s=u^e\). & Prop.~\ref{prop:full-marked-formal-morsification}\\[2pt]
\(Z(A),\ HH_0(A)\) & Centre of \(A\) and its cocenter \(A/[A,A]\), where \([A,A]\) is the linear span of commutators. & Lem.~\ref{lem:centre-cocenter-duality}\\[2pt]
\(\mathcal C_H,\ \mathsf M_H\) & Quotient of \(\widehat{\mathcal A}_\ell\) by the ideal generated by \(e_0\), and its relative cocenter \(HH_0(\mathcal C_H/R_H)\). & \S\,\ref{subsubsec:marked-quotient-downstairs}\\[2pt]
\(\mathcal C_0,\ \mathsf M_0\) & Pullback of \(\mathcal C_H\) to the original arc and its cocenter over \(R_1\). The marked comparison identifies \(\mathcal C_0\) with \(A\). & Lem.~\ref{lem:arc-quotient-supported}\\[2pt]
\(\Pi_\lambda,\ e_p\) & Deformed preprojective algebra over \(\widetilde R\) and the vertex idempotent for the retained vertex \(p\). & \S\,\ref{subsec:universal-marked-cocenter}\\[2pt]
\(E_p,\ \mathcal D_p,\ z\) & Corner of the central extension at \(p\), the deformed preprojective corner \(\mathcal D_p=E_p/(z)\), and the central element \(z\). & \S\,\ref{subsec:universal-marked-cocenter}\\[2pt]
\(\mathsf M_p,\ N_p\) & Free cocenter \(HH_0(E_p/\widetilde R)\). Its rank \(N_p\) counts positive roots with \(m_p(\alpha)>0\). & \eqref{eq:central-extension-corner}, \eqref{eq:N_p-definition}\\[2pt]
\(\Psi_p,\ \Theta\) & Matrix of multiplication by \(z\) on \(\mathsf M_p\), and a presentation matrix for \(\mathsf M_H\) over \(R_H\), respectively. & \eqref{eq:universal-cocenter-presentation}, \eqref{eq:downstairs-presentation}\\[2pt]
\(\Theta_A\) & Presentation matrix for \(\mathcal C_H\). Up to a unit its determinant is \(\prod_j P_j^{j^2}\), whereas \(\det\Theta\) contains each \(P_j\) once. & \eqref{eq:weighted-contraction-determinant}, Thm.~\ref{thm:marked-cocenter}\\[2pt]
\end{notationgroup}

\begin{notationgroup}{Semisimple fibres and associative laws}
\(S_{\mathrm{GV}},\ N_{\mathrm{GV}}\) & Split algebra \(\prod_{j=1}^\ell M_j(\C)^{\times n_j}\) and its number of matrix factors \(N_{\mathrm{GV}}=\sum_j n_j\). & \eqref{eq:general-gv-target}\\[2pt]
\(d,\ X_d,\ G\) & For \(d=\dim A\), the scheme of associative laws on a fixed \(d\)-dimensional space \(V\) with identity vector \(1\), and its change-of-basis group \(G=\Stab_{\GL(V)}(1)\). & \S\,\ref{sec:orbit-criterion}\\[2pt]
\(HH^i(A),\ A^e\) & Hochschild cohomology \(\operatorname{Ext}^i_{A^e}(A,A)\), where \(A^e=A\otimes_\C A^{\mathrm{op}}\). & \S\,\ref{subsec:structural-reduction}\\[2pt]
\end{notationgroup}

\begin{notationgroup}{The length-five \(E_8\) calculation}
\(B_5,\ Q_\rho,\ h_m\) & Closed corner \(B_5\) at the retained vertex, the algebra \(Q_\rho\) defined by the three relations with a graded surjection onto \(B_5\), and the graded dimension \(h_m=\dim_\C HH_0(B_5)_m\). & \eqref{eq:e8-presentation}, Prop.~\ref{prop:e8-exact-modular}\\[2pt]
\end{notationgroup}

\clearpage
\begin{notationgroup}{Flops with several retained curves}
\(J,\ I,\ W_I\) & Vertices retained in the partial resolution, the complementary set of deleted vertices, and the Weyl subgroup generated by their reflections \(s_i\) for \(i\in I\). & \S\,\ref{subsec:multicurve-marked-quotient}\\[2pt]
\(S,\ R_J\) & Completed Cartan ring \(S=\widehat{\mathcal O}_{\mathfrak h,0}\) and invariant ring \(R_J=S^{W_I}\) of the marked base. & \S\,\ref{subsec:multicurve-marked-quotient}\\[2pt]
\(\pi_I,\ \mathcal B_J\) & Restriction to the simple-root coefficients at vertices in \(J\), and the set of distinct nonzero restrictions of positive roots. & \S\,\ref{subsec:multicurve-marked-quotient}\\[2pt]
\(P_\beta,\ \delta_J\) & Reduced equation of the curve divisor indexed by \(\beta\in\mathcal B_J\), and the product \(\delta_J=\prod_\beta P_\beta\). & Lem.~\ref{lem:multicurve-root-products}, Prop.~\ref{prop:multicurve-reduced-annihilator}\\[2pt]
\(s_\beta,\ n_\beta\) & Matrix size \(s_\beta=\sum_{i\in J}\beta_i\) and genus-zero GV multiplicity of the class \(\sum_{i\in J}\beta_i[C_i]\). & \S\,\ref{sec:marked-semisimple}\\[2pt]
\(C_J,\ D_J\) & Universal contraction quotient over \(R_J\) and its extension to \(S\), the preprojective corner \(D_J=e_J\Pi_{\Delta,\lambda}e_J\). & \eqref{eq:multicurve-cartan-quotient}\\[2pt]
\(\mathscr E_J,\ M_J,\ \mathcal H_J\) & Retained corner of the central extension, its free cocenter \(M_J=HH_0(\mathscr E_J/S)\), and the cocenter \(\mathcal H_J=HH_0(C_J/R_J)\) of the contraction quotient. & Thm.~\ref{thm:multicurve-universal-cocenter}\\[2pt]
\(N_J\) & Number of positive roots with nonzero restriction to \(J\). This is the rank of \(M_J\) over \(S\). & Thm.~\ref{thm:all-corner-trace-excision}, Thm.~\ref{thm:multicurve-universal-cocenter}\\[2pt]
\(\Theta_J,\ \Theta_{\mathrm{tr},J}\) & Presentation matrices over \(R_J\) for the algebra \(C_J\) and its cocenter \(\mathcal H_J\), respectively. & \eqref{eq:multicurve-weighted-resolution}, \eqref{eq:multicurve-cocenter-resolution}\\[2pt]
\(E,\ V_\beta,\ S_E\) & Algebra \(E=\bigoplus_{i\in J}\C e_i\) of labelled projectors, block space \(V_\beta=\bigoplus_{i\in J}\C^{\beta_i}\), and marked semisimple algebra \(S_E=\prod_\beta\End_\C(V_\beta)^{\times n_\beta}\) with \(N_{\mathrm{GV}}=\sum_\beta n_\beta\) factors. & \eqref{eq:marked-target-definition}\\[2pt]
\(V_{ij},\ c_{ij}\) & Peirce space \(e_iVe_j\) of a fixed \(E\)-bimodule \(V\), and its dimension. & \S\,\ref{subsec:marked-orbit-criterion}\\[2pt]
\(X_E(V),\ G_E\) & Scheme of associative laws with the prescribed \(E\)-actions and idempotents, and the group of \(E\)-bimodule automorphisms fixing every \(e_i\). & \S\,\ref{subsec:marked-orbit-criterion}\\[2pt]
\end{notationgroup}

\begin{notationgroup}{Centres and their deformations}
\(\Com_c,\ \Sm_c\) & Scheme of commutative associative laws with fixed identity on \(\C^c\), and its smoothable component, the reduced orbit closure of the split algebra \(\C^c\). & \S\,\ref{subsec:centre-morphism-setup}\\[2pt]
\(\mathsf Z_U,\ \mathscr F_U,\ \cent\) & Centre bundle on \(U\subset X_d\), its scheme of unit-preserving frames, and the map \(\mathscr F_U\to\Com_c\) recording multiplication in the chosen frame. & Def.~\ref{def:centre-morphism}\\[2pt]
\(T^1_{\mathrm{com}}(Z)\) & First-order commutative deformations of \(Z\) modulo infinitesimal changes of basis preserving the unit. & \S\,\ref{subsec:centre-morphism-setup}\\[2pt]
\(\chi_B\) & For \(B\in U(\C)\), the map \(HH^2(B)\to T^1_{\mathrm{com}}(Z(B))\) taking a first-order associative deformation class to its induced centre deformation class. & \S\,\ref{subsec:centre-tangent-map}\\[2pt]
\(A_n,\ Z_n\) & Type \(D\) algebra \(A_n=D_{n,2}\) for \(n\ge2\) and its centre \(Z_n=Z(A_n)\). The case \(n=2\) is computed in Appendix~\ref{app:d22-computation}. & \eqref{eq:dn2-presentation}\\[2pt]
\(\su,\ \sv,\ \sw\) & Central elements \(x^2,\ x^{2n-1},\ x^{2n-2}y\) in \(A_n\), respectively. The sans-serif letters distinguish them from deformation parameters. & Lem.~\ref{lem:dn2-structure}\\[2pt]
\end{notationgroup}
\endgroup
\clearpage



\begin{thebibliography}{99}
\addcontentsline{toc}{section}{References}
\frenchspacing

\bibitem[Amiot(09)]{Amiot}
C.~Amiot, \emph{Cluster categories for algebras of global dimension 2 and quivers with potential}, Ann. Inst. Fourier \textbf{59} (2009), 2525--2590. \href{https://arxiv.org/abs/0805.1035v2}{arXiv:0805.1035v2}.

\bibitem[August(20)]{August}
J.~August, \emph{The tilting theory of contraction algebras}, Adv. Math. \textbf{374} (2020), Article 107372. \href{https://doi.org/10.1016/j.aim.2020.107372}{doi:10.1016/\allowbreak j.aim.2020.107372}. \href{https://arxiv.org/abs/1802.10366v2}{arXiv:1802.10366v2}.

\bibitem[Bowman--Pulji\'c--Smoktunowicz(24)]{BowmanPuljicSmoktunowicz}
D.~Bowman, D.~Pulji\'c, and A.~Smoktunowicz, \emph{A construction of deformations to general algebras}, Int. Math. Res. Not. IMRN (2024), no.~12, 9827--9848. \href{https://doi.org/10.1093/imrn/rnae077}{doi:10.1093/imrn/rnae077}.

\bibitem[Brion(10)]{Brion}
M.~Brion, \emph{Introduction to actions of algebraic groups}, Les cours du C.I.R.M. \textbf{1} (2010), no.~1, 1--22.

\bibitem[Brown--Wemyss(25)]{BrownWemyss}
G.~Brown and M.~Wemyss, \emph{Local normal forms of noncommutative functions}, Forum Math. Pi \textbf{13} (2025), Paper e8. \href{https://doi.org/10.1017/fmp.2025.2}{doi:10.1017/fmp.2025.2}.

\bibitem[Bryan--Katz--Leung(01)]{BKL}
J.~Bryan, S.~Katz, and N.~C. Leung, \emph{Multiple covers and the integrality conjecture for rational curves in Calabi--Yau threefolds}, J. Algebraic Geom. \textbf{10} (2001), no.~3, 549--568. Preprint: \href{https://arxiv.org/abs/math/9911056v1}{arXiv:math/9911056v1}.

\bibitem[Chevalley(55)]{Chevalley1955}
C.~Chevalley, \emph{Invariants of finite groups generated by reflections}, Amer. J. Math. \textbf{77} (1955), no.~4, 778--782. \href{https://doi.org/10.2307/2372597}{doi:10.2307/2372597}.

\bibitem[Donovan--Wemyss(16)]{DonovanWemyss}
W.~Donovan and M.~Wemyss, \emph{Noncommutative deformations and flops}, Duke Math. J. \textbf{165} (2016), no.~8, 1397--1474.

\bibitem[Dugas(12)]{Dugas}
A.~Dugas, \emph{Periodicity of \(d\)-cluster-tilted algebras}, J. Algebra \textbf{368} (2012), 40--52. \href{https://arxiv.org/abs/1007.2811v2}{arXiv:\allowbreak1007.2811v2}.

\bibitem[Etingof--Rains(06)]{EtingofRains}
P.~Etingof and E.~Rains, \emph{Central extensions of preprojective algebras, the quantum Heisenberg algebra, and 2-dimensional complex reflection groups}, J. Algebra \textbf{299} (2006), no.~2, 570--588. \href{https://arxiv.org/abs/math/0503393v3}{arXiv:math/0503393v3}.

\bibitem[Etingof--Latour--Rains(07)]{EtingofLatourRains}
P.~Etingof, F.~Latour, and E.~Rains, \emph{On central extensions of preprojective algebras}, J. Algebra \textbf{313} (2007), no.~1, 165--175. \href{https://arxiv.org/abs/math/0606403v1}{arXiv:math/0606403v1}.

\bibitem[Gabriel(75)]{Gabriel}
P.~Gabriel, \emph{Finite representation type is open}, in V.~Dlab and P.~Gabriel (eds.), \emph{Representations of Algebras}, Lecture Notes in Math. \textbf{488}, Springer, 1975, pp.~132--155. \href{https://doi.org/10.1007/BFb0081219}{doi:10.1007/BFb0081219}.

\bibitem[Gerstenhaber(64)]{Gerstenhaber}
M.~Gerstenhaber, \emph{On the deformation of rings and algebras}, Ann. of Math. (2) \textbf{79} (1964), 59--103.

\bibitem[Gille--Parimala--Suresh(21)]{GilleParimalaSuresh}
P.~Gille, R.~Parimala, and V.~Suresh, \emph{Local triviality for \(G\)-torsors}, \href{https://arxiv.org/abs/1901.04722v5}{arXiv:1901.04722v5} (2021).

\bibitem[Harrison(62)]{Harrison}
D.~K. Harrison, \emph{Commutative algebras and cohomology}, Trans. Amer. Math. Soc. \textbf{104} (1962), 191--204.

\bibitem[Herschend--Iyama(11)]{HerschendIyama}
M.~Herschend and O.~Iyama, \emph{Selfinjective quivers with potential and 2-representation-finite algebras}, Compos. Math. \textbf{147} (2011), no.~6, 1885--1920. \href{https://doi.org/10.1112/S0010437X11005367}{doi:10.1112/S0010437X11005367}. Preprint: \href{https://arxiv.org/abs/1006.1917v2}{arXiv:1006.1917v2}.

\bibitem[Hua--Toda(18)]{HuaToda}
Z.~Hua and Y.~Toda, \emph{Contraction algebra and invariants of singularities}, Int. Math. Res. Not. IMRN (2018), no.~10, 3173--3198. Preprint: \href{https://arxiv.org/abs/1601.04881v2}{arXiv:1601.04881v2}. The theorem numbering used here follows this preprint.

\bibitem[Iyama--Kalck--Wemyss--Yang(15)]{IKWY}
O.~Iyama, M.~Kalck, M.~Wemyss, and D.~Yang, \emph{Frobenius categories, Gorenstein algebras and rational surface singularities}, Compos. Math. \textbf{151} (2015), no.~3, 502--534. \href{https://doi.org/10.1112/S0010437X14007647}{doi:10.1112/S0010437X14007647}. Preprint: \href{https://arxiv.org/abs/1209.4215v3}{arXiv:1209.4215v3}.

\bibitem[Jasso--Keller--Muro(24)]{JKM}
G.~Jasso, B.~Keller, and F.~Muro, \emph{The Donovan--Wemyss conjecture via the derived Auslander--Iyama correspondence}, in \emph{Triangulated Categories in Representation Theory and Beyond}, Abel Symposia \textbf{17}, Springer, 2024, 105--140. \href{https://doi.org/10.1007/978-3-031-57789-5_4}{doi:10.1007/978-3-031-57789-5\_4}. Preprint: \href{https://arxiv.org/abs/2301.11593v3}{arXiv:2301.11593v3}. Theorem numbering follows this preprint.

\bibitem[Jelisiejew--Smoktunowicz(26)]{JS}
J.~Jelisiejew and A.~Smoktunowicz, \emph{On a conjecture by Michael Wemyss regarding the calculation of GV invariants}, \href{https://arxiv.org/abs/2602.05919v1}{arXiv:2602.05919v1} (2026).

\bibitem[Karmazyn(19)]{Karmazyn}
J.~Karmazyn, \emph{The length classification of threefold flops via noncommutative algebras}, Adv. Math. \textbf{343} (2019), 393--447. Preprint page references and the convention comparison in Remark~\ref{rem:karmazyn-sign-conventions} use \href{https://arxiv.org/abs/1709.02720v1}{arXiv:1709.02720v1} (2017).

\bibitem[Karmazyn(18)]{KarmazynDeformations}
J.~Karmazyn, \emph{Deformations of algebras defined by tilting bundles}, J. Algebra \textbf{513} (2018), 388--434.

\bibitem[Karmazyn--Pan(26)]{KP26}
J.~Karmazyn and S.~Pan, \emph{Centres and semisimple deformations of contraction algebras}, \href{https://arxiv.org/abs/2611.09220v1}{arXiv:2611.09220v1} (2026).

\bibitem[Kawamata(24)]{KawamataNC}
Y.~Kawamata, \emph{Deformations over non-commutative base}, C. R. Math. Acad. Sci. Paris \textbf{362} (2024), special issue S1, 159--169. \href{https://doi.org/10.5802/crmath.622}{doi:10.5802/crmath.622}. Preprint: \href{https://arxiv.org/abs/2303.06323v7}{arXiv:2303.06323v7}. Theorem numbering follows this preprint; its Theorem~3.1 is Theorem~4 in the published version.

\bibitem[Keller(11)]{KellerCY}
B.~Keller, with an appendix by M.~Van den Bergh, \emph{Deformed Calabi--Yau completions}, J. Reine Angew. Math. \textbf{654} (2011), 125--180. \href{https://arxiv.org/abs/0908.3499v5}{arXiv:0908.3499v5}.

\bibitem[Keller(18)]{KellerCYErratum}
B.~Keller, \emph{Erratum to ``Deformed Calabi--Yau completions''}, \href{https://arxiv.org/abs/1809.01126}{arXiv:1809.01126} (2018).

\bibitem[Keller(06)]{KellerDG}
B.~Keller, \emph{On differential graded categories}, \href{https://arxiv.org/abs/math/0601185v5}{arXiv:math/0601185v5} (2006).

\bibitem[Katz--Morrison(92)]{KatzMorrison}
S.~Katz and D.~R. Morrison, \emph{Gorenstein threefold singularities with small resolutions via invariant theory for Weyl groups}, J. Algebraic Geom. \textbf{1} (1992), no.~3, 449--530. Preprint: \href{https://arxiv.org/abs/alg-geom/9202002v1}{arXiv:alg-geom/9202002v1}.

\bibitem[Linckelmann(13)]{Linckelmann}
M.~Linckelmann, \emph{Tate duality and transfer in Hochschild cohomology}, J. Pure Appl. Algebra \textbf{217} (2013), no.~12, 2387--2399. \href{https://arxiv.org/abs/1211.5999v1}{arXiv:1211.5999v1}.

\bibitem[Nabijou--Wemyss(24)]{NabijouWemyss}
N.~Nabijou and M.~Wemyss, \emph{GV and GW invariants via the enhanced movable cone}, Moduli \textbf{1} (2024), Paper e8, 1--38. \href{https://doi.org/10.1112/mod.2024.4}{doi:10.1112/mod.2024.4}.

\bibitem[Namikawa(08)]{NamikawaFlops}
Y.~Namikawa, \emph{Flops and Poisson deformations of symplectic varieties}, Publ. Res. Inst. Math. Sci. \textbf{44} (2008), no.~2, 259--314. \href{https://doi.org/10.2977/prims/1210167328}{doi:10.2977/prims/1210167328}.

\bibitem[Poonen(08)]{Poonen}
B.~Poonen, \emph{The moduli space of commutative algebras of finite rank}, J. Eur. Math. Soc. \textbf{10} (2008), no.~3, 817--836. \href{https://doi.org/10.4171/JEMS/131}{doi:10.4171/JEMS/131}.

\bibitem[Riemenschneider(n.d.)]{RiemenschneiderComplexSpaces}
O.~Riemenschneider, \emph{Supplement: A short introduction to complex analytic spaces}, Chapter~19 of online lecture notes, Section~19.22, Theorems~19.57 and~19.59. \href{https://www.math.uni-hamburg.de/home/riemenschneider/bookchap30.pdf}{Author's online notes}.

\bibitem[Sard(42)]{Sard}
A.~Sard, \emph{The measure of the critical values of differentiable maps}, Bull. Amer. Math. Soc. \textbf{48} (1942), 883--890. \href{https://webhomes.maths.ed.ac.uk/~v1ranick/papers/sard.pdf}{Original article}.

\bibitem[Smoktunowicz(25)]{SmoktunowiczFlat}
A.~Smoktunowicz, \emph{On flat deformations and their applications}, \href{https://arxiv.org/abs/2509.10121v2}{arXiv:2509.10121v2} (2025).

\bibitem[Smoktunowicz--Smoktunowicz(26)]{SmoktunowiczSmoktunowicz}
A.~Smoktunowicz and A.~Smoktunowicz, \emph{On computing deformations of contraction algebras of type E}, \href{https://arxiv.org/abs/2610.10157v1}{arXiv:2610.10157v1} (2026).

\bibitem[Stacks Project()]{Stacks}
The Stacks Project Authors, \emph{The Stacks Project}, \href{https://stacks.math.columbia.edu}{stacks.math.columbia.edu}, 2026; Tags~\href{https://stacks.math.columbia.edu/tag/00MP}{00MP} and~\href{https://stacks.math.columbia.edu/tag/08IB}{08IB}.

\bibitem[Toda(15)]{Toda}
Y.~Toda, \emph{Non-commutative width and Gopakumar--Vafa invariants}, Manuscripta Math. \textbf{148} (2015), no.~3--4, 521--533.

\bibitem[Van den Bergh(04)]{VanDenBergh}
M.~Van den Bergh, \emph{Three-dimensional flops and non-commutative rings}, Duke Math. J. \textbf{122} (2004), no.~3, 423--455.

\bibitem[Katz(08)]{KatzGV}
S.~Katz, \emph{Genus zero Gopakumar--Vafa invariants of contractible curves}, J. Differential Geom. \textbf{79} (2008), no.~2, 185--195. \href{https://arxiv.org/abs/math/0601193}{arXiv:math/0601193}.

\bibitem[Iyama--Oppermann(13)]{IyamaOppermann}
O.~Iyama and S.~Oppermann, \emph{Stable categories of higher preprojective algebras}, Adv. Math. \textbf{244} (2013), 23--68. \href{https://arxiv.org/abs/0912.3412v2}{arXiv:0912.3412v2}.

\bibitem[Crawley-Boevey(01)]{CrawleyBoeveyMomentMap}
W.~Crawley-Boevey, \emph{Geometry of the moment map for representations of quivers}, Compos. Math. \textbf{126} (2001), 257--293. \href{https://www.math.uni-bielefeld.de/~wcrawley/moment4.pdf}{Author's version}.

\bibitem[P\'erez Rodr\'iguez(16)]{PerezRodriguezFormalDeformations}
M.~P\'erez Rodr\'iguez, \emph{Deformation of formal schemes through local homology}, J.~Algebra \textbf{445} (2016), 78--102. \href{https://doi.org/10.1016/j.jalgebra.2015.07.041}{doi:10.1016/j.jalgebra.2015.07.041}. Version checked: \href{https://arxiv.org/abs/2402.02879v1}{arXiv:2402.02879v1}. The cited locators follow this preprint.

\bibitem[Serre(56)]{SerreGAGA}
J.-P.~Serre, \emph{G\'eom\'etrie alg\'ebrique et g\'eom\'etrie analytique}, Ann. Inst. Fourier \textbf{6} (1956), 1--42, nos.~6, 9 and~12. \href{https://doi.org/10.5802/aif.59}{doi:10.5802/aif.59}.

\end{thebibliography}
\end{document}